\documentclass[a4paper,12pt]{article}
\usepackage{dsfont}
\usepackage{amssymb}
\usepackage{latexsym}
\usepackage{amsmath}
\usepackage{xcolor}
\usepackage{comment}
\usepackage{amsthm}
\usepackage[dvips]{graphicx}
\usepackage{layout} 
\usepackage{ascmac} 
\usepackage{ulem}
\usepackage{enumerate}
\usepackage{bm}
\usepackage{bbm} 
\usepackage[bbgreekl]{mathbbol} 
\usepackage{authblk}
\usepackage{setspace} 
\usepackage{eufrak}
\newif\ifrs
\rstrue
\ifrs \usepackage{mathrsfs} \fi  
\newif\ifcol
\coltrue 

\newtheorem{theorem}{Theorem}[section]

\newtheorem{lemma}[theorem]{Lemma}

\newtheorem{proposition}[theorem]{Proposition}
\newtheorem{corollary}[theorem]{Corollary}
\newtheorem{remark}[theorem]{Remark}

\numberwithin{equation}{section}
\newtheorem{theorem*}{Theorem}
\newtheorem{ass*}[theorem*]{Assumption}
\newtheorem{note*}[theorem*]{Note}
\newtheorem{lemma*}[theorem*]{Lemma}
\newtheorem{definition*}[theorem*]{Definition}
\newtheorem{proposition*}[theorem*]{Proposition}
\newtheorem{corollary*}[theorem*]{Corollary}
\newtheorem{remark*}[theorem*]{Remark}
\newtheorem{example*}[theorem*]{Example}
\numberwithin{equation}{section}
\newif\ifcol
\colfalse
\ifcol
\newcommand{\colorr}{\color[rgb]{0.8,0,0}}

\newcommand{\colorb}{\color[rgb]{0,0,0.8}}

\newcommand{\colorn}{\color[rgb]{1,1,1}}

\newcommand{\coloro}{\color[rgb]{1,0.4,0}}
\newcommand{\coloroy}{\color[rgb]{1,0.95,0}}

\else
\newcommand{\colorb}{\color{black}}

\newcommand{\colorr}{\color{black}}

\newcommand{\colorn}{\color{black}}
\newcommand{\coloro}{\color{black}}

\newcommand{\coloroy}{\color{black}}
\fi
\newcommand{\cred}{\color{black}}
\newif\ifcol
\colfalse
\ifcol
\newcommand{\sred}{\color[rgb]{0.8,0,0}}

\else
\newcommand{\sred}{\color{black}}
\fi
\excludecomment{en-text}
\includecomment{jp-text}
\includecomment{comment}
\def\mH{{\mathfrak H}}

\def\infm{{\infty\text{--}}}
\def\inftym{\infm}
\def\koko{{\coloroy{koko}}}
\def\bd{\begin{description}}
\def\ed{\end{description}}

\def\D2{\bbD_{2,\infty-}}

\def\tj{{t_j}}
\def\tjm{{t_{j-1}}}

\def\n{{\bf n}}

\def\D{{\bf D}}

\def\R{{\bf R}}

\def\V{{\bf V}}

\def\calb{{\cal B}}
\def\calc{{\cal C}}

\def\cale{{\cal E}}
\def\calf{{\cal F}}

\def\calh{{\cal H}}
\def\cali{{\cal I}}
\def\calj{{\cal J}}
\def\calk{{\cal K}}
\def\call{{\cal L}}
\def\calm{{\cal M}}

\def\calq{{\cal Q}}

\def\cals{{\cal S}}

\def\calv{{\cal V}}

\def\calx{{\cal X}}
\def\caly{{\cal Y}}

\def\ds{\displaystyle}
\def\yeq{\>=\>}
\def\yleq{\>\leq\>}
\def\ygeq{\>\geq\>}

\def\sfm{{\sf m}}

\def\sfd{{\sf d}}

\def\simleq{\ \raisebox{-.7ex}{$\stackrel{{\textstyle <}}{\sim}$}\ }

\def\ep{\epsilon}
\def\half{\frac{1}{2}}

\def\Iku{\Rightarrow}

\def\down{\downarrow}

\def\y{\vspace*{3mm}\\}
\def\halflineskip{\vspace*{3mm}}
\def\nn{\nonumber}
\def\be{\begin{equation}}
\def\ee{\end{equation}}
\def\bea{\begin{eqnarray}}
\def\eea{\end{eqnarray}}
\def\beas{\begin{eqnarray*}}
\def\eeas{\end{eqnarray*}}
\def\bi{\begin{itemize}}
\def\ei{\end{itemize}}
\def\im{\item}
\def\bd{\begin{description}}
\def\ed{\end{description}}
\def\r{\right}

\def\dotc{\stackrel{\circ}{C}}

\newcommand{\bbD}{{\mathbb D}}
\newcommand{\bbE}{{\mathbb E}}

\newcommand{\bbH}{{\mathbb H}}
\newcommand{\bbI}{{\mathbb I}}
\newcommand{\bbJ}{{\mathbb J}}
\newcommand{\bbK}{{\mathbb K}}
\newcommand{\bbL}{{\mathbb L}}

\newcommand{\bbN}{{\mathbb N}}

\newcommand{\bbR}{{\mathbb R}}
\newcommand{\bbS}{{\mathbb S}}
\newcommand{\bbT}{{\mathbb T}}

\newcommand{\bbV}{{\mathbb V}}
\newcommand{\bbW}{{\mathbb W}}
\newcommand{\bbX}{{\mathbb X}}
\newcommand{\bbY}{{\mathbb Y}}
\newcommand{\bbZ}{{\mathbb Z}}

\def\mba{{\mathbb a}}
\def\mbbb{{\mathbb b}}

\def\mbs{\mathbb s}
\def\mbt{\mathbb t}

\def\mbv{\mathbb v}

\def\mbx{\mathbb x}

\def\mfh{\mH}
\def\tk{{t_k}}
\def\tkm{{t_{k-1}}}
\def\ol{\overline}
\def\wt{\widetilde}
\def\dotx{\stackrel{\circ}{X}}

\def\mfh{{\EuFrak H}}
\def\dota{\stackrel{\circ}{a}\!} 
\def\ddota{\stackrel{\circ\circ}{a}\!}
\def\tti{{\tt i}}
\def\onelineskip{\halflineskip\halflineskip}

\newcommand{\sfx}{{\sf x}}
\newcommand{\sfz}{{\sf z}}

\def\sfd{{\sf d}}
\def\sfm{{\sf m}}

\begin{document}

\title{Asymptotic expansion of a variation with anticipative weights 
\footnote{
This work was in part supported by 
Japan Science and Technology Agency CREST JPMJCR14D7; 
Japan Society for the Promotion of Science Grants-in-Aid for Scientific Research 
No. 17H01702 (Scientific Research);  
and by a Cooperative Research Program of the Institute of Statistical Mathematics. 
}
}
\author{Nakahiro Yoshida}
\affil{Graduate School of Mathematical Sciences, University of Tokyo
\footnote{Graduate School of Mathematical Sciences, University of Tokyo: 3-8-1 Komaba, Meguro-ku, Tokyo 153-8914, Japan. e-mail: nakahiro@ms.u-tokyo.ac.jp}}
\affil{CREST, Japan Science and Technology Agency
}
\date{First version: May 19, 2020, \\
Second version: December 29, 2020
}
\maketitle
\ \\
{\it Summary} \hspace{5pt}
Asymptotic expansion of a variation with anticipative weights is derived 
by the theory of asymptotic expansion for Skorohod integrals having a mixed normal limit. 
The expansion formula is expressed with the quasi-torsion, quasi-tangent and other random symbols. 
To specify these random symbols, 
it is necessary to classify the level of the effect of each term appearing in 
the stochastic expansion of the variable in question.  
To solve this problem, we consider a class $\call$ of certain sequences $(\cali_n)_{n\in\bbN}$ of Wiener functionals 
and we give a systematic way of estimation of the order of $(\cali_n)_{n\in\bbN}$. 
Based on this method, we introduce a notion of {\it exponent} of 
the sequence $(\cali_n)_{n\in\bbN}$, and 
investigate the stability and contraction effect of the operators $D_{u_n}$ and $D$ on $\call$, 
{\sred where $u_n$ is the integrand of a Skorohod integral.} 
After constructed these machineries, we derive asymptotic expansion of the variation having 
anticipative weights. An application to robust volatility estimation is mentioned. 

\section{Introduction}
Let $\bbW=\{\bbW(h)\}_{h\in\mH}$ be an isonormal Gaussian process on a real separable Hilbert space $\mH$ 
with inner product $\langle\cdot,\cdot\rangle$. 
That is, $\bbW$ is a system of centered Gaussian variables on a probability space 
$(\Omega,\calf,P)$ with 
\bea\label{202003130048} 
E\big[\bbW(h_1)\bbW(h_2)\big] &=& \langle h_1,h_2\rangle
\eea
for $h_1,h_2\in\mH$. 
Then we have the Malliavin derivative $D$ and the divergence operator $\delta$, i.e., 
the adjoint operator $D^*$ of $D$. The Sobolev space with differentiability index $s$ 
and integrability index $p$ is denoted by $\bbD^{s,p}$. 
The $q$-fold multiple Wiener integral is denoted by $I_q(f)=\delta^k(f)$ for $f\in\mH^{\odot q}$, 
the space of symmetric elements of $\mH^{\otimes q}$. 
For the Malliavin calculus, we refer the reader to 
Watanabe \cite{watanabe1984lectures}. 
Ikeda and Watanabe \cite{IkedaWatanabe1989} and 
Nualart \cite{Nualart2006}. 
Nourdin and Pecati \cite{nourdin2012normal} is an excellent exposition of applications 
of the Malliavin calculus to central limit theorems by Stein's method.

In this paper, we are specially interested in asymptotic expansion of a variation of a Wiener process 
with random weights. This problem is {\sred deeply} related to the expansion of 
the power variation of a diffusion process time-discretely observed; 
see Podolskij and Yoshida \cite{podolskij2016edgeworth} and Podolskij et al. \cite{podolskij2017edgeworth}. 
In what follows, we suppose that $\mH=L^2({\colorr \bbR_+})$, the $L^2$-space of 
$\bbR$-valued functions on ${\sred\bbR_+=[0,\infty)}$, and that 
\beas 
\langle h_1,h_2\rangle &=& \int_0^{{\colorr\infty}} h_1(t) h_2(t) dt
\eeas
for $h_1,h_2\in\mH$. 
Let $w_t=\bbW(1_{[0,t]})$ for $t\in{\sred\bbR_+}$. 
Then process $w=(w_t)_{t\in{\sred\bbR_+}}$ is a one-dimensional standard Wiener process. 

Let $t_j=t^n_j=j/n$ $(j=0,1,...,n)$ and let denote by $1_j$ the indicator function 
$1_{I_j}$, $I_j=[\tjm,\tj]$. 
{\colorb Suppose that $\calq$ is a finite subset of $\{2,3,4,....\}$. 
Let us consider the variation
\bea\label{202004021410}
V_n
&=&
\sum_{q\in\calq}
n^{2^{-1}(q-1)}\sum_{j=1}^na_j(q) I_q(1_j^{\otimes q})
\eea
with random weights $a_j(q)=a_j^n(q)$ depending on $j\in\bbJ_n=\{1,...,n\}$ and $n\in\bbN$. 
The weight $a_j(q)$ is not necessarily predictable with respect to the Brownian filtration. 
}

For each $q\in\calq$, we consider an $\mH$-valued functional
\bea\label{202003130046}
u_n &=& \sum_{q\in\calq}u_n(q)
\eea
given by
\bea\label{202003130047}
u_n(q)
&=& 
n^{-1/2}\sum_{j=1}^n a_j(q) I_{q-1}((n^{1/2}1_j)^{\otimes(q-1)})n^{1/2}1_j 
\nn\\&=& 
n^{2^{-1}(q-1)}\sum_{j=1}^na_j(q) I_{q-1}(1_j^{\otimes (q-1)})1_j. 
\eea
For example, $a_j(q)=a_q(w_\tj,w_{1-\tj})$ for a function $a_q$ and a Wiener process $w$, 
however we will treat more general functionals $a_j(q)$. 
Suppose for the meantime that 
\beas
a_j(q)\in\bbD^{1,2+}=\cup_{p>2}\bbD^{1,p}
\quad (j\in\bbJ_n;\>q\in\calq)
\eeas
\begin{en-text}
Suppose that 
\bea\label{0112121022}
a_j(q)\in\bbD^{1,\infty-}=\cap_{p>1}\bbD^{1,p}
\quad (j\in\bbJ_n;\>q\in\calq)
\eea
\end{en-text}
In particular $u_n(q)\in \text{Dom}(\delta)$ and 
\bea\label{0112121101}
M_n
&:=&
\delta(u_n)
\eea
is well defined and in $L^2$. 
Equivalently, 
\bea\label{202003171321} 
M_n &=& \sum_{q\in\calq}M_n(q)
\eea
where
\bea\label{202003181159} 
M_n(q) &=& \delta(u_n(q))
\nn\\&=&
n^{2^{-1}(q-1)}\sum_{j=1}^na_j(q) I_q(1_j^{\otimes q})
-n^{2^{-1}(q-1)}\sum_{j=1}^n(D_{1_j}a_j(q) )I_{q-1}(1_j^{\otimes (q-1)})
\eea
Here $D_{1_j}a_j(q)=\langle Da_j(q),1_j\rangle$. 

Consider a sequence of random variables 
\bea\label{202003130001}
Z_n &=& M_n+n^{-1/2}N_n\quad(n\in\bbN)
\eea
where $(N_n)_{n\in\bbN}$ is a sequence of random variables. 
Variable $N_n$ is of order one and satisfying some regularity conditions, 
as specified more clearly later. 
Therefore $Z_n$ is a perturbation of $M_n$ by $n^{-1/2}N_n$. 
{\colorb In particular, $V_n$ is an example of $Z_n$ if one sets 
\bea\label{202004151750} 
\sum_{q\in\calq}n^{q/2}\sum_{j=1}^n(D_{1_j}a_j(q) )I_{q-1}(1_j^{\otimes (q-1)})
\eea
to $N_n$. 
}\noindent
The variable $Z_n$ of the form (\ref{202003130001}) is standard in analysis of 
the variation of a process driven by $w$ and statistics for volatility as well. 
The perturbation (\ref{202003130001}) appeared in 
\begin{en-text}
Yoshida \cite{Yoshida1997,yoshida2001malliavin,yoshida2013martingale, yoshida2012asymptotic}, 
Dalalyan and Yoshida \cite{DalalyanYoshida2011}, 
Podolskij and Yoshida \cite{podolskij2016edgeworth}, 
Podolskij et al. \cite{podolskij2018edgeworth,podolskij2018edgeworth, podolskij2018edgeworthToappear} 
for a martingale $M_n$, and Nualart and Yoshida \cite{nualart2019asymptotic} 
for a Skorohod integral $M_n$. 
\end{en-text}
\cite{Yoshida1997}, \cite{yoshida2001malliavin}, \cite{yoshida2013martingale}, \cite{yoshida2012asymptotic}, 
\cite{DalalyanYoshida2011}, 
\cite{podolskij2016edgeworth}, 
\cite{podolskij2018edgeworthToappear},  
\cite{nualart2019asymptotic} etc.

Statistics appearing here is the so-called non-ergodic statistics, 
where the limit distribution of $M_n$ 
is typically a mixed normal distribution. 
In the statistical context, the Fisher information becomes random even in the limit, and 
it causes a variance mixture of normal distributions. 
In Section \ref{202003141649}, we will provide asymptotic expansion of the distribution of $(V_n,X_\infty)$ 
for a reference variable $X_\infty$. 
Section \ref{202004021507} treats the general functional $Z_n$ of (\ref{202003130001}) without the specific expression (\ref{202004151750}) for $N_n$, and derives asymptotic expansion of 
the joint distribution of $(Z_n,X_n)$ for a sequence of reference variables $X_n$. 

Asymptotic expansion for a sequence of random variables asymptotically having 
a mixed normal distribution was derived in Yoshida \cite{yoshida2013martingale} 
for continuous martingales 
by introducing the notion of random symbols 
called {\it adaptive random symbol} for {\it tangent} 
and {\it anticipative random symbol} for {\it torsion}, and applied to 
the power variation by Podolskij and Yoshida \cite{podolskij2016edgeworth}, 
the quadratic variation under microstructure noise by 
Podolskij et al. \cite{podolskij2017edgeworth}, and to 
a precise error estimate for the Euler-Maruyama scheme 
by Podolskij et al. \cite{podolskij2018edgeworthToappear}. 
Even if one starts with the martingale framework, asymptotic expansion is not closed in it, 
and the Malliavin derivative appears in the anticipative random symbol. 
Nualart and Yoshida \cite{nualart2019asymptotic} introduced 
new random symbols called {\it quasi-tangent} ({\it q-tangent}) and {\it quasi-torsion} ({\it q-torsion}), 
and provided 
asymptotic expansion of the Skorohod integral with these random symbols.\footnote{We remark that the correspondence between the two pairs of random symbols 
is not one to one. In particular, the quasi-tangent does not correspond to the tangent 
expressed by the adaptive random symbol. }
They applied this scheme to some examples including a randomly weighted 
quadratic form of the increments of a fractional Brownian motion. 
Section \ref{202003141634} recalls the theory of asymptotic expansion for Skorohod integrals.

As summarized in Section \ref{202003141634},  
it is necessary to identify each limit of the q-tangent and the q-torsion 
for deriving the asymptotic expansion formula. 
When we apply the theory to $V_n$, 
it is not easy to distinguish the effect of each term in the limit since 
the expression of random symbols becomes fairly complicated, and 
really it is more difficult 
for $M_n$ of (\ref{0112121101}) than the quadratic form treated 
in Nualart and Yoshida \cite{nualart2019asymptotic}. 
To solve this problem, we consider a class $\call$ of sequences $(\cali_n)_{n\in\bbN}$ of Wiener functionals; 
the definition of $\call$ is in Section \ref{202004270705}. 
In Section \ref{202003141640}, we present estimates to $(\cali_n)_{n\in\bbN}$ in $\call$ 
in order to make the classifications easier in derivation of the expansion formula. 
The author hopes this section is of interest even apart from asymptotic expansion. 
Based on the estimate in Section \ref{202003141640}, 
we introduce a notion of {\it exponent} of 
the sequence $(\cali_n)_{n\in\bbN}$ in Section \ref{202004270637}, and 
investigate the effect of the operators $D_{u_n(q)}$ and $D$ 
{\sred applied} on $\call$. 
Section \ref{202003141657} gives a proof to the results in Section \ref{202003141649} 
by applying 
the theory of asymptotic expansion for Skorohod integrals recalled in Section \ref{202003141634}, while 
Section \ref{202004270654} validates the results in Section \ref{202004021507}. 
Our results are general. There are many applications though we only mention two examples in 
Sections \ref{202005180614} and \ref{202004261806}. 
The latter section treats asymptotic expansion for an anticipative functional 
appearing in robust volatility estimation. 

Concluding Introduction, we mention developments 
of the theory of asymptotic expansion (Edgeworth expansion) 
in the central limit cases, listing 
some basic literature in the early days (because giving a complete list in a paragraph is impossible) 
to guide the reader to the history of asymptotic expansion and to indicate where we are now. 
When the limit of a sequence of random variables is exactly Gaussian, 
as it is ordinarily true under a sufficiently fast mixing property, 
asymptotic expansion was proved 
in various scenarios. 
Basic references are 
Cram\'er \cite{cramer1928composition, cramer2004random}, 
Gnedenko and Kolmogorov \cite{gnedenko1954limit}, 
Bhattacharya \cite{Bhattacharya1971}, 
Petrov \cite{Petrov1975}, 
Bhattacharya and Ranga Rao \cite{BhattacharyaRanga1976} and 
Bhattacharya and Ghosh \cite{BhattacharyaGhosh1978} among huge historical literature 
for independent observations. 
Asymptotic expansion has been established as a basic tool supporting the modern statistical theory. 
Higher-order asymptotic theory cannot exist without asymptotic expansion: 
Akahira and Takeuchi \cite{AkahiraTakeuchi1981}, Pfanzagl \cite{Pfanzagl1985} and Ghosh \cite{Ghosh1994}.  
%
Historically, 
multivariate analysis is one of the fields where asymptotic expansion techniques have been well developed 
(Okamoto \cite{okamoto1963asymptotic}). 
Informative textbooks for asymptotic expansion in the 
multivariate analysis are 
Anderson \cite{anderson1962introduction} and 
Fujikoshi, Ulyanov and Shimizu \cite{fujikoshi2011multivariate}. 
Theory of Bootstrap methods (Efron \cite{efron1979bootstrap}) 
has a theoretical basis on asymptotic expansion (Hall \cite{Hall1992}). 
The reader is referred to 
Barndorff-Nielsen \cite{barndorff2012parametric}, 
Jensen \cite{jensen1995saddlepoint} and 
Pace and Salvan \cite{PaceSalvan1997} 
for saddle-point approximations. 
Combining asymptotic expansion with the notion of the $\alpha$-connection, 
information geometry 
gave an interpretation of the higher-order efficiency of estimation  
in terms of the curvature of the fibre associated with the estimator (Amari \cite{Amari1985}). 
Asymptotic expansion was also applied to construction of information criteria for 
model selection as well as prediction problems; 
e.g., Konishi and Kitagawa \cite{KonishiKitagawa1996}, 
Uchida and Yoshida \cite{UchidaYoshida2001,UchidaYoshida2004}, Komaki \cite{Komaki1996}. 
Bhattacharya et al. \cite{bhattacharya2016course} is an excellent textbook including 
a concise exposition of asymptotic expansion. 
Dependent cases were studied by 
G\"otze and Hipp \cite{GotzeHipp1983, GotzeHipp1994} for (approximate) Markov chains having a mixing property. 
It should be commented that in this short paragraph, we can not list a large number of studies behind the epoch-making paper by G\"otze and Hipp \cite{GotzeHipp1983}, but 
a few related papers to be mentioned are 
Nagaev \cite{nagaev1957some} and Tikhomirov \cite{tikhomirov1981convergence},  
and Lahiri \cite{lahiri1996asymptotic} after \cite{GotzeHipp1983}. 
Specialized treatments are possible for Gaussian processes (Taniguchi \cite{Taniguchi1991}). 
With stochastic analysis, asymptotic expansion was studied 
by Kusuoka and Yoshida \cite{KusuokaYoshida2000} and Yoshida \cite{Yoshida2004} 
for continuous-time mixing Markov processes, and by Mykland \cite{Mykland1992,Mykland1993} and 
Yoshida \cite{Yoshida1997, yoshida2001malliavin} for martingales 
{\sred about moments and distributions, respectively.} 
The martingale expansion with mixed normal limit distribution (Yoshida \cite{yoshida2013martingale}) 
has been applied to power variations and the Euler-Maruyama approximation 
of a stochastic differential equation, as already mentioned. 
The Malliavin calculus played an essential role in the above studies 
\cite{Yoshida1997, yoshida2001malliavin, KusuokaYoshida2000, Yoshida2004,yoshida2013martingale}. 
Asymptotic expansion of the heat kernel was an important application of the Malliavin calculus 
(Watanabe \cite{Watanabe1987}, Ikeda and Watanabe \cite{IkedaWatanabe1989}), though 
purely distributional treatment is not necessary, differently from the limit theorems we are now aiming at. 
In statistics and finance, asymptotic expansion has various applications 
even if limited to the objects in stochastic analysis: 
for example, 
Yoshida \cite{Yoshida1992, Yoshida1993, Yoshida1996a, yoshida2003conditional} 
based on the Watanabe theory \cite{Watanabe1987} in the Malliavin calculus, 
Dermoune and Kutoyants \cite{DermouneKutoyants1995}, Kutoyants \cite{Kutoyants2004, Kutoyants1998}, 
Sakamoto and Yoshida \cite{
SakamotoYoshida1996, 
SakamotoYoshida1998a, 
SakamotoYoshida2003, 
SakamotoYoshida2004, 
SakamotoYoshida2009, 
SakamotoYoshida2010} together with an application 
(Sakamoto et al. \cite{sakamoto2004expansions}) 
to an inadmissibility problem (Stein's phenomenon), 
Masuda and Yoshida \cite{
MasudaYoshida2005}, 
Kutoyants and Yoshida \cite{KutoyantsYoshida2007}, 
Dalalyan and Yoshida \cite{DalalyanYoshida2011}, 
Li \cite{li2013maximum} 
in statistics for stochastic processes, 
and in finance 
Yoshida \cite{Yoshida1992a}, 
Kunitomo and Takahashi \cite{KunitomoTakahashi2001}, 
Takahashi and Yoshida \cite{TakahashiYoshida2004, TakahashiYoshida2005}, 
Uchida and Yoshida \cite{UchidaYoshida2004b}, 
Li \cite{li2014closed, li2016estimating} and Cai et al. \cite{cai2014closed} 
for option pricing among plenty of 
literature; 
a related paper is Hayashi and Ishikawa \cite{hayashi2012composition}. 
YUIMA has implemented asymptotic expansion techniques 
for automatic fast computation of the expectation of a generic functional defined 
by a system of stochastic differential equations 
(\cite{Yuima2014}, \cite{iacus2017simulation}). 
%
We refer the reader to Yoshida \cite{yoshida2016asymptotic} for a short exposition of  
asymptotic expansion and additional references therein. 
Today asymptotic expansion has been motivated once again 
in the surge of studies on central limit theorems for Wiener functionals after the forth moment theorem 
(Nualart and Peccati \cite{nualart2005central}). 
Tudor and Yoshida \cite{tudor2019asymptotic} provided 
asymptotic expansion for general Wiener functionals. 
Arbitrary order of asymptotic expansion for Wiener functionals was obtained by 
Tudor and Yoshida \cite{tudor2019high} in the ergodic case and applied to a stochastic wave equation. 
{\sred This scheme was also applied to a mixed fractional Brownian motion 
in Tudor and Yoshida \cite{tudor2020asymptotic}.} 
Yoshida \cite{yoshida2013martingale} and Nualart and Yoshida \cite{nualart2019asymptotic} derived asymptotic expansion having a mixed normal limit. 
The last paper is the starting point of this article.

\section{Asymptotic expansion of {\colorb$V_n$} 
}\label{202003141649}
{\colorb
We will present asymptotic expansion to the joint law $\call\{(V_n,X_\infty)\}$ 
for $V_n$ defined in (\ref{202004021410}) and a reference variable $X_\infty$. 
The random asymptotic variance of $V_n$ is an example of $X_\infty$ but 
a generic variable $X_\infty$ will be treated. 
In this section, we give priority to simplicity of presentation though 
we can generalize the results in this section by replacing $X_\infty$ 
by a sequence of variables $X_n$ ($n\in\bbN$), as will do 
it for $Z_n$ of (\ref{202003130001}) with a generic $N_n$ together with $X_n$ 
in Section \ref{202004021507}. 
}

\subsection{Asymptotic expansion of $E[f(V_n,X_\infty)]$ for a smooth function $f$}
The following set of conditions is concerning certain regularity of $a_j(q)$ 
and limits of some related variables, in particular, it involves $a(t,q)$ for $(t,q)\in[0,1]\times\calq$. 
We will write 
\bea\label{202004201141}
G_\infty &=& \sum_{q\in\calq}q! \int_0^1 a(t,q)^2dt.
\eea
To write the asymptotic expansion formula, we use measurable random fields 
$(a(t,q))_{t\in[0,1]}$, $(\dota(t,s,q))_{(t,s)\in[0,1]^2}$, 
$(\dot{a}(t,2))_{t\in[0,1]}$ (if $2\in\calq$), 
$(\ddota(t,s,q))_{(t,s)\in[0,1]^2}$ 
$(q\in\calq)$
{\colorr and $(\ddot{a}(t,2))_{t\in[0,1]}$ (if $2\in\calq$)}.

\begin{en-text}
\bd
\im[[A\!\!]] 
\begin{enumerate}[(i)]
\im\label{la} a
\im b
\im c 
\end{enumerate}
(\ref{la})
\ed
\end{en-text}

\bd
\im[[A\!\!]] The following properties hold.  
\begin{enumerate}[(i)] 
\im\label{ai}
$a_j(q)\in\bbD^{\infty}$ for $j\in{\sred\bbJ_n}$, 
$n\in\bbN$ and $q\in\calq$, 
and that 
\bea\label{2020031201810a} 
\sup_{n\in\bbN}\sup_{j\in\bbJ_n}
\sup_{t_1,...,t_i\in[0,1]}
\big\|D_{t_i}\cdots D_{t_1}a_j(q)\big\|_{p}&<&\infty
\eea
for every $i\in\bbZ_+=\{0,1,...\}$, $p>1$ and $q\in\calq$. 
(Recall that $a_j(q)$ depends on $n$.)
%

\begin{en-text}
\im\label{aii} 
For every $p>1$, 
\bea\label{202003251201} &&
\max_{q\in\calq}\sup_{t\in[0,1]}\big\|a(t,q)\big\|_p
+\max_{q\in\calq}\sup_{(t,s)\in[0,1]^2}\big\|\dota(t,s,q)\big\|_p
\nn\\&&
+\sup_{t\in[0,1]}\big\|\dot{a}(t,2)\big\|_p
+\max_{q\in\calq}\sup_{(t,s)\in[0,1]^2}\big\|\ddota(t,s,q)\big\|_p
\><\> 
\infty.
\eea
\end{en-text}

\im\label{aiii} 
For every $q\in\calq$ and every $\ep>0$, 
\bea\label{202003251202}
E\bigg[\int_{[0,1]}1_{\big\{\big|\sum_{j=1}^n 1_{I_j}(s)a_j(q) - a(s,q) \big|>\ep\big\}}ds\bigg]
&\to&
0
\eea
as $n\to\infty$. 
\im\label{aiv} 
For every $q\in\calq$ and every $\ep>0$, 
\bea\label{202003251203}
E\bigg[\int_{[0,1]^2}
1_{\big\{\big|\sum_{j,k=1}^n1_{I_k\times I_j}(t,s)nD_{1_k}a_j(q) - \dota\>(t,s,q) \big|>\ep\big\}} dsdt 
\bigg]
&\to&
0
\eea
as $n\to\infty$. 
\im\label{av} 
For every $\ep>0$, 
\bea\label{202003251204}
E\bigg[ \int_{[0,1]}
1_{\big\{\big|\sum_{k=1}^n1_{I_k}(t)nD_{1_k}a_k(2) - \dot{a}(t,2) \big|>\ep\big\}} dt \bigg]
&\to&
0
\eea
as $n\to\infty$ if $2\in\calq$. 
\im\label{avi} 
For every $q\in\calq$ and every $\ep>0$, 
\bea\label{202003251205}
E\bigg[ \int_{[0,1]^2}
1_{\big\{
\big|\sum_{j,k=1}^n1_{I_k\times I_j}(t,s) n^2D_{1_k}D_{1_k}a_j(q)-\ddota(t,s,q)\big|>\ep\big\}}dsdt\bigg]
&\to&
0
\eea
as $n\to\infty$. 
\im\label{avii} 
{\colorr If $2\in\calq$, then 
} 
\bea\label{202003251205ddot}
E\bigg[ \int_{[0,1]^2}
1_{\big\{
\big|\sum_{j=1}^n1_{ I_j}(t) n^2D_{1_j}D_{1_j}a_j({\colorr2})-\ddot{a}(t,2)\big|>\ep\big\}}dsdt\bigg]
&\to&
0
\eea
as $n\to\infty$ for every $\ep>0$.

\im\label{aviii} 
$G_\infty\in\bbD^{\infty}$ and 
\beas
\sup_{t_1,...,t_i\in[0,1]}\|D_{t_i}\cdots D_{t_1}G_\infty\|_p<\infty
\eeas
for every $i\in\bbN$ and $p>1$. 
Moreover, 
\bea\label{202003251920}
n^{-1}\sum_j\sum_{q\in\calq}q!a_j(q)^2-G_\infty
&=& 
o_p(n^{-1/2})
\eea
as $n\to\infty$. 

\im\label{aix} 
$X_\infty\in\bbD^{\infty}$ and 
\beas
\sup_{t_1,...,t_i\in[0,1]}\|D_{t_i}\cdots D_{t_1}X_\infty\|_p<\infty
\eeas
for every $i\in\bbN$ and $p>1$. 
Moreover, there exists a measurable process $\ddot{X}_\infty=(\ddot{X}_\infty(t))_{t\in[0,1]}$ 
such that 
\begin{en-text}
\beas 
\sup_{t\in[0,1]}\|\ddot{X}_\infty(t)\|_p &<& \infty
\eeas
for every $p>1$ and 
\end{en-text}
\bea\label{202003251205x}
E\bigg[ \int_{[0,1]}
1_{\big\{
\big|\sum_{j=1}^n1_{I_j}(t) 
n^2D_{1_j}D_{1_j}X_\infty-\ddot{X}_\infty(t)\big|>\ep\big\}}dt\bigg]
&\to&
0
\eea
as $n\to\infty$. 
\end{enumerate}
\ed 

\begin{remark}\rm
{\sred 
$D_{t_i}\cdots D_{t_1}F$ denotes a density that represents 
the $i$-th Malliavin derivative $D^iF$ of a functional $F$. 
$D_{t_i}\cdots D_{t_1}F$ is also denoted by $D_{t_1,...,t_i}F$.
In the condition, the}
infinite order of smoothness of $a_j(q)$ is for $L^\inftym$-estimate, 
{\sred as explained later.} 
$X_\infty\in\bbD^\infty$ is assumed for simplicity in presentation 
but the differentiability can be relaxed. 
Condition $[A]$ (\ref{avii}) is for handling the specific $N_n$ with expression (\ref{202004151750}) 
in the case of $V_n$. 
\end{remark}
\halflineskip

In statistics of volatility, special interest is in even polynomials of the increments of the process, 
i.e., the case $\calq\subset2\bbN$. 
However, we will treat a more general case where 
\bea\label{202004111810}
\calq\cap(\calq+1) &=& \phi, 
\eea
in other words, 
\beas
\min\{|q_1-q_2|;\>q_1,q_2\in\calq,\>q_1\not=q_2\}\geq2.
\eeas
This condition is satisfied when one considers a variation with an even power, e.g. 
The Hermite power variation with anticipative weights is a simple example of the case $\#\calq=1$, that satisfies 
(\ref{202004111810}). 
It is possible to remove Condition (\ref{202004111810}) but we do not pursue it here because 
the results and the presentation will be much more complicated. 

{\sred Let} 
\bea\label{202003181607} 
a^{(3,0)}(t,s,q_1,q_2)
&=&
\ddota(t,s,q_1)a(s,q_2)a(t,2)
+\dota(t,s,q_1)\dota(t,s,q_2)a(t,2)
\nn\\&&
+\dota(t,s,q_1)a(s,q_2)\dot{a}(t,2). 
\eea
We need the notion of random symbols (Yoshida \cite{yoshida2013martingale}). 
In this paper, we will consider only polynomial random symbols. 
A polynomial random symbols is a random polynomial of $(\sfz,\sfx)$ 
the coefficients of which are integrable random variables. 
We introduce three random symbols. 
{\sred For $\bar{q}(3)=q_1+q_2+q_3$, let} 
\bea\label{202003181606}
{\mathfrak S}^{(3,0)}(\tti\sfz,\tti\sfx)
&=& 
\frac{1}{3}\sum_{q_1,q_2,q_3\in\calq}{\colorb 1_{\{\bar{q}(3):even\}}}
(q_1+q_2)(q_1-1)c_{0.5\bar{q}(3)-2}(q_1-2,q_2-1,q_3-1)
\nn\\&&\hspace{50pt}\times
\int_0^1 a(q_1)_ta(q_2)_ta(q_3)_tdt\>(\tti\sfz)^3
\nn\\&&
+\sum_{q_1\in\calq}{\colorr1_{\{2\in\calq\}}}{\colorr q_1!}
\int_{[0,1]^2}a^{(3,0)}(t,s,q_1,{\colorr q_1})dsdt\>(\tti\sfz)^3
\nn\\&&
+\sum_{q_1\in\calq}{\colorr1_{\{2\in\calq\}}}{\colorr q_1!}
\int_{t=0}^1
\bigg[\bigg(
{\sred 2^{-1}}D_tG_\infty
(\tti \sfz)^5+D_tX_\infty(\tti\sfz)^3(\tti\sfx)\bigg)
\nn\\&&\hspace{120pt}\times
\int_0^1\dota(t,s,q_1)a(s,{\colorr q_1})ds \bigg]a(t,2)dt, 
\eea
where the random field $a^{(3,0)}(t,s,q_1,q_2)$ 
in (\ref{202003181606}) is given by (\ref{202003181607}) 
and $c_\nu(q_1,q_2,q_3)$ is defined by (\ref{202004111042}). 
The integrals appearing in the expression (\ref{202003181606}) exist a.s.
{\sred under $[A]$, since e.g. $\dota(\cdot,\cdot,q)\in L^\inftym(dPdtds)$ 
by Fatou's lemma. 
}
Let 
\bea\label{202004161636}
{\mathfrak S}^{(1,1)}(\tti\sfz,\tti\sfx)
&=& 
{\colorr 1_{\{2\in\calq\}}\bigg\{}
\int_0^1\ddot{X}_\infty(t)a(t,2)dt(\tti\sfz)(\tti\sfx)
\nn\\&&\hspace{50pt}
+\int_0^1D_tX_\infty\dot{a}(t,2)dt(\tti\sfz)(\tti\sfx)
\nn\\&&\hspace{50pt}
+\half\int_0^1(D_tG_\infty)(D_tX_\infty)a(t,2)dt(\tti\sfz)^3(\tti\sfx)
\nn\\&&\hspace{50pt}
+\int_0^1(D_tX_\infty)^2a(t,2)dt(\tti\sfz)(\tti\sfx)^2
{\colorr\bigg\}}
\nn\\&&
\eea
and let 
\bea\label{202004161637}
{\mathfrak S}^{(1,0)}(\tti\sfz)
&=&
1_{\{2\in\calq\}}\bigg\{
\half\int_0^1(D_t{\sred G}_\infty)\dot{a}(t,2)dt(\tti\sfz)^3
\nn\\&&\hspace{50pt}
+\int_0^1(D_tX_\infty)\dot{a}(t,2)dt(\tti\sfz)(\tti\sfx)
+\int_0^1\ddot{a}(t,2)dt(\tti\sfz)\bigg\}. 
\eea
Then the full random symbol for Theorem \ref{0112121225} is given by 
\bea\label{202004161621}
\mathfrak{S}(\tti\sfz,\tti,\sfx)
&=&
\mathfrak{S}^{(3,0)}(\tti\sfz,\tti\sfx) 
+\mathfrak{S}^{(1,1)}(\tti\sfz,\tti\sfx) 
+\mathfrak{S}^{(1,0)}(\tti\sfz,\tti\sfx) 
\eea
where the random symbols 
$\mathfrak{S}^{(3,0)}(\tti\sfz,\tti\sfx)$, 
$\mathfrak{S}^{(1,1)}(\tti\sfz,\tti\sfx)$ and 
$\mathfrak{S}^{(1,0)}(\tti\sfz,\tti\sfx)$ are given by 
(\ref{202003181606}), (\ref{202004161636}) and 
(\ref{202004161637}), respectively. 

The space of bounded functions on $\bbR^2$ having bounded derivatives up to order 
{\sred $K$} 
is denoted by $C^{{\sred K}}_b(\bbR^2)$, 
equipped with the norm 
$\|f\|=\sum_{k=0}^{{\sred K}}\sup_{(z,x)\in\bbR^2}|\partial_{(z,x)}^kf(z,x)|$. 
The following is the first result, the proof of which is given in Section \ref{202003141657}. 
\begin{theorem}\label{0112121225}
Suppse that Condition $[A]$ and (\ref{202004111810}) are fulfilled. Then 
\beas 
\sup_{f\in \calb}\bigg|
E\big[f(V_n,X_\infty)\big] - 
\bigg\{E\big[f(G_\infty^{1/2}\zeta,X_\infty)\big]
+n^{-1/2}E\big[\mathfrak{S}(\partial_z,\partial_x)f(G_\infty^{1/2}\zeta,X_\infty)\big]\bigg\}\bigg|
&=& 
o(n^{-1/2})
\eeas
as $n\to\infty$ 
for any bounded set $\calb$ in $C^{8}_b(\bbR^2)$, 
where 
$\mathfrak{S}$ is given by (\ref{202004161621}) and 
$\zeta\sim N(0,1)$ independent of $\calf$. 
\end{theorem}

\begin{remark}\rm 
Theorem \ref{0112121225} gives stable convergence of $V_n$ to $G_\infty^{1/2}\zeta$, in particular. 
Of course, it is possible to reduce the conditions for the asymptotic expansion 
only to prove the stable convergence, but it is an exercise and so details are omitted. 
\end{remark}

\subsection{Expansion of $E[f(V_n,X_\infty)]$ for a measurable function $f$}
If $(V_n,X_\infty)$ is asymptotically non-degenerate in Malliavin's sense, 
the differentiability of $f$ can be removed. 
Application of the Malliavin calculus to the non-degeneracy issue 
requires more smoothness of functionals. 
We say $U_n=O_{\bbD^\infty}(b_n)$ as $n\to\infty$ for 
a sequence of variables $U_n$ and a sequence of positive numbers $b_n$ 
if $\|U_n\|_{s,p}=O(b_n)$ as $n\to\infty$ 
for all $(s,p)\in\bbR\times(1,\infty)$. 
For a multi-dimensional differentiable functional $F=(F^i)_{i=1,...,\bar{i}}$, we write 
$\Delta_F=\det(\langle DF^{i_1},DF^{i_2}\rangle)_{i_1,i_2=1,...,\bar{i}}$. 
\bd\im[[A$^\sharp$\!\!]]
{\bf (I)} $[A]$ (\ref{ai}), (\ref{aiii}), (\ref{aiv}), (\ref{av}), (\ref{avi}), (\ref{aix}) and 
the following three conditions hold. 
\bd\im
\bd
\im[(i$^\sharp$)] For every $t,s\in[0,1]$ and $q\in\calq$, 
$a(t,q)\in\bbD^{5,\infty}$, $\dota(t,s,q)\in\bbD^{5,\infty}$, $\dot{a}(t,2)\in\bbD^{5,\infty}$ (if $2\in\calq$) and $\ddota(t,s,q)\in\bbD^{5,\infty}$, and 
\bea\label{202003251201} &&
\max_{q\in\calq}\sup_{t\in[0,1]}\big\|a(t,q)\big\|_{5,p}
+\max_{q\in\calq}\sup_{(t,s)\in[0,1]^2}\big\|\dota(t,s,q)\big\|_{5,p}
\nn\\&&
+1_{\{2\in\calq\}}\sup_{t\in[0,1]}\big\|\dot{a}(t,2)\big\|_{5,p}
+\max_{q\in\calq}\sup_{(t,s)\in[0,1]^2}\big\|\ddota(t,s,q)\big\|_{5,p}
\><\> 
\infty
\eea
for every $p>1$. 
\ed
\ed
\bd\im
\bd
\im[(\ref{avii}$^\sharp$)] 
{\colorr If $2\in\calq$, then $\sup_{t\in[0,1]}\big\|\ddot{a}(t,2)\big\|_{5,p}<\infty$ 
for every $p>1$, and } 
\bea\label{202003251205sharp}
E\bigg[ \int_{[0,1]^2}
1_{\big\{
\big|\sum_{j=1}^n1_{ I_j}(t) n^2D_{1_j}D_{1_j}a_j({\colorr2})-\ddot{a}(t,2)\big|>\ep\big\}}dsdt\bigg]
&\to&
0
\eea
as $n\to\infty$ for every $\ep>0$. 
\ed
\ed
\bd\im
\bd
\im[(\ref{aviii}$^\sharp$)] 
$G_\infty\in\bbD^{\infty}$ and 
\beas
\sup_{t_1,...,t_i\in[0,1]}\|D_{t_i}\cdots D_{t_1}G_\infty\|_p<\infty
\eeas
for every $i\in\bbN$ and $p>1$. 
Moreover, there exists a positive number $\kappa$ such that 
\bea\label{2020041446}
n^{-1}\sum_{j=1}^n\sum_{q\in\calq}q!a_j(q)^2-G_\infty
&=& 
O_{\bbD^{\infty}}(n^{-(\half+\kappa)})
\eea
as $n\to\infty$. 
\ed
\ed
\bd\im[(II)] 
$G_\infty^{-1}, \Delta_{X_\infty}^{-1}\in L^{\infty-}$. 
\begin{en-text}
\bd
\im[(a)] $G_\infty^{-1}\in L^{\infty-}$. 
\im[(b)] There exists a positive number $\kappa$ such that 
\beas 
P\big[\Delta_{(M_n,X_\infty)}<s_n\big] &=& O(n^{-(\half+\kappa)})
\eeas
for some positive random variables $s_n\in\bbD^{8,\infty}$ satisfying 
\beas 
\sup_{n\in\bbN}\big(\|s_n^{-1}\|_p+\|s_n\|_{8,p}\big) &<& \infty
\eeas 
for every $p>1$. 
\ed
\end{en-text}
\ed
\ed

Condition (\ref{2020041446}) is realistic. An example will be given in Section \ref{202004201231}. 

For a polynomial random symbol $\varsigma(\tti\sfz,\tti\sfx)=\sum_\alpha c_\alpha(\tti\sfz,\tti\sfx)^\alpha$, 
where each $c_\alpha$ is a random variable and the multi-index $\alpha$ is in $\bbZ_+^2$, 
the expected action of the adjoint  $\varsigma(\partial_z,\partial_x)^*$ 
to $\phi(z;0,G_\infty)\delta_x(X_\infty)$ is defined by 
\bea\label{202004171910}
E\bigg[\varsigma(\partial_z,\partial_x)^*\bigg\{\phi(z;0,G_\infty)\delta_x(X_\infty)\bigg\}\bigg]
&=& 
\sum_\alpha(-\partial_z,-\partial_x)^\alpha 
E\bigg[c_\alpha\phi(z;0,G_\infty)\delta_x(X_\infty)\bigg]
\nn\\&=&
\sum_\alpha(-\partial_z,-\partial_x)^\alpha 
\bigg\{E\big[c_\alpha\phi(z;0,G_\infty)\big|X_\infty=x\big]p^{X_\infty}(x)\bigg\}
\nn\\&&
\eea
Here $\delta_x(X_\infty)$ is Watanabe's delta function, that is, 
the pull-back of the delta function $\delta_x$ by $X_\infty$, 
that is a generalized Wiener functional in ${\widetilde{\bbD}}^{-\infty}=\cup_{s\in\bbR}\cap_{p>1}\bbD^{s,p}$, suppose that $X_\infty$ is non-degenerate in Malliavin's sense, as well as 
smoothness of functionals appearing. 
$p^{X_\infty}(x)$ is a good version of density function of the distribution $\call\{X_\infty\}$. 
It and its continuous derivatives exist so that the differentiations in (\ref{202004171910}) make 
the ordinary sense. 
See Section \ref{202005171900}. 

Define the (polynomial) random symbol $\mathfrak{S}_n$ by 
\beas 
\mathfrak{S}_n
&=&
1+n^{-1/2}\mathfrak{S}, 
\eeas
where the random symbol $\mathfrak{S}$ is given in (\ref{202004161621})
for Theorem \ref{202004171845} forthcoming. 
The asymptotic expansion formula is expressed by the density 
\bea\label{202004191710} 
p_n(z,x) &=& 
E\bigg[\mathfrak{S}_n(\partial_z,\partial_x)^*
\bigg\{\phi(z;0,G_\infty)\delta_x(X_\infty)\bigg\}\bigg]. 
\eea
The density function $p_n$ is well defined under Condition $[A^\sharp]$.

Let $\cale(M,\gamma)$ denote the set of measurable functions 
$f:\bbR^2\to\bbR$ satisfying 
$|f(z,x)|\leq M(1+|z|+|x|)^\gamma$ for all $(z,x)\in\bbR^2$. 
\begin{theorem}\label{202004171845}
Suppose that Condition $[A^\sharp]$ is satisfied. Then 
\beas
\sup_{f\in\cale(M,\gamma)}
\bigg| E\big[f(V_n,X_\infty)\big] 
-\int_{\bbR^2}f(z,x)p_n(z,x)dzdx\bigg|
&=& 
o(r_n)
\eeas
as $n\to\infty$ for every $(M,\gamma)\in(0,\infty)^2$. 
Here the density $p_n(z,x)$ is defined by (\ref{202004191710}) 
for the random symbol $\mathfrak{S}$ given in (\ref{202004161621}). 
\end{theorem}
\halflineskip

Proof of Theorem \ref{202004171845} is in Section \ref{202004171914}. 
In Theorems \ref{0112121225} and \ref{202004171845}, 
$N_n$ and $X_n$ were specific. 
We will generalize them in Section \ref{202004021507} after clarifying what kind conditions 
should be posed by the general theory recalled in Section \ref{202003141634}.

\section{Basic results in asymptotic expansion}

\subsection{Donsker-Watanabe's delta function}\label{202005171900}
Denote by $\cals(\bbR^\sfd)$ the space of rapidly decreasing smooth functions on $\bbR^\sfd$. 
The dual of $\cals(\bbR^\sfd)$ is the space of Schwartz distributions and denoted by $\cals'(\bbR^\sfd)$. 
Let $A=1+|x|^2-\half\Delta$, where $\Delta$ is the Laplacian on $\bbR^\sfd$. 
Define the norm $\|u\|_{2k}=\sup_{x\in\bbR\sfd}|A^ku(x)|$ for $u\in\cals(\bbR^\sfd)$ and 
$k\in\bbZ$. 
For negative $k$, $A^k$ is well defined as an integral operator. 

Denote $\hat{C}(\bbR^\sfd)$ the space of continuous functions $g$ on $\bbR^\sfd$ 
such that $\lim_{|x|\to\infty}g(x)=0$. 
The spadce $\hat{C}(\bbR^\sfd)$ is equipped with the sup-norm $\|\cdot\|_\infty$. 
For $m\in\bbZ_+$, defined $\hat{C}^{-2m}(\bbR^\sfd)$ by 
\beas 
\hat{C}^{-2m}(\bbR^\sfd) &=& 
\big\{g\in\cals'(\bbR^\sfd);\>A^{-m}g\in\hat{C}(\bbR^\sfd)\big\}. 
\eeas
For $g\in\hat{C}^{-2m}(\bbR^\sfd)$, we defined $\|g\|_{-2m}$ by 
\beas 
\|g\|_{-2m} &=& \|A^{-m}g\|_\infty. 
\eeas
By the representation theorem of $\cals'(\bbR^\sfd)$, we know 
\beas 
\cals'(\bbR^\sfd) &=& \bigcup_{m\in\bbZ_+}\hat{C}^{-2m}(\bbR^\sfd). 
\eeas
More precisely, 
every $g\in\cals'(\bbR^\sfd)$ has a representation 
\beas 
g(x) &=& (1+|x|^2)^{s/2}\partial^\alpha f(x)
\eeas
for some $\alpha\in\bbZ_+^\sfd$, $s\geq0$ and some bounded measurable function $f$, 
where the partial derivatives is in the sense of the generalized derivatives. 
Then $\partial^\nu A^{-m}\partial^\lambda g \in \hat{C}(\bbR^\sfd)$ whenever 
$2m>|\nu|+|\lambda|+|\alpha|+s$. 
Moreover, 
\beas 
\|\partial^\nu A^{-m}\partial^\lambda g\|_\infty &\leq & C \|f\|_\infty
\eeas
for some constant $C$ depending on $m$, $|\nu|$, $|\lambda|$, $|\alpha|$, $s$ and $\sfd$, 
but independent of $f$. 
In the case of the Dirac delta function $\delta_a$ on $\bbR^\sfd$, 
\beas 
\partial_x^{\bf n}\delta_a(x) &\in& \hat{C}^{-2m}(\bbR^\sfd)
\eeas
for $m\in\bbZ_P$ satisfying $m>(|\n|+\sfd)/2$. 
Therefore, 
$\partial_x^{\n}\delta_a(F)$ is well defined for $F\in\bbD^{2m+1,\infty}(\bbR^\sfd)$ 
and $m>(|\n|+\sfd)/2$ if $\Delta_F^{-1}\in L^\inftym$. 
The composite functional $\partial_a^{\n}\delta_a(F)$ requires the same regularity. 
In this paper, we use $\delta_x(X_\infty)$ and its derivatives of limited order. 
For example, the differentiability requested to $X_\infty$ in Condition $[C]$ 
in Section \ref{202003141634} is sufficient for the use. 
\begin{en-text}
\koko
For $k\in\bbZ$, we define $C_{2k}$ as the completion of $\cals(\bbR^\sfd)$ with respect to the norm $\|\cdot\|_{2k}$. 
Then by the representation theorem of $\cals'(\bbR^\sfd)$, 
it is known that 
\beas 
\cals'(\bbR^\sfd) &=& \bigcup_{k\in\bbZ}C_{2k}. 
\eeas
Suppose that $F=(F^i)_{i=1}^\sfd$ 
\end{en-text}
We refer the reader to 
Watanbe \cite{watanabe1984lectures} and Ikeda and Watanabe \cite{IkedaWatanabe1989} for 
details of the notion of the generalized Wiener functionals.

\subsection{Asymptotic expansion of Skorohod integrals}\label{202003141634}
This section overviews the machinery given by Nualart and Yoshida \cite{nualart2019asymptotic} 
for asymptotic expansion of Skorohod integrals, 
and extends their results {\sred by a perturbation method}. 
Since there are many other applications than this paper gives 
and the {\sred extension is} general, in this section, we keep the multi-dimensional setting 
in \cite{nualart2019asymptotic}. 
In this section, we do not need the structure of the Hilbert space for a standard Brownian motion. 
So we will work on a probability space $(\Omega,\calf,P)$ equipped with an isonormal Gaussian process 
$\bbW=\{W(h)\}_{h\in\mH}$ on a real separable Hilbert space $\mH$ having an inner product 
$\langle\cdot,\cdot\rangle$ and the norm $|\cdot|_\mH=\langle\cdot,\cdot\rangle^{1/2}$. 

We will consider a sequence of $\sfd$-dimensional random variables $Z_n$ having a decomposition
\bea\label{202003291147}
Z_n &=& M_n+r_nN_n
\eea
where 
$M_n=\delta(u_n)=\big(\delta(u_n^i)\big)_{i=1,...,\sfd}$ is the Skorohod integral 
for {\sred a} 
generic $u_n=(u_n^i)_{i=1,...,\sfd}\in{\mathcal Dom}(\delta)^\sfd$, 
$N_n$ is a $\sfd$-dimensional random variable and 
$(r_n)_{n\in\bbN}$ is a sequence of positive numbers such that $\lim_{n\to\infty}r_n=0$. 
{\sred Later,} we will set $r_n=n^{-1/2}$ as in (\ref{202003130001}), 
when using the asymptotic expansion for Skorohod integrals to the applications of this paper. 
%
A reference variable will be denoted by a $\sfd_1$-dimensional random variable $X_n$. 
Our interest is in deriving asymptotic expansion for the joint distribution $\call\{(Z_n,X_n)\}$. 
To be considered here is the situation where $M_n\to^dG_\infty^{1/2}\zeta$ 
as $n\to\infty$, where $G_\infty$ is a $\sfd\times\sfd$ positive symmetric random matrix 
and $\zeta$ is a $\sfd$-dimensional standard normal random vector independent of $G_\infty$ and possibly defined on an extension of $(\Omega,\calf,P)$. 
See Nourdin, Nualart and Pecati \cite{nourdin2016quantitative} for 
the stable limit theorems on the Wiener space. 

The $r$ times tensor product of a vector $v$ is denoted by $v^{\otimes r}$. 
For a tensor $T=(T_{i_1,..,i_k})_{i_1,...,i_k}$ and vectors 
$v_1=(v^{i_1}_1)_{i_1}$,...,$v_k=(v_k^{i_k})_{i_k}$, we write 
\beas 
T[v_1,...,v_k] 
&=& 
T[v_1\otimes\cdots \otimes v_k]
\yeq
\sum_{i_1,...,i_k}T_{i_1,...,i_k}v^{i_1}_1\cdots v_k^{i_k}. 
\eeas
The brackets $[\ ]$ stands for a multilinear mapping. 
For $\mH$-valued tensors $S=(S_i)$ and $T=(T_j)$, 
the tensor with components $(\langle S_i,T_j\rangle)_{i,j}$ is denoted by $\langle S,T\rangle$. 
This rule is {\sred also} applied to {\sred tensor-valued} 
tensors. 
The Malliavin derivative of a matrix-valued functional $F=(F^{ij})$ is denoted by $DF$.

Nualart and Yoshida \cite{nualart2019asymptotic} defined three random symbols analogous to the adaptive random symbol (tangent) and the anticipative random symbol (torsion) introduced in Yoshida \cite{yoshida2013martingale}, though these two sets of random symbols are not corresponding one-to-one. 
The {\bf quasi-tangent}  (q-tangent) is defined by 
\begin{en-text}
\beas 
{\sf qTan}(\tti\sfz)^2
&=& 
r_n^{-1}\big(D_{u_n\tti\sfz}M_n\tti\sfz-G_\infty(\tti\sfz)^2\big). 
\eeas
\end{en-text}
\beas 
{\sf qTan}[{\cred\tti\sfz_1,\tti\sfz_2}]&=&
r_n^{-1}\bigg(\big\langle DM_n[\tti{\cred\sfz_1}],u_n[\tti\sfz_2]\big\rangle
-G_\infty[\tti\sfz_1,\tti\sfz_2]\bigg)
\eeas
for $(\tti\sfz_1,\tti\sfz_2)\in
(\tti\bbR^{\sfd})^2$. 
The {\bf quasi-torsion} (q-torsion) and {\bf modified quasi-torsion} (modified q-torsion) are defined by 
\begin{en-text}
\beas
{\sf qTor}(\tti\sfz)^3
&=& 
r_n^{-1}D_{u_n\tti\sfz}D_{u_n\tti\sfz}M_n\tti\sfz
\eeas
\end{en-text}
\beas
{\sf qTor}[{\cred\tti\sfz_1,\tti\sfz_2,\tti\sfz_3}]
&=&
r_n^{-1}\bigg\langle D\big\langle DM_n[\tti{\cred\sfz_1}],u_n[\tti{\cred\sfz_2}]\big\rangle_\mfh,u_n[\tti{\cred\sfz_3}]\bigg\rangle
\eeas
and 
\beas
{\sf mqTor}[{\cred\tti\sfz_1,\tti\sfz_2,\tti\sfz_3}]&=&
r_n^{-1}\big\langle DG_\infty[{\cred\tti\sfz_1,\tti\sfz_2}],u_n[\tti{\cred\sfz_3}]\big\rangle
\eeas
for $(\tti\sfz_1,\tti\sfz_2,\tti\sfz_3)\in(\tti\bbR^{\sfd})^3$, respectively. 
\begin{en-text}
\beas 
{\sf mqTor}(\tti\sfz)^3
&=&
r_n^{-1}D_{u_n\tti\sfz}G_\infty(\tti\sfz)^2,
\eeas
respectively. 
\end{en-text}
These random symbols are well defined under suitable differentiability of variables. 
They are related by the formula
\beas 
\big\langle D\>{\sf qTan}[(\tti\sfz)^{\otimes2}],u_n[\tti\sfz]\big\rangle_\mfh 
&=& 
{\sf qTor}[(\tti\sfz)^{\otimes3}]-{\sf mqTor}[(\tti\sfz)^{\otimes3}], 
\eeas
or symbolically, 
\beas 
D_{u_n}{\sf qTan} &=& {\sf qTor}-{\sf mqTor}. 
\eeas
{\sred We could} replace $D_{{\sred u_n}}M_n=\big\langle DM_n,u_n\big\rangle$ by its symmetric version 
\beas 
\half\big(
\big\langle DM_n,u_n\big\rangle+\big\langle u_n, DM_n\big\rangle\big)
\eeas
{\sred when defining ${\sf qTan}$, but we go with the above definition.} 

When $(Z_n,X_n)\to^d(G_\infty^{1/2}\zeta,X_\infty)$ as $n\to\infty$, 
we have 
\bea\label{202013291441}
E\big[f(Z_n,X_n)\big] - E\big[f(G_\infty^{1/2}\zeta,X_\infty)\big] &=& o(1)
\eea
as $n\to\infty$ for any bounded continuous function 
$f:\bbR^{{\sred\check{\sfd}}}\to\bbR$, {\sred$\check{\sfd}=\sfd+\sfd_1$,} 
by definition of the convergence 
in distribution. 
Asymptotic expansion gives more precise approximation to $E\big[f(Z_n,X_n)\big]$ than 
(\ref{202013291441}) as 
\bea\label{202013291445}
E\big[f(Z_n,X_n)\big] - E\big[f(G_\infty^{1/2}\zeta,X_\infty)\big]
-r_n\bbL(f) &=& o(r_n), 
\eea
where $\bbL(f)$ is a linear functional of $f$. 
The approximation (\ref{202013291445}) holds for smooth function $f$ in general. 
However, if the distribution $\call\{(Z_n,X_n)\}$ is regular to some extent, then 
(\ref{202013291445}) is valid for non-differentiable functions $f$, and moreover 
uniformly in $f$ in a certain class of measurable functions. 
Then the notion of {\it local density} is necessary and a truncation technique will be required 
to carry out the validation. 
For this purpose, we introduce a random variable $\psi_n:\Omega\to[0,1]$ depending on $n\in\bbN$. 

Write $\dotx_n=r_n^{-1}(X_n-X_\infty)$, 
\beas 
G^{(2)}_n &=& r_n{\sf qTan}\yeq D_{u_n}M_n-G_\infty
\eeas
and 
\beas 
G^{(3)}_n &=& r_n{\sf mqTor}\yeq D_{u_n}G_\infty.
\eeas

We write 
\beas 
D_{u_n}^k=\overbrace{D_{u_n}\cdots D_{u_n}}^{k}. 
\eeas
Define the following random symbols: 
\beas 
\mathfrak{S}^{(3,0)}_n(\tti\sfz) 
&=& 
\mathfrak{S}^{(3,0)}_n(\tti\sfz,\tti\sfx) 
\yeq
\frac{1}{3}r_n^{-1}D_{u_n{\sred[\tti\sfz]}}^2M_n{\sred[\tti\sfz]}
\yeq
\frac{1}{3}{\sf qTor}\>[(\tti\sfz)^{\otimes3}],
\\
\mathfrak{S}^{(2,0)}_{0,n}(\tti\sfz) 
&=&
\mathfrak{S}^{(2,0)}_{0,n}(\tti\sfz,\tti\sfx) 
\yeq 
\half r_n^{-1}G_n^{(2)}(\tti\sfz)
\\&=&
\half r_n^{-1}\big(D_{u_n[\tti\sfz]}M_n[\tti\sfz]-G_\infty[(\tti\sfz)^{\otimes2}]\big)
\yeq \half{\sf qTan}\>[(\tti\sfz)^{\otimes2}],
\\
\mathfrak{S}^{(1,1)}_n(\tti\sfz,\tti\sfx) 
&=&
r_n^{-1}D_{u_n[\tti\sfz]}X_\infty[\tti\sfx],
\\
\mathfrak{S}^{(1,0)}_n(\tti\sfz) 
&=&
\mathfrak{S}^{(1,0)}_n(\tti\sfz,\tti\sfx) 
\yeq N_n[\tti\sfz],
\\
\mathfrak{S}^{(0,1)}_n(\tti\sfx)
&=&
\mathfrak{S}^{(0,1)}_n(\tti\sfz,\tti\sfx)
\yeq 
\dotx_n[\tti\sfx],
\\
\mathfrak{S}^{(2,0)}_{1,n}(\tti\sfz)
&=&
\mathfrak{S}^{(2,0)}_{1,n}(\tti\sfz,\tti\sfx)
\yeq 
D_{u_n[\tti\sfz]}N_n[\tti\sfz],
\\
\mathfrak{S}^{(1,1)}_{1,n}(\tti\sfz,\tti\sfx)
&=&
D_{u_n[\tti\sfz]}\dotx_n[\tti\sfx]
\eeas
for $\sfz\in\bbR^{{\sred\sfd}}$ and $\sfx\in\bbR^{{\sred\sfd_1}}$. 

Let 
\beas 
\Psi(\sfz,\sfx) &=& \exp\big(2^{-1}G_\infty{\sred [(\tti\sfz)^{\otimes2}]}
+X_\infty{\sred [\tti\sfx]}\big)
\eeas
for $\sfz\in\bbR^{{\sred \sfd}}$ and $\sfx\in\bbR^{{\sred \sfd_1}}$. 
Suppose that there are some (polynomial) random symbols 
$\mathfrak{S}^{(3,0)}(\tti\sfz,\tti\sfx)$, 
$\mathfrak{S}^{(2,0)}_{0}(\tti\sfz,\tti\sfx)$, 
$\mathfrak{S}^{(1,1)}(\tti\sfz,\tti\sfx)$, 
$\mathfrak{S}^{(1,0)}(\tti\sfz,\tti\sfx)$, 
$\mathfrak{S}^{(0,1)}(\tti\sfz,\tti\sfx)$, 
$\mathfrak{S}^{(2,0)}_{1}(\tti\sfz,\tti\sfx)$
and 
$\mathfrak{S}^{(1,1)}_{1}(\tti\sfz,\tti\sfx)$. 
Condition $[B]$ in Section 3.2 of Nualart and Yoshida \cite{nualart2019asymptotic} 
is rephrased as follows. 
\bd\im[[B\!\!]]{\bf (i)} $u_n\in\bbD^{4,p}(\mH{\sred \otimes\bbR^\sfd})$, $G_\infty\in\bbD^{3,p}{\sred (\bbR^\sfd\otimes\bbR^\sfd)}$, 
$N_n{\sred \in\bbD^{3,p}(\bbR^\sfd)}$, $X_n,X_\infty\in\bbD^{3,p}{\sred (\bbR^{\sfd_1})}$ and $\psi_n\in\bbD^{2,p_1}$ for some $p$ and $p_1$ satisfying 
$5p^{-1}+p_1^{-1}\leq1$. 
\bd
\im[(ii)] As $n\to\infty$, 
\bea\label{202003301101}
\|u_n\|_{1,p} &=& O(1)
\eea
\bea\label{202003301102}
\sum_{k=2,3}\|G^{(k)}_n\|_{p/2} &=& O(r_n)
\eea
\bea\label{202003301103}
\|D_{u_n}G^{(2)}_n\|_{p/3} &=& O(r_n)
\eea
\bea\label{202003301104}
\|D_{u_n}G^{(3)}_n\|_{p/3} &=& o(r_n)
\eea
\bea\label{202003301105}
\sum_{k=2,3}\|D_{u_n}^2G^{(k)}_n\|_{p/4} &=& o(r_n)
\eea
\bea\label{202003301106}
\|D_{u_n}X_\infty\|_p &=& O(r_n)
\eea
\bea\label{202003301107}
\|D_{u_n}^2X_\infty\|_{\sred p/3} &=& o(r_n)
\eea
\bea\label{202003301108}
\|D_{u_n}^3X_\infty\|_{\sred p/4} &=& o(r_n)
\eea
\bea\label{202003301109}
\|N_n\|_{3,p}+\|\dotx_n\|_{3,p}&=& O(1)
\eea
%
\bea\label{202003301110}
\|D_{u_n}^2N_n\|_{p/3}+\|D_{u_n}^2\dotx_n\|_{p/3} &=& o(1)
\eea
\bea\label{202003301111}
\|D_{u_n}^3N_n\|_{{\sred p/4}}+\|D_{u_n}^3\dotx_n\|_{{\sred p/4}} &=& o(1)
\eea
%
\bea\label{202003301112}
\|1-\psi_n\|_{2,p_1} &=& o(r_n). 
\eea
\ed
\im[(iii)] 
For every $\sfz\in\bbR$ and $\sfx\in\bbR$, 
\bea\label{202004220244}
\lim_{n\to\infty} E\big[\Psi(\sfz,\sfx)\mathfrak{T}_n(\tti\sfz,\tti\sfx)\psi_n\big] 
&=& E\big[\Psi(\sfz,\sfx)\mathfrak{T}(\tti\sfz,\tti\sfx)\big]
\eea
for the pairs of random symbols $(\mathfrak{T}_n,\mathfrak{T})$$=$
$(\mathfrak{S}^{(3,0)}_n,\mathfrak{S}^{(3,0)})$, 
$(\mathfrak{S}^{(2,0)}_{0,n},\mathfrak{S}^{(2,0)}_0)$, 
$(\mathfrak{S}^{(1,1)}_n,\mathfrak{S}^{(1,1)})$, 
$(\mathfrak{S}^{(1,0)}_n,\mathfrak{S}^{(1,0)})$, 
$(\mathfrak{S}^{(0,1)}_n,\mathfrak{S}^{(0,1)})$, 
$(\mathfrak{S}^{(2,0)}_{1,n},\mathfrak{S}^{(2,0)}_1)$, 
and 
$(\mathfrak{S}^{(1,1)}_{1,n},\mathfrak{S}^{(1,1)}_1)$. 
\ed
\onelineskip

Here we are assuming that each $\mathfrak{T}$ is a random polynomial of $(\tti\sfz,\tti\sfx)$ 
with integrable coefficients. 
The convergence in (iii) of Condition $[B]$ does not necessarily mean 
the $L^1$-convergence of the random symbol $\mathfrak{T}_n$ to $\mathfrak{T}$. 
It only requires ``weak'' convergence in that the projection of $\mathfrak{T}_n$ onto $\Psi$
converges to that of $\mathfrak{T}$ onto $\Psi$. 
Indeed, the degree of the random polynomial $\mathfrak{T}$ may differ from that of $\mathfrak{T}_n$. 

\begin{en-text}
{\colorr We can extend the results in \cite{nualart2019asymptotic} 
by removing their conditions (3.17) and (3.18), i.e., 
\bea\label{202003301110}
\|D_{u_n}^2N_n\|_{p/3}+\|D_{u_n}^2\dotx_n\|_{p/3} &=& o(1)
\eea
\bea\label{202003301111}
\|D_{u_n}^3N_n\|_{p/4}+\|D_{u_n}^3\dotx_n\|_{p/4} &=& o(1)
\eea
in the present case of their $W_n=0$. 
This extension can be done by using a perturbation method in Sakamoto and Yoshida \cite{SakamotoYoshida2003}. 
We will discuss it in Section \koko
}
\end{en-text}
{\colorr
In some cases, $N_n$ admits a decomposition $N_n=N_n'+N_n''$ such that 
Conditions (\ref{202003301110}) and (\ref{202003301111}) can be verified easily for $N_n'$. 
For $N_n''$, we can try the perturbation method with the help of a limit theorem. 
The perturbation method was used in Yoshida \cite{Yoshida1997, yoshida2001malliavin, yoshida2013martingale} 
and 
Sakamoto and Yoshida \cite{SakamotoYoshida2003}. 
}

The random symbol $\mathfrak{S}(\tti\sfz,\tti\sfx)$ is defined by 
\bea\label{202004161620}
\mathfrak{S}(\tti\sfz,\tti\sfx)
&=&
\mathfrak{S}^{(3,0)}(\tti\sfz,\tti\sfx) 
+\mathfrak{S}^{(2,0)}_{0}(\tti\sfz,\tti\sfx) 
+\mathfrak{S}^{(1,1)}(\tti\sfz,\tti\sfx) 
+\mathfrak{S}^{(1,0)}(\tti\sfz,\tti\sfx) 
\nn\\&&
+\mathfrak{S}^{(0,1)}(\tti\sfz,\tti\sfx)
+\mathfrak{S}^{(2,0)}_{1}(\tti\sfz,\tti\sfx)
+\mathfrak{S}^{(1,1)}_{1}(\tti\sfz,\tti\sfx). 
\eea
Denote by $\beta_0$ the degree in $(\sfz,\sfx)$ of the random symbol $\mathfrak{S}$. 
Let $\beta=\max\{7,\beta_0\}$. 
Now Theorem 3.3 of Nualart and Yoshida \cite{nualart2019asymptotic} is rephrased 
by the following theorem. 
\begin{theorem}\label{202004011550}
Suppose that Condition $[B]$ is satisfied. 
Then 
\beas 
\sup_{f\in \calb}\bigg|
E\big[f(Z_n,X_n)\big] - 
\bigg\{E\big[f(G_\infty^{1/2}\zeta,X_\infty)\big]
+r_nE\big[\mathfrak{S}(\partial_z,\partial_x)f(G_\infty^{1/2}\zeta,X_\infty)\big]\bigg\}\bigg|
&=& 
o(r_n)
\eeas
as $n\to\infty$ 
for any bounded set $\calb$ in $C^{\beta+1}_b(\bbR^{\sred\check{\sfd}})$, 
where $\zeta\sim N_{\sred\sfd}(0,{\sred I_\sfd})$ independent of $\calf$. 
\end{theorem}

The role of $\psi_n$ is implicit in Theorem \ref{202004011550}. 
Simplification by setting $\psi_n=1$ is obvious if the other conditions are satisfied. 
However, the truncation functional $\psi_n$ is indispensable when 
asymptotic expansion of $E[f(Z_n,X_n)]$ for measurable functions $f$, as discussed below.

Asymptotic expansion of $E[f(Z_n,X_n)]$ is possible for measurable functions 
if the law $\call\{(Z_n,X_n)\}$ is sufficiently regular. 
In the present context, it is natural to asses the regularity 
by the local non-degeneracy of the Malliavin covariance of the variables in question. 
Here we recall a set of conditions among others considered in 
Nualart and Yoshida \cite{nualart2019asymptotic}. 
It involves a parameter $\ell{\sred =\check{\sfd}+8}$ concerning differentiability. 
We denote by $\beta_x$ the maximum degree in $x$ of $\mathfrak{S}$. 
Let $d_2={\sred(\ell+\beta_x-7)\vee\big(2\lfloor(\sfd_1+2)/2\rfloor}
+2\lfloor(\beta_x+1)/2\rfloor\big)$. 
Let $\bbD^{s,\infty}=\cap_{p>1}\bbD^{s,p}$. 
For a multi-dimensional functional $F$, we write $\Delta_F=\det\sigma_{F}$, 
$\sigma_F$ being the Malliavin covariance matrix of $F$. 
The following condition rephrases Condition $[C^\natural]$ of Section 7 
of Nualart and Yoshida \cite{nualart2019asymptotic}. 

\bd
\im[[C\!\!]] \hspace{7pt}{\bf (i)} 
$u_n\in{\sred\bbD^{\ell+1,\infty}(\mH\otimes\bbR^\sfd)}$, 
$G_\infty\in{\sred\bbD^{(\ell+1)\vee d_2,\infty}(\bbR^\sfd\otimes\bbR^\sfd)}$, 
$N_n\in{\sred\bbD^{\ell,\infty}(\bbR^\sfd)}$, 
$X_n\in{\sred\bbD^{\ell,\infty}(\bbR^{\sfd_1})}$, 
$X_\infty\in{\sred\bbD^{\ell\vee(d_2+1),\infty}(\bbR^{\sfd_1})}$. 
\bd
\im[(ii)] There exists some positive number $\kappa$ such that 
the following estimates hold for every $p>1$: 
\bea\label{202004011421}
\|u_n\|_{{\sred\ell},p} &=& O(1)
\eea
\bea\label{202004011422}
\|G^{(2)}_n\|_{{\sred\ell-2},p} &=& O(r_n)
\eea
\bea\label{202004011423}
\|G^{(3)}_n\|_{{\sred\ell-2},p} &=& O(r_n)
\eea
\bea\label{202004011424}
\|D_{u_n}G^{(3)}_n\|_{{\sred\ell-1},p} &=& O(r_n^{1+\kappa})
\eea
\bea\label{202004011425}
\|D_{u_n}^2G^{(2)}_n\|_{{\sred\ell-3},p} &=& O(r_n^{1+\kappa})
\eea
\bea\label{202004011426}
\|D_{u_n}X_\infty\|_{{\sred\ell-3},p} &=& O(r_n)
\eea
\bea\label{202004011427}
\|D_{u_n}^2X_\infty\|_{{\sred\ell-2},p} &=& O(r_n^{1+\kappa})
\eea
\bea\label{202004011428}
\|N_n\|_{{\sred\ell-1},p}+\|\dotx_n\|_{{\sred\ell-1},p} &=& O(1)
\eea
\bea\label{202004011429}
\|D_{u_n}^2N_n\|_{{\sred\ell-2},p}+\|D_{u_n}^2\dotx_n\|_{{\sred\ell-2},p} &=& O(r_n^{{\colorr \kappa}}). 
\eea
%
\im[(iii)] 
For each pair of random symbols $(\mathfrak{T}_n,\mathfrak{T})$$=$
$(\mathfrak{S}^{(3,0)}_n,\mathfrak{S}^{(3,0)})$, 
$(\mathfrak{S}^{(2,0)}_{0,n},\mathfrak{S}^{(2,0)}_0)$, 
$(\mathfrak{S}^{(1,1)}_n,\mathfrak{S}^{(1,1)})$, 
$(\mathfrak{S}^{(1,0)}_n,\mathfrak{S}^{(1,0)})$, 
$(\mathfrak{S}^{(0,1)}_n,\mathfrak{S}^{(0,1)})$, 
$(\mathfrak{S}^{(2,0)}_{1,n},\mathfrak{S}^{(2,0)_1})$, 
and 
$(\mathfrak{S}^{(1,1)}_{1,n},\mathfrak{S}^{(1,1)}_1)$, 
the following conditions are fulfilled: 
\bd
\im[(a)] $\mathfrak{T}$ is a (polynomial) random symbol the coefficients of which are in 
$\bbD^{{\sred\check{\sfd}}+\beta_x{\sred+1},1+}=\cup_{p>1}\bbD^{{\sred\check{\sfd}}+\beta_x{\sred+1},p}$. 
\im[(b)] For some $p>1$, there exists a (polynomial) random symbol $\overline{\mathfrak{T}}_n$ 
that has $L^p$ coefficients and the same degree as $\mathfrak{T}$, 
such that 
\beas 
E\big[\Psi(\sfz,\sfx)\mathfrak{T}_n(\tti\sfz,\tti\sfx)\big] 
&=& 
E\big[\Psi(\sfz,\sfx)\overline{\mathfrak{T}}_n(\tti\sfz,\tti\sfx)\big] 
\eeas
and 
$\overline{\mathfrak{T}}_n\to\mathfrak{T}$ in $L^p$ as $n\to\infty$. 
\ed

\im[(iv)] 
\bd
\im[(a)] $G_\infty^{-1}\in L^{\infty-}$. 
\im[(b)] There exists a positive number $\kappa$ such that 
\beas 
P\big[\Delta_{(M_n,X_\infty)}<s_n\big] &=& O(r_n^{1+\kappa})
\eeas
for some positive random variables $s_n\in\bbD^{{\sred\ell-2},\infty}$ satisfying 
\beas 
\sup_{n\in\bbN}\big(\|s_n^{-1}\|_p+\|s_n\|_{{\sred\ell-2},p}\big) &<& \infty
\eeas 
for every $p>1$. 
\ed
\ed
\ed
\halflineskip

{\sred In the case $\sfd=\sfd_1=1$, 
the indices appearing in $[C]$ are }
\beas {\sred 
\ell=10\quad\text{and}\quad
d_2=(3+\beta_x)\vee\big(2+2\lfloor(\beta_x+1)/2\rfloor\big). 
}
\eeas
\noindent
In Condition $[C]$, the convergence 
$\overline{\mathfrak{T}}_n\to\overline{\mathfrak{T}}$ in $L^p$ means 
each coefficient of $\overline{\mathfrak{T}}_n$ converges to the corresponding coefficient of $\overline{\mathfrak{T}}$ in $L^p$. 
The functional $s_n$ in Condition $[C]$ (iv) (b)
makes a suitable truncation functional $\psi_n$, 
though we just mention it here; 
see Nualart and Yoshida \cite{nualart2019asymptotic} for details of a construction of $\psi_n$. 
The truncation method is essential because $(M_n,X_\infty)$ is not necessarily 
uniformly non-degenerate, but locally asymptotically non-degenerate. 
In statistics, variables are not always defined as a smooth Wiener functional, 
and then there is no way the uniform non-degeneracy can be obtained. 
Condition $[C]$ is much stronger than Condition $[B]$ {\sred even} with $\psi_n=1$. 

\begin{en-text}
{\colorr 
Here we extended Theorem 7.6 by removing the condition 
\bea\label{202004011429}
\|D_{u_n}^2N_n\|_{{\colorr8},p}+\|D_{u_n}^2\dotx_n\|_{8,p} &=& O(r_n^{{\colorr \kappa}}). 
\eea
See Section \koko for details. 
}
\end{en-text}

Define the (polynomial) random symbol $\mathfrak{S}_n$ by 
\bea\label{202004210745}
\mathfrak{S}_n
&=&
1+r_n\mathfrak{S},
\eea
where $\mathfrak{S}$ is the general random symbol defined by (\ref{202004161620}). 
In this general setting, 
the asymptotic expansion formula is given by the density 
\bea\label{202004211625}
p_n(z,x) &=& 
E\bigg[\mathfrak{S}_n(\partial_z,\partial_x)^*
\bigg\{\phi(z;0,G_\infty)\delta_x(X_\infty)\bigg\}\bigg]. 
\eea
The density function $p_n$ is well defined under Condition $[C]$.

The following theorem rephrases Theorem 7.6 of of Nualart and Yoshida \cite{nualart2019asymptotic}. 
{\sred 
In the present context, we denote by $\cale(M,\gamma)$ 
the set of measurable functions 
$f:\bbR^{\check{\sfd}}\to\bbR$ satisfying 
$|f(z,s)|\leq M(1+|z|+|x|)^\gamma$ for all $(z,x)\in\bbR^2$. 
}
\begin{theorem}\label{202004011551}
Suppose that Condition $[C]$ is satisfied. Then 
\beas
\sup_{f\in\cale(M,\gamma)}
\bigg| E\big[f(Z_n,X_n)\big] 
-\int_{{\sred\bbR^{\check{\sfd}}}}f(z,x)p_n(z,x)dzdx\bigg|
&=& 
o(r_n)
\eeas
as $n\to\infty$ for every $(M,\gamma)\in(0,\infty)^2$. 
\end{theorem}
\halflineskip

Proofs of Theorems \ref{202004011551}and \ref{202004011551} are in 
Nualart and Yoshida \cite{nualart2019asymptotic}.

\subsection{Perturbation method} 
We recall the perturbation method used in 
Yoshida \cite{Yoshida1997, yoshida2001malliavin, yoshida2013martingale} 
and in Sakamoto and Yoshida \cite{SakamotoYoshida2003}. 
Consider a $\sfd$-dimensional random variable $\bbS_T$ that has the decomposition
\beas 
\bbS_T &=& \bbX_T+r_T\bbY_T
\eeas
where $\bbX_T$ and $\bbY_T$ are $\sfd$-dimensional random variables 
{\sred for $T>0$}, and 
$r_T$ is a positive number such that $\lim_{T\to\infty}r_T=0$. 
Let $m\in\bbN$. 
Suppose that positive numbers $M$, $\gamma$, $\ell_1$ and $\ell_2$, and 
a functional $\xi_T$ 
satisfy the following conditions:
\bd
\im[(i)] 
$\ds 
m\>>\>\half\gamma+1,\quad 
\ell_1\ygeq\ell_2-1,\quad
\ell_2\ygeq(2m+1)\vee(\sfd+3).
$
\im[(ii)] 
$\ds 
\sup_T\|\bbX_T\|_{\ell_2,p}+\sup_T\|\bbY_T\|_{\ell_2,p}\><\>\infty
$
for every $p>1$. 
\im[(iii)]
$\ds (\bbX_T,\bbY_T) \>\to^d\> (\bbX_\infty,\bbY_\infty)$
as $T\to\infty$ for some random variables $\bbX_\infty$ and $\bbY_\infty$. 
\im[(iv)] 
$\ds \sup_T\|\xi_T\|_{\ell_1,p}\><\>\infty$ for every $p>1$. 
\im[(v)] 
$\ds P[|\xi_T|>1/2] \yeq O(r_T^\alpha)$ for some $\alpha>1$. 
\im[(vi)] 
$\ds\sup_TE[1_{\{|\xi_T|<1\}}\Delta_{\bbX_T}^{-p}]\><\>\infty$ 
for every $p>1$. 
\im[(vii)] 
There exists a signed measure $\Psi_T$ on $\bbR^\sfd$ such that 
\beas 
\sup_{f\in\cale(M,\gamma)}\bigg|E\big[f(\bbX_T)\big]
-\int_{\bbR^\sfd}f(z){\sred\Psi_T}(dz)\bigg|
&=& 
o(r_T)
\eeas
as $n\to\infty$, where $\cale(M,\gamma)$ is the set of measurable functions $f:\bbR^\sfd\to\bbR$ 
such that $|f(z)|\leq M(1+|x|)^\gamma$. 
\ed
Under these conditions, Sakamoto and Yoshida \cite{SakamotoYoshida2003} showed the following result. 
\begin{theorem}\label{202005171835}
It holds that 
\beas
\sup_{f\in\cale(M,\gamma)}\bigg|E\big[f(\bbS_T)\big] 
-\int_{\bbR^\sfd}f(z)\big\{{\sred\Psi_T}(dz)+{\sred r_T}g_\infty(z)dz\big\}\bigg|
&=& 
o(r_T)
\eeas
as ${\sred T}\to\infty$, where 
\beas 
g_\infty(z) &=& -\partial_z\cdot\big\{E\big[\bbY_\infty|\>\bbX_\infty=z\big]p^{\bbX_\infty}(z)\big\}.
\eeas
\end{theorem}
Under the assumptions, there exist a differentiable density of $\bbX_\infty$ 
and a version of the conditional expectation appearing above, and 
also the measure $|g_\infty(z)|dz$ has any order of moment.

\section{Generalization
}\label{202004021507}
We will consider the variable $Z_n$ having the expression (\ref{202003130001}). 
The functional $M_n$ is defined by (\ref{0112121101}) with the kernel $u_n$  
given by (\ref{202003130046}) and (\ref{202003130047}). 
The principal part $M_n$ of $Z_n$ is the same as before, 
however, 
\underline{the variables $N_n$ and $X_n$ {\sred ($n\in\bbN$)} are general} in this section. 
We do not assume $N_n$ is specified by (\ref{202004151750}). 
Such a generalization is motivated by asymptotic expansion of a parametric estimator 
for a stochastic differential equation, 
for example, where $N_n$ involves derivatives (with respect to the statistical parameter) 
of the quasi-log likelihood function of the model and is 
utterly different from and more complicated than (\ref{202004151750}). 
{\sred 
A multi-dimensional extension is straightforward but we will only consider 
one-dimensional variables for notational simplicity. 
}

A sufficient condition for asymptotic expansion of the joint distribution $\call\{(Z_n,X_n)\}$ 
will be in part similar to the one presented in the general theory in Section \ref{202003141634} 
since we do not assume any specific form of $N_n$ or $X_n$. 
We keep $G_\infty$ defined by (\ref{202004201141}) as well as 
measurable random fields 
$(a(t,q))_{t\in[0,1]}$, $(\dota(t,s,q))_{(t,s)\in[0,1]^2}$, 
$(\dot{a}(t,2))_{t\in[0,1]}$ (if $2\in\calq$) 
and 
$(\ddota(t,s,q))_{(t,s)\in[0,1]^2}$ 
$(q\in\calq)$. 
{\sred 
Let $X_\infty$ be a one-dimensional random variable in $\bbD^\infty$. According to the theory in Section \ref{202003141634}, 
$\Psi(\sfz,\sfx)$ is defined by 
\beas 
\Psi(\sfz,\sfx) &=& \exp\big(2^{-1}G_\infty(\tti z)^2+X_\infty\tti x\big)
\eeas
for $\sfz\in\bbR$ and $\sfx\in\bbR$, in the present situation.}
Denote $\dotx_n=n^{1/2}(X_n-X_\infty)$. 
Following the theory, let 
\beas 
{\mathfrak S}_n^{(1,0)}(\tti\sfz) &=& {\mathfrak S}_n^{(1,0)}(\tti\sfz,\tti\sfx) \yeq N_n\tti\sfz,
\eeas
\beas 
{\mathfrak S}_n^{(0,1)}(\tti\sfx) &=& {\mathfrak S}_n^{(0,1)}(\tti\sfz,\tti\sfx) \yeq \dotx_n\tti\sfx,
\eeas
\beas 
{\mathfrak S}_{1,n}^{(2,0)}(\tti\sfz)
&=& 
{\mathfrak S}_{1,n}^{(2,0)}(\tti\sfz,\tti\sfx) 
\yeq 
D_{u_n}N_n(\tti\sfz)^2.
\eeas
and let 
\beas 
{\mathfrak S}_{1,n}^{(1,1)}(\tti\sfz,\tti\sfx) 
&=& 
D_{u_n}\dotx_n(\tti\sfz)(\tti\sfx).
\eeas
We suppose that (polynomial) random symbols $\mathfrak{S}^{(1,0)}(\tti\sfz,\tti\sfx)$, 
$\mathfrak{S}^{(0,1)}(\tti\sfz,\tti\sfx)$, $\mathfrak{S}^{(2,0)}_1(\tti\sfz,\tti\sfx)$ 
and $\mathfrak{S}^{(1,1)}_1(\tti\sfz,\tti\sfx)$ are given. 
Condition (\ref{202004111810}) is assumed in this section.

\subsection{Asymptotic expansion of $E[f(Z_n,X_n)]$ for a differentiable function $f$}\label{202004211610}
We derive asymptotic expansion of the expectation $E[f(Z_n,X_n)]$ for differentiable functions $f$. 
\bd
\im[[D\!\!]] 
{\bf(I)} Condition $[A]$ is satisfied except for $[A]$ (\ref{avii}). 
\bd
\im[(II)] $N_n, X_n\in\bbD^\infty$ and 
the properties (\ref{202003301109}), (\ref{202003301110}) and (\ref{202003301111}) hold for every $p>1$. 
\im[(III)] 
For every $\sfz\in\bbR$ and $\sfx\in\bbR$, 
\beas 
\lim_{n\to\infty} E\big[\Psi(\sfz,\sfx)\mathfrak{T}_n(\tti\sfz,\tti\sfx)
\big] 
&=& E\big[\Psi(\sfz,\sfx)\mathfrak{T}(\tti\sfz,\tti\sfx)\big]
\eeas
for the pairs of random symbols $(\mathfrak{T}_n,\mathfrak{T})$$=$
$(\mathfrak{S}^{(1,0)}_n,\mathfrak{S}^{(1,0)})$, 
$(\mathfrak{S}^{(0,1)}_n,\mathfrak{S}^{(0,1)})$, 
$(\mathfrak{S}^{(2,0)}_{1,n},\mathfrak{S}^{(2,0)}_1)$, 
and 
$(\mathfrak{S}^{(1,1)}_{1,n},\mathfrak{S}^{(1,1)}_1)$. 
\ed
\ed
\onelineskip

The random symbol ${\mathfrak S}$ is defined by 
\bea\label{202004220342}
\mathfrak{S}(\tti\sfz,\tti\sfx)
&=&
\mathfrak{S}^{(3,0)}(\tti\sfz,\tti\sfx) 
+\mathfrak{S}^{(1,1)}(\tti\sfz,\tti\sfx) 
+\mathfrak{S}^{(1,0)}(\tti\sfz,\tti\sfx) 
\nn\\&&
+\mathfrak{S}^{(0,1)}(\tti\sfz,\tti\sfx)
+\mathfrak{S}^{(2,0)}_{1}(\tti\sfz,\tti\sfx)
+\mathfrak{S}^{(1,1)}_{1}(\tti\sfz,\tti\sfx),
\eea
{\sred where} the components 
$\mathfrak{S}^{(3,0)}(\tti\sfz,\tti\sfx)$  
and $\mathfrak{S}^{(1,1)}(\tti\sfz,\tti\sfx)$  
are gven by 
(\ref{202003181606}) 
and 
{\sred (\ref{202004161636})}, 
respectively, 
and 
$\mathfrak{S}^{(1,0)}(\tti\sfz,\tti\sfx)$, $\mathfrak{S}^{(0,1)}(\tti\sfz,\tti\sfx)$,  
$\mathfrak{S}^{(2,0)}_{1}(\tti\sfz,\tti\sfx)$ and $\mathfrak{S}^{(1,1)}_{1}(\tti\sfz,\tti\sfx)$ 
are characterized by Condition $[D]$ (III). 
Remark that $\mathfrak{S}^{(2,0)}_{0}(\tti\sfz,\tti\sfx)=0$ 
{\sred in the general expression of $\mathfrak{S}$ in (\ref{202004161620}).}
We define the random symbol ${\mathfrak S}_n$ by 
\bea\label{202004211634}
\mathfrak{S}_n
&=&
1+n^{-1/2}\mathfrak{S},
\eea
{\sred Let $\beta_1$ is the maximum degree of 
the (polynomial) random symbols $\mathfrak{S}^{(1,0)}(\tti\sfz,\tti\sfx)$, 
$\mathfrak{S}^{(0,1)}(\tti\sfz,\tti\sfx)$, $\mathfrak{S}^{(2,0)}_1(\tti\sfz,\tti\sfx)$ 
and $\mathfrak{S}^{(1,1)}_1(\tti\sfz,\tti\sfx)$. 
Remark that the possible maximum degree of $\mathfrak{S}^{(3,0)}(\tti\sfz,\tti\sfx)$ and 
$\mathfrak{S}^{(1,1)}(\tti\sfz,\tti\sfx)$ is $5$.}
The proof of the following theorem is given in Section \ref{202004210757}. 
\begin{theorem}\label{202004210802}
Suppose that Condition $[D]$ is satisfied. Then 
\beas 
\sup_{f\in \calb}\bigg|
E\big[f(Z_n,X_n)\big] - 
E\big[\mathfrak{S}_n(\partial_z,\partial_x)f(G_\infty^{1/2}\zeta,X_\infty)\big]
\bigg|
&=& 
o(n^{-1/2})
\eeas
as $n\to\infty$ 
for any bounded set $\calb$ in 
$C^{{\sred\beta+1}}_b(\bbR^2)$, 
where 
{\sred$\beta=\max\{7,\beta_1\}$ and}
$\zeta\sim N(0,1)$ independent of $\calf$. 
\end{theorem}


\subsection{Expansion of $E[f(Z_n,X_n)]$ for a measurable function $f$}
The index $\beta_x$ is the maximum degree in $\sfx$ of ${\mathfrak S}$ in Section \ref{202004211610}. 
The condition we will work in this section is 
\bd\im[[D$^\sharp$\!\!]]
{\bf (I)} Satisfied are 
$[A]$ (\ref{ai}), (\ref{aiii}), (\ref{aiv}), (\ref{av}), (\ref{avi}), (\ref{aix}), 
$[A^\sharp]$ (\ref{aviii}$^\sharp$), $[A^\sharp]$ (II) and 
%
\bd\im
\bd
\im[(\ref{ai}$^{\sf x}$)] For every $t,s\in[0,1]$ and $q\in\calq$, 
$a(t,q)\in\bbD^{3+\beta_x,\infty}$, $\dota(t,s,q)\in\bbD^{3+\beta_x,\infty}$, $\dot{a}(t,2)\in\bbD^{3+\beta_x,\infty}$ (if $2\in\calq$) and $\ddota(t,s,q)\in\bbD^{3+\beta_x,\infty}$, and 
\bea\label{202003251201} &&
\max_{q\in\calq}\sup_{t\in[0,1]}\big\|a(t,q)\big\|_{3+\beta_x,p}
+\max_{q\in\calq}\sup_{(t,s)\in[0,1]^2}\big\|\dota(t,s,q)\big\|_{3+\beta_x,p}
\nn\\&&
+1_{\{2\in\calq\}}\sup_{t\in[0,1]}\big\|\dot{a}(t,2)\big\|_{3+\beta_x,p}
+\max_{q\in\calq}\sup_{(t,s)\in[0,1]^2}\big\|\ddota(t,s,q)\big\|_{3+\beta_x,p}
\><\> 
\infty
\nn\\&&
\eea
for every $p>1$. 
\ed
\ed
\bd
\im[(II)] $N_n, X_n\in\bbD^\infty$ and 
the properties (\ref{202004011428}) and (\ref{202004011429}) {\sred hold} 
for some positive nuber $\kappa$. 
\im[(III)] 
The random symbols $\mathfrak{S}^{(1,0)}$, $\mathfrak{S}^{(0,1)}$, $\mathfrak{S}^{(2,0)}_1$ 
and $\mathfrak{S}^{(1,1)}_1$ are a (polynomial) random symbol 
the coefficients of which are in $\bbD^{3+\beta_x,1+}$. 
Moreover, 
for each pairs of random symbols $(\mathfrak{T}_n,\mathfrak{T})$$=$
$(\mathfrak{S}^{(1,0)}_n,\mathfrak{S}^{(1,0)})$, 
$(\mathfrak{S}^{(0,1)}_n,\mathfrak{S}^{(0,1)})$, 
$(\mathfrak{S}^{(2,0)}_{1,n},\mathfrak{S}^{(2,0)}_1)$, 
and 
$(\mathfrak{S}^{(1,1)}_{1,n},\mathfrak{S}^{(1,1)}_1)$, and for some $p>1$, 
there exists a (polynomial) random symbol $\overline{\mathfrak T}_n$ 
that has $L^p$ coefficients and the same degree as ${\mathfrak T}$ such that 
\beas 
E\big[\Psi(\sfz,\sfx)\mathfrak{T}_n(\tti\sfz,\tti\sfx)
\big] 
&=& E\big[\Psi(\sfz,\sfx)\overline{\mathfrak{T}}_n(\tti\sfz,\tti\sfx)\big]
\eeas
and $\overline{\mathfrak T}_n\to{\mathfrak T}$ in $L^p$ as $n\to\infty$ 
for every $\sfz\in\bbR$ and $\sfx\in\bbR$. 
\ed
\ed
\onelineskip

The random symbol ${\mathfrak S}_n$ is defined by (\ref{202004211634}) 
for the present ${\mathfrak S}$ described in (\ref{202004220342}). 
We consider the asymptotic expansion formula (\ref{202004211625}) 
for the present ${\mathfrak S}_n$. 
This function is well defined under Condition $[D^\sharp]$. 
The error of the asymptotic expansion of $E[f(Z_n,X_n)]$ for a measurable function $f$ 
is estimated as follows. Proof is given in Section \ref{202004211638}. 
\begin{theorem}\label{202004211628}
Suppose that Condition $[D^\sharp]$ is satisfied. Then 
\bea\label{202005010703}
\sup_{f\in\cale(M,\gamma)}
\bigg| E\big[f(Z_n,X_n)\big] 
-\int_{\bbR^2}f(z,x)p_n(z,x)dzdx\bigg|
&=& 
o(r_n)
\eea
as $n\to\infty$ for every $(M,\gamma)\in(0,\infty)^2$. 
\end{theorem}
\halflineskip

Theorem \ref{202004211628} can provide an asymptotic expansion formula of $\call\{Z_n\}$ 
without a reference variable $X_n$. 
Indeed, even when $X_n$ is not given, we can apply Theorem \ref{202004211628} 
to the joint distribution $\call\{(Z_n,X_\infty)\}$ by adding an independent variable $X_n=X_\infty\sim N(0,1)$ 
and extending the Malliavin calculus for $Z_n$ to incorporate $X_\infty$. 
Though this is a possible way, reconstruction of the proof without $X_n$ 
has a merit that it reduces the differentiability index. 
So we directly apply Theorem 7.7 of Nualart and Yoshida \cite{nualart2019asymptotic} 
to the case without the reference variable $X_n$ nor $X_\infty$. 
Let $\widehat{\Psi}(\sfz)=\exp(2^{-1}G_\infty(\tti\sfz)^2)$ for $\sfz\in\bbR$. 

\bd\im[[D$^\natural$\!\!]]
{\bf (I)} Satisfied are 
$[A]$ (\ref{ai}), (\ref{aiii}), (\ref{aiv}), (\ref{av}), (\ref{avi}), 
$[A^\sharp]$ (\ref{aviii}$^\sharp$) and 
%
\bd\im
\bd
\im[(\ref{ai}$^\natural$)] For every $t,s\in[0,1]$ and $q\in\calq$, 
$a(t,q)\in\bbD^{3,\infty}$, $\dota(t,s,q)\in\bbD^{3,\infty}$, $\dot{a}(t,2)\in\bbD^{3,\infty}$ (if $2\in\calq$) and $\ddota(t,s,q)\in\bbD^{3,\infty}$, and 
\bea\label{202003251201} &&
\max_{q\in\calq}\sup_{t\in[0,1]}\big\|a(t,q)\big\|_{3,p}
+\max_{q\in\calq}\sup_{(t,s)\in[0,1]^2}\big\|\dota(t,s,q)\big\|_{3,p}
\nn\\&&
+1_{\{2\in\calq\}}\sup_{t\in[0,1]}\big\|\dot{a}(t,2)\big\|_{3,p}
+\max_{q\in\calq}\sup_{(t,s)\in[0,1]^2}\big\|\ddota(t,s,q)\big\|_{3,p}
\><\> 
\infty
\nn\\&&
\eea
for every $p>1$. 
\ed
\ed
\bd
\im[(II)] $N_n\in\bbD^\infty$, 
$\|N_n\|_{8,p}\yeq O(1)$ and $\|D_{u_n}^2N_n\|_{7,p}\yeq O(r_n^\kappa)$ ($p>1$) 
as $n\to\infty$ for some positive nuber $\kappa$. 
\im[(III)] 
The random symbols $\mathfrak{S}^{(1,0)}$, 
and $\mathfrak{S}^{(2,0)}_1$ 
are a (polynomial) random symbol 
the coefficients of which are in $\bbD^{2,1+}$. 
Moreover, 
for each pairs of random symbols $(\mathfrak{T}_n,\mathfrak{T})$$=$
$(\mathfrak{S}^{(1,0)}_n,\mathfrak{S}^{(1,0)})$, 
and $(\mathfrak{S}^{(2,0)}_{1,n},\mathfrak{S}^{(2,0)}_1)$,  
and for some $p>1$, 
there exists a (polynomial) random symbol $\overline{\mathfrak T}_n$ 
that has $L^p$ coefficients and the same degree as ${\mathfrak T}$ such that 
\beas 
E\big[\widehat{\Psi}(\sfz)\mathfrak{T}_n(\tti\sfz)\big] 
&=& 
E\big[\widehat{\Psi}(\sfz)\overline{\mathfrak{T}}_n(\tti\sfz)\big]
\eeas
and $\overline{\mathfrak T}_n\to{\mathfrak T}$ in $L^p$ as $n\to\infty$ 
for every $\sfz\in\bbR$.

\im[(IV)] $G_\infty^{-1}\in L^\inftym$. 
\ed
\ed
\onelineskip

In the present situation with no reference variable $X_n$, we define the random symbol ${\mathfrak S}$ 
by 
\beas 
{\mathfrak S}(\tti\sfz)
&=& 
{\mathfrak S}^{(3,0)}(\tti\sfz)+{\mathfrak S}^{(2,0)}_0(\tti\sfz)
+{\mathfrak S}^{(1,0)}(\tti\sfz)+{\mathfrak S}^{(2,0)}_1(\tti\sfz).
\eeas
Let 
\beas
\hat{p}_n(z) &=& 
E\big[\widehat{\mathfrak{S}}_n(\partial_z)^*\phi(z;0,G_\infty)\big]
\eeas
for $\widehat{\mathfrak{S}}_n(\tti\sfz)=1+n^{-1/2}\widehat{\mathfrak{S}}(\tti\sfz)$. 
Denote by $\cale_1(M,\gamma)$ the set  of measurable functions $f:\bbR\to\bbR$ 
satisfying $|f(z)|\leq M(1+|z|)^\gamma$ for $M,\gamma>0$. 

\begin{theorem}\label{202005120214}
Suppose that Condition $[D^\natural]$ is satisfied. 
Let $(M,\gamma)\in(0,\infty)^2$. Then 
\bea\label{202005010707}
\sup_{f\in\cale_1(M,\gamma)}
\bigg| E\big[f(Z_n)\big] 
-\int_{\bbR}f(z)\hat{p}_n(z)dz\bigg|
&=& 
o(r_n)
\eea
as $n\to\infty$. 
\end{theorem}

\begin{en-text}
Here the density function $\hat{p}$ is defined by 
\beas 
\hat{p}(z) 
&=& 
\int_\bbR p_n(z,x)dx
\eeas
for $p_n(z,x)$ in (\ref{202004211625}) with ${\mathfrak S}$ described by 
\beas
\mathfrak{S}(\tti\sfz,\tti\sfx) 
&=&
\mathfrak{S}^{(3,0)}(\tti\sfz,\tti\sfx) 
+\mathfrak{S}^{(1,0)}(\tti\sfz,\tti\sfx) 
+\mathfrak{S}^{(2,0)}_{1}(\tti\sfz,\tti\sfx).
\eeas
We remark that $D_{u_n\tti\sfz}X_\infty=0$ by orthogonality of $X_\infty$ to the original variables. 
The random symbol ${\mathfrak S}_n(\tti\sfz,\tti\sfx)$ can be written by ${\mathfrak S}_n(\tti\sfz)$ 
since it does not depend on $\sfx$, and then

By the blowing-up, $\mathfrak{S}_n(\partial_z,\partial_x)$ involves $\tti\sfx$ but 
the effect of it vanishes after all under independency of $X_\infty$ and the fact that $DX_\infty$ is deterministic. 
We refer the reader to Nualart and Yoshida \cite{nualart2019asymptotic} 
for rigorous validation of the asymptotic expansion by $\hat{p}_n(z)$. 
In this situation, Condition $[D^\sharp]$ (III) should be as follows. \koko
The random symbols $\mathfrak{S}^{(1,0)}$ and 
$\mathfrak{S}^{(2,0)}_1$ 
are a (polynomial) random symbol 
the coefficients of which are in $\bbD^{3+\beta_x,1+}$. 
Moreover, 
for every $\sfz\in\bbR$ and $\sfx\in\bbR$, 
\beas 
\lim_{n\to\infty} E\big[\hat{\Psi}(\sfz)\mathfrak{T}_n(\tti\sfz)\psi_n\big] 
&=& E\big[\Psi(\sfz,\sfx)\mathfrak{T}(\tti\sfz)\big]
\eeas
for the pairs of random symbols $(\mathfrak{T}_n,\mathfrak{T})$$=$
$(\mathfrak{S}^{(1,0)}_n,\mathfrak{S}^{(1,0)})$, 
$(\mathfrak{S}^{(0,1)}_n,\mathfrak{S}^{(0,1)})$, 
$(\mathfrak{S}^{(2,0)}_{1,n},\mathfrak{S}^{(2,0)}_1)$, 
and 
$(\mathfrak{S}^{(1,1)}_{1,n},\mathfrak{S}^{(1,1)}_1)$, 
where $\hat{\Psi}(\sfz)=\exp(2^{-1}G_\infty[(\tti\sfz)^{\otimes2}]$. 
Indeed, $\Psi(\sfz,\sfx)=\hat{\Psi}(\sfz)e^{-x^2/2}$ and 
the general ${\mathfrak S}_n$ in \koko
{\colorr Write Theorem.}
\end{en-text}

\section{An example}\label{202005180614}
We will consider an example in this section. 
For a Wiener process $w=(w_t)_{t\in[0,1]}$ realized by an isonormal Gaussian process over 
a real separable Hilbert space $\mH$. 
We will consider a quadratic variation with anticipative weights 
\beas 
\V_n &=& \sum_{j=1}^na(w_{1-\tj})(\Delta_jw)^2
\eeas
where $a\in C^\infty_p(\bbR)$ and $\tj=j/n$. 
Let us simply denote $f_t=f(w_t)$ as before. 
It is easy to see 
\beas 
\V_n &\to^p& \V_\infty
\eeas
as $n\to\infty$ for  
$\V_\infty = \int_0^1 a_{1-t}dt$. 
Let $Z_n= \sqrt{n}(V_n-V_\infty)$. Then 
\beas 
Z_n&:=&\sqrt{N}\big(\V_n-\V_\infty \big)
\yeq 
M_n+n^{-1/2}N_n
\eeas
where 
$M_n =\delta(u_n)$ for 
\bea\label{202005180753}
u_n &=& n^{1/2}\sum_{j=1}^na_{1-\tj}I_1(1_j)1_j
\eea
and 
\bea\label{202005181317}
N_n 
&=& 
n\sum_{j=1}^n(D_{1_j}a_{1-\tj})I_1(1_j)+n\sum_{j=1}^n
\int_{{\sred \tjm}}^{{\sred \tj}}(a_{1-\tj}-a_{1-t})dt
\nn\\&=:&
N_n^{(1)}+N_n^{(2)}.
\eea
This $Z_n$ is an example of that in (\ref{202003130001}) if one sets 
$\calq=\{2\}$ and $a_j(2)=a_{1-\tj}$. 
For a reference variable, we shall consider $X_\infty=w_1$. 
Therefore the asymptotic expansion of the joint law $\call\{(Z_n,X_\infty)\}$ 
is given by Theorem \ref{202004211628}. 

We will verify Condition $[D^\sharp]$. 
Condition $[A]$ (\ref{ai}) is obvious. 
For 
\beas 
a(s,2)=a_{1-s},
\eeas 
Condition $[A]$ (\ref{aiii}) holds. 
Since 
\beas 
D_{1_k}a_{1-\tj} &=& a'(w_{1-\tj})n^{-1}1_{\{t_k\leq 1-t_j\}},
\eeas
Condition $[A]$ (\ref{aiv}) is satisfied with 
\beas
\dota(t,s,2)=a'(w_{1-s})1_{\{t\leq 1-s\}}, 
\eeas
and Condition $[A]$ (\ref{av}) with
\beas 
\dot{a}(t,2)=a'(w_{1-t})1_{\{t\leq 1/2\}}.
\eeas
We have 
\beas 
D_{1_k}D_{1_k}a(w_{1-\tj}) 
&=& 
a''(w_{1-\tj})n^{-2}1_{\{t_k\leq 1-\tj\}}
\eeas
and hence Condition $[A]$ (\ref{avi}) is satisfied by
\beas 
\ddota(t,s,2) 
&=& 
a''(w_{1-s})1_{\{t\leq 1-s\}}.
\eeas
Condition $[A]$ (\ref{aix}) is met for $X_\infty=w_1$ and $\ddot{X}_\infty(t)=0$. 

In the present situation.  
\beas
G_\infty 
&=& 
2\int_0^1a_{1-s}^2ds\yeq 2\int_0^1 a_s^2ds. 
\eeas
Then Condition $[A^\sharp]$ (\ref{aviii}$^\sharp$) is satisfied. 
In fact, Section \ref{202004201231} shows (\ref{2020041446}). 
We assume that $G_\infty^{-1}\in L^\inftym$. Then 
Condition $[A^\sharp]$ (II) is fulfilled. 
Condition $[D^\sharp]$ (i$^{\sf x}$) is easy to check. 

The variable $X_\infty=w_1$ satisfies the conditions in $[D^\sharp]$ (II) concerning it, so we verify it for $N_n$. 
{\sred According to Section \ref{202004270705},} 
the exponent of the first component of (\ref{202005181317}) is 
$e\big(N_n^{(1)}\big)=0$, 
which means $N_n^{(1)}=O_{\bbD^\infty}(1)$. 
Moreover, 
Proposition \ref{202004091035} applied to $N_n^{(1)}$ represented by (\ref{202005181317}) 
ensures $D_{u_n}N_n^{(1)}=O_{\bbD^\infty}(n^{-1/2})$, and hence $D_{u_n}^2N_n^{(1)}=O_{\bbD^\infty}(n^{-1/2})$. Thus $N_n^{(1)}$ satisfies the related part of Condition $[D^\sharp]$ (II). 

For $N_n^{(2)}$, we have 
\beas 
N_n^{(2)}
&=&
-n\sum_{j=1}^n\int_\tjm^\tj(a_t-a_\tjm)dt
\eeas
Let 
\beas 
A_n(j;x) &=& -\int_0^1a'\big(w_x+v(w_x-w_\tjm)\big)dv.
\eeas
for $x\in\bbT$, the one-dimensional torus identified with $[0,1]$, 
and let 
\beas 
\dot{\bbH}_n^{(0)}(x)
&=&
n\sum_{j=1}^n A_n(j;x)I_1\big(1_{[\tjm,x]}1_j(x)\big)
\eeas
Then 
\beas
N_n^{(2)}
&=&
\int_\bbT\dot{\bbH}_n^{(0)}(x)dx,
\eeas
where $dx$ is the Haar measure on the compact group $\bbT$ compatible with the Lebesgue measure. 
The functional ${\stackrel{\circ}{\bbH}}_n(\mbx)$ of (\ref{202005190511}) is 
\beas 
{\stackrel{\circ}{\bbH^{(0)}_n}}(x)
&=& 
n\sum_{j=1}^n A_n(j;\tjm+x)I_1\big(1_{[\tjm,\tjm+x]}1_j(\tjm+x)\big).
\eeas
for $\dot{\bbH}_n^{(0)}(x)$. 
The exponent 
\beas 
e\bigg({\stackrel{\circ}{\bbH^{(0)}_n}}(x)\bigg)
&=&
1
\eeas
and hence 
\beas 
{\stackrel{\circ}{\bbH^{(0)}_n}}(x)
&=& 
O_{\bbD^\infty}(n)
\eeas
uniformly in $x\in\bbT$ because this estimate follows from Theorem \ref{0110221242}. 
Therefore we obtain 
\beas 
N_n^{(2)} &=& O_{L^\inftym}(1)
\eeas
by Theorem \ref{202005190527} applied to the case where 
\beas
e^{(1)}_{n,j}(x)
&=& 
I_1\big(1_{[\tjm,x]}1_j(x)\big)
\yeq 
I_1\big(1_{[\tjm,x]}1_1(x-\tjm)\big)1_{[0,n^{-1}]}(x-\tjm),
\\
g^{(1)}_{n,j} &=& -\tjm,
\\
\chi^{(1)}_n &=& [0,n^{-1}].
\eeas
The $i$-th derivative $D^iN_n^{(2)}$ can be estimated in a similar manner. 
Indeed, from (\ref{202005190606}), 
\begin{en-text}
\beas 
D^i_{\mbv}N_n^{(2)}
&=&
\int_\bbT D^i_{|v|_{\mH^{\otimes i}}^{-1}\mbv}\>
{\stackrel{\circ}{\bbH^{(0)}_n}}(x)1_{[0,n^{-1}]}(x)
dx\>|v|_{\mH^{\otimes i}}
\eeas
for non-zero $\mbv\in\mH^{\otimes i}$. 
\end{en-text}
\bea\label{202005190625}  
D^iN_n^{(2)}
&=&
\int_\bbT D^i\>
{\stackrel{\circ}{\bbH^{(0)}_n}}(x)1_{[0,n^{-1}]}(x)dx
\eea
for $i\in\bbZ_+$. 
The (uniform-in-$x$ version of) stability of the operation $D^i$ proved by Proposition \ref{202004121624} 
shows 
\beas
\big\|D^iN_n^{(2)}\big\|_p
&\leq&
\int_\bbT \>\bigg\|D^i{\stackrel{\circ}{\bbH^{(0)}_n}}(x)\bigg\|_p1_{[0,n^{-1}]}(x)dx
\yeq 
O(1)
\eeas
for every $p>1$. 
We obtained the estimate $N_n^{(2)}=O_{\bbD^\infty}(1)$ as $n\to\infty$. 

Operating $D_{u_n}$ to (\ref{202005190625}) in $i=0$, we have 
\bea\label{202005190639} 
D^kD_{u_n}N_n^{(2)}
&=&
\int_\bbT D^kD_{u_n}{\stackrel{\circ}{\bbH^{(0)}_n}}(x)1_{[0,n^{-1}]}(x)dx
\eea
for $k\in\bbZ_+$. 
Since ${\stackrel{\circ}{\bbH^{(0)}_n}}(x)$ involved only the first chaos, 
Proposition \ref{202004091035} (ii) undertakes the estimate 
\beas
e\bigg(D_{u_n}{\stackrel{\circ}{\bbH^{(0)}_n}}(x)\bigg)
&\leq& 
0.5, 
\eeas
and the stability Proposition \ref{202004121624} gives 
\beas
\bigg\|D^iD_{u_n}{\stackrel{\circ}{\bbH^{(0)}_n}}(x)\bigg\|_p
&=& 
O(n^{0.5})
\eeas
for every $p>1$. 
Combining this and the representation (\ref{202005190639}), we conclude 
$D_{u_n}N_n^{(2)}=O_{\bbD^\infty}(n^{-0.5})$. 
We can also say $D_{u_n}^kN_n^{(2)}=O_{\bbD^\infty}(n^{-0.5})$ for $k\geq1$ 
since $u_n=O_{\bbD^\infty}(1)$. 
Thus, the part concerning $N_n^{(2)}$ in $[D^\sharp]$ (II) has been verified. 
$[D^\sharp]$ (II) has been checked. 
In particular, $D_{u_n}N_n=O_{\bbD^\infty}(n^{-0.5})$.

The random symbol ${\mathfrak S}^{(3,0)}$ defined by (\ref{202003181606}) is in this case 
\bea\label{202005191239}
{\mathfrak S}^{(3,0)}(\tti\sfz,\tti\sfx)
&=&
\begin{en-text}
\frac{4}{3}\int_0^1 a_{1-t}^3dt\>(\tti\sfz)^3
+\int_{[0,1]^2} a^{(3,0)}(t,s,2,2)dsdt(\tti\sfz)^3
\nn\\&&
+\int_0^1\big(D_tG_\infty(\tti\sfz)^5+(\tti\sfz)^3(\tti\sfx)\big)
\bigg[\int_0^{1-t} a'_{1-s}
a_{1-s}ds\bigg]\>a_{1-t}\>dt
\nn\\&=&
\end{en-text}
\frac{4}{3}\int_0^1 a_{1-t}^3dt\>(\tti\sfz)^3
+{\sred 2}\int_{[0,1]^2} a^{(3,0)}(t,s,2,2)dsdt(\tti\sfz)^3
\nn\\&&
+\int_0^1D_tG_\infty
\bigg[\int_0^{1-t} a'_{1-s}
a_{1-s}ds\bigg]\>a_{1-t}\>dt\>(\tti\sfz)^5
\nn\\&&
+{\sred 2}\int_0^1
\bigg[\int_0^{1-t} a'_{1-s}
a_{1-s}ds\bigg]\>a_{1-t}\>dt\>(\tti\sfz)^3(\tti\sfx)
\eea
where 
$a^{(3,0)}(t,s,q_1,q_2)$ given by (\ref{202003181607}) is 
\beas
a^{(3,0)}(t,s,2,2)
&=&
{\sred a''_{1-s}}
a_{1-s}a_{1-t}1_{\{t\leq 1-s\}}
+(a'_{1-s})^2a_{1-t}1_{\{t\leq 1-s\}}
\nn\\&&
+a'_{1-s}a_{1-s}a'_{1-t}1_{\{t\leq 1-s\}}1_{\{t\leq 1/2\}}. 
\eeas

The random symbol ${\mathfrak S}^{(1,1)}$ defined in (\ref{202004161636}) is now 
\bea\label{202005191240}
{\mathfrak S}^{(1,1)}(\tti\sfz,\tti\sfx)
&=& 
 \int_0^\half a'(w_{1-t})dt(\tti\sfz)(\tti\sfx)
+\half\int_0^1(D_tG_\infty)a_{1-t}dt\>(\tti\sfz)^3(\tti\sfx)
\nn\\&&
+\int_0^1a_{1-t}dt(\tti\sfz)(\tti\sfx)^2. 
\nn\\&&
\eea

It is now obvious that 
${\mathfrak S}^{(0,1)}=0$, 
${\mathfrak S}^{(2,0)}_1=0$ and 
${\mathfrak S}^{(1,1)}_1=0$. 
%
We have 
\beas 
N_n^{(2)}
&=&
{\sred 
-n\sum_{j=1}^n\int_\tjm^\tj\big(a_t-a_\tjm\big)dt
}
\nn\\&=&
-\sum_{j=1}^n a'_\tjm I_1(g_j)
-n\sum_{j=1}^n \int_\tjm^\tj\int_\tjm^t 2^{-1}a''_sds dt
+O_{L^\inftym}(n^{-0.5})
\eeas
where 
$g_j$ is defined in (\ref{202005060441}), and 
we used the esimate 
\beas
-n\sum_{j=1}^n a''_\tjm\int_\tjm^\tj\int_\tjm^t\int_\tjm^sdw_rdw_sdt
&=&
O_{L^\inftym}(n^{-0.5}). 
\eeas
Then, for ${\mathfrak S}^{(1,0)}_n(\tti\sfz)=N_n\tti\sfz$, we have 
\beas 
E\big[\Psi(\sfz,\sfx){\mathfrak S}^{(1,0)}_n(\tti\sfz)\big]
&=&
E\big[\Psi(\sfz,\sfx)N_n\big]\tti\sfz
\nn\\&=&
E\bigg[\Psi(\sfz,\sfx)n\sum_{j=1}^n(D_{1_j}a_{1-\tj})I_1(1_j)\bigg]\tti\sfz
\nn\\&&
-E\bigg[\Psi(\sfz,\sfx)
\sum_{j=1}^n a'_\tjm I_1(g_j)\bigg]\tti\sfz
\nn\\&&
-E\bigg[\Psi(\sfz,\sfx)
n\sum_{j=1}^n \int_\tjm^\tj\int_\tjm^t 2^{-1}a''_sds dt\bigg]\tti\sfz
+o(1)
\nn\\&=&
\sum_{j=1}^n E\bigg[D_{1_j}\big(\Psi(\sfz,\sfx)a'_{1-\tj}1_{\{\tj\leq1/2\}}\big)\bigg]\tti\sfz
\nn\\&&
-\sum_{j=1}^nE\bigg[D_{g_j}\big(\Psi(\sfz,\sfx) a'_\tjm \big)\bigg]\tti\sfz
\nn\\&&
-\sum_{j=1}^nE\bigg[\Psi(\sfz,\sfx)
n \int_\tjm^\tj\int_\tjm^t 2^{-1}a''_sds dt\bigg]\tti\sfz
+o(1)
\eeas
In the present problem, 
\beas 
\Psi(\sfz,\sfx)
&=&
\exp\big(2^{-1}G_\infty(\tti\sfz)^2+w_1\tti\sfx\big).
\eeas
Therefore, we obtain 
\beas &&
\sum_{j=1}^n E\bigg[D_{1_j}\big(\Psi(\sfz,\sfx)a'_{1-\tj}1_{\{\tj\leq1/2\}}\big)\bigg]\tti\sfz
\nn\\&=&
\sum_{j=1}^n E\bigg[\Psi(\sfz,\sfx)\big\{
2^{-1}D_{1_j}G_\infty(\tti\sfz)^2
+n^{-1}\tti\sfx\big\}
a'_{1-\tj}1_{\{\tj\leq1/2\}}\bigg]\tti\sfz
\nn\\&&
+\sum_{j=1}^n E\bigg[\Psi(\sfz,\sfx)a''_{1-\tj}1_{\{\tj\leq1/2\}}n^{-1}\bigg]\tti\sfz
\nn\\&\to&
E\bigg[\Psi(\sfz,\sfx)\bigg\{
\int_0^\half2^{-1}D_tG_\infty a'_{1-t}dt\>(\tti\sfz)^3
+\int_0^\half a'_{1-t}\>dt(\tti\sfz)(\tti\sfx)
+\int_0^\half a''_{1-t}dt\>(\tti\sfz)\bigg\}\bigg], 
\eeas
\beas &&
-\sum_{j=1}^nE\bigg[D_{g_j}\big(\Psi(\sfz,\sfx) a'_\tjm \big)\bigg]\tti\sfz
\nn\\&=&
-\sum_{j=1}^nE\bigg[\Psi(\sfz,\sfx) \big\{2^{-1}D_{g_j}G_\infty(\tti\sfz)^2+2^{-1}n^{-1}\tti\sfx\big\}a'_\tjm
\bigg]\tti\sfz
\nn\\&&
-\sum_{j=1}^nE\bigg[\Psi(\sfz,\sfx) D_{g_j}a'_\tjm \bigg]\tti\sfz
\nn\\&\to&
{\sred 
E\bigg[\Psi(\sfz,\sfx)\bigg\{
-\int_0^1 4^{-1}D_tG_\infty a'_t dt\>(\tti\sfz)^3-\int_0^12^{-1}a'_tdt(\tti\sfz)(\tti\sfx)
\bigg\}\bigg]
}
\eeas
and 
\beas 
-\sum_{j=1}^nE\bigg[\Psi(\sfz,\sfx)
n \int_\tjm^\tj\int_\tjm^t 2^{-1}a''_sds dt\bigg]\tti\sfz
&\to&
E\bigg[\Psi(\sfz,\sfx)(-4^{-1})\int_0^1 a''_tdt\tti\sfz\bigg]
\eeas
Thus, we conclude 
\beas 
E\big[\Psi(\sfz,\sfx){\mathfrak S}^{(1,0)}_n(\tti\sfz)\big]
&\to&
E\big[\Psi(\sfz,\sfx){\mathfrak S}^{(1,0)}(\tti\sfz)\big]
\eeas
as $n\to\infty$ for 
\bea\label{202005191238}
{\mathfrak S}^{(1,0)}(\tti\sfz)
&=& 
\begin{en-text}
\int_0^\half2^{-1}D_tG_\infty a'_{1-t}dt\>(\tti\sfz)^3
+\int_0^\half a'_{1-t}\>dt(\tti\sfz)(\tti\sfx)
+\int_0^\half a''_{1-t}dt\>(\tti\sfz)
\nn\\&&
-\int_0^1 4^{-1}D_tG_\infty a'_t dt\>(\tti\sfz)^3-\int_0^12^{-1}a'_tdt(\tti\sfz)(\tti\sfx)
-\int_0^12^{-1}a'_tdt\>(\tti\sfz)
\nn\\&&
-4^{-1}\int_0^1 a''_tdt(\tti\sfz). 
\nn\\&=&
\end{en-text}
\int_0^\half2^{-1}D_tG_\infty a'_{1-t}dt\>(\tti\sfz)^3
-\int_0^1 4^{-1}D_tG_\infty a'_t dt\>(\tti\sfz)^3
\nn\\&&
+\int_0^\half a'_{1-t}\>dt(\tti\sfz)(\tti\sfx)
-\int_0^12^{-1}a'_tdt(\tti\sfz)(\tti\sfx)
\nn\\&&
{\sred 
+\int_0^\half a''_{1-t}dt\>(\tti\sfz)
-4^{-1}\int_0^1 a''_tdt(\tti\sfz). 
}
\eea

In the present case, the full random symbol 
\beas 
{\mathfrak S} 
&=& 
{\mathfrak S}^{(3,0)}+{\mathfrak S}^{(1,1)}+{\mathfrak S}^{(1,0)}
\eeas 
is expressed more precisely with 
(\ref{202005191239}), (\ref{202005191240}) and (\ref{202005191238}). 
For the density function $p_n(z,x)$ for 
{\sred 
${\mathfrak S}_n$
defined by (\ref{202004211634}), 
}
the asymptotic expansion {\sred of $\call\{(Z_n,w_1)\}$} 
can be obtained {\sred from Therome \ref{202004211628}}
under $G_\infty^{-1}\in L^\inftym$ for measurable functions.

\section{
Robust realized volatility}\label{202004261806}
Robustness is an established notion in statistics 
(Huber \cite{huber2004robust}, Hampel et al. \cite{hampel2011robust}). 
According to the concept of robust statistics, 
the most effective part of the estimator should be removed or diminished to attain robustness  
because that part is most sensitive to contamination of the data. 
In volatility estimation, it is known that the realized volatility is weak against jumps, and 
multi-power variations and 
(predictable) thresholding methods for example have been studied to approach this issue. 
Recently, Inatsugu and Yoshida \cite{inatsugu2018global} proposed 
global jump filters to realize robust estimation for volatility of an It\^o process. 
For correct rejection of contaminated observations, their filters use a selector depending on the whole data. 
The resulting estimator obviously destroys the martingale structure, 
however resists jumps well, especially consecutive jumps, for which 
traditional methods like by-power variation do not work. 
Related to such construction, we will consider another type of robust variation with a global filter 
and investigate asymptotic expansion as an example of our scheme.

Consider a unique strong solution $X=(X_t)_{t\in\bbR_+}$ to the stochastic differential equation 
\beas 
\left\{
\begin{array}{ccl}
dX_t &=& \sigma(X_t)dw_t+b(X_t)dt,\quad t\in\bbR_+\\
X_0 &=& x_0.
\end{array}
\right.
\eeas
Here $w=(w_t)$ is a Brownian motion, 
$x_0\in\bbR$, $\sigma,b\in C^\infty$ and we suppose 
every positive order of derivative of $\sigma$ and $b$ is bounded. 
The uniform nondegeneracy $\inf_x\sigma(x)>0$ is assumed. 

The most popular volatility estimator in finance is 
the realized volatility $U_n$ defined by 
\beas 
U_n  &=& \sum_{j=1}^n(\Delta_jX)^2
\eeas
on the interval $[0,1]$, where $\Delta_jX=X_\tj-X_\tjm$ with $\tj=j/n$. 
It is well known that $U_n$ is sensitive to jumps. 
A robust variation of the realized volatility we will consider is 
\bea\label{202005041049}
\bbV_n &=& \sum_{j=1}^n\Phi\left(U_n,L_{n,j}\right)(\Delta_jX)^2
\eea
where $\Phi\in C^\infty_p(\bbR^2)$, the space of smooth functions with derivatives of at most polynomial growth.  
The functional $L_{n,j}$ is defined by 
\beas 
L_{n,j} &=& {\colorb \eta_{n,j}^{-1}}\sum_{{\colorb k\in K_j}}(\Delta_kX)^2
\eeas
{\colorb with 
$
\eta_{n,j}= \# K_j/n
$, 
where 
\beas 
K_j &=& K_{n,j} \yeq \big\{k\in\bbN;\> 
\underline{K}_j\leq k\leq \ol{K}_j\big\}
\eeas
and 
\beas 
\underline{K}_j\yeq\underline{K}_{n,j}\yeq( j-\lfloor n\lambda\rfloor+1)\vee 1,
&&
\ol{K}_j\yeq\ol{K}_{n,j}\yeq(j+\lfloor n\lambda\rfloor-1)\wedge n
\eeas
{\sred for the locality parameter $\lambda\in(0,1)$}. 
}

An example of robust volatility measures is  
\beas 
\bbV_n^* &=& \sum_{j=1}^n\varphi\left(U_n^{-1}L_{n,j}\right)(\Delta_jX)^2
\eeas
where $\varphi\in C_b^\infty(\bbR)$, the space of smooth bounded functions having bounded derivatives. 
The statistic $\bbV_n^*$ selects increments $\Delta_jX$ according to 
the relative magnitude $L_{n,j}$ of the local volatility to the estimated average spot volatility. 
For example, $\varphi$ is a smooth function taking $1$ in $\{x;|x|\leq C\}$ and 
$0$ in $\{x;|x|\geq C\}$ for some positive constant $C$. 
Then, especially when ${\sred \lambda}$ is small, an increment $\Delta_jX$ will be removed from the sum 
if an increment in either $[\tjm,\tj]$ or intervals
{\sred near around} 
is contaminated with a jump. 
Let $\psi:\bbR_+\to[0,1]$ be a smooth function such that $\psi(x)=1$ for $x\in[0,1/c_0]$ and 
$\psi(x)=0$ for $x\in[2/c_0,\infty)$ for $c_0=\inf_{x\in\bbR}\sigma(x)^2/2$. 
Then it is easy to see the asymptotic expansion of $\bbV_n^*$ coincides with that of 
\beas 
\bbV_n^{**} &=& \sum_{j=1}^n\psi(U_n^{-1})\varphi\left(U_n^{-1}L_{n,j}\right)(\Delta_jX)^2
\eeas
The statistic $\bbV_n^{**}$ is exactly an example of $\bbV_n$ 
for $\Phi(x,y)=\psi(|x|^{-1})\varphi(|x|^{-1}|y|)$. 

There is no mathematical difficulty even if 
one adopts {\sred various local estimators} 
but we consider $L_{n,j}$ for notational simplicity. 
{\sred Indeed,} 
as for $U_n$, it was possible to choose a multi-power variation or another global volatility estimator. 

Simply denote $f_t=f(X_t)$ for a function $f$. 
Let 
\beas 
U_\infty \yeq \int_0^1 \sigma_t^2dt,&& 
{\colorb 
L_{\infty,t}\yeq \eta_{\infty,t}^{-1}\int_{(t-\lambda)\vee0}^{(t+\lambda)\wedge1}\sigma_s^2ds
}
\eeas
for
{\colorb
\beas 
\eta_{\infty,t}&=&\{(t+\lambda)\wedge1\}-\{(t-\lambda)\vee0\}. 
\eeas 
}\noindent

{\sred 
We will consider estimation of the filtered integrated volatility 
\bea\label{202101010110}
\ol{\bbV}_n
&=& 
\int_0^1 \Phi\big(U_n,L_{n,j}\big)\sigma_t^2dt,
\eea
while estimation of the limit  
\beas
\bbV_\infty 
&=& 
\int_0^1 \Phi\big(U_\infty,L_{\infty,t}\big)\sigma_t^2dt. 
\eeas
of $\ol{\bbV}_n$ also makes sense though the description will be more involved. 
}

\subsection{Stochastic expansion of the error}
We introduce some systematic notation. 
We write 
$b^{[1]}=\sigma$ and $b^{[2]}=b$. 
For a $C^2$ function $f$, denote 
$f^{[1]}=\sigma f' = b^{[1]}f'$, $f^{[2]}=\half \sigma^2 f''+bf'=\half (b^{[1]})^2f''+b^{[2]}f'$. 
It\^o's formula is then written by 
\beas 
f_t=f_s+\int_s^t f^{[1]}_rdw_r+\int_s^t f^{[2]}_rdr
\eeas
for $t>s\geq0$. 
For $b^{[1]}=\sigma$ and $b^{[2]}=b$, 
define $b^{[i_1,i_2,...,i_m]}$ inductively by 
\beas 
b^{[i_1,i_2,...,i_m]} &=& \big(b^{[i_1,i_2,...,i_{m-1}]}\big)^{[i_m]}
\eeas
for $i_1,...,i_m\in\{1,2\}$, $m\geq2$.

Let 
\beas
\ol{\theta}_t &=& \Phi\big(U_\infty,L_{\infty,t}\big).
\eeas 
\begin{en-text}
We consider $M_n=\delta(u_n)$ for $u_n$ defined by 
\bea\label{202005042346}
u_n 
&=& 
n^{1/2}\sum_{j=1}^n\ol{\theta}_\tjm(b^{[1]}_\tjm)^2I_1(1_j)1_j.
\eea
\end{en-text}
\begin{en-text}
In what follows, $\calx_n\equiv^a\caly_n$ means $\calx_n$ and $\caly_n$ are 
asymptotically equivalent, i.e., 
$\calx_n-\caly_n\to^p0$ as $n\to\infty$ 
for two sequences of random variables $\calx_n$ and $\caly_n$. 
\end{en-text}
To simplify notation, we will often use the symbol $\beta=(b^{[1]})^2$. Then 
\beas 
\beta^{[1]} &=& 2(b^{[1]})^2(b^{[1]})' \yeq 2b^{[1]}b^{[1,1]}
\eeas
and 
\beas 
\beta^{[2]} &=& (b^{[1,1]})^2+2b^{[1]}b^{[1,2]}.
\eeas
Define $S_j^{(\ref{202005041155})}$ by 
\bea\label{202005041155}
S_j^{(\ref{202005041155})}
&=&
{\colorr b^{[1]}_\tjm b^{[1,1]}_\tjm  I_3(1_j^{\otimes3})
+2hb^{[1]}_\tjm b^{[1,1]}_\tjm  I_1(1_j)}
\nn\\&&
+2hb^{[1]}_\tjm b^{[2]}_\tjm I_1(1_j) 
+h^2b^{[1]}_\tjm b^{[1,2]}_\tjm 
+h^2b^{[1]}_\tjm b^{[2,1]}_\tjm 
\nn\\&&
+h^2(b^{[2]}_\tjm)^2 
+2^{-1}h^2(b^{[1,1]}_\tjm)^2
-h\beta^{[1]}_\tjm I_1(g_j)- 2^{-1}h^2{\colorr \beta^{[2]}_\tjm}.
\eea

Let 
\bea\label{202005030358}
F_j^{(\ref{202005030358})}
&=&
(b^{[1]}_\tjm)^2 I_2(1_j^{\otimes2})
\eea
and let 
\begin{en-text}
\bea\label{202005110334}
\Theta_j
&=&
{\colorb 
\Phi(U_\infty,L_{\infty,\tjm })
+\sum_{k=1}^n\partial_1\Phi(U_\infty,L_{\infty,\tkm})\int_\tkm^\tk\beta_rdr 
}
\nn\\&&\hspace{10pt}
{\colorb
+\sum_{k=1}^n1_{\{ j\in K_k\}}
\partial_2\Phi(U_\infty,L_{\infty,\tkm})\eta_{n,k}^{-1}\int_\tkm^\tk\beta_rdr.
}
\eea
\end{en-text}
\bea\label{202005110334}
\Theta_j
&=&
{\sred 
\Phi(U_\infty,L_{\infty,\tjm })
}
\eea
Let 
\bea\label{202005120042}
\Psi_{j,k}(x,y) 
&=& 
\Psi_{n,j,k}(x,y) 
\yeq
{\colorb1_{\{k\in\bbJ_n\}}}\partial_1\Phi(x,y)
+{\colorb1_{\{k\in K_j\}}}\partial_2\Phi(x,y){\colorb\eta_{n,j}}^{-1}.
\eea
Let 
\bea\label{202005131723}
\Xi_{j,k,\ell}(x,y) 
&=& 
\Xi_{n,j,k,\ell}(x,y) 
\nn\\&=&
\half{\colorb1_{\{k,\ell\in\bbJ_n\}}}\partial_1^2\Phi(x,y)
+\partial_1\partial_2\Phi(x,y){\colorb1_{\{k\in\bbJ_n\}}1_{\{\ell\in K_j\}}\eta_{n,j}^{-1}}
\nn\\&&
+\half\partial_2^2\Phi(x,y)
{\colorb1_{\{k,\ell\in K_j\}}\eta_{n,j}^{-2}}. 
\eea
We introduce a process $(\mba_s)_{s\in[0,1]}$ defined by 
{\sred 
\bea\label{202005120000} 
\mba_s 
&=&
\Phi(U_\infty,L_{\infty,s})\beta_s
.
\eea
}
{\colorb 
\begin{en-text}
\bea\label{202005120000} 
\mba_s 
&=&
\bigg\{\Phi(U_\infty,L_s)
+\int_0^1\big(\partial_1\Phi(U_\infty,L_{\infty,r})\big)\beta_rdr
\nn\\&&\hspace{10pt}
+\int_{(s-\lambda)\vee0}^{(s+\lambda)\wedge1}\big(\partial_2\Phi(U_\infty,L_{\infty,r})\big)\eta_{\infty,r}^{-1}\beta_rdr\bigg\}\beta_s.
\eea
\end{en-text}
Let 
\bea\label{202005141652} 
G_\infty
=
\int_0^12\>\mba_t^2\>dt. 
\eea

Now the scaled error 
\beas 
Z_n 
&=&
n^{1/2}\big(\bbV_n - {\sred \ol{\bbV}_n}
\big)
\eeas
admits the following stochastic expansion. 
\begin{lemma}\label{202005051550}
\bea\label{202005050001}
Z_n 
&=&
M_n+n^{-1/2}N_n
\eea
for 
{\colorr 
$M_n=\delta(u_n)$ for $u_n$ defined in (\ref{202005042346}) to be mentioned subsequently, and 
}
{\sred
\bea\label{202005061538}
N_n 
&=&
n\sum_{j=1}^n \big(D_{1_j}({\colorr\Theta_j}
(b^{[1]}_\tjm)^2)\big) I_1(1_j)
+n\sum_{j=1}^n\Theta_j S_j^{(\ref{202005041155})}
\nn\\&& 
+n\sum_{j,k=1}^n\Psi_{j,k}(U_\infty,L_{\infty,\tjm})F_j^{(\ref{202005030358})}F_k^{(\ref{202005030358})}
+R_n^{(\ref{202005110436})},
\eea
}
\begin{en-text}
{\colorr  
\bea\label{202005061538}
N_n 
&=&
n\sum_{j=1}^n \big(D_{1_j}({\colorr\Theta_j}
(b^{[1]}_\tjm)^2)\big) I_1(1_j)
+n\sum_{j=1}^n\Theta_j S_j^{(\ref{202005041155})}
\nn\\&& 
+n\sum_{j,k=1}^n\Psi_{j,k}(U_\infty,L_{\infty,\tjm})F_j^{(\ref{202005030358})}F_k^{(\ref{202005030358})}
\nn\\&& 
+n\sum_{j,k,\ell=1}^n\Xi_{j,k,\ell}(U_\infty,L_{\infty,\tjm})\int_\tjm^\tj\beta_tdt F_k^{(\ref{202005030358})}F_\ell^{(\ref{202005030358})}
+R_n^{(\ref{202005110436})},
\eea
}
\end{en-text}
where 
$R_n^{(\ref{202005110436})}\in\bbD^\infty$ satisfying 
\bea\label{202005110436}
R_n^{(\ref{202005110436})}
&=&
O_{\bbD^\infty}(n^{-0.5}). 
\eea
\begin{en-text}
\beas 
N_n 
&=&
n\sum_j \big(D_{1_j}(\ol{\theta}_\tjm(b^{[1]}_\tjm)^2)\big) I_1(1_j)
\nn\\&& 
+n\sum_{j=1}^n\ol{\theta}_\tjm S_j^{(\ref{202005041155})}
\nn\\&& 
+n\sum_{j,k=1}^n\Psi_{j,k}(U_\infty,L_{\infty,\tjm})F_j^{(\ref{202005030358})}F_k^{(\ref{202005030358})}
+n\sum_{j,k=1}^n\Psi_{j,k}(U_\infty,L_{\infty,\tjm})F_j^{(\ref{202005030358})}S_k^{(\ref{202005041155})}
\nn\\&&
+n\sum_{j,k,\ell=1}^n\Xi_{j,k,\ell}(U_\infty,L_{\infty,\tjm})F_j^{(\ref{202005030358})}F_k^{(\ref{202005030358})}F_\ell^{(\ref{202005030358})}
\nn\\&& 
+n\sum_{j,k=1}^n\Psi_{j,k}(U_\infty,L_{\infty,\tjm})\int_\tjm^\tj\beta_tdt F_k^{(\ref{202005030358})}
+n\sum_{j,k=1}^n\Psi_{j,k}(U_\infty,L_{\infty,\tjm})\int_\tjm^\tj\beta_tdt S_k^{(\ref{202005041155})}
\nn\\&&
+n\sum_{j,k,\ell=1}^n\Xi_{j,k,\ell}(U_\infty,L_{\infty,\tjm})\int_\tjm^\tj\beta_tdt F_k^{(\ref{202005030358})}F_\ell^{(\ref{202005030358})}
+nR_n^{(\ref{202005050041})}.
\eeas
\end{en-text}
\end{lemma}
\halflineskip
Proof of Lemma \ref{202005051550} is in Section \ref{202005141700}. 

\subsection{Asymptotic expansion of $Z_n$}
We can now apply the results in Section \ref{202004021507}. 
However, the asymptotic expansion term $N_n$ consists of two terms of different nature. 
Because of this, it is {\sred useful} to treat them in different ways. 
%
%
Based on the stochastic expansion (\ref{202005050001}), 
we consider the following decomposition of $Z_n$: 
\bea\label{202005120237}
Z_n &=& Z_n^{\sf o} +n^{-1/2}N_n^{\sf x}, 
\eea
with
\bea\label{202005120238}
Z_n^{\sf o} &=& M_n+n^{-1/2}N_n^{\sf o}
\eea
where
\bea\label{202005120741} 
N_n^{\sf o}
&=&
n\sum_{j=1}^n \big(D_{1_j}({\colorr\Theta_j}
(b^{[1]}_\tjm)^2)\big) I_1(1_j)
+n\sum_{j=1}^n\Theta_j S_j^{(\ref{202005041155})}
+R_n^{(\ref{202005110436})}
\eea
and 
{\sred 
\bea\label{202005120742} 
N_n^{\sf x}
&=&
n\sum_{j,k=1}^n\Psi_{j,k}(U_\infty,L_{\infty,\tjm})F_j^{(\ref{202005030358})}F_k^{(\ref{202005030358})}
\eea
}
\begin{en-text}
\bea\label{202005120742} 
N_n^{\sf x}
&=&
n\sum_{j,k=1}^n\Psi_{j,k}(U_\infty,L_{\infty,\tjm})F_j^{(\ref{202005030358})}F_k^{(\ref{202005030358})}
\nn\\&& 
+n\sum_{j,k,\ell=1}^n\Xi_{j,k,\ell}(U_\infty,L_{\infty,\tjm})\int_\tjm^\tj\beta_tdt F_k^{(\ref{202005030358})}F_\ell^{(\ref{202005030358})}.
\eea
\end{en-text}
We will apply the results in Section \ref{202004021507} to $Z_n^{\sf o}$ of (\ref{202005120238}) 
to obtain asymptotic expansion for the law $\call\{(Z_n^{\sf o},X_1)\}$. 
Then the variable $Z_n$ is a perturbation of $Z_n^{\sf o}$ by ${\sred n^{-1/2}}N_n^{\sf x}$, as in (\ref{202005120237}). 
To incorporate the effect of $N_n^{\sf x}$, the perturbation method will be used.

Among several possible exposition, we only give a result on the asymptotic expansion for a measurable function $f$ 
in the existence of the reference variable $X_1$ though the situation can be easily generalized. 
A non-degeneracy condition for $G_\infty$ will be assumed because it can degenerate in general. 
On the other hand, the non-degeneracy of $X_1$ is quite standard. 
Proof of the following theorem is in Section \ref{202005141659}.
\begin{theorem}\label{202005110453}
Suppose that $G_\infty^{-1}$ and $\Delta_{X_1}^{-1}$ are in $L^\inftym$. Then 
there exists a (polynomial) random symbol ${\mathfrak A}(\tti\sfz,\tti\sfx)$ 
with coefficients in $\bbD^\infty$ 
such that 
\beas\label{202005141414}
\sup_{f\in\cale(M,\gamma)}
\bigg| E\big[f(Z_n,X_1)\big] 
-\int_{\bbR^2}f(z,x){\mathfrak p}_n(z,x)dzdx\bigg|
&=& 
o(r_n)
\eeas
as $n\to\infty$ for every $(M,\gamma)\in(0,\infty)^2$, 
where 
\beas
{\mathfrak p}_n(z,x) &=& 
E\bigg[\mathfrak{A}_n(\partial_z,\partial_x)^*
\bigg\{\phi(z;0,G_\infty)\delta_x(X_1)\bigg\}\bigg]. 
\eeas
and 
\beas 
{\mathfrak A}_n(\tti\sfz,\tti\sfx)
&=& 
1+n^{-1/2}{\mathfrak A}(\tti\sfz,\tti\sfx). 
\eeas
An expression of ${\mathfrak A}$ is presented in Section \ref{2020051437}.
\end{theorem}
}

\section{Action of the operators $D_{u_n(q)}$ and $D^i$}\label{202004270637}
\subsection{Exponent of a functional}\label{202004270705}
We will consider 
the functional $\cali_n$ defined by
\bea\label{202004251257}
\cali_n 
&=&
n^{\alpha}
\sum_{j_1,...,j_m=1}^nA_n(j_1,...,j_m)
I_{q_1}(f_{j_1}^{(1)})\cdots I_{q_m}(f_{j_m}^{(m)})
\eea
where $n\in\bbN$, $\alpha\in\bbR$, $m\in\bbN$ and $(q_1,...,q_m)\in\bbZ^m$. 
Suppose that 
$f_{j_i}^{(i)}=f_{j_i}^{(i),n}\in\mH^{\tilde{\otimes}q_j}$ such that 
$\text{supp}(f_{j_i}^{(i)})\subset I_{j_i}^{\>q_i}$ and that 
\bea\label{2020042559}
\max_{i=1,...,m}\>\sup_{n\in\bbN}\sup_{j_i=1,...,n}
\sup_{t_1,...,t_{q_i}\in[0,{\sred1}]}|f_{j_i}^{(i)}(t_1,...,t_{q_i})|&\leq&C_1
\eea
for some finite constant $C_1$. 
The following condition will be assumed for the coefficients $A_n(j_1,...,j_m)$ in this section. 
We are writing $\bbJ_n=\{1,...,n\}$. 
\bd\im[[S\!\!]] \label{202005131639}
$A_n(j_1,...,j_m)\in\bbD^\infty$ for all $j_1,...,j_m\in\bbJ_n$, $n\in\bbN$, and 
\beas 
 C(i,p) \>:=\>
\sup_{n\in\bbN}\sup_{j_1,...,j_m\in\bbJ_n}\sup_{t_1,...,t_i\in[0,1]}
\big\|D^i_{t_1,...,t_i}A_n(j_1,...,j_m)\big\|_p &<&\infty
\eeas
for every $i\in\bbZ_+$ and $p>1$.  
\ed

Let $\bar{q}(m)=q_1+\cdots+q_m$. 
For $(q_1,...,q_m)\in\bbZ^m$, define $\sfm_1(q_1,...,q_m)$ and $\sfm_0(q_1,...,q_m)$ by 
\beas 
\sfm_1(q_1,...,q_m) \yeq \#\big\{m'\in\{1,...,m\};\>q_{m'}>0\big\}
&\text{and}& 
\sfm_0(q_1,...,q_m) \yeq \#\big\{m'\in\{1,...,m\};\>q_{m'}=0\big\},
\eeas
respecrively. 
We define {\it the exponent} $e(\cali_n)$ of $\cali_n$ admitting an expression (\ref{202004051400}) by 
\beas 
 e(\cali_n) 
 &=& 
\alpha-\half \bigg(\bar{q}(m)- \sfm_1(q_1,...,q_m)\bigg)
+\sfm_0(q_1,...,q_m)-\infty1_{\{\min\{q_1,...,q_m\}<0\}}
\\&=&
\alpha-\half \bar{q}(m)+m-\half\sfm_1(q_1,...,q_m)-\infty1_{\{\min\{q_1,...,q_m\}<0\}}.
\eeas
The exponent is defined for the sequence $(\cali_n)_{n\in\bbN}$; 
in this sense, it is 
more rigorous to write $e((\cali_n)_{n\in\bbN})$ for $e(\cali_n)$, but we adopt simpler notation. 
For example, 
$
e(\cali_n) = 0
$
for $\cali_n$ admitting the representation (\ref{202004051400}) 
if $\alpha=2^{-1}(\bar{q}(m)-m)$ and $(q_1,...,q_m)\in\bbN^m$. 
The exponent $e(\cali_n)$ depends on an expression (\ref{202004251257}) of $\cali_n$. 
Often we will not write explicitly which expression is considered when writing an exponent 
if there is no fear of confusion. 
Theorem \ref{0110221242} ensures 
\bea\label{202004051415}
\cali_n &=& O_{L^\inftym}(n^{e(\cali_n)})
\eea
as $n\to\infty$ under $[S]$ 
if $(q_1,...,q_m)\in\bbZ_+^m$. 
In fact, 
if $q_1,...,q_k>0=q_{k+1}=\cdots=q_m$, for example, then 
\beas 
\|\cali_n\|_p 
&\leq&
n^\alpha\sum_{j_{k+1},...,j_m}
\bigg\|\sum_{j_1,...,j_k}
A_n(j_1,...,j_m)
I_{q_1}(f_{j_1}^{(1)})\cdots I_{q_m}(f_{j_m}^{(m)})
\bigg\|_p
\\&\leq&
n^{\alpha+(m-k)-2^{-1}(q_1+\cdots+q_k-k)}
\\&=&
O(n^{e(\cali_n)})
\eeas
for $p>1$. If $(q_1,...,q_m)$ involves a negative component, then 
$e(\cali_n)=-\infty$ and (\ref{202004051415}) is valid in the sense that 
\beas 
0 &=& O_{L^\inftym}(n^{-\infty}). 
\eeas
Thus, (\ref{202004051415}) holds for all $(q_1,...,q_m)\in\bbZ^m$. 

Let $\cals$ be the set of sequences $(\cali_n)_{n\in\bbN}$ 
such that each $\cali_n$ admits a representation {\sred (\ref{202004251257})} 
for some $\alpha\in\bbR$, 
$m\in\bbN$ and $(q_1,..,q_m)\in\bbZ^m$ (independent of $n$) 
and 
coefficients $A_n(j_1,...,j_m)$ satisfying $[S]$. 
Let $\call$ be the linear space generated by all finite linear combination of elements of $\cals$. 
For $\calj_n\in\call$, we define the exponent $e(\calj_n)$ of $\calj_n$ by 
the maximum exponent of the summands in $\calj_n$. 
We remark that each summand of $\calj_n$ has a specific representation to determine its exponent, 
and hence $e(\calj_n)$ depends on these representations.
Furthermore, we define the exponent $e(\calk_n)$ 
by 
\beas 
e(\calk_n)&=&\max_{i=1,...,\bar{i}}e(\calj_n^{(i)})
\eeas
for a sequence $(\calk_n)_{n\in\bbN}$ of variables that satisfies 
\beas 
|\calk_n|\leq\sum_{i=1}^{\bar{i}}|\calj_n^{(i)}|
\eeas 
for some sequences $(\calj_n^{(i)})_{n\in\bbN}\in\call$, $i\in\{1,...,\bar{i}\}$. 
The exponent $e(\calk_n)$ depends on the dominating sequences $(\calj_n^{(i)})_{n\in\bbN,\>i\in\{1,...,{\bar{i}}\}}$ 
and their representation. 
The set of such sequences $(\calk_n)_{n\in\bbN}$ is denoted by $\calm$. 
By definition, 
\bea\label{202004121326}
\calk_n &=& O_{L^\inftym}(n^{e(\calk_n)})
\eea
for $(\calk_n)_{n\in\bbN}\in\calm$.

\subsection{Operator $D_{u_n(q)}$}
We are interested in the effect of the action of the projector $D_{u_n}$ on the functional $\cali_n$. 
Though the following argument is easily applied to more general functions $f_j$ but 
we will specially pay attention to $\cali_n$ re-defined by 
\bea\label{202004051400}
\cali_n 
&=&
n^{\alpha}
\sum_{j_1,...,j_m=1}^nA_n(j_1,...,j_m)I_{q_1}(1_{j_1}^{\otimes q_1})\cdots I_{q_m}(1_{j_m}^{\otimes q_m})
\eea
where $n\in\bbN$, $\alpha\in\bbR$, $m\in\bbN$ and $(q_1,...,q_m)\in\bbZ^m$. 
By convention,
\beas 
I_{q}(1_j^{\otimes q}) 
&=&
\left\{\begin{array}{cl}
1&(q=0)\y
0&(q<0)
\end{array}\right.
\eeas

In the present problem, the product formula has the expression 
\bea\label{202003151619} 
I_{q_1}(1_j^{\otimes(q_1)}) I_{q_2}(1_j^{\otimes(q_2)}) \cdots I_{q_\mu}(1_j^{\otimes(q_\mu)})
&=&
\sum_{\nu\in\bbZ_+}c_\nu(q_1,q_2,...,q_\mu)
I_{\sum_{i=1}^\mu q_i-2\nu}\big(1_j^{\otimes(\sum_{i=1}^\mu q_i-2\nu)}\big)n^{-\nu}
\nn\\&&
\eea
for some constants $c_\nu(q_1,q_2,...,q_\mu)$ depending on $\nu\in\bbZ_+$ 
and $(q_1,q_2,...,q_\mu)\in\bbZ_+^\mu$. 
The function $(q_1,q_2,...,q_\mu)\mapsto c_\nu(q_1,q_2,...,q_\mu)$ is symmetric, and 
may take $0$; e.g., 
$c_\nu(q_1,q_2,...,q_\mu)=0$ whenever $\sum_{i=1}^\mu q_i-2\nu<0$. 
We have 
\beas 
c_\nu(q_1,q_2) &=& c_\nu(q_1,q_2,0)\yeq
1_{\{0\leq\nu\leq (q_1\wedge q_2)\}}\nu!
\left(\begin{array}{c}q_1\\ \nu\end{array}\right)
\left(\begin{array}{c}q_2\\ \nu\end{array}\right)
\eeas
for $q_1,q_2\in\bbN$, and 
\bea\label{202004111042}
c_\nu(q_1,q_2,q_3) 
&=&
\sum_{\nu_1\in\bbZ_+}
\bigg[
1_{\{0\leq\nu_1\leq (q_2\wedge q_3)\}}1_{\{\nu_1\leq\nu\leq((q_1+\nu_1)\wedge(q_2+q_3-\nu_1))\}}
\nn\\&&\hspace{50pt}\times
\nu_1!
\left(\begin{array}{c}q_2\\ \nu_1\end{array}\right)
\left(\begin{array}{c}q_3\\ \nu_1\end{array}\right)
(\nu-\nu_1)!
\left(\begin{array}{c}q_1\\ \nu-\nu_1\end{array}\right)
\left(\begin{array}{c}q_2+q_3-2\nu_1\\  \nu-\nu_1\end{array}\right)\bigg]
\nn\\&&
\eea
for $q_1,q_2,q_3\in\bbZ_+$ (it is easy to verify that this three-factor product formula is valid for indices 
possibly being $0$). 
Let $c_x=0$ when $x\not\in\bbZ$.

If $c_\nu(q_1,q_2,q_3)\not=0$ for $\nu=\bar{q}(3)/2$, then 
$\bar{q}(3)$ is even and 
\bea\label{202004111043}
\nu\yeq \frac{q_1+q_2+q_3}{2} \leq \min\{(q_1+\nu_1),(q_2+q_3-\nu_1)\}
\eea
for some $\nu_1$ such that 
\bea\label{202004111045}
0\leq\nu_1\leq\min\{q_2,q_3\}.
\eea
Remark that if $(a+b)/2\leq\min\{a,b\}$ for two real numbers $a$ and $b$, then $a=b=(a+b)/2$. 
Therefore, 
\bea\label{202004111046} 
q_1+\nu_1\yeq q_2+q_3-\nu_1\yeq \frac{q_1+q_2+q_3}{2}, 
\eea
which implies 
\bea\label{202004111044}
\nu_1=\frac{-q_1+q_2+q_3}{2}.
\eea
From (\ref{202004111045}) and (\ref{202004111046}), we obtain 
\bea\label{202004111047}
q_1\yleq q_2+q_3,\quad q_2\yleq q_1+q_3,\quad q_3\yleq q_1+q_2. 
\eea
Let $\bbT$ be the set of $(q_1,q_2,q_3)\in\bbZ_+^3$ satisfying  
the ``triangular condition'' (\ref{202004111047}). 
For $\nu=\bar{q}(3)$ and $\nu_1$ given by (\ref{202004111044}), simple calculus gives 
\bea\label{202004111101}
c_{\bar{q}(3)/2}(q_1,q_2,q_3)
&=&
\frac{q_1!q_2!q_3!}
{\left(\frac{q_1+q_2-q_3}{2}\right)!\left(\frac{q_2+q_3-q_1}{2}\right)!\left(\frac{q_1+q_3-q_2}{2}\right)!}
\eea
from (\ref{202004111042}). 
In conclusion, $c_{\bar{q}(3)/2}(q_1,q_2,q_3)$ is possibly non-zero
only when $\bar{q}(3)$ is even and $(q_1,q_2,q_3)\in\bbT$, and then 
$c_{\bar{q}(3)/2}(q_1,q_2,q_3)$ is given by (\ref{202004111101}) 
(for $\nu=\bar{q}(3)/2$ and $\nu_1$ in (\ref{202004111044})).

The exponent $e(D_{u_n(q)}\cali_n)$ is estimated as follows. 
\begin{proposition}\label{202004091035}
Suppose that Conditions $[S]$ and (\ref{2020031201810a}) of $[A]$ are fulfilled. 
Suppose that $u_n(q)$ is given by (\ref{202003130047}) with some $q\in\{2,3,...\}$. 
Then, for $\cali_n$ with the representation (\ref{202004051400}), 
\bd
\im[(i)] $(D_{u_n(q)}\cali_n)_{n\in\bbN}\in\call$ and 
\beas
e(D_{u_n(q)}\cali_n) \yleq e(\cali_n).
\eeas
\im[(ii)] If there is no $m'\in\{1,...,m\}$ such that $q_{m'}=q$, then 
\beas 
e(D_{u_n(q)}\cali_n) &\leq& e(\cali_n)-\half.
\eeas
\ed
\end{proposition}
\proof 
Recall $u_n(q)$ defined in (\ref{202003130047}): 
\beas
u_n(q)
&=& 
n^{2^{-1}(q-1)}\sum_{j=1}^na_j(q) I_{q-1}(1_j^{\otimes (q-1)})1_j
\eeas
where $q\in\{2,3,...\}$. 
Let $\check{q}(m)=\bar{q}(m)+q$. 
We have
\bea&&
D_{u_n(q)}\cali_n
\nn\\&=&
n^{\alpha+0.5q-1.5}
\sum_{j_1,...,j_m,j_{m+1}=1}^n
\big(nD_{1_{j_{m+1}}}A_n(j_1,...,j_m)\big)a_{j_{m+1}}(q)
I_{q_1}(1_{j_1}^{\otimes q_1})\cdots I_{q_m}(1_{j_m}^{\otimes q_m})
I_{q-1}(1_{j_{m+1}}^{\otimes(q-1)})
\nn\\&&
+
\sum_{m'=1}^m
n^{\alpha+0.5q-1.5}\sum_{j_1,...,j_m}^nA_n(j_1,...,j_m)a_{j_{m'}}(q)
q_{m'}I_{q_1}(1_{j_1}^{\otimes q_1})\cdots I_{q_{m'-1}}(1_{j_{m'-1}}^{\otimes q_{m'-1}})
\nn\\&&\hspace{150pt}\times
I_{q_{m'}-1}(1_{j_{m'}}^{\otimes(q_{m'}-1)})I_{q-1}(1_{j_{m'}}^{\otimes(q-1)})
\nn\\&&\hspace{150pt}\times
I_{q_{m'+1}}(1_{j_{m'+1}}^{\otimes q_{m'+1}})\cdots I_{q_m}(1_{j_m}^{\otimes q_m})
\nn\\&=:&
\cali_n(1)+\sum_{m'=1}^m\cali_n(2,m').
\eea
We notice that 
\beas 
\big\|nD_{1_{j_{m+1}}}A_n(j_1,...,j_m)\big\|_p
&=&
\bigg\|n\int_0^1 D_tA_n(j_1,...,j_m)1_{j_{m+1}}(t)dt\bigg\|_p
\\&\leq&
n\int_0^1 \big\|D_tA_n(j_1,...,j_m)\big\|_p1_{j_{m+1}}(t)dt
\\&\leq&
\sup_{t\in[0,1]}\big\|D_tA_n(j_1,...,j_m)\big\|_p
\eeas
for every $p>1$. 
Therefore $(\cali_n(1))_{n\in\bbN}\in\cals$ and 
\bea\label{202004051850}
e(\cali_n(1)) 
&=& 
(\alpha+0.5q-1.5)-0.5(\bar{q}(m)+q-1)+(m+1)
\nn\\&&
-0.5\sfm_1(q_1,...,q_m,q-1)-\infty1_{\{\min\{q_1,..,q_m,q-1\}<0\}}
\nn\\&=& 
\alpha-0.5\bar{q}(m)+m-0.5\sfm_1(q_1,...,q_m)-\infty1_{\{\min\{q_1,..,q_m\}<0\}}-0.5
\nn\\&=&
e(\cali_n)-0.5.
\eea
\begin{en-text}
and hence 
\bea\label{202004051452}
\cali_n(1) &=& O_{L^\inftym}(n^{e(\cali_n)-0.5})
\eea
by (\ref{202004051415}). 
\end{en-text}

If $q_{m'}\leq0$, then $e(\cali_n(2,m'))=-\infty$. 
If $q_{m'}=1$, then 
for the expression 
\beas 
\cali_n(2,m')
&=&
n^{\alpha+0.5q-1.5}\sum_{j_1,...,j_m}A_n(j_1,...,j_m)
q_{m'}I_{q_1}(1_{j_1}^{\otimes q_1})\cdots I_{q_{m'-1}}(1_{j_{m'-1}}^{\otimes q_{m'-1}})
\nn\\&&\hspace{150pt}\times
I_{q-1}(1_{j_{m'}}^{\otimes(q-1)})
I_{q_{m'+1}}(1_{j_{m'+1}}^{\otimes q_{m'+1}})\cdots I_{q_m}(1_{j_m}^{\otimes q_m})
\eeas
of $\cali_n(2,m')$, we may say 
\bea
e(\cali_n(2,m'))
&=& 
(\alpha+0.5q-1.5)-0.5\big(\bar{q}(m)-1+q-1\big)
+m
\nn\\&&
-0.5\sfm_1(q_1,...,q_{m'-1},q-1,q_{m'+1},...,q_m)
-\infty1_{\{\min\{q_1,...,q_{m'-1},q-1,q_{m'+1},...,q_m)\}<0\}}
\nn\\&=&
(\alpha+0.5q-1.5)-0.5\big(\bar{q}(m)-1+q-1\big)
+m
\nn\\&&
-0.5\sfm_1(q_1,...,q_{m'-1},q_{m'},q_{m'+1},...,q_m)
-\infty1_{\{\min\{q_1,...,q_{m'-1},q_{m'},q_{m'+1},...,q_m)\}<0\}}
\nn\\&=&
e(\cali_n)-0.5.
\eea

If $q_{m'}\geq2$, then we apply the product formula (\ref{202003151619}) on p.\pageref{202003151619} 
to obtain 
\beas 
I_{q_{m'}-1}(1_{j_{m'}}^{\otimes(q_{m'}-1)})I_{q-1}(1_{j_{m'}}^{\otimes(q-1)})
&=&
\sum_{\nu\in\bbZ_+}c_\nu(q_{m'}-1,q-1)I_{q_{m'}+q-2-2\nu}\big(1_{j_{m'}}^{\otimes(q_{m'}+q-2-2\nu)}\big)
n^{-\nu}
\eeas
Therefore, 
\bea\label{202004051452}
\cali_n(2,m')
&=& 
\sum_{\nu\in\bbZ_+}c_\nu(q_{m'}-1,q-1)\cali_n(2,m',\nu)
\eea
where 
\bea\label{202004051455}
\cali_n(2,m',\nu)
&=&
n^{\alpha+0.5q-1.5-\nu}\sum_{j_1,...,j_m}A_n(j_1,...,j_m)
{\sred a_{j_{m'}}(q)}
q_{m'}I_{q_1}(1_{j_1}^{\otimes q_1})\cdots I_{q_{m'-1}}(1_{j_{m'-1}}^{\otimes q_{m'-1}})
\nn\\&&\hspace{150pt}\times
I_{q_{m'}+q-2-2\nu}\big(1_{j_{m'}}^{\otimes(q_{m'}+q-2-2\nu)}\big)
\nn\\&&\hspace{150pt}\times
I_{q_{m'+1}}(1_{j_{m'+1}}^{\otimes q_{m'+1}})\cdots I_{q_m}(1_{j_m}^{\otimes q_m})
\eea

If $q_{m'}\geq2$ and $q_{m'}=q$, then 
\bea
e\big(\cali_n(2,m',\nu)\big)
&=&
(\alpha+0.5q-1.5-\nu)-0.5\big(\bar{q}(m)+q-2-2\nu\big)+m
\nn\\&& 
-0.5\sfm_1(q_1,...,q_{m'-1},q_{m'}+q-2-2\nu,q_{m'+1},...,q_m)
\nn\\&&
-\infty1_{\{\min\{q_1,...,q_{m'-1},q_{m'}+q-2-2\nu,q_{m'+1},...,q_m\}<0\}}. 
\eea
In this case, for $\nu=q-1$, 
\bea
e\big(\cali_n\big(2,m',2^{-1}(q_{m'}+q)-1\big)\big)
&=&
(\alpha-0.5q-0.5)-0.5\big(\bar{q}(m)-q\big)+m
\nn\\&& 
-0.5\big(\sfm_1(q_1,...,q_{m'-1},q_{m'},q_{m'+1},...,q_m)-1\big)
\nn\\&&
-\infty1_{\{\min\{q_1,...,q_{m'-1},q_{m'},q_{m'+1},...,q_m\}<0\}}
\nn\\&=&
e(\cali_n)
\eea
So, if a term satisfying $q_{m'}\geq2$ and $q_{m'}=q$ appears, then 
in general 
$e(D_{u_n}\cali_n))\leq e(\cali_n)$, according to the other estimates (\ref{202004061607}) below.

If $q_{m'}\geq2$ and $q_{m'}\not=q$, then $q_{m'}+q-2-2\nu$ never attains $0$. 
Indeed, the terms for $\nu\leq\min\{q_{m'}-1,q-1\}$ can only remain 
in the product formula, and this implies $\nu$ can be $2^{-1}(q_{m'}+q)-1$ 
(equivalently, $q_{m'}+q-2-2\nu=0$) 
among such $\nu$s only when $q_{m'}=q$. 
We may consider $\nu$ such that $q_{m'}+q-2-2\nu>0$ for the least favorable case; 
otherwise $q_{m'}+q-2-2\nu<0$ and then the exponent is $-\infty$. 
In this situation, 
\bea\label{202004061607}
\cali_n(2,m',\nu)
&=&
(\alpha+0.5q-1.5-\nu)-0.5\big(\bar{q}(m)+q-2-2\nu\big)+m
\nn\\&& 
-0.5\sfm_1(q_1,...,q_{m'-1},q_{m'}+q-2-2\nu,q_{m'+1},...,q_m)
\nn\\&&
-\infty1_{\{\min\{q_1,...,q_{m'-1},q_{m'}+q-2-2\nu,q_{m'+1},...,q_m\}<0\}}
\nn\\&=&
(\alpha+0.5q-1.5-\nu)-0.5\big(\bar{q}(m)+q-2-2\nu\big)+m
\nn\\&& 
-0.5\sfm_1(q_1,...,q_{m'-1},q_{m'},q_{m'+1},...,q_m)
\nn\\&&
-\infty1_{\{\min\{q_1,...,q_{m'-1},q_{m'},q_{m'+1},...,q_m\}<0\}}
\nn\\&=&
e(\cali_n)-0.5.
\eea
\qed\halflineskip

\begin{en-text}
\bea\label{202004051553}
\cali_n(2,m',2^{-1}(q_{m'}+q)-1)
&=& 
\sum_{j_{m'}} \cali_n(2,m',2^{-1}(q_{m'}+q)-1)_{j_{m'}}
\eea
where 
\bea\label{202004051500}&&
\cali_n(2,m',2^{-1}(q_{m'}+q)-1))_{j_{m'}}
\nn\\&=&
n^{2^{-1}\check{q}(m)-2-(2^{-1}(q_{m'}+q)-1)}\sum_{j_1,...,j_{m'-1},j_{m'+1},...,j_m}^nA_n(j_1,...,j_m)
q_{m'}I_{q_1}(1_{j_1}^{\otimes q_1})\cdots I_{q_{m'-1}}(1_{j_{m'-1}}^{\otimes q_{m'-1}})
\nn\\&&\hspace{150pt}\times
I_{q_{m'+1}}(1_{j_{m'+1}}^{\otimes q_{m'+1}})\cdots I_{q_m}(1_{j_m}^{\otimes q_m}), 
\eea
and hence 
\bea\label{202004051503}
e\big(\cali_n(2,m',2^{-1}(q_{m'}+q)-1)_{j_{m'}}\big)
&=&
0.5m-1.5. 
\eea
for $j_{m'}\in\bbJ_n$. 

If $q_{m'}\geq2$ and $q_{m'}\not=q$, then 
$q_{m'}+q-2-2\nu>0$ if 
$I_{q_{m'}+q-2-2\nu}\big(1_{j_{m'}}^{\otimes(q_{m'}+q-2-2\nu)}\big)$ remain. 
Indeed, the terms for $\nu\leq\min\{q_{m'}-1,q-1\}$ can only remain 
in the product formula, and this implies $\nu$ can be $2^{-1}(q_{m'}+q)-1$ 
(equivalently, $q_{m'}+q-2-2\nu=0$) 
among such $\nu$s only when $q_{m'}=q$. 
Therefore, we conclude 
\bea\label{202004051536}
e\big(\cali_n(2,m',\nu)\big)
&=&
\big(0.5\check{q}(m)-2-\nu\big)-0.5\big(\check{q}(m)-2-2\nu-m)
\nn\\&=&
0.5m-1
\eea
from the exression (\ref{202004051455}). 
\end{en-text}

\subsection{Operator $D^i$}

Recall the sequence of functionals $\cali_n$ admitting the respresentation (\ref{202004051400}): 
\beas
\cali_n 
&=&
n^{\alpha}
\sum_{j_1,...,j_m=1}^nA_n(j_1,...,j_m)I_{q_1}(1_{j_1}^{\otimes q_1})\cdots I_{q_m}(1_{j_m}^{\otimes q_m})
\eeas
where $n\in\bbN$, $\alpha\in\bbR$, $m\in\bbN$ and $(q_1,...,q_m)\in\bbZ^m$. 
In this section, we will asses the effect of the operator $D^i$ applied to $\cali_n$. 
\begin{proposition}\label{202004121624}
Suppose that Conditions $[S]$ and (\ref{2020031201810a}) of $[A]$ are fulfilled. 
Then, for $\cali_n$ with the representation (\ref{202004051400}) and $i\in\bbZ_+$, 
\beas 
\big|D^i\cali_n\big|_{\mH^{\otimes i}} &=& O_{L^\inftym}(n^{e(\cali_n)})
\eeas
as $n\to\infty$. 
\end{proposition}
\proof
{\sred 
We may assume that $q_1,...,q_m\geq0$ 
since otherwise the claim is true due to $\cali_n=0$. 
}
We have 
\bea\label{202004061653}
|D\cali_n|_\mH^2
&\leq& 
\calj^{(1)}_n+\calj^{(2)}_n
\eea
where
\beas
\calj^{(1)}_n
&=&
2^mn^{2\alpha}
\sum_{j_1,...,j_m,\atop j_1',...,j_m'=1}^n
\langle DA_n(j_1,...,j_m),DA_n(j_1',...,j_m')\rangle
\nn\\&&\hspace{80pt}\times
I_{q_1}(1_{j_1}^{\otimes q_1})\cdots I_{q_m}(1_{j_m}^{\otimes q_m})
I_{q_1}(1_{j_1'}^{\otimes q_1})\cdots I_{q_m}(1_{j_m'}^{\otimes q_m})
\eeas
and 
\beas 
\calj^{(2)}_n
&=&
2^mn^{2\alpha-1}\sum_{m'=1}^m q_{m'}^2
\sum_{{j_1,...,j_m,\atop j_1',...,j_m'=1}\atop j_{m'}=j_{m'}'}^n
A_n(j_1,...,j_m)A_n(j_1',...,j_m')
\nn\\&&\hspace{150pt}\times
I_{q_1}(1_{j_1}^{\otimes q_1})\cdots I_{q_{m'-1}}(1_{j_{m'-1}}^{\otimes q_m})
I_{q_{m'+1}}(1_{j_{m'+1}}^{\otimes q_1})\cdots I_{q_m}(1_{j_m}^{\otimes q_m})
\nn\\&&\hspace{150pt}\times
I_{q_1}(1_{j_1'}^{\otimes q_1})\cdots I_{q_{m'-1}}(1_{j_{m'-1}'}^{\otimes q_m})
I_{q_{m'+1}}(1_{j_{m'+1}'}^{\otimes q_1})\cdots I_{q_m}(1_{j_m'}^{\otimes q_m})
\nn\\&&\hspace{150pt}\times
I_{q_{m'}-1}(1_{j_{m'}}^{\otimes(q_{m'}-1)})^2
\eeas
Then 
\bea\label{202004121456}
e(\calj^{(1)}_n) 
&=& 
2\alpha-\bar{q}(m)+2m-\sfm_1(q_1,...,q_m)-\infty1_{\{\min\{q_1,...,q_m\}<0\}}
\yeq 2e(\cali_n).
\eea
For $\calj_n^{(2)}$, we apply the product formula 
\beas 
I_{q_{m'}-1}(1_{j_{m'}}^{\otimes(q_{m'}-1)})^2
&=& 
\sum_\nu c_\nu(q_{m'}-1,q_{m'}-1)I_{2q_{m'}-2-2\nu}\big(1_{j_{m'}}^{\otimes(2q_{m'}-2-2\nu)}\big)n^{-\nu}
\eeas
to obtain 
\bea\label{202004121457}
e(\calj^{(2)}_n) 
&\leq&
2e(\cali_n)
\eea
since
\begin{en-text}
\beas &&
\big(2\alpha-1-\nu\big)-0.5\big(2\bar{q}(m)-2-2\nu\big)+(2m-1)
\nn\\&&
-0.5\times2\sfm_1(q_1,...,q_{m'-1},q_{m'}-1,q_{m'+1},...,q_m)
-\infty1_{\{\min\{q_1,...,q_{m'-1},q_{m'}-1,q_{m'+1},...,q_m\}<0\}}
\nn\\&=&
2e(\cali_n)+1_{\{q_{m'}=1\}}-\infty1_{\{q_{m'}\leq0\}}-1
\nn\\&\leq&
2e(\cali_n). 
\eeas
\end{en-text}
{\sred 
\beas &&
\big(2\alpha-1-\nu\big)-0.5\big(2\bar{q}(m)-2-2\nu\big)+(2m-1)
\nn\\&&
-0.5\times\sfm_1(q_1,...,q_{m'-1},q_{m'+1},...,q_m,q_1,...,q_{m'-1},q_{m'+1},...,q_m,
2q_{m'}-2-2\nu
)
\nn\\&&
-\infty1_{\{\min\{q_1,...,q_{m'-1},2q_{m'}-2-2\nu,q_{m'+1},...,q_m\}<0\}}
\nn\\&\leq&
2\alpha-2\times0.5\bar{q}(m)+2m-1
\nn\\&&
-0.5\times\sfm_1(q_1,...,q_{m'-1},q_{m'+1},...,q_m,q_1,...,q_{m'-1},q_{m'+1},...,q_m)
\nn\\&=&
2\alpha-2\times0.5\bar{q}(m)+2m-1
-0.5\times\big(2\sfm_1(q_1,...,q_m)-21_{\{q_{m'}>0\}}\big)
\nn\\&\leq&
2e(\cali_n)
\eeas
due to $\min\{q_1,...,q_m\}\geq0$ by assumption. 
}
From (\ref{202004061653}), (\ref{202004121456}) and (\ref{202004121457}), 
we obtain
\beas 
e\big(|D\cali_n|_\mH^2\big) &\leq& 2e(\cali_n)
\eeas
and hence 
\beas
|D\cali_n|_\mH &=& O(n^{e(\cali_n)}). 
\eeas

The derivative $D^s\cali_n$ for ${\sred s}\geq2$ can be estimated in a similar manner. 
We consider 
\beas 
\calk_n &=& 
\bigg|
n^{\alpha}
\sum_{j_1,...,j_m=1}^nD^kA_n(j_1,...,j_m)
I_{p_1}(1_{j_1}^{\otimes p_1})\cdots I_{p_m}(1_{j_m}^{\otimes p_m})
\odot1_{j_1}^{\otimes(q_1-p_1)}\odot\cdots\odot1_{j_m}^{\otimes(q_m-p_m)}
\bigg|_{\mH^{\otimes s}}^2
\eeas
for $k,p_1,...,p_m\in\bbZ_+$ such that $p_i\leq q_i$ ($i=1,...,m$) and 
$k+\sum_{i=1}^m(q_i-p_i)=s$. 
In general, $p_{\ell+1}=\cdots=p_m=0$, and in this situation, we may consider 
\begin{en-text}
\beas
\calk_n' &=& 
\bigg|
n^{\alpha}
\sum_{j_1,...,j_m=1}^nD^kA_n(j_1,...,j_m)I_{p_1}(1_{j_1}^{\otimes p_1})\cdots I_{p_\ell}(1_{j_\ell}^{\otimes p_\ell})
\nn\\&&\hspace{150pt}\otimes
1_{j_1}^{\otimes(q_1-p_1)}\otimes\cdots\otimes
1_{j_1}^{\otimes(q_\ell-p_\ell)}\otimes1_{j_1}^{\otimes q_{\ell+1} }
\otimes\cdots\otimes1_{j_m}^{\otimes q_m}
\bigg|_{\mH^{\otimes s}}^2
\eeas
where $p_1,...,p_\ell>0$. 
Now 
\bea\label{202004121617} 
\calk_n' &=& 
n^{2\alpha+(m-\ell)-(\bar{q}(\ell)-\bar{p}(\ell))-(q_{\ell+1}+\cdots+q_m)}
\nn\\&&\times
\sum_{j_1,...,j_\ell,j_1',...,j_\ell'}
\bigg\{n^{-(m-\ell)}\sum_{j_{\ell+1},...,j_m}
\langle D^kA_n(j_1,...,j_m),D^kA_n(j_1',...,j_m')\rangle_{\mH^{\otimes k}}\bigg\}
\nn\\&&\times
I_{p_1}(1_{j_1}^{\otimes p_1})\cdots I_{p_\ell}(1_{j_\ell}^{\otimes p_\ell})
I_{p_1}(1_{j_1'}^{\otimes p_1})\cdots I_{p_\ell}(1_{j_\ell'}^{\otimes p_\ell})
\eea
For the representation (\ref{202004121617}) of $\calk_n'$, we have 
\beas
e(\calk_n')
&=& 
\big\{2\alpha+(m-\ell)-(\bar{q}(\ell)-\bar{p}(\ell))-(q_{\ell+1}+\cdots+q_m)\big\}
-0.5\times 2\bar{p}(\ell)+2\ell-0.5\times 2\ell
\\&=&
2\alpha+m-\bar{q}(m)
\yeq 
2e(\cali_n)
\eeas
\end{en-text}
{\sred %
\beas
\calk_n' &=& 
\bigg|
n^{\alpha}
\sum_{j_1,...,j_m=1}^nD^kA_n(j_1,...,j_m)I_{p_1}(1_{j_1}^{\otimes p_1})\cdots I_{p_\ell}(1_{j_\ell}^{\otimes p_\ell})
I_{q_{\ell+1}}(1_{j_{\ell+1}}^{\otimes {q_{\ell+1}}})\cdots
I_{q_\kappa}(1_{j_\kappa}^{\otimes {q_\kappa}})
\nn\\&&\hspace{150pt}\otimes
1_{j_1}^{\otimes(q_1-p_1)}\otimes\cdots\otimes
1_{j_\ell}^{\otimes(q_\ell-p_\ell)}
\bigg|_{\mH^{\otimes s}}^2
\eeas
where $p_1,...,p_\ell\geq0$. 
This expression causes, when $q_1,...,q_\ell,q_{\ell+1},...,q_\kappa>0$ 
and $q_{\kappa+1}=\cdots=q_m=0$, 
from the $q_i-p_i(>0)$ times derivative of 
$I_{p_1}(1_{j_1}^{\otimes p_1}),....,I_{p_m}(1_{j_m}^{\otimes p_\ell})$ 
and no derivative of 
$I_{q_{\ell+1}}(1_{j_{\ell+1}}^{\otimes q_{\ell+1}}),....,I_{q_\kappa}(1_{j_\kappa}^{\otimes q_\kappa})$ with $q_{\ell+1},...,q_\kappa>0$. 
%
Now 
\bea\label{202004121617} 
\calk_n' &=& 
n^{2\alpha-\sum_{i=1}^\ell(q_i-p_i)+2(m-\kappa)}
\nn\\&&\times
\sum_{j_1,...,j_\ell\atop
j_{\ell+1},...,j_{\kappa},j_{\ell+1}',...,j_{\kappa}'}
\bigg\{n^{-2(m-\kappa)}\sum_{j_{\kappa+1},...,j_m\atop j_{\kappa+1}',...,j_m'}
\langle D^kA_n(j_1,...,j_m),D^kA_n(j_1',...,j_m')\rangle_{\mH^{\otimes k}}\bigg\}
\nn\\&&\times
I_{p_1}(1_{j_1}^{\otimes p_1})^2\cdots I_{p_\ell}(1_{j_\ell}^{\otimes p_\ell})^2\>
I_{q_{\ell+1}}(1_{j_{\ell+1}}^{\otimes {q_{\ell+1}}})\cdots I_{q_\kappa}(1_{j_\kappa}^{\otimes q_\kappa})
I_{q_{\ell+1}}(1_{j_{\ell+1}'}^{\otimes q_1})\cdots I_{q_\ell}(1_{j_\kappa'}^{\otimes q_\kappa})
\eea
where $j_1=j_1'$, ..., $j_\ell=j_\ell'$. 
The factor $I_{p_1}(1_{j_1}^{\otimes p_1})^2\cdots I_{p_\ell}(1_{j_\ell}^{\otimes p_\ell})^2$ 
is a weighted sum of the terms 
\bea\label{202012291041}
I_{2p_1-2\nu_1}(1_{j_1}^{2p_1-2\nu_1})\cdots I_{2p_\ell-2\nu_\ell}(1_{j_\ell}^{2p_\ell-2\nu_\ell})
n^{-(\nu_1+\cdots+\nu_\ell)}
\eea
with 
$\nu_i\in\{0,...,p_i\}$ for $i=1,...,\ell$. 
Therefore, for the representation (\ref{202004121617}) of $\calk_n'$
expanded with (\ref{202012291041}), we have $e(\calk_n')$ 
as the maximum (in $\nu_i$'s) of 
\beas 
&& 
2\alpha-\sum_{i=1}^\ell(q_i-p_i)+2(m-\kappa)
-\sum_{i=1}^\ell\nu_i
\nn\\&&
-\half\bigg\{\sum_{i=1}^\ell(2p_i-2\nu_i)+\sum_{i=\ell+1}^\kappa 2q_i\bigg\}
+\big(\ell+2(\kappa-\ell)\big)
\nn\\&&
-\half\sfm_1\big(2p_1-2\nu_1,...,2p_\ell-2\nu_\ell,q_{\ell+1},....,q_\kappa,q_{\ell+1},....,q_\kappa\big)
\nn\\&\leq&
2\alpha-\sum_{i=1}^\kappa q_i+2m-\kappa
\nn\\&=&
2\alpha-\sum_{i=1}^m q_i+2m-\sfm_1(q_1,...,q_m)
\nn\\&=&
2e(\cali_n)
\eeas
}
Consequently, 
\beas
|D^i\cali_n|_{\mH^{\otimes i}}^2 &=& O(n^{2e(\cali_n)})
\eeas
for $i\in\bbZ_+$. 
\qed\halflineskip
%

{\sred 
\subsection{Stability of the class $[S]$}
\begin{proposition}\label{202012292110}
Suppose that $\cali_n$ admits an expression (\ref{202004251257}) 
with $A_n(j_1,...,j_m)$ satisfying Condition $[S]$. 
Then a version of density $D^i_{t_1,...,t_i}\cali_n$ satisfies 
\beas 
\sup_{n\in\bbN}\sup_{t_1,...,t_i\in[0,1]}\big\|D^i_{t_1,...,t_i}\cali_n\big\|_p
&\leq& 
K(i,p)n^{e(\cali_n)}
\eeas
for every $i\in\bbZ$ and $p>1$, 
where $K(i,p)$ is a constant determined by a finite number of $C(i,p)$'s in Condition $[S]$. 
\end{proposition}
\proof 
According to the discussion below, we may show the claim for $i=1$. 
We have 
\beas 
D_s\cali_n 
&=&
n^{\alpha}
\sum_{j_1,...,j_m=1}^nD_sA_n(j_1,...,j_m)
I_{q_1}(f_{j_1}^{(1)})\cdots I_{q_m}(f_{j_m}^{(m)})
\nn\\&&
+
n^{\alpha}
\sum_{j_2,...,j_m=1}^nA_n(j_1,...,j_m)q_1I_{q_1-1}(f_{j_1(s)}^{(1)})
I_{q_2}(f_{j_2}^{(2)})\cdots I_{q_m}(f_{j_m}^{(m)})+\cdots
\nn\\&=:&
\bbK_n^{(0)}+\bbK_n^{(1)}+\cdots+\bbK_n^{(m)}
\eeas
where $j_1(s)$ is $j_1$ such that $s\in I_{j_1}$. 
Obviously, $\sup_{t\in[0,1]}\sup_{n\in\bbN}\big\|\bbK_n^{(0)}\big\|_p=O(n^{e(\cali_n)})$ 
for every $p>1$ because $e(\cali_n)=e(\bbK_n^{(0)})$. 
For $\bbK_n^{(1)}$, the $L^p$-norm $\|\bbK_n^{(1)}\|_p$ is bounded by 
the product of 
$q_1\|I_{q_1-1}(f_{j_1(s)}^{(1)})\|_{2p}=O(n^{-0.5(q_1-1)})$ and 
\beas &&
\sup_{n\in\bbN}\sup_{j\in\bbJ_n}\bigg\|
n^{\alpha}
\sum_{j_2,...,j_m=1}^nA_n(j_1,...,j_m)
I_{q_2}(f_{j_2}^{(2)})\cdots I_{q_m}(f_{j_m}^{(m)})
\bigg\|_{2p}
\yeq
O(n^{\xi})
\eeas
for 
\beas 
\xi 
&=& 
\alpha-0.5(q_2+\cdots+q_m)+(m-1)-0.5\sfm_1(q_2,...,q_m).
\eeas
Therefore, 
\beas 
\big\|\bbK_n^{(1)}\big\|_p
&=& 
O(n^{\xi'})
\eeas
for 
\beas 
\xi' 
&=& 
\alpha-0.5(q_2+\cdots+q_m)+(m-1)-0.5\sfm_1(q_2,...,q_m)
-0.5(q_1-1)
\nn\\&\leq&
\alpha-0.5(q_1+\cdots+q_m)+m-0.5\sfm_1(q_1,...,q_m)
\nn\\&\leq&
e(\cali_n)
\eeas
We obtain the result since similar estimates hold for other terms and 
it is easy to see each estimate does not depend on $s\in[0,1]$. 
\qed\halflineskip

We can strengthen Proposition \ref{202004091035} (ii). 

\begin{corollary}\label{202012292154}
Under the conditions in Proposition \ref{202004091035}, 
if there is one $m'\in\{1,...m\}$ such that $q_{m'}\not=q$, then 
\beas 
e(D_{u_n(q)}\cali_n) &\leq& e(\cali_n)-\half.
\eeas
\end{corollary}
\proof 
We may suppose that $q_1,...,q_m\geq0$ and $q_1\not=q$. 
Then 
\bea\label{202012292207}
\cali_n &=& 
n^{\alpha-{\sf e}}\sum_{j_1}A'_n(j_1)I_{q_1}(1_{j_1}^{\otimes q_1})
\eea
for 
\beas 
A'_n(j_1)
&=&
n^{{\sf e}}
\sum_{j_2,...,j_m}A_n(j_1,j_2,...,j_m)
I_{q_2}(1_{j_2}^{\otimes q_2})\cdots I_{q_m}(1_{j_m}^{\otimes q_m}),
\eeas
where 
\beas 
{\sf e}
&=& 
e(A'_n(j_1))
\yeq 
-0.5(q_2+\cdots q_m)+(m-1)-0.5\sfm_1(q_2,...,q_m).
\eeas
By Proposition \ref{202012292110}, 
we see $A_n'(j_1)$ satisfies Condition $[S]$. 
Then, by Proposition \ref{202004091035} (ii), the representation (\ref{202012292207}) gives 
\beas 
e(D_{u_n(q)}\cali_n) 
&\leq& 
\big(\alpha+{\sf e}-0.5q_1+1-0.5\sfm_1(q_1)\big)-0.5
\nn\\&=&
e(\cali_n)-0.5.
\eeas
\qed\halflineskip
}

\section{Proof of Theorem \ref{0112121225}}\label{202003141657}
\begin{en-text}
\begin{comment}
\begin{itembox}{Comments}
$a_j(q)\in\bbD_{3,\inftym}$
\end{itembox}
\halflineskip
\end{comment}
\end{en-text}

\subsection{The first projection $\boldsymbol{D_{u_n}M_n}$}\label{202004211724}
For the first projection $D_{u_n(q_2)}M_n(q_1)=\langle DM_n(q_1),u_n(q_2)\rangle$, we have
\bea\label{20200317316}
D_{u_n(q_2)}M_n(q_1)
&=&
D_{u_n(q_2)}\bigg(n^{2^{-1}(q_1-1)}\sum_{j=1}^na_j(q_1) I_{q_1}(1_j^{\otimes q_1})\bigg)
\nn\\&&
-D_{u_n(q_2)}\bigg(n^{2^{-1}(q_1-1)}\sum_{j=1}^n(D_{1_j}a_j(q_1)) I_{q_1-1}(1_j^{\otimes (q_1-1)})\bigg)
\nn\\&=&
\bbI_1+\bbI_2+\bbI_3+\bbI_4
\eea
where 
\bea\label{202004041122} 
\bbI_1
&=&
n^{2^{-1}(q_1+q_2)-1}\sum_{j,k}q_1a_j(q_1)a_k(q_2)\langle 1_j,1_k\rangle I_{q_1-1}(1_j^{\otimes(q_1-1)}) I_{q_2-1}(1_k^{\otimes(q_2-1)}),
\eea
\bea\label{202004041123} 
\bbI_2
&=&
n^{2^{-1}(q_1+q_2)-1}\sum_{j,k}(D_{1_k}a_j(q_1)) a_k(q_2) I_{q_1}(1_j^{\otimes q_1})I_{q_2-1}(1_k^{\otimes(q_2-1)}),
\eea
\bea\label{202004041124} 
\bbI_3
&=&
-n^{2^{-1}(q_1+q_2)-1}\sum_{j,k}(q_1-1)(D_{1_j}a_j(q_1)) a_k(q_2) \langle 1_j,1_k\rangle 
I_{q_1-2}(1_j^{\otimes(q_1-2)})I_{q_2-1}(1_k^{\otimes(q_2-1)})
\eea
and
\bea\label{202004041125} 
\bbI_4
&=&
-n^{2^{-1}(q_1+q_2)-1}\sum_{j,k}(D_{1_k}D_{1_j}a_j(q_1)) a_k(q_2) I_{q_1-1}(1_j^{\otimes (q_1-1)})I_{q_2-1}(1_k^{\otimes (q_2-1)}).
\eea

\subsection{The second projection $\boldsymbol{D_{u_n}^2M_n}$
and the quasi-torsion}\label{202004220553}
In this section, we will investigate the quasi-torsion. 
The second projection $D_{u_n}^2M_n=D_{u_n}D_{u_n}M_n=\langle D\langle DM_n,u_n\rangle,u_n\rangle$ 
involves various patterns of terms. 
Let 
\bea\label{202003151547} 
\beta_{j,k} &=& \langle 1_j,1_k\rangle\yeq n^{-1}1_{\{j=k\}}, 
\eea
the Kronecker's delta, and simply denote 
\beas 
\bar{q}(3)&=&q_1+q_2+q_3
\eeas 
for $(q_1,q_2,q_3)\in\calq^3$. 
For random variables $V_n$ and positive numbers $b_n$, we write $V_n=O_{L^\inftym}(b_n)$ 
if $\|V_n\|_p=O(b_n)$ as $n\to\infty$ for every $p>1$. 
%
Let $\sfz\in\bbR$ and $\sfx\in\bbR$.

\subsubsection{$\boldsymbol{D_{u_n(q_3)}\bbI_4}$}
Recall the expression equivalent to (\ref{202004041125}): 
\beas 
\bbI_4&=&
-n^{2^{-1}(q_1+q_2)-3}\sum_{j,k}(n^2D_{1_k}D_{1_j}a_j(q_1)) a_k(q_2) I_{q_1-1}(1_j^{\otimes (q_1-1)})I_{q_2-1}(1_k^{\otimes (q_2-1)}). 
\eeas
The exponent of $\bbI_4$ by this representation is 
\bea\label{202004121725}
e(\bbI_4) 
&=& 
\big(0.5(q_1+q_2)-3\big)-0.5(q_1+q_2-2)+2-0.5\times2
\yeq
-1
\eea
since $q_1,q_2\geq2$. 
\begin{en-text}
We have 
\beas
D_{u_n(q_3)}\bbI_4
&=&
-n^{0.5(\bar{q}(3)-3)}\sum_{j,k,\ell}(q_1-1)(D_{1_k}D_{1_j}a_j(q_1))a_k(q_2)
\beta_{j,\ell}a_\ell(q_3)
\\&&\hspace{80pt}\times
I_{q_1-2}(1_j^{\otimes (q_1-2)})I_{q_2-1}(1_k^{\otimes (q_2-1)})I_{q_3-1}(1_\ell^{\otimes(q_3-1)})
\\&&
-n^{0.5(\bar{q}(3)-3)}\sum_{j,k,\ell}(q_2-1)(D_{1_k}D_{1_j}a_j(q_1))a_k(q_2)\beta_{k,\ell}a_\ell(q_3)
\\&&\hspace{80pt}\times
I_{q_1-1}(1_j^{\otimes (q_1-1)})I_{q_2-2}(1_k^{\otimes (q_2-2)})I_{q_3-1}(1_\ell^{\otimes(q_3-1)})
\\&&
-n^{0.5(\bar{q}(3)-3)}\sum_{j,k,\ell}(D_{1_\ell}\{(D_{1_k}D_{1_j}a_j(q_1)) a_k(q_2)\})a_\ell(q_3)
\\&&\hspace{80pt}\times
I_{q_1-1}(1_j^{\otimes (q_1-1)})I_{q_2-1}(1_k^{\otimes (q_2-1)})I_{q_3-1}(1_\ell^{\otimes(q_3-1)})
\\&=:&
\bbI_{4,1}+\bbI_{4,1}'+\bbI_{4,2}.
\eeas
Then, for the Brownian motion $w$, the exponent of $n$ of $\|\bbI_{4,1}\|_p$ is bounded by 
\beas 
0.5(\bar{q}(3)-3)-0.5(\bar{q}(3)-7)-3&=&-1
\eeas
for every $p>1$, 
where $0.5(3q(3)-7)="2^{-1}(\bar{q}-\sfm)"$ from Theorem \ref{0110221242}, 
 and the number $3$ comes from 
the maximum order of $(D_{1_k}D_{1_j}a_j)\beta_{j,\ell}$. 
Thus $\bbI_{4,1}=O_{L^\inftym}(n^{-1})$, and by the same reason $\bbI_{4,1}'=O_{L^\inftym}(n^{-1})$. 
Similarly, $\bbI_{4,2}=O_{L^\inftym}(n^{-1.5})$. 
\end{en-text}
Then Proposition \ref{202004091035} gives 
\bea\label{202004131414}
e(D_{u_n}\bbI_4) 
&\leq& 
-1
\eea
Therefore,  
\bea\label{202004041615} 
D_{u_n(q_3)}\bbI_4 &=& O_{L^\inftym}(n^{-1})
\eea
by (\ref{202004051415}). 
In particular, 
\bea\label{202003161330}
E[\Psi(\sfz,\sfx)D_{u_n(q_3)}\bbI_4]
&=&
O(n^{-1})
\eea
as $n\to\infty$ for every $\sfz\in\bbR$ and $\sfx\in\bbR$.

\subsubsection{$\boldsymbol{D_{u_n(q_3)}\bbI_1}$}
Recall 
\beas
\bbI_1 &=& 
n^{2^{-1}(q_1+q_2)-1}\sum_{j,k}q_1a_j(q_1)a_k(q_2)\beta_{j,k} I_{q_1-1}(1_j^{\otimes(q_1-1)}) I_{q_2-1}(1_k^{\otimes(q_2-1)})
\eeas
Then 
\bea\label{202003161223}
D_{u_n(q_3)}\bbI_1
&=&
n^{0.5(\bar{q}(3)-3)}\sum_{j,k,\ell}q_1(q_1-1)a_j(q_1)a_k(q_2)a_\ell(q_3)\beta_{j,k} \beta_{j,\ell}
\nn\\&&\hspace{80pt}\times
I_{q_1-2}(1_j^{\otimes(q_1-2)}) I_{q_2-1}(1_k^{\otimes(q_2-1)}) I_{q_3-1}(1_\ell^{\otimes(q_3-1)})
\nn\\&&
+n^{0.5(\bar{q}(3)-3)}\sum_{j,k,\ell}q_1(q_2-1)a_j(q_1)a_k(q_2)a_\ell(q_3)\beta_{j,k} \beta_{k,\ell}
\nn\\&&\hspace{80pt}\times
I_{q_1-1}(1_j^{\otimes(q_1-1)}) I_{q_2-2}(1_k^{\otimes(q_2-2)}) I_{q_3-1}(1_\ell^{\otimes(q_3-1)})
\nn\\&&
+n^{0.5(\bar{q}(3)-3)}\sum_{j,k,\ell}q_1(D_{1_\ell}(a_j(q_1)a_k(q_2)))a_\ell(q_3)\beta_{j,k}
\nn\\&&\hspace{80pt}\times
 I_{q_1-1}(1_j^{\otimes(q_1-1)}) I_{q_2-1}(1_k^{\otimes(q_2-1)}) I_{q_3-1}(1_\ell^{\otimes(q_3-1)})
\nn \\&=:&
 \bbI_{1,1}+  \bbI_{1,1}'+\bbI_{1,2}
\eea

In the case (\ref{202003151547}), 
\beas
\bbI_{1,1}
&=&
n^{0.5\bar{q}(3)-3.5}\sum_{j}q_1(q_1-1)a_j(q_1)a_j(q_2)a_j(q_3)
\\&&\hspace{80pt}\times
I_{q_1-2}(1_j^{\otimes(q_1-2)}) I_{q_2-1}(1_j^{\otimes(q_2-1)}) I_{q_3-1}(1_j^{\otimes(q_3-1)}).
\eeas
By (\ref{202003151619}), 
\bea\label{202004041456}
 \bbI_{1,1}
&=&
\bbI_{1,1,1}+\bbI_{1,1,2},
\eea
where 
\bea\label{202004131435}
\bbI_{1,1,1}
&=& 
n^{-1.5}\sum_jq_1(q_1-1)a_j(q_1)a_j(q_2)a_j(q_3) 
c_{0.5\bar{q}(3)-2}(q_1-2,q_2-1,q_3-1)
\eea
and 
\beas \bbI_{1,1,2}
&=&
\sum_{\nu:\>\bar{q}(3)-4-2\nu>0} c_{\nu}(q_1-2,q_2-1,q_3-1)
\nn\\&&\hspace{50pt}\times
n^{0.5\bar{q}(3)-3.5-\nu}\sum_jq_1(q_1-1)a_j(q_1)a_j(q_2)a_j(q_3) 
I_{\bar{q}(3)-4-2\nu}(1_j^{\otimes (\bar{q}(3)-4-2\nu)}).
\eeas
We have 
\bea\label{202004131436}
e(\bbI_{1,1,1})
&=& -1.5-0.5\times0+1-0.5\times0\yeq -0.5
\eea
and 
\beas 
e(\bbI_{1,1,2})
&=& 
\big(0.5\bar{q}(3)-3.5-\nu\big)-0.5\big(\bar{q}(3)-4-2\nu\big)
+1-0.5
\yeq -1.
\eeas
Thus, $\bbI_{1,1,1}$ determines the order of the symbol. 
\begin{en-text}
\beas 
E[\Psi(\sfz)\bbI_{1,1}]
&=&
q_1(q_1-1)c_{0.5\bar{q}(3)-2}(q_1-2,q_2-1,q_3-1)
\\&&\hspace{30pt}\times
n^{-1.5+0.5\ep(q)}\sum_{j=1}^n
E\big[D^{\ep(\bar{q}(3))}_{1_j}(\Psi(\sfz)a_j(q_1)a_j(q_2)a_j(q_3))\big]+O(n^{-1})
\eeas
where $\ep(q)=1_{\{q:\text{odd}\}}$. 
When $\bar{q}(3)$ is odd, this term has no effect. When $\bar{q}(3)$ is even, this term remains. 
Under the present assumption that all $q_1,q_2,q_3$ are even, 
this term remains. 
\end{en-text}
More precisely, 
\bea\label{202003161225}
E[\Psi(\sfz,\sfx)\bbI_{1,1}]
&=&
{\colorb 1_{\{\bar{q}(3):even\}}}
q_1(q_1-1)c_{0.5\bar{q}(3)-2}(q_1-2,q_2-1,q_3-1)
\nn\\&&\hspace{30pt}\times
n^{-1.5}\sum_{j=1}^n
E\big[\Psi(\sfz,\sfx)a_j(q_1)a_j(q_2)a_j(q_3)\big]+O(n^{-1}). 
\eea

Similarly, 
\beas 
\bbI_{1,1}'
&=&
\bbI_{1,1,1}'+\bbI_{1,1,2}'
\eeas
where
\bea\label{202004131443} 
\bbI_{1,1,1}'
&=&
n^{-1.5}\sum_jq_1(q_2-1)a_j(q_1)a_j(q_2)a_j(q_3) 
c_{0.5\bar{q}(3)-2}(q_1-1,q_2-2,q_3-1)
\eea
and 
$\bbI_{1,1,2}'\in\call$. 
Moreover, 
\bea\label{202004131444}
e(\bbI_{1,1,1}') &=& -0.5
\eea
and 
\beas 
e(\bbI_{1,1,2}') &=& -1
\eeas
In particular, since $e(\bbI_{1,1,2}')=O_{L^\inftym}(n^{-1/2})$, 
\bea\label{202003161226}
E[\Psi(\sfz,\sfx)\bbI_{1,1}']
&=&
{\colorb 1_{\{\bar{q}(3):even\}}}
q_1(q_2-1)c_{0.5\bar{q}(3)-2}(q_1-1,q_2-2,q_3-1)
\nn\\&&\hspace{30pt}\times
n^{-1.5}\sum_{j=1}^n
E\big[\Psi(\sfz,\sfx)a_j(q_1)a_j(q_2)a_j(q_3)\big]+O(n^{-1}). 
\eea

To estimate $E[\Psi(\sfz,\sfx)\bbI_{1,2}]$, 
we apply the product formula to $I_{q-1}(1_j^{\otimes(q-1)})^2$ first and obtain 
\begin{en-text}
\beas
\bbI_{1,2}
&=&
n^{0.5(\bar{q}(3)-3)}\sum_{j,\ell}q_1(D_{1_\ell}(a_j(q_1)a_j(q_2)))a_\ell(q_3)n^{-1}
\\&&\hspace{80pt}\times
{\colorr I_{q_1-1}(1_j^{\otimes(q_1-1)}) I_{q_2-1}(1_j^{\otimes(q_2-1)}) }
I_{q_3-1}(1_\ell^{\otimes(q_3-1)})
\\&=&
n^{0.5(\bar{q}(3)-3)}\sum_{j,\ell}q_1(D_{1_\ell}(a_j(q_1)a_j(q_2)))a_\ell(q_3)n^{-1}
\\&&\hspace{50pt}\times
{\colorr \sum_\nu c_\nu(q_1-1,q_2-1) I_{q_1+q_2-2-2\nu}(1_j^{\otimes (q_1+q_2-2-2\nu)}) n^{-\nu}}
I_{q_3-1}(1_\ell^{\otimes(q_3-1)})
\eeas
and hence 
\end{en-text}
\bea\label{202004041510}
\bbI_{1,2}
&=&
\bbI_{1,2,1}+\bbI_{1,2,2}
\eea
where 
\bea\label{202004131428}  
\bbI_{1,2,1}
&=& 
n^{0.5q_3-1.5}\sum_\ell\bigg\{n^{-1}\sum_j
{\colorr 1_{\{q_1=q_2\}}q_1!}(nD_{1_\ell}(a_j(q_1)a_j(q_2)))a_\ell(q_3)\bigg\}
I_{q_3-1}(1_\ell^{\otimes(q_3-1)})
\nn\\&&
\eea
and 
\beas 
\bbI_{1,2,2}
&=&
\sum_{\nu:\>q_1+q_2-2-2\nu>0}
c_\nu(q_1-1,q_2-1) 
n^{0.5\bar{q}(3)-{\sred 3.5}-\nu}\sum_{j,\ell}{\sred n}
q_1(D_{1_\ell}(a_j(q_1)a_j(q_2)))a_\ell(q_3)
\\&&\hspace{50pt}\times
{\colorr I_{q_1+q_2-2-2\nu}(1_j^{\otimes (q_1+q_2-2-2\nu)})}
I_{q_3-1}(1_\ell^{\otimes(q_3-1)})
\eeas
Then 
\bea\label{202004131429} 
e(\bbI_{1,2,1}) 
&=& 
-0.5
\eea
and 
\beas 
e(\bbI_{1,2,2}) 
&=& 
{\sred \max_{\nu:q_1+q_2-2-2\nu>0}\bigg\{}
\big(0.5\bar{q}(3)-{\sred 3.5}-\nu\big)-0.5\big(\bar{q}(3)-3-2\nu\big)+{\sred 2}-0.5\times2
{\sred \bigg\}}
\nn\\&=&
-1. 
\eeas
Thus, the IBP gives 
\beas 
E[\Psi(\sfz,\sfx)\bbI_{1,2}]
&=&
n^{0.5q_3-1.5}
{\colorr 1_{\{q_1=q_2\}}q_1!}
\\&&\times
\sum_{j,\ell}E\big[\Psi(\sfz,\sfx)(D_{1_\ell}(a_j(q_1)a_j(q_2)))a_\ell(q_3)
I_{q_3-1}(1_\ell^{\otimes(q_3-1)})\big]
+O(n^{-1})
\\&=&
n^{0.5q_3-1.5}
{\colorr 1_{\{q_1=q_2\}}q_1!}
\\&&\times
\sum_{j,\ell}E\big[D^{q_3-1}_{1_\ell^{\otimes(q_3-1)}}\big\{
\Psi(\sfz,\sfx)(D_{1_\ell}(a_j(q_1)a_j(q_2)))a_\ell(q_3)\big\}\big]
+O(n^{-1})
\\&=&
{\colorr 1_{\{q_1+q_2:\text{ even}\}}}O(n^{-0.5(q_3-1)})+O(n^{-1})
\eeas
\begin{en-text}
Here the exponent is obtained by 
\beas 
0.5(q-1)-1+2-q&=&-0.5q+0.5
\eeas
\end{en-text}
{\colorr
Only when $q_1=q_2$ and $q_3=2$, this term has effect to the asymptotic expansion. 
Therefore 
\bea\label{202003161227}
E[\Psi(\sfz,\sfx)\bbI_{1,2}]
&=&
{\colorr1_{\{q_1=q_2,\>q_3=2\}}}
n^{-0.5}{\colorr  q_1!}
\sum_{j,\ell}E\big[D_{1_\ell}\big\{
\Psi(\sfz,\sfx)(D_{1_\ell}(a_j(q_1)^2)a_\ell(2)\big\}\big]
\nn\\&&
+O(n^{-1}).
\eea
}

From (\ref{202003161223}), (\ref{202003161225}), (\ref{202003161226}) and (\ref{202003161227}), 
it follows that 
\bea\label{202003161326}
E[\Psi(\sfz,\sfx)D_{u_n(q_3)}\bbI_1]
&=&
{\colorb 1_{\{\bar{q}(3):even\}}}
q_1(q_1-1)c_{0.5\bar{q}(3)-2}(q_1-2,q_2-1,q_3-1)
\nn\\&&\hspace{30pt}\times
n^{-1.5}\sum_{j=1}^n
E\big[\Psi(\sfz,\sfx)a_j(q_1)a_j(q_2)a_j(q_3)\big]
\nn\\&&
+{\colorb 1_{\{\bar{q}(3):even\}}}
q_1(q_2-1)c_{0.5\bar{q}(3)-2}(q_1-1,q_2-2,q_3-1)
\nn\\&&\hspace{30pt}\times
n^{-1.5}\sum_{j=1}^n
E\big[\Psi(\sfz,\sfx)a_j(q_1)a_j(q_2)a_j(q_3)\big]+O(n^{-1})
\nn\\&&
+{\colorr1_{\{q_1=q_2,\>q_3=2\}}}
n^{-0.5}{\colorr  q_1!}
\nn\\&&\times
\sum_{j,\ell}E\big[D_{1_\ell}\big\{
\Psi(\sfz,\sfx)(D_{1_\ell}({\colorr a_j(q_1)^2}))a_\ell(2)\big\}\big]
+O(n^{-1})
\eea
We remark that, on the right-hand side of (\ref{202003161326}), 
 the contribution possibly occurs from the first term only when $(q_1-2,q_2-1,q_3-1)\in\bbT$, 
 and 
 from the second term only when $(q_1-1,q_2-2,q_3-1)\in\bbT$. 

From the above discussion, we know 
\bea\label{202004131405} 
D_{u_n(q_3)}\bbI_1 &=& \check{\bbI}_1+ \hat{\bbI}_1
\eea
where 
\bea\label{202004131406}
\check{\bbI}_1 &=& \bbI_{1,1,1}+\bbI_{1,1,1}'+\bbI_{1,2,1}
\eea
with 
\beas 
e(\check{\bbI}_1) &=& -0.5
\eeas
and 
\bea\label{202004131408} 
e( \hat{\bbI}_1) &=& -1.
\eea

\begin{en-text}
When $q_1=2$, the sum $\bbI_{1,1}$ has the expression 
\beas
\bbI_{1,1}
&=&
n^{0.5\bar{q}(3)-3.5}\sum_{j}2a_j(2)a_j(q_2)a_j(q_3)
I_{q_2-1}(1_j^{\otimes(q_2-1)}) I_{q_3-1}(1_j^{\otimes(q_3-1)}).
\eeas
By the product formula, 
\beas 
\bbI_{1,1}
&=&
\sum_\nu c_\nu(q_2-1,q_3-1)
n^{0.5\bar{q}(3)-3.5-\nu}\sum_{j}2a_j(2)a_j(q_2)a_j(q_3)
I_{q_2+q_3-2-2\nu}(1_j^{\otimes(q_2+q_3-2-2\nu)}).
\eeas
Therefore, the summands have the exponent 
\beas 
&&
\big(0.5\bar{q}(3)-3.5-\nu\big)-0.5(q_2+q_3-2-2\nu)+1-0.51_{\{q_2+q_3-2-2\nu\}}
\\&=&
-0.5-0.51_{\{q_2+q_3-2-2\nu\}}.
\eeas
That $q_2+q_3-2-2\nu=0$ is equivalent to that $\nu=q_2-1=q_3-1$, since 
$\nu\leq\min\{q_2-1,q_3-1\}$. 
Therefore, 
\beas 
e(\bbI_{1,1})
&=&
\left\{\begin{array}{ll}
-0.5&(q_1=2,\>q_2=q_3)\y
-1&(q_1=2,\>q_2\not=q_3)
\end{array}\right.
\eeas

When $q_1>2$, we apply the product formula twice to the sum $\bbI_{1,1}$ to obtain 
\beas
\bbI_{1,1}
&=&
\bbI_{1,1,1}+\bbI_{1,1,2}
\eeas
where
\beas
\bbI_{1,1,1}
&=&
1_{\{q_1+q_2:\text{ even}\}}
n^{(0.5\bar{q}(3)-3.5)-0.5(q_1+q_2-3)}
\nn\\&&\times
\sum_{j}q_1(q_1-1)a_j(q_1)a_j(q_2)a_j(q_3)
c_{0.5(q_1+q_2-3)} I_{q_3-1}(1_j^{\otimes(q_3-1)})
\eeas
and 
\beas 
\bbI_{1,1,2}
&=&
\sum_{\nu_1:q_1+q_2-3-2\nu_1>0}
n^{(0.5\bar{q}(3)-3.5)-\nu_1}
\sum_{j}q_1(q_1-1)a_j(q_1)a_j(q_2)a_j(q_3)c_{\nu_1}(q_1-2,q_2-1)
\\&&\hspace{80pt}\times
 I_{q_1+q_2-3-2\nu_1}(1_j^{\otimes(q_1+q_2-3-2\nu_1)}) I_{q_3-1}(1_j^{\otimes(q_3-1)})
\nn\\&=&
\sum_{\nu_1:q_1+q_2-3-2\nu_1>0}\sum_{\nu_2}
n^{(0.5\bar{q}(3)-3.5)-\nu_1-\nu_2}
\nn\\&&\times
\sum_{j}q_1(q_1-1)a_j(q_1)a_j(q_2)a_j(q_3)c_{\nu_1}(q_1-2,q_2-1)
c_{\nu_2}(q_1+q_2-3-2\nu_1,q_3-1)
\\&&\hspace{80pt}\times
 I_{\bar{q}(3)-4-2\nu_1-2\nu_2}(1_j^{\otimes(\bar{q}(3)-4-2\nu_1-2\nu_2)}). 
\eeas
The exponent 
\beas 
e(\bbI_{1,1,1})
&=& 
\big\{(0.5\bar{q}(3)-3.5)-0.5(q_1+q_2-3)\big\}
-0.5(q_3-1)+1-0.5
\nn\\&=&
-1. 
\eeas
About the exponent of $\bbI_{1,1,2}$, 
\beas 
e(\bbI_{1,1,2})
&\leq&
\max_{\nu_1,\nu_2}\bigg\{
\big\{(0.5\bar{q}(3)-3.5)-\nu_1-\nu_2\big\}-0.5\big(\bar{q}(3)-4-2\nu_1-2\nu_2\big)\bigg\}
\nn\\&&
+1-0.5\>1_{\{\bar{q}(3)-4-2\nu_1-2\nu_2>0\}}
\nn\\&=&
\max_{\nu_1,\nu_2}\bigg\{-0.5-0.5\>1_{\{\bar{q}(3)-4-2\nu_1-2\nu_2>0\}}\bigg\}. 
\eeas
The case $\bar{q}(3)-4-2\nu_1-2\nu_2=0$ occurs only when 
\beas 
\nu_2\yeq 0.5(\bar{q}(3)-4-3\nu_1)\yeq q_1+q_2-3-2\nu_1 \yeq q_3-1,
\eeas
since $\nu_2\leq\min\{q_1+q_2-3-2\nu_1,q_3-1\}$. 
Then $\nu_1=0.5(q_1+q_2-q_3)-1$.

\koko

By (\ref{202003151619}), 
\bea\label{202004041456}
 \bbI_{1,1}
&=&
n^{0.5(\bar{q}(3)-3)}\sum_jq_1(q_1-1)a_j(q_1)a_j(q_2)a_j(q_3) n^{-2}
\nn\\&&\hspace{50pt}\times
\sum_{\nu\in\bbZ_+} c_{\nu}(q_1-2,q_2-1,q_3-1)I_{\bar{q}(3)-4-2|\nu|}(1_j^{\otimes (\bar{q}(3)-4-2\nu)})
n^{-\nu}
\eea
Therefore, by the integration-by-parts (IBP) 
(use differentiability of $X_\infty$ also), we can say  
\beas 
|E[\Psi(\sfz,\sfx)\bbI_{1,1}]|
&=&
O(n^{-\xi^{(\ref{202003151625})}})
\eeas
as $n\to\infty$ for every $(\sfz,\sfx)\in\bbR^2$. 
where 
\bea\label{202003151625}
\xi^{(\ref{202003151625})} &=& 0.5+0.5 \times{\bf1}_{\{\bar{q}(3):odd\}}.
\eea
Here the exponent of $n$ was obtained by 
\beas &&
0.5(\bar{q}(3)-3)+1-2-(\bar{q}(3)-4-2\nu)-\nu
\\&=&
-0.5(\bar{q}(3)-3)+\nu
\quad(\nu:\>\bar{q}(3)-4-2\nu\geq0)
\\&\leq&
-0.5(\bar{q}(3)-3)+0.5\bar{q}(3)-2-0.5 \times{\bf1}_{\{\bar{q}(3):odd\}}
\\&=&
-0.5-0.5 \times{\bf1}_{\{\bar{q}(3):odd\}}=-\xi^{(\ref{202003151625})}.
\eeas
%
Thus, the last term in the product formula determines the order of the symbol. 
\end{en-text}
\begin{en-text}
\beas 
E[\Psi(\sfz)\bbI_{1,1}]
&=&
q_1(q_1-1)c_{0.5\bar{q}(3)-2}(q_1-2,q_2-1,q_3-1)
\\&&\hspace{30pt}\times
n^{-1.5+0.5\ep(q)}\sum_{j=1}^n
E\big[D^{\ep(\bar{q}(3))}_{1_j}(\Psi(\sfz)a_j(q_1)a_j(q_2)a_j(q_3))\big]+O(n^{-1})
\eeas
where $\ep(q)=1_{\{q:\text{odd}\}}$. 
When $\bar{q}(3)$ is odd, this term has no effect. When $\bar{q}(3)$ is even, this term remains. 
Under the present assumption that all $q_1,q_2,q_3$ are even, 
this term remains. 
\end{en-text}
\begin{en-text}
More precisely, 
\bea\label{202003161225}
E[\Psi(\sfz,\sfx)\bbI_{1,1}]
&=&
{\colorb 1_{\{\bar{q}(3):even\}}}
q_1(q_1-1)c_{0.5\bar{q}(3)-2}(q_1-2,q_2-1,q_3-1)
\nn\\&&\hspace{30pt}\times
n^{-1.5}\sum_{j=1}^n
E\big[\Psi(\sfz,\sfx)a_j(q_1)a_j(q_2)a_j(q_3)\big]+O(n^{-1}). 
\eea

Similarly, 
\bea\label{202003161226}
E[\Psi(\sfz,\sfx)\bbI_{1,1}']
&=&
{\colorb 1_{\{\bar{q}(3):even\}}}
q_1(q_2-1)c_{0.5\bar{q}(3)-2}(q_1-1,q_2-2,q_3-1)
\nn\\&&\hspace{30pt}\times
n^{-1.5}\sum_{j=1}^n
E\big[\Psi(\sfz,\sfx)a_j(q_1)a_j(q_2)a_j(q_3)\big]+O(n^{-1}). 
\eea

To estimate $E[\Psi(\sfz,\sfx)\bbI_{1,2}]$, 
we apply the product formula to $I_{q-1}(1_j^{\otimes(q-1)})^2$ first and obtain 
\beas
\bbI_{1,2}
&=&
n^{0.5(\bar{q}(3)-3)}\sum_{j,\ell}q_1(D_{1_\ell}(a_j(q_1)a_j(q_2)))a_\ell(q_3)n^{-1}
\\&&\hspace{80pt}\times
{\colorr I_{q_1-1}(1_j^{\otimes(q_1-1)}) I_{q_2-1}(1_j^{\otimes(q_2-1)}) }
I_{q_3-1}(1_\ell^{\otimes(q_3-1)})
\\&=&
n^{0.5(\bar{q}(3)-3)}\sum_{j,\ell}q_1(D_{1_\ell}(a_j(q_1)a_j(q_2)))a_\ell(q_3)n^{-1}
\\&&\hspace{50pt}\times
{\colorr \sum_\nu c_\nu(q_1-1,q_2-1) I_{q_1+q_2-2-2\nu}(1_j^{\otimes (q_1+q_2-2-2\nu)}) n^{-\nu}}
I_{q_3-1}(1_\ell^{\otimes(q_3-1)})
\eeas
and hence 
\bea\label{202004041510}
\bbI_{1,2}
&=&
n^{0.5(\bar{q}(3)-3)}\sum_{j,\ell}q_1(D_{1_\ell}(a_j(q_1)a_j(q_2)))a_\ell(q_3)n^{-1}
{\colorr 1_{\{q_1+q_2:\text{ even}\}}}
\nn\\&&\hspace{50pt}\times
{\colorr  c_{0.5(q_1+q_2-2)}(q_1-1,q_2-1) n^{-0.5(q_1+q_2-2)}}
I_{q_3-1}(1_\ell^{\otimes(q_3-1)})
\nn\\&&+O_{L^\inftym}(n^{-1})
\quad(\text{by }L^p\text{ estimate})
\eea
Here the last term of $O_{L^\inftym}(n^{-1})$ is obtained 
by the $L^p$-estimate using Theorem \ref{0110221242}; the exponent is counted by 
\beas 
0.5(\bar{q}(3)-3)-1(D_{1_\ell})-1
-0.5\big((q_1+q_2-2-2\nu)+(q_3-1)-2\big)-\nu
&=&
-1
\eeas
whenever $q_1+q_2-2-2\nu\geq1$. 
\end{en-text}

\subsubsection{$\boldsymbol{D_{u_n(q_3)}\bbI_2}$}
Recall 
\beas
\bbI_2
&=&
n^{2^{-1}(q_1+q_2)-1}\sum_{j,k}(D_{1_k}a_j(q_1)) a_k(q_2) I_{q_1}(1_j^{\otimes q_1})I_{q_2-1}(1_k^{\otimes(q_2-1)}).
\eeas
Then the projection $D_{u_n(q_3)}\bbI_2$ admits the decomposition
\bea\label{202003161439}
D_{u_n(q_3)}\bbI_2
&=&
n^{0.5(\bar{q}(3)-3)}\sum_{j,k,\ell}(D_{1_k}a_j(q_1)) a_k(q_2) q_1I_{q_1-1}(1_j^{\otimes (q_1-1)})\beta_{j,\ell}a_\ell(q_3)
I_{q_2-1}(1_k^{\otimes(q_2-1)})I_{q_3-1}(1_\ell^{\otimes(q_3-1)})
\nn\\&&
+n^{0.5(\bar{q}(3)-3)}\sum_{j,k,\ell}(D_{1_k}a_j(q_1)) a_k(q_2) I_q(1_j^{\otimes q_1}) (q_2-1)I_{q_2-2}(1_k^{\otimes(q_2-2)})\beta_{k,\ell}
\nn\\&&\hspace{50pt}\times
a_\ell(q_3) I_{q_3-1}(1_\ell^{\otimes(q_3-1)})
\nn\\&&
+n^{0.5\bar{q}(3)-3.5}\sum_{j,k,\ell}\big\{n^2D_{1_\ell}\big((D_{1_k}a_j(q_1)) a_k(q_2) \big)\big\}a_\ell (q_3)
I_q(1_j^{\otimes q_1})I_{q_2-1}(1_k^{\otimes(q_2-1)})I_{q_3-1}(1_\ell^{\otimes(q_3-1)})
\nn\\&=:&
\bbI_{2,1}+\bbI_{2,2}+\bbI_{2,3}.
\eea

The sum $\bbI_{2,1}$ has the same pattern as $\bbI_{1,2}$. Thus, 
\bea\label{202004041539}
\bbI_{2,1}
&=&
n^{0.5(\bar{q}(3)-3)-1}\sum_{j,k}(D_{1_k}a_j(q_1)) a_k(q_2)a_j(q_3)  q_1I_{q_1-1}(1_j^{\otimes (q_1-1)})I_{q_3-1}(1_j^{\otimes (q_3-1)})
I_{q_2-1}(1_k^{\otimes(q_2-1)})
\nn\\&=&
n^{0.5(\bar{q}(3)-3)-1}\sum_{j,k}(D_{1_k}a_j(q_1)) a_k(q_2)a_j(q_3)  q_1
\nn\\&&\times
\sum_{\nu\in\bbZ_+
}
c_\nu(q_1-1,q_3-1)I_{q_1+q_3-2-2\nu}(1_j^{\otimes(q_1+q_3-2-2\nu)})n^{-\nu}
I_{q_2-1}(1_k^{\otimes(q_2-1)})
\nn\\&=&
\bbI_{2,1,1}+\bbI_{2,1,2}
\eea
where 
\beas 
\bbI_{2,1,1}
&=&
n^{0.5q_2-1.5}\sum_k\bigg\{n^{-1}\sum_j(nD_{1_k}a_j(q_1)) a_k(q_2)a_j({\colorr q_1})  \bigg\}
\nn\\&&\hspace{50pt}\times
{\colorr 1_{\{q_1=q_3\}}q_1!}
I_{q_2-1}(1_k^{\otimes(q_2-1)})
\eeas
and 
\beas
\bbI_{2,1,2}
&=&
\sum_{\nu\in\bbZ_+,\atop q_1+q_3-2-2\nu>0
}
n^{0.5\bar{q}(3)-3.5-\nu}\sum_{j,k}(nD_{1_k}a_j(q_1)) a_k(q_2)a_j(q_3)  q_1
\nn\\&&\hspace{30pt}\times
c_\nu(q_1-1,q_3-1)I_{q_1+q_3-2-2\nu}(1_j^{\otimes(q_1+q_3-2-2\nu)})
I_{q_2-1}(1_k^{\otimes(q_2-1)}).
\nn\\&&
\eeas
\begin{en-text}
If Theorem \ref{0110221242} is applied, 
the exponent of the $L^p$-norm of the first term in the last expression is 
\beas 
0.5(\bar{q}(3)-3)-1+1(\sum_j)-1(D_{1_k})-(q_1+q_3-2)/2-0.5(q_2-2)(\sum_k)&=&-0.5,
\eeas
therefore this sum can remain and is treated below. 
The $L^p$-exponent of the second term is 
\beas 
0.5(\bar{q}(3)-3)-1-1(D_{1_k})-0.5(\bar{q}(3)-3-2\nu-2)(\sum_{j,k})-\nu&=&-1, 
\eeas
therefore the second term is negligible. 
\end{en-text}
Then 
\beas 
e(\bbI_{2,1,1})
&=& 
(0.5q_2-1.5)-0.5(q_2-1)+1-0.5 
\yeq
-0.5,
\eeas
therefore this sum can remain and is treated below. 
On the other hand, 
\beas 
e(\bbI_{2,1,2})
&=& 
\max_{\nu;\>q_1+q_3-2-2\nu>0}\bigg\{
(0.5\bar{q}(3)-3.5-\nu)-0.5(\bar{q}(3)-3-2\nu)+2-1\bigg\}
\yeq
-1, 
\eeas
therefore this term is negligible. 

By these observations, we compute the projection of $\bbI_{2,1}$ to $\Psi(\sfz,\sfx)$: 
\beas 
E[\Psi(\sfz,\sfx)\bbI_{2,1}] 
&=& 
{\colorr q_1! 1_{\{q_1=q_3\}}}
n^{0.5q_2-1.5}
\\&&\times
\sum_{j,k}
E\bigg[\Psi(\sfz,\sfx)
(D_{1_k}a_j(q_1)) a_k(q_2)a_j({\colorr q_1})  
I_{q_2-1}(1_k^{\otimes(q_2-1)})
\bigg]
+O(n^{-1})
\\&=& 
{\colorr q_1! 1_{\{q_1=q_3\}}}
n^{0.5q_2-1.5}
\\&&\times
\sum_{j,k}
E\bigg[D^{q_2-1}_{1_k^{\otimes(q_2-1)}}\big(\Psi(\sfz,\sfx)
(D_{1_k}a_j(q_1)) a_k(q_2)a_j({\colorr q_1}) \big)
\bigg]
+O(n^{-1})
\eeas
by the IBP. 
The order of the order of this term is therefore 
\beas
(0.5q_2-1.5)+2(\sum_{j,k})-q_2&=&-0.5q_2+0.5,
\eeas
and 
{\colorr 
only the case $q_2=2$ remains. 
Consequently, 
we obtain 
\bea\label{202003161430} 
E[\Psi(\sfz,\sfx)\bbI_{2,1}] 
&=& 
{\colorr 1_{\{q_1=q_3,\>q_2=2\}}}
n^{-0.5}q_1!
\nn\\&&\times
E\bigg[
\sum_{j,k}\big\{D_{1_k}\big(\Psi(\sfz)(D_{1_k}a_j(q_1)) a_k(2)a_j(q_1))\big)\big\}  \bigg]
+O(n^{-1})
\eea
}

The sum $\bbI_{2,2}$ is evaluated as follow: 
\bea\label{20200404552}
\bbI_{2,2} 
&=&
n^{0.5(\bar{q}(3)-3)-1}\sum_{j,k}(D_{1_k}a_j(q_1)) a_k(q_2)a_k(q_3) 
I_{q_1}(1_j^{\otimes q_1}) (q_2-1)I_{q_2-2}(1_k^{\otimes(q_2-2)})
 I_{q_3-1}(1_k^{\otimes(q_3-1)})
\nn\\&=&
n^{0.5(\bar{q}(3)-3)-1}\sum_{j,k}(D_{1_k}a_j(q_1)) a_k(q_2)a_k(q_3) (q_2-1) I_{q_1}(1_j^{\otimes q_1})
\nn\\&&\hspace{50pt}\times
{\colorr \sum_{\nu\in\bbZ_+}c_\nu(q_2-2,q_3-1)I_{q_2+q_3-3-2\nu}(1_k^{\otimes(q_2+q_3-3-2\nu)})n^{-\nu}}
\nn\\&=&
\bbI_{2,2,1}+\bbI_{2,2,2}
\eea
where
\bea\label{202004101846}
\bbI_{2,2,1} 
&=&
n^{0.5q_1-1}\sum_j\bigg\{n^{-1}\sum_k(nD_{1_k}a_j(q_1)) a_k(q_2)a_k({\colorr q_2-1}) 
{\colorr (q_2-1)!}
 \bigg\}
I_{q_1}(1_j^{\otimes q_1})
\nn\\&&\hspace{50pt}\times
{\colorr1_{\{q_2=q_3+1\}} }
\nn\\&&
\eea
and 
\bea\label{202004101847}
\bbI_{2,2,2}
&=&
\sum_{\nu:q_2+q_3-3-2\nu>0}c_\nu(q_2-2,q_3-1)
n^{0.5\bar{q}(3)-3.5-\nu}
\nn\\&&\hspace{50pt}\times
\sum_{j,k}(nD_{1_k}a_j(q_1)) a_k(q_2)a_k(q_3) (q_2-1) I_{q_1}(1_j^{\otimes q_1})
I_{q_2+q_3-3-2\nu}(1_k^{\otimes(q_2+q_3-3-2\nu)}).
\nn\\&&
\eea

From the representation (\ref{202004101846}), 
\bea\label{202004101848}
e(\bbI_{2,2,1} ) 
&=&
(0.5q_1-1)-0.5q_1+1-0.5 \yeq -0.5
\eea
and 
\bea\label{202004101849}
e(\bbI_{2,2,2} ) 
&=&
\max_{\nu;\>q_2+q_3-3-2\nu>0}\bigg\{
(0.5\bar{q}(3)-3.5-\nu)-0.5(\bar{q}(3)-3-2\nu)+2-0.5\times2
\bigg\}
\nn\\&=&
-1. 
\eea
Then
\beas 
E[\Psi(\sfz,\sfx)\bbI_{2,2}]
&=&
{\colorr1_{\{q_2=q_3+1\}} }
n^{0.5q_1-1} \sum_{j,k} {\colorr (q_2-1)!}
\\&&\times
E\bigg[\Psi(\sfz,\sfx) (D_{1_k}a_j(q_1)) a_k(q_2)a_k({\colorr q_2-1}) I_{q_1}(1_j^{\otimes q_1})
\bigg]+O(n^{-1})
\\&=&
{\colorr1_{\{q_2=q_3+1\}} }
n^{0.5q_1-1} {\colorr (q_2-1)!} 
\\&&\times
\sum_{j,k}E\bigg[D^{q_1}_{1_j^{\otimes q_1}}\big\{\Psi(\sfz,\sfx) (D_{1_k}a_j(q_1)) a_k(q_2)
a_k({\colorr q_2-1})\big\}
\bigg]+O(n^{-1})
\eeas
by IBP. The exponent to $n$ of this expression is 
\beas 
(0.5q_1-1)+2(\sum_{j,k})-q_1(D_{1_j^{\otimes q_1}})-1(D_{1_k})
&=& -0.5q_1\leq-1\quad(q_1\geq2)
\eeas
Thus, $\bbI_{2,2}$ asymptotically has no effect to the random symbol 
(even without Condition (\ref{202004111810})): 
\bea\label{202004111223}
E[\Psi(\sfz,\sfx)\bbI_{2,2}]
&=&
O(n^{-1}). 
\eea

{\sred We remark that,} 
under Condition (\ref{202004111810}), $\bbI_{2,2,1}=0$ since $q_2\not=q_3+1$ for $q_2,q_3\in\calq$, 
hence 
\bea\label{202004121242}
e(\bbI_{2,2}) &=& -1,
\eea
and in particular, (\ref{202004111223}) follows directly.

Based on the representation of $\bbI_{2,3}$ in (\ref{202003161439}), 
\beas 
e(\bbI_{2,3})
&=&
\big(0.5\bar{q}(3)-3.5\big)-0.5\big(\bar{q}(3)-2\big)+3-0.5\times3
\yeq-1
\eeas
and we conclude 
\bea\label{202004041609}
\bbI_{2,3}&=&O_{L^\inftym}(n^{-1})
\eea
and 
\bea\label{202004111224}
E\big[\Psi(\sfz,\sfx)\bbI_{2,3}\big] &=& O(n^{-1}). 
\eea

Consequently, 
\bea\label{202003171201}
E[\Psi(\sfz,\sfx)D_{u_n(q_3)}\bbI_2] 
&=& 
{\colorr 1_{\{q_1=q_3,\>q_2=2\}}}
n^{-0.5} {\colorr q_1!}
\nn\\&&\times
E\bigg[
\sum_{j,k}\big\{D_{1_k}\big(\Psi(\sfz,\sfx)(D_{1_k}a_j(q_1)) a_k(2)a_j({\colorr q_1})\big)\big\}  \bigg]
+O(n^{-1})
\eea
by (\ref{202003161439}), (\ref{202003161430}), (\ref{202004111223}) and (\ref{202004111224}).

\subsubsection{$\boldsymbol{D_{u_n(q_3)}\bbI_3}$}
In the present case of a Brownian motion, 
\bea\label{202004111300} 
\bbI_3
&=&
\bbJ_{3,1}+\bbJ_{3,2}
\eea
where 
\bea\label{202004111301} 
\bbJ_{3,1}
&=&
-n^{-1.5}\sum_{j}1_{\{q_1=q_2+1\}}(q_1-1)!(nD_{1_j}a_j(q_1)) a_j(q_2)
I_0(1_j^{\otimes0})
\eea
and 
\bea\label{202004111302}
\bbJ_{3,2}
&=&
-
\sum_{\{\nu:\>q_1+q_2-3-2\nu>0\}}c_\nu(q_1-2,q_2-1)(q_1-1)
\nn\\&&\hspace{50pt}\times
n^{0.5(q_1+q_2)-3-\nu}\sum_{j}(nD_{1_j}a_j(q_1)) a_j(q_2)
I_{q_1+q_2-3-2\nu}(1_j^{\otimes(q_1+q_2-3-2\nu)})
\eea
Then 
\bea\label{202004111303}
e(\bbJ_{3,1}) &=& -0.5
\eea
and 
\bea\label{202004111304}
e(\bbJ_{3,2}) &=& -1.
\eea
Applying Proposition \ref{202004091035} to $\bbJ_{3,1}$ having the representation (\ref{202004111301}) and 
the exponent (\ref{202004111303}), we see $D_{u_n(q_3)}\bbJ_{3,1}\in\call$ and 
\bea\label{202004111305}
e(D_{u_n(q_3)}\bbJ_{3,1}) &\leq& -1
\eea
since $0\not=q_3\in\calq$. 
Simple application of Proposition \ref{202004091035} gives 
\bea\label{202004111306}
e(D_{u_n(q_3)}\bbJ_{3,2}) &\leq& -1
\eea
Thus, from (\ref{202004111300}), (\ref{202004111305}) and (\ref{202004111306}), 
we conclude $D_{u_n(q_3)}\bbI_3\in\call$ and 
\bea\label{202004111307}
e(D_{u_n(q_3)}\bbI_3) &\leq& -1. 
\eea
In particular, 
\begin{en-text}
Then 
\beas &&
D_{u_n(q_3)}\bbI_3 
\\&=&
-n^{0.5\bar{q}(3)-2.5}\sum_{j,\ell}(q_1-1)(D_{1_j}a_j(q_1)) a_j(q_2)
\sum_{\nu\in\bbZ_+}
c_\nu(q_1-2,q_2-1)(q_1+q_2-3-2\nu)
\\&&\hspace{50pt}\times
I_{q_1+q_2-4-2\nu}(1_j^{\otimes(q_1+q_2-4-2\nu)})\beta_{j,\ell}
a_\ell(q_3) I_{q_3-1}(1_\ell^{\otimes(q_3-1)})n^{-\nu}
\\&&
-n^{0.5\bar{q}(3)-2.5}\sum_{j,\ell}(q_1-1)\{D_{1_\ell}((D_{1_j}a_j(q_1)) a_j(q_2) )\}a_\ell(q_3)
\\&&\hspace{50pt}\times
\sum_{\nu\in\bbZ_+}
c_\nu(q_1-2,q_2-1)
I_{q_1+q_2-3-2\nu}(1_j^{\otimes(q_1+q_2-3-2\nu)})I_{q_3-1}(1_\ell^{\otimes(q_3-1)})n^{-\nu}
\\&=&
-n^{0.5\bar{q}(3)-3.5}\sum_{j}(q_1-1)(D_{1_j}a_j(q_1)) a_j(q_2)a_\ell(q_3)
\sum_{\nu\in\bbZ_+}
c_\nu(q_1-2,q_2-1)(q_1+q_2-3-2\nu)
\\&&\hspace{50pt}\times
I_{q_1+q_2-4-2\nu}(1_j^{\otimes(q_1+q_2-4-2\nu)})
 I_{q_3-1}(1_j^{\otimes(q_3-1)})n^{-\nu}
\\&&
-n^{0.5\bar{q}(3)-2.5}\sum_{j,\ell}(q_1-1)\{D_{1_\ell}((D_{1_j}a_j(q_1)) a_j(q_2) )\}a_\ell(q_3)
\\&&\hspace{50pt}\times
\sum_{\nu\in\bbZ_+}
c_\nu(q_1-2,q_2-1)
I_{q_1+q_2-3-2\nu}(1_j^{\otimes(q_1+q_2-3-2\nu)})I_{q_3-1}(1_\ell^{\otimes(q_3-1)})n^{-\nu}
\\&=:&
\bbI_{3,1}+\bbI_{3,2}.
\eeas

Simple $L^p$-estimate gives an estimate for the exponent of $\bbI_{3,1}$ in $L^\inftym$: 
\beas 
\big(0.5\bar{q}(3)-3.5\big)-1(D_{1_j})+1(\sum_j)
-{\colorb0.5(q_1+q_2-3-2\nu)-0.5(q_3-1)}-\nu{\colorb\yeq -1}
\eeas
and hence $\bbI_{3,1}=O_{L^\inftym}(n^{{\colorb-1}})$. 
Theorem \ref{0110221242} gives the following estimate of the exponent of $\bbI_{3,2}$: 
\beas 
\big(0.5\bar{q}(3)-2.5\big)-2(D_{1_\ell}D_{1_j})-0.5(\bar{q}(3)-4-2\nu-2)-\nu=-1.5
\eeas
and hence $\bbI_{3,2}=O_{L^\inftym}(n^{-1.5})$. 
Thus, 
\end{en-text}
\bea\label{202003171314}
D_{u_n(q_3)}\bbI_3 =O_{L^\inftym}(n^{-1})
\eea
has no effect asymptotically.

\subsubsection{{\sf qTor}}
From (\ref{202003130046}) and (\ref{202003171321}), we have 
\beas 
D_{u_n}D_{u_n}M_n 
&=& 
\sum_{q_1,q_2,q_3\in\calq}D_{u_n(q_3)}D_{u_n(q_2)}M_n(q_1)
\eeas
and, from (\ref{20200317316}), 
\beas 
D_{u_n(q_3)}D_{u_n(q_2)}M_n(q_1)
&=& 
D_{u_n(q_3)}\bbI_1+D_{u_n(q_3)}\bbI_2+D_{u_n(q_3)}\bbI_3+D_{u_n(q_3)}\bbI_4.
\eeas
Then, by the estimates (\ref{202003161326}), (\ref{202003171201}), (\ref{202003171314}) 
and (\ref{202003161330}), 
and by the symmetry of the functions $c_\nu(\cdots)$, 
we obtain 
\beas 
E\big[\Psi(\sfz,\sfx)D_{u_n}D_{u_n}M_n\big]
&=&
\sum_{q_1,q_2,q_3\in\calq}
E\big[\Psi(\sfz,\sfx)D_{u_n(q_3)}D_{u_n(q_2)}M_n(q_1)\big]
\\&=&
\sum_{q_1,q_2,q_3\in\calq}{\colorb 1_{\{\bar{q}(3):even\}}}
(q_1+q_2)(q_1-1)c_{0.5\bar{q}(3)-2}(q_1-2,q_2-1,q_3-1)
\\&&\hspace{50pt}\times
n^{-1.5}\sum_{j=1}^n
E\big[\Psi(\sfz,\sfx)a_j(q_1)a_j(q_2)a_j(q_3)\big]
\\&&
+\sum_{q_1,q_2,q_3\in\calq}{\colorr1_{\{q_1=q_2,\>q_3=2\}}}
n^{-0.5}  {\colorr q_1!}
\\&&\hspace{50pt}\times
\sum_{j,k=1}^nE\big[D_{1_k}\big\{
\Psi(\sfz,\sfx)(D_{1_k}({\colorr a_j(q_1)^2}))a_k(2)\big\}\big]
\\&&
+\sum_{q_1,q_2,q_3\in\calq}{\colorr 1_{\{q_1=q_3,\>q_2=2\}}}
n^{-0.5}{\colorr q_1!}
\\&&\hspace{50pt}\times
\sum_{j,k=1}^nE\big[
D_{1_k}\big(\Psi(\sfz,\sfx)(D_{1_k}a_j(q_1)) a_k(2)a_j({\colorr q_1})\big)  \big]
\\&&+O(n^{-1})
\eeas
or equivalently, 
\bea\label{202003171503} &&
n^{-0.5}E\big[\Psi(\sfz,\sfx){\sf qTor}(\tti\sfz)^3\big]
\nn\\&\equiv&
E\big[\Psi(\sfz,\sfx)D_{u_n}D_{u_n}M_n(\tti\sfz)^3\big]
\nn\\&=&
\sum_{q_1,q_2,q_3\in\calq}{\colorb 1_{\{\bar{q}(3):even\}}}
(q_1+q_2)(q_1-1)c_{0.5\bar{q}(3)-2}(q_1-2,q_2-1,q_3-1)
\nn\\&&\hspace{50pt}\times
n^{-1.5}\sum_{j=1}^n
E\big[\Psi(\sfz,\sfx)a_j(q_1)a_j(q_2)a_j(q_3)\big](\tti\sfz)^3
\nn\\&&
+3\sum_{q_1\in\calq}{\colorr1_{\{2\in\calq\}}}
n^{-0.5}  {\colorr q_1!}
\sum_{j,k=1}^nE\big[\Psi(\sfz,\sfx)D_{1_k}\big\{
(D_{1_k}a_j(q_1))a_j({\colorr q_1})a_k(2)\big\}\big](\tti\sfz)^3
\nn\\&&
+3\sum_{q_1\in\calq}{\colorr1_{\{2\in\calq\}}}
n^{-0.5}{\colorr q_1!}
\sum_{j,k=1}^nE\bigg[
\Psi(\sfz,\sfx)
\big(2^{-1}D_{1_k}G_\infty(\tti \sfz)^5+D_{1_k}X_\infty(\tti\sfz)^3(\tti\sfx)\big)
\nn\\&&\hspace{150pt}\times
(D_{1_k}a_j(q_1)) a_j({\colorr q_1})a_k(2)\big)  \bigg]
\nn\\&&+O(n^{-1})
\eea

\begin{en-text}
\beas 
\bbI_3
&=&
-n^{q-1}\sum_{j,k}(q-1)(D_{1_j}a_j) a_k \langle 1_j,1_k\rangle 
I_{q-2}(1_j^{\otimes(q-2)})I_{q-1}(1_k^{\otimes(q-1)})
\\&=&
-n^{q-2}\sum_{j}(q-1)(D_{1_j}a_j) a_j 
I_{q-2}(1_j^{\otimes(q-2)})I_{q-1}(1_j^{\otimes(q-1)})
\\&=&
-n^{q-2}\sum_{j}(q-1)(D_{1_j}a_j) a_j 
\sum_{\nu\in\calc(q-2,q-1)}
I_{3q-3-2\nu}(1_j^{\otimes(3q-3-2\nu)})n^{-\nu}
\eeas
\end{en-text}
\begin{en-text}
\bi
\im After all that, the following terms in $E\big[\Psi(\sfz)D_{u_n(q_3)}D_{u_n(q_2)}M_n(q_1)\big]$ 
have effect:
\beas 
E[\Psi(\sfz)\bbI_{1,1}]
&=&
q_1(q_1-1)c_{0.5\bar{q}(3)-2}(q_1-2,q_2-1,q_3-1)
\\&&\hspace{30pt}\times
n^{-1.5}\sum_{j=1}^n
E\big[\Psi(\sfz)a_j(q_1)a_j(q_2)a_j(q_3)\big]+O(n^{-1})
\eeas

\beas 
E[\Psi(\sfz)\bbI_{1,1}']
&=&
q_1(q_2-1)c_{0.5\bar{q}(3)-2}(q_1-1,q_2-2,q_3-1)
\\&&\hspace{30pt}\times
n^{-1.5}\sum_{j=1}^n
E\big[\Psi(\sfz)a_j(q_1)a_j(q_2)a_j(q_3)\big]+O(n^{-1})
\eeas
\beas 
E[\Psi(\sfz)\bbI_{1,2}]
&=&
n^{-0.5}  c_{0.5(q_1+q_2-2)}(q_1-1,q_2-1) q_1
\\&&\times
\sum_{j,\ell}E\big[D_{1_\ell}\big\{
\Psi(\sfz)(D_{1_\ell}(a_j(q_1)a_j(q_2)))a_\ell(2)\big\}\big]1_{\{q_3=2\}}
+O_{L^\inftym}(n^{-1}).
\eeas
\beas 
E[\Psi(\sfz)\bbI_{2,1}] 
&=& 
n^{-0.5}c_{0.5(q_1+q_3-2)}(q_1-1,q_3-1)q_1
\\&&\times
\sum_{j,k}E\big[
D_{1_k}\big(\Psi(\sfz)(D_{1_k}a_j(q_1)) a_k(2)a_j(q_3))\big)  \big]1_{\{q_2=2\}}
\\&&+O(n^{-1})
\eeas

\im symmetry in $c_\nu(...)$
\im 
\im Equivalently, {\sf q-Tor} is obtained from 
\beas &&
E\big[\Psi(\sfz)D_{u_n}D_{u_n}M_n\big]
\\&=&
\sum_{q_1,q_2,q_3\in\calq}
E\big[\Psi(\sfz)D_{u_n(q_3)}D_{u_n(q_2)}M_n(q_1)\big]
\\&=&
\sum_{q_1,q_2,q_3\in\calq}
(q_1+q_2)(q_1-1)c_{0.5\bar{q}(3)-2}(q_1-2,q_2-1,q_3-1)
n^{-1.5}\sum_{j=1}^n
E\big[\Psi(\sfz)a_j(q_1)a_j(q_2)a_j(q_3)\big]
\\&&
+\sum_{q_1,q_2\in\calq}n^{-0.5}  c_{0.5(q_1+q_2-2)}(q_1-1,q_2-1) q_1
\sum_{j,k=1}^nE\big[D_{1_k}\big\{
\Psi(\sfz)(D_{1_k}(a_j(q_1)a_j(q_2)))a_k(2)\big\}\big]
\\&&
+\sum_{q_1,q_2\in\calq}
n^{-0.5}c_{0.5(q_1+q_2-2)}(q_1-1,q_2-1)q_1
\sum_{j,k=1}^nE\big[
D_{1_k}\big(\Psi(\sfz)(D_{1_k}a_j(q_1)) a_k(2)a_j(q_2))\big)  \big]
\\&&+O(n^{-1})
\eeas
\ei
\end{en-text}

\subsubsection{Approximate projection of {\sf qTor}}
Define the random symbol ${\mathfrak S}^{(3,0)}_n(\tti\sfz)$ by 
\beas 
{\mathfrak S}^{(3,0)}_n(\tti\sfz)\yeq{\mathfrak S}^{(3,0)}_n(\tti\sfz,\tti\sfx)
\yeq \frac{1}{3}\>{\sf qTor}\>(\tti\sfz)^3. 
\eeas
Then 
\bea\label{202003181604}
E\big[\Psi(\sfz,\sfx){\mathfrak S}^{(3,0)}_n(\tti\sfz)\big]
&=& 
E\big[\Psi(\sfz,\sfx)\widetilde{\mathfrak S}^{(3,0)}_n(\sfz,\sfx)\big]+O(n^{-0.5})
\eea
for an approximate projection, that is a random symbol $\widetilde{\mathfrak S}^{(3,0)}_n(\sfz,\sfx)$ for ${\mathfrak S}^{(3,0)}_n(\sfz)$  defined by 
\bea\label{202003181605}
\widetilde{\mathfrak S}^{(3,0)}_n(\tti\sfz,\tti\sfx)
&=&
\frac{1}{3}\sum_{q_1,q_2,q_3\in\calq}{\colorb 1_{\{\bar{q}(3):even\}}}
(q_1+q_2)(q_1-1)c_{0.5\bar{q}(3)-2}(q_1-2,q_2-1,q_3-1)
\nn\\&&\hspace{50pt}\times
n^{-1}\sum_{j=1}^n
a_j(q_1)a_j(q_2)a_j(q_3)(\tti\sfz)^3
\nn\\&&
+\sum_{q_1\in\calq}{\colorr1_{\{2\in\calq\}}}{\colorr q_1!}
\sum_{j,k=1}^n
D_{1_k}\big\{(D_{1_k}a_j(q_1))a_j({\colorr q_1})a_k(2)\big\}(\tti\sfz)^3
\nn\\&&
+\sum_{q_1\in\calq}{\colorr1_{\{2\in\calq\}}}{\colorr q_1!}
\sum_{j,k=1}^n
\bigg[\big(2^{-1}D_{1_k}G_\infty(\tti \sfz)^5+D_{1_k}X_\infty(\tti\sfz)^3(\tti\sfx)\big)
\nn\\&&\hspace{120pt}\times
(D_{1_k}a_j(q_1)) a_j({\colorr q_1})a_k(2)  \bigg].
\eea

{\sred Recall (\ref{202003181607}): }
\beas
a^{(3,0)}(t,s,q_1,q_2)
&=&
\ddota(t,s,q_1)a(s,q_2)a(t,2)
+\dota(t,s,q_1)\dota(t,s,q_2)a(t,2)
\nn\\&&
+\dota(t,s,q_1)a(s,q_2)\dot{a}(t,2). 
\eeas
The random symbol ${\mathfrak S}^{(3,0)}(\tti\sfz,\tti\sfx)$ is defined by (\ref{202003181606}) on p.\pageref{202003181606}. 
\begin{en-text}
\beas 
D_{1_k}\big\{(D_{1_k}a_j(q_1))a_j(q_2)a_k(2)\big\}
&=&
n^{-2}a^{(3,0)}_1(t,s,q_1,q_2)
+o_{L^\inftym}(n^{-2})
\eeas
where the term $o_{L^\inftym}(n^{-2})$ can be estimated uniformly in 
$t\in I_k$, $s\in I_j$ and $(k,j)\in\{1,...,n\}^2$. 
\end{en-text}

Let $r>0$. 
\bd
\im[[A$'$\!\!]]$_r$ 
Conditions (\ref{ai}){\sred,}
(\ref{aviii}) and (\ref{aix}) of $[A]$ are satisfied. 
Moreover, the following conditions are fulfilled (in place of (\ref{aiii})-(\ref{avii}) of $[A]$). 
\begin{en-text}
{\bf(i)} 
$a_j(q)\in\bbD^{\infty}$ for $j\in\{1,...,n\}$, $n\in\bbN$ and $q\in\calq$, 
and that 
\bea\label{2020031201810} 
\sup_{n\in\bbN}\sup_{j\in\bbJ_n}
\sup_{t_1,...,t_i\in[0,1]}
\big\|D_{t_i}\cdots D_{t_1}a_j(q)\big\|_{p}&<&\infty
\eea
for every $i\in\bbZ_+=\{0,1,...\}$, $p>1$ and $q\in\calq$. 
\end{en-text}
%
%
\begin{en-text}
Moreover, there exist measurable random fields 
$(a(t,q))_{t\in[0,1]}$, $(\dota(t,s,q))_{(t,s)\in[0,1]^2}$, 
$(\dot{a}(t,2))_{t\in[0,1]}$ 
and 
$(\ddota(t,s,q))_{(t,s)\in[0,1]^2}$ 
$(q\in\calq)$ 
such that 
\end{en-text}
\begin{en-text}
For every $p>1$, 
\bea\label{202003201813} &&
\max_{q\in\calq}\sup_{t\in[0,1]}\big\|a(t,q)\big\|_p
+\max_{q\in\calq}\sup_{(t,s)\in[0,1]^2}\big\|\dota(t,s,q)\big\|_p
\nn\\&&
+\sup_{t\in[0,1]}\big\|\dot{a}(t,2)\big\|_p
+\max_{q\in\calq}\sup_{(t,s)\in[0,1]^2}\big\|\ddota(t,s,q)\big\|_p
\><\> 
\infty.
\eea
\end{en-text}
\bd
\im[(\ref{aiii})] For every $q\in\calq$, 
\bea\label{202003201816}
\sum_{j=1}^n \int_0^1
1_{I_j}(s)\big|
a_j(q) - a(s,q) \big|^r ds 
&\to^p&
0
\eea
as $n\to\infty$. 
\im[(\ref{aiv})] For every $q\in\calq$, 
\bea\label{2020031201805}
\sum_{j,k=1}^n \int_0^1\int_0^1
1_{I_k\times I_j}(t,s)\big|
nD_{1_k}a_j(q) - \dota(t,s,q) \big|^r dsdt 
&\to^p&
0
\eea
as $n\to\infty$. 
\im[(\ref{av})] 
As $n\to\infty$, 
\bea\label{202003241611}
\sum_{k=1}^n \int_0^1
1_{I_k}(t)\big|
nD_{1_k}a_k(2) - \dot{a}(t,2) \big|^r {\sred dt}
&\to^p&
0
\eea
if $2\in\calq$. 
\im[(\ref{avi})] 
For every $q\in\calq$, 
\bea\label{202003201508}
\sum_{j,k=1}^n \int_{[0,1]^2}
1_{I_k\times I_j}(t,s)
\big|n^2D_{1_k}D_{1_k}a_j(q_1)-\ddota(t,s,q)\big|^rdsdt
&\to^p&
0
\eea
as $n\to\infty$. 
\im[(\ref{avii})]
An $n\to\infty$, 
\bea\label{202004170843}
\sum_{j=1}^n\int_0^11_j(t)\big|n^2D_{1_j}D_{1_j}a_j(2)-\ddot{a}(2)\big|^rdt
&\to^p& 
0
\eea
if $2\in\calq$. 


\begin{en-text}
$G_\infty\in\bbD^\infty$ and $\sup_{t\in[0,1]}\|D_tG_\infty\|_p<\infty$ for every $p>1$. 
\im[(vi)] 
$X_\infty\in\bbD^\infty$ and $\sup_{t\in[0,1]}\|D_tX_\infty\|_p<\infty$ for every $p>1$. 
\end{en-text}
\ed
\ed 
\halflineskip
\begin{lemma} The following properties are equivalent: 
\bd\im[(a)] $[A]$ holds. 
\im[(b)] $[A']_r$ holds for some $r>0$. 
\im[(c)] $[A']_r$ holds for all $r>0$. 
\ed
\end{lemma}
\proof 
Let 
\beas 
\varphi_n(t,s)
&=&
\sum_{j,k}1_{I_k\times I_j}(t,s)n^2D_{1_k}D_{1_k}a_j(q_1)-\ddota(t,s,q_1). 
\eeas
Under (\ref{2020031201810a}) and (\ref{202003251205}), it follows from 
Fatou's lemma applied to the measure $dPdsdt$ that 
\bea\label{202005021306}
E\bigg[\int_{[0,1]^2}\big|\ddota(t,s,q_1)\big|^rdsdt\bigg] &<& \infty
\eea
for every $r>1$. 
We have
\beas &&
E\bigg[\sum_{j,k}\int_{[0,1]^2}1_{I_k\times I_j}(t,s)
\big|n^2D_{1_k}D_{1_k}a_j(q_1)-\ddota(t,s,q_1)\big|^rdsdt\bigg]
\\&=&
E\bigg[\int_{[0,1]^2}
\big|\varphi_n(t,s)\big|^rdsdt\bigg]
\\&\leq&
E\bigg[\int_{[0,1]^2}\big|\varphi_n(t,s)\big|^r1_{\{|\varphi_n(t,s)|> A\}}dsdt\bigg]
+
E\bigg[\int_{[0,1]^2}\big|\varphi_n(t,s)\big|^r1_{\{\ep<|\varphi_n(t,s)|\leq A\}}dsdt\bigg]
\\&&
+
E\bigg[\int_{[0,1]^2}\big|\varphi_n(t,s)\big|^r1_{\{|\varphi_n(t,s)|\leq \ep\}}dsdt\bigg]
\\&\leq&
A^{-1}\int_{[0,1]^2}E\big[\big|\varphi_n(t,s)\big|^{r+1}\big]dsdt
+
A^rE\bigg[\int_{[0,1]^2}1_{\{|\varphi_n(t,s)|>\ep\}}dsdt\bigg]
+\ep^r
\eeas
for constants $A>\ep>0$. 
This shows (\ref{202003251205}) implies (\ref{202003201508}) for any $r>0$ 
under (\ref{2020031201810a}) and (\ref{202005021306}). 
Other implications are verified in the same way. 
Therefore we proved $(a)\Iku(c)$. 
It is easy to show $(c)\Iku(b)\Iku(a)$. 
\qed\halflineskip

\begin{en-text}
Suppose that $a_j(q)\in\bbD^{\infty}$ for $j\in\{1,...,n\}$, $n\in\bbN$ and $q\in\calq$, 
and that 
\beas
\sup_{n\in\bbN}\sup_{j\in\bbJ_n}
\sup_{t_1,...,t_i\in[0,1]}
\big\|D^i_{t_1,...,t_i}a_j(q)\big\|_{p}&<&\infty
\eeas
for every $i\in\bbZ_+=\{0,1,...\}$, $p>1$ and $q\in\calq$. 
Recall that $a_j(q)$ depends on $n$. 

Let 
\beas
G_\infty &=& \sum_{q\in\calq}q! \int_0^1 a(q)_t^2dt.
\eeas

Characterization of $a(q)_t$ \koko

Let 
\beas 
\Psi(\sfz,\sfx) &=& \exp\big(2^{-1}G_\infty(\tti\sfz)^2+X_\infty\tti\sfx\big)
\eeas
for $\sfz\in\bbR$ and $\sfx\in\bbR$. 
\end{en-text}

We will show the convergence of $\sum_{j,k}(D_{1_k}D_{1_k}a_j(q_1))a_j(q_2)a_k(2)$ {\colorr when $2\in\calq$}. 
We have 
\beas
\Phi_1&:=&
\bigg|\sum_{j,k}(D_{1_k}D_{1_k}a_j(q_1))a_j(q_2)a_k(2)
- \int_{[0,1]^2}\ddota(t,s,q_1)a(s,q_2)a(t,2)dsdt\bigg| 
\\&\leq& 
\Phi_{1,1}+\Phi_{1,2}+\Phi_{1,3}
\eeas
for 
\beas 
\Phi_{1,1}
&=&
\bigg|\sum_{j,k}\int_{[0,1]^2}1_{I_k\times I_j}(t,s)
\big\{n^2D_{1_k}D_{1_k}a_j(q_1)-\ddota(t,s,q_1)\big\}
a_j(q_2)a_k(2)dsdt\bigg|,
\eeas
\beas 
\Phi_{1,2}
&=&
\bigg|\sum_{j,k}\int_{[0,1]^2}1_{I_k\times I_j}(t,s)
\ddota(t,s,q_1)\big\{a_j(q_2)-a(s,q_2)\big\}
a_k(2)dsdt\bigg|
\eeas
\beas 
\Phi_{1,3}
&=&
\bigg|\sum_{j,k}\int_{[0,1]^2}1_{I_k\times I_j}(t,s)
\ddota(t,s,q_1)a(s,q_2)\big\{a_k(2)-a(t,2)\big\}
dsdt\bigg|.
\eeas
H\"older's inequality gives
\beas
\Phi_{1,1} 
&\leq&
\bigg[\sum_{j,k}\int_{[0,1]^2}1_{I_k\times I_j}(t,s)
\big|n^2D_{1_k}D_{1_k}a_j(q_1)-\ddota(t,s,q_1)\big|^r
dsdt\bigg]^{1/r}
\\&&\times
\bigg[\sum_{j,k}\int_{[0,1]^2}1_{I_k\times I_j}(t,s)
|a_j(q_2)a_k(2)|^{r'}dsdt\bigg]^{1/r'}
\eeas
where $r'=r/(r-1)$ for any $r>1$.

Then, 
by (\ref{202003201508}) and (\ref{2020031201810a}), 
we see 
$\Phi_{1,1}\to^p0$ as $n\to\infty$. 
Similarly, the convergences $\Phi_{1,i}\to^p0$ ($i=2,3$) follow from 
(\ref{202003201816}), 
(\ref{2020031201810a}) and (\ref{202005021306}). 
Therefore, $\Phi_1\to^p0$ as $n\to\infty$. 
Condition (\ref{2020031201810a}) 
ensures boundedness of 
the variables $\{\sum_{j,k}(D_{1_k}D_{1_k}a_j(q_1))a_j(q_2)a_k(2)\}_{n\in\bbN}$ in $L^\inftym$, and so we obtain 
\bea\label{202003241641}
\sum_{j,k}(D_{1_k}D_{1_k}a_j(q_1))a_j(q_2)a_k(2)
&\to& 
\int_{[0,1]^2}\ddota(t,s,q_1)a(s,q_2)a(t,2)dsdt\quad\text{in } L^\inftym
\eea
as $n\to\infty$; 
the limit is in $L^\inftym$ since  $E\big[\int_{[0,1]}|a(s,q_2)|^rds\big]<\infty$ for every $r>1$, 
which follows from Fatou's lemma and the assumed convergence in measure. 
Similarly we obtain 
\bea\label{202003241646}
\sum_{j,k}(D_{1_k}a_j(q_1))(D_{1_k}a_j(q_2))a_k(2)
&\to& 
\int_{[0,1]^2}\dota(t,s,q_1)\dota(t,s,q_2)a(t,2)dsdt\quad\text{in } L^\inftym
\eea
from (\ref{202003201816}) and (\ref{2020031201805}), 
and also 
\bea\label{202003241648}
\sum_{j,k}(D_{1_k}a_j(q_1))a_j(q_2)D_{1_k}a_k(2)
&\to& 
\int_{[0,1]^2}\dota(t,s,q_1)a(s,q_2)\dot{a}(t,2)dsdt\quad\text{in } L^\inftym
\eea
as $n\to\infty$ 
from 
(\ref{202003201816}), (\ref{2020031201805}) and (\ref{202003241611}). 
%
Since
\beas 
D_{1_k}\big\{(D_{1_k}a_j(q_1))a_j(q_2)a_k(2)\big\}
&=&
(D_{1_k}D_{1_k}a_j(q_1))a_j(q_2)a_k(2)
\\&&
+(D_{1_k}a_j(q_1))(D_{1_k}a_j(q_2))a_k(2)
\\&&
+(D_{1_k}a_j(q_1))a_j(q_2)D_{1_k}a_k(2), 
\eeas
it follows from (\ref{202003241641}), (\ref{202003241646}), (\ref{202003241648}) 
and (\ref{202003181607}) that 
\bea\label{202003241656}
\sum_{j,k=1}^n
D_{1_k}\big\{(D_{1_k}a_j(q_1))a_j(q_2)a_k(2)\big\}
&\to&
\int_{[0,1]^2}a^{(3,0)}(t,s,q_1,q_2)dsdt\quad\text{in }L^\inftym
\eea
as $n\to\infty$. 

Similarly, it is possible to prove 
\bea\label{202003251710}
n^{-1}\sum_j a_j(q_1)a_j(q_2)a_j(q_3)
&\to&
\int_0^1 a(q_1)_ta(q_2)_ta(q_3)_tdt\quad\text{in }L^\inftym.
\eea
\begin{en-text}
We apply the H\"older inequality to the integral with respect to 
the probability measure 
\beas
\sum_{j,k}n1_k(t_2)1_k(t_1)1_j(s)dsdt_1dt_2
\eeas
on $[0,1]^3$ to obtain 
\bea\label{2020031201807} 
\Phi_{11}
&\leq&
\bigg[\sum_{j,k}\int_{(s,t_2)\in[0,1]^2}1_k(t_2)1_j(s)
\big|(D_{t_2}a_j(q_2))-\dota(t_2,s,q_2)\big|^r dsdt_2\bigg]^{1/r}
\nn\\&&\times
\bigg[\sum_{j,k}\int_{t_1=0}^1\int_{(s,t_2)\in[0,1]^2}n1_k(t_2)1_k(t_1)1_j(s)
|D_{t_1}a_j(q_1)|^{r'}|a_k(2)|^{r'}dsdt_2dt_1\bigg]^{1/r'}
\eea
where $r'=r/(r-1)$. 
Since the second factor on the right-hand side of (\ref{2020031201807}) is tight uniformly in $n\in\bbN$ 
thanks to the boundedness of the moment ensured by (\ref{2020031201810}) of $[A]$, 
Condition (\ref{2020031201805}) of $[A2]$ implies $\Phi_{1,1}\to^p0$ as $n\to\infty$. 
Similarly, 
with the help of  (\ref{2020031201810}) and (\ref{202003201813}), 
Condition (\ref{2020031201805}) gives $\Phi_{1,2}\to^p0$. 
In quite the same way, under (\ref{202003201813}), 
Condition (\ref{202003201816}) gives $\Phi_{1,3}\to^p0$. 
Therefore, we obtain $\Phi_1\to^p0$ as $n\to\infty$. 
\end{en-text}

Since the mapping $t\to D_tG_\infty(\omega)$ is in $L^2([0,1],dt)$ for almost all $\omega$, 
with approximation in $L^2([0,1],dt)$ by a continuous function if necessary, we see 
\bea\label{202003251650}
\int_{[0,1]}\bigg|\sum_{j=1}^n1_{I_j}(t)nD_{1_j}G_\infty - D_tG_\infty\bigg|dt &\to& 0\quad a.s. 
\eea
as $n\to\infty$; in fact, for a continuous process $g=(g(t))_{t\in[0,1]}$, one has the estimate 
\beas 
\int_{[0,1]}\bigg|\sum_j1_j(t)\big(nD_{1_j}G_\infty-g(t)\big)\bigg|dt
&\leq&
\int_{[0,1]}\big|D_sG_\infty-g(s)\big|ds
+\sup_{s,t\in[0,1]\atop |s-t|\leq 1/n}\big|g(s)-g(t)\big|.
\eeas
From (\ref{202003251650}) in particular,  
\bea\label{202003181406}
E\bigg[\int_0^1 
1_{\big\{\big|\sum_{j=1}^n1_{I_j}(t)nD_{1_j}G_\infty-D_tG_\infty\big|>\ep\big\}} dt\bigg]&\to& 0
\eea
as $n\to\infty$ for every $\ep>0$. 
Therefore, in the present situation, we get 
\bea\label{202003181406}
\sum_{j=1}^n\int_0^1 1_{I_j}(t)\big|nD_{1_j}G_\infty-D_tG_\infty\big|^r dt&\to^p& 0
\eea
as $n\to\infty$ for any $r>0$. 
In the same way, 
\bea\label{202003181406}
E\bigg[\int_0^1 
1_{\big\{\big|\sum_{j=1}^n1_{I_j}(t)nD_{1_j}X_\infty-D_tX_\infty\big|>\ep\big\}} dt\bigg]&\to& 0
\eea
as $n\to\infty$ for every $\ep>0$, 
and in particular, 
\bea\label{202003181406}
\sum_{j=1}^n\int_0^1 1_{I_j}(t)\big|nD_{1_j}X_\infty-D_tX_\infty\big|^r dt&\to^p& 0
\eea
as $n\to\infty$ for any $r>0$. 
Thus, we obtain 
\bea\label{202003251711}
&&
\sum_{j,k=1}^n
\bigg[\big(2^{-1}D_{1_k}G_\infty(\tti \sfz)^5+D_{1_k}X_\infty(\tti\sfz)^3(\tti\sfx)\big)
(D_{1_k}a_j(q_1)) a_j({\colorr q_1})a_k(2)  \bigg]
\nn\\&\to&
\int_{0}^1
\bigg[\bigg(
{\sred 2^{-1}}D_tG_\infty
(\tti \sfz)^5+D_tX_\infty(\tti\sfz)^3(\tti\sfx)\bigg)
\int_0^1\dota(t,s,q_1)a(s,q_2)ds \bigg]a(t,2)dt.
\quad\text{in }L^\inftym
\nn\\&&
\eea
as $n\to\infty$. 
[Though likely $D_tG_\infty=\int_{s=0}^1\sum_{q\in\calq}q!\dota(t,s,q)a(s,q)ds$, 
we do not assume it. ]

Therefore, from (\ref{202003181604}), (\ref{202003181605}), (\ref{202003181606}), 
(\ref{202003241656}), (\ref{202003251710}) and (\ref{202003251711}), 
we obtain 
\bea\label{202003181608}
E\big[\Psi(\sfz,\sfx){\mathfrak S}^{(3,0)}_n(\tti\sfz)\big]
&\to&
E\big[\Psi(\sfz,\sfx){\mathfrak S}^{(3,0)}(\tti\sfz,\tti\sfx)\big]
\eea
as $n\to\infty$, where ${\mathfrak S}^{(3,0)}(\tti\sfz,\tti\sfx)$ is defined by (\ref{202003181606}).

\begin{en-text}
\bd
\im[[A2$'$\!\!]] 
There exist measurable random fields 
$(a(t,q))_{t\in[0,1]}$, $(\dota(t,s,q))_{(t,s)\in[0,1]^2}$, 
$(\dot{a}(t,2))_{t\in[0,1]}$ 
and 
$(\ddota(t,s,q,q'))_{(t,s)\in[0,1]^2}$ 
$(q,q'\in\calq)$ 
such that 
\beas &&
\max_{q\in\calq}\bigg\|\sup_{t\in[0,1]}|a(t,q)|\bigg\|_p
+\max_{q\in\calq}\bigg\|\sup_{(t,s)\in[0,1]^2}|\dota(t,s,q)|\bigg\|_p
\\&&
+\bigg\|\sup_{t\in[0,1]}|\dot{a}(t,2)|\bigg\|_p
+\max_{q,q'\in\calq}\sup_{(t,s)\in[0,1]^2}\big\|\ddota(t,s,q,q')\big\|_p
\><\> 
\infty
\eeas 
 for every $p>1$, and the following conditions hold. 
\bd\im[(i)] For every $q\in\calq$, 
\beas
\sum_{j=1}^n \int_0^1
1_{I_j}(s)\big|
a_j(q) - a(s,q) \big| ds 
&\to^p&
0
\eeas
as $n\to\infty$.
\im[(ii)] For every $q\in\calq$, 
\beas
\sum_{j,k=1}^n \int_0^1\int_0^1
1_{I_k\times I_j}(t,s)\big|
D_ta_j(q) - \dota(t,s,q) \big| dsdt 
&\to^p&
0
\eeas
as $n\to\infty$.
\im[(iii)] 
\beas
\sum_{k=1}^n \int_0^1
1_{I_k}(t)\big|
D_ta_k(2) - \dot{a}(t,2) \big| dsdt 
&\to^p&
0
\eeas
as $n\to\infty$.
\im[(iv)] For every pair $(q_1,q_2)\in\calq^2$ such that $q_1+q_2$ is even, 
\bea\label{202003201508}
\sum_{j,k=1}^n \int_{[0,1]^2}
1_{I_k\times I_k}(t',t)
\big(D_{t'}D_{t}a_j(q_1)\big)a_j(q_2)a_k(2)dtdt' 
&\to^p&
\int_{[0,1]^2}\ddota(t,s,q_1,q_2) dt 
\nn\\&&
\eea
as $n\to\infty$. 
\ed
\ed 
\end{en-text}
\begin{en-text}
The convergence (\ref{202003201508}) has a special flavor than other convergences.  
It can be written as 
\bea\label{202003201517}
\sum_{j,k=1}^n \int_{[0,1]^3}
n1_{I_k\times I_k\times I_j}(t',t,s)
\big(D_{t'}D_{t}a_j(q_1)\big)a_j(q_2)a_k(2)dsdtdt' 
&\to^p&
\int_{[0,1]^2}\ddota(t,s,q_1,q_2)  dsdt . 
\nn\\&&
\eea

The expression (\ref{202003201517}) is more intuitive. 
This suggests a convergence of measures on $[0,1]^3$ to a limit degenerate on 
the set $\{(t,t);t\in[0,1]\}\times[0,1]$. 
In particular, $(t',t)$ does not correspond to $(t,s)$ in the expression (\ref{202003201508}). 
\bd
\im[[X\!\!]] 
$X_\infty\in\bbD^\infty$ and 
\bea\label{202003181406}
\sup_{j\in\{1,...,n\}}\sup_{t\in I_j}\big\|D_{1_j}X_\infty-n^{-1}D_tX_\infty\big\|_p
&=&
o(n^{-1})
\eea
as $n\to\infty$ for every $p>1$. 
\ed
In (\ref{202003181406}), $D_t$ denotes a density function 
representing an element in $\mH$ given by the $D$-derivative of a functional. 
\end{en-text}
\subsubsection{Example}
Let $\calq=\{2\}$, $X_\infty=0$ and let $a\in C^\infty(\bbR)$ 
such that any order of derivative of $a$ is of at most polynomial growth. 
Then $q_1=q_2=q_3=2$ and $\bar{q}(3)=6$ for $q_1,q_2,q_3\in\calq$. 
Let $a_t=a(B_t)$ and $a_j(2)=a(w_\tjm)$ for $j\in\{1,...,n\}$. 
This is the case of a quadratic form of the increments of a Brownian motion 
with predictable weights
\beas 
M_n &=& n^{1/2}\sum_{j=1}^na(w_\tjm)I_2(1_j^{\otimes2}).
\eeas
We remark that the second term of the expression (\ref{202003181159}) vanishes 
in this example because of the strong predictability of $a_j(q)$, while 
such a predictability condition is not assumed in the general context 
in this paper. 

Let 
\beas 
G_\infty &=& \int_0^12a_t^2dt. 
\eeas
Then 
\beas 
\Psi(\sfz,\sfx)\yeq \exp\bigg(-\int_0^1a_t^2dt\bigg)
\yeq \exp\bigg(\int_0^1a_t^2dt(\tti\sfz)^2\bigg).
\eeas 

In this special case, 
\bea\label{202003171656} 
E\big[\Psi(\sfz,\sfx){\mathfrak S}^{(3,0)}_n(\tti\sfz)\big]
&=&
\frac{1}{3}
E\big[\Psi(\sfz,\sfx){\sf qTor}(\tti\sfz)^3\big]
\nn\\&=&
\frac{4}{3}n^{-1}\sum_{j=1}^n
E\big[\Psi(\sfz,\sfx)a_j(2)a_j(2)a_j(2)\big](\tti\sfz)^3
\nn\\&&
+2
\sum_{j,k=1}^nE\big[\Psi(\sfz,\sfx)D_{1_k}\big\{
(D_{1_k}a_j(2))a_j(2)a_k(2)\big\}\big](\tti\sfz)^3
\nn\\&&
+
\sum_{j,k=1}^nE\bigg[
\Psi(\sfz,\sfx)
\big(D_{1_k}G_\infty(\tti \sfz)^5\big)
(D_{1_k}a_j(2)) a_j(2)a_k(2)\big)  \bigg]
\nn\\&&+O(n^{-0.5})
\eea
We can use an approximate projection 
\beas 
\widetilde{\mathfrak S}^{(3,0)}_n(\tti\sfz)
&=& 
\frac{4}{3}n^{-1}\sum_{j=1}^n
a_j(2)a_j(2)a_j(2)(\tti\sfz)^3
\nn\\&&
+2
\sum_{j,k=1}^nD_{1_k}\big\{
(D_{1_k}a_j(2))a_j(2)a_k(2)\big\}(\tti\sfz)^3
\nn\\&&
+
\sum_{j,k=1}^n\bigg[
\big(D_{1_k}G_\infty(\tti \sfz)^5\big)
(D_{1_k}a_j(2)) a_j(2)a_k(2)\big)  \bigg]
\eeas
Then it is not difficult to show from (\ref{202003171656}) that 
\bea\label{202003171724}
E\big[\Psi(\sfz,\sfx){\mathfrak S}^{(3,0)}_n(\tti\sfz)\big]
\yeq
E\big[\Psi(\sfz,\sfx)\widetilde{\mathfrak S}^{(3,0)}_n(\tti\sfz)\big]
+O(n^{-0.5})
&\to&
E\big[\Psi(\sfz,\sfx)\sigma(\tti\sfz) \big]
\eea
where 
\beas 
\sigma(\tti\sfz) &=& 
\frac{4}{3}\int_0^1a_t^3dt(\tti\sfz)^3
+2\int_0^1\int_t^1(aa')'_sa_tdsdt (\tti\sfz)^3
+4\int_0^{{\sred 1}}\bigg(\int_t^1a_sa_s'ds\bigg)^2a_t dt(\tti\sfz)^5.
\eeas
Here $a'_s=a'(w_s)$ and $(aa')'_s=(aa')'(w_s)$. 

The random symbol $\sigma(\tti\sfz)$ coincides 
with the full random symbol $\sigma$ 
(without $\tti\sfx$) of Yoshida (\cite{yoshida2013martingale}, p.918) and also with 
${\mathfrak S}^{3,0}$ 
of Nualart and Yoshida (\cite{nualart2019asymptotic}, p.37). 
In this example, it is known that the full random symbol, which is the sum of 
the adaptive random symbol (tangent) and the anticipative random symbol (torsion), in the martingale expansion 
corresponds to only the quasi-torsion in the theory of asymptotic expansion for Skorohod integrals. 
In other words, the quasi-tangent disappears in this example; see the above papers or Section \ref{202003171756}.

\begin{en-text}
The formula for $3\times(\text{q-}{\sf Tor})$ is:  
\beas
E[\Psi(\sfz)D_{u_n}^2M_n]
&=&
4n^{-0.5}\sum_{j=1}^nn^{-1}
E\big[\Psi(\sfz)a_j^3\big]
\\&&
+
n^{-0.5}2\sum_{j,\ell} E\big[D_{1_\ell}\{\Psi(\sfz)(D_{1_\ell}(a_j^2))a_\ell \}\big]
\\&&
+
n^{-0.5}2\sum_{j,k}E\big[
D_{1_k}\big(\Psi(\sfz)(D_{1_k}a_j) a_ka_j)\big) \big]+O(n^{-1}).
\eeas
We have 
\beas 
n^{1/2}\sum_{j,\ell} E\big[D_{1_\ell}\{\Psi(\sfz)(D_{1_\ell}(a_j^2))a_\ell \}\big]
&=&
\sum_{j,\ell} E\big[D_{1_\ell}\{\Psi(\sfz)2a_ja'_j\alpha_{\tjm,\ell}a_\ell \}\big]
\\&=&
4\sum_{j,\ell} E\bigg[\Psi(\sfz)\int_0^1a_ta'_t\alpha_{t,\ell}dt(\tti\sfz)^2
a_ja'_j\alpha_{\tjm,\ell}a_\ell \bigg]
\\&&
+2\sum_{j,\ell} E\big[\Psi(\sfz)\{(aa')'_j\alpha_{\tjm,\ell}^2a_\ell \}\big]
\\&\to&
4E\bigg[\Psi(\sfz)\int_0^t\bigg(\int_\ell^1a_sa_s'ds\bigg)^2a_t dt\bigg](\tti\sfz)^2
\\&&
+2\sum_{j,\ell} E\bigg[\Psi(\sfz)\int_0^1\int_t^1(aa')'_sa_tdsdt \bigg]
\eeas
where $\alpha_{t,\ell}=\langle 1_{[0,t]},1_\ell\rangle$. 
Also, 
\beas 
n^{1/2}\sum_{j,k}E\big[D_{1_k}\big(\Psi(\sfz)(D_{1_k}a_j) a_ka_j)\big) \big]
&=&
\sum_{j,k}E\big[D_{1_k}\big(\Psi(\sfz)a'_j\alpha_{\tjm,k} a_ka_j)\big) \big]
\\&=&
\sum_{j,k}E\bigg[\Psi(\sfz)\int_0^12a_ta'_t\alpha_{t,k}dt(\tti\sfz)^2a'_j\alpha_{\tjm,k} a_ka_j\bigg](\tti\sfz)^2
\\&&
+\sum_{j,k}E\big[\Psi(\sfz)(a'a)'_j\alpha_{\tjm,k}\alpha_{\tjm,k} a_k \big]
\\&&
+\sum_{j,k}E\big[\Psi(\sfz)\int_0^1a'_ja_j\alpha_{\tjm,k} a'_k\alpha_{t_{k-1},k}\big]
\\&\to&
2E\bigg[\Psi(\sfz)\bigg(\int_t^1a_sa'_sds\bigg)^2a_tdt\bigg](\tti\sfz)^2
\\&&
+E\bigg[\Psi(\sfz)\int_0^1\int_t^1(a'a)'_s a_t dsdt\bigg].
\eeas
Therefore, 
\beas 
\frac{1}{3}n^{1/2}E[\Psi(\sfz)D_{u_n}^2M_n]
&\to&
\frac{4}{3}E\bigg[\Psi(\sfz)\int_0^1a_t^3dt\bigg]
\\&&
+\frac{8}{3}
E\bigg[\Psi(\sfz)\int_0^t\bigg(\int_\ell^1a_sa_s'ds\bigg)^2a_t dt\bigg](\tti\sfz)^2
\\&&
+\frac{4}{3}\sum_{j,\ell} E\bigg[\Psi(\sfz)\int_0^1\int_t^1(aa')'_sa_tdsdt \bigg]
\\&&
+\frac{4}{3}E\bigg[\Psi(\sfz)\bigg(\int_t^1a_sa'_sds\bigg)^2a_tdt\bigg](\tti\sfz)^2
+\frac{2}{3}E\bigg[\Psi(\sfz)\int_0^1\int_t^1(a'a)'_s a_t dsdt\bigg]
\\&=&
\frac{4}{3}E\bigg[\Psi(\sfz)\int_0^1a_t^3dt\bigg]
\\&&
+4
E\bigg[\Psi(\sfz)\int_0^t\bigg(\int_\ell^1a_sa_s'ds\bigg)^2a_t dt\bigg](\tti\sfz)^2
+2E\bigg[\Psi(\sfz)\int_0^1\int_t^1(aa')'_sa_tdsdt \bigg]
\eeas
which coincides 
with the full random symbol $\sigma$  (without $\tti\sfx$) of Yoshida (p.918, SPA2013) and also with 
${\mathfrak S}^{3,0}$ (without $\tti\sfx$) of Nualart and Yoshida (p.37, EJP2019). 
\end{en-text}

\subsection{The quasi-tangent}\label{202003171756}
In this section, we will consider the quasi-tangent. 
For the meantime, we work without Condition (\ref{202004111810}) to clarify the necessity of the assumption. 
Recall the decomposition (\ref{20200317316}) on p.\pageref{20200317316} 
of $D_{u_n(q_2)}M_n(q_1)$. 
\begin{en-text}
Recall $D_{u_n(q_2)}M_n(q_1)=\langle DM_n(q_1),u_n(q_2)\rangle$
\beas 
D_{u_n(q_2)}M_n(q_1)
&=&
D_{u_n(q_2)}\bigg(n^{2^{-1}(q_1-1)}\sum_{j=1}^na_j(q_1) I_{q_1}(1_j^{\otimes q_1})\bigg)
\\&&
-D_{u_n(q_2)}\bigg(n^{2^{-1}(q_1-1)}\sum_{j=1}^n(D_{1_j}a_j(q_1)) I_{q_1-1}(1_j^{\otimes (q_1-1)})\bigg)
\\&=&
n^{2^{-1}(q_1+q_2)-1}\sum_{j,k}q_1a_j(q_1)a_k(q_2)\langle 1_j,1_k\rangle I_{q_1-1}(1_j^{\otimes(q_1-1)}) I_{q_2-1}(1_k^{\otimes(q_2-1)})
\\&&
+n^{2^{-1}(q_1+q_2)-1}\sum_{j,k}(D_{1_k}a_j(q_1)) a_k(q_2) I_{q_1}(1_j^{\otimes q_1})I_{q_2-1}(1_k^{\otimes(q_2-1)})
\\&&
-n^{2^{-1}(q_1+q_2)-1}\sum_{j,k}(q_1-1)(D_{1_j}a_j(q_1)) a_k(q_2) \langle 1_j,1_k\rangle 
I_{q_1-2}(1_j^{\otimes(q_1-2)})I_{q_2-1}(1_k^{\otimes(q_2-1)})
\\&&
-n^{2^{-1}(q_1+q_2)-1}\sum_{j,k}(D_{1_k}D_{1_j}a_j(q_1)) a_k(q_2) I_{q_1-1}(1_j^{\otimes (q_1-1)})I_{q_2-1}(1_k^{\otimes (q_2-1)})
\\&=:&
\bbI_1+\bbI_2+\bbI_3+\bbI_4.
\eeas
\end{en-text}
By the product formula (\ref{202003151619}) on p.\pageref{202003151619}, 
\bea\label{202004121705}
\bbI_1
&=&
n^{2^{-1}(q_1+q_2)-1}\sum_{j,k}q_1a_j(q_1)a_k(q_2)\langle 1_j,1_k\rangle I_{q_1-1}(1_j^{\otimes(q_1-1)}) I_{q_2-1}(1_k^{\otimes(q_2-1)})
\nn\\&=&
\bbJ_{1,1}+\bbJ_{1,2},
\eea
where 
\bea\label{202004121706}
\bbJ_{1,1}
&=&
n^{-1}\sum_{j}1_{\{q_1=q_2\}}q_1!a_j(q_1)^2
\eea
and 
\bea\label{202004121707}
\bbJ_{1,2}
&=&
\sum_{\nu:\bar{q}(2)-2-2\nu>0} c_\nu(q_1-1,q_2-1)
n^{0.5(q_1+q_2)-2-\nu}\sum_{j}q_1a_j(q_1)a_j(q_2)
I_{\bar{q}(2)-2-2\nu}(1_j^{\otimes(\bar{q}(2)-2-2\nu)}).
\nn\\&&
\eea
Since
\bea
e\big(\bbI_1-\bbJ_{1,1}\big)\yeq e(\bbJ_{1,2})\yeq -0.5,
\eea
the effect of $\bbJ_{1,2}$ to the random symbol can remain. 
After IBP, the exponent of $n$ for $E[\Psi(\sfz,\sfx)\bbJ_{1,2}]$ becomes 
\beas &&
\max_{\nu:\bar{q}(2)-2-2\nu>0,\atop\nu\leq\min\{q_1-1,q_2-1\}}\bigg\{\big(0.5\bar{q}(2)-2-\nu\big)+1(\sum_j)-(\bar{q}(2)-2-2\nu)
\bigg\}
\\&=&
\max_{\nu:\bar{q}(2)-2-2\nu>0,\atop\nu\leq\min\{q_1-1,q_2-1\}}\big\{-0.5\bar{q}(2)+1+\nu\big\}
\\&=&
-0.5|q_1-q_2| -1_{\{q_1=q_2\}}
\eeas
The term is negligible when $|q_1-q_2|>1$ or $q_1=q_2$. 
More precisely, remaining is 
Case $[|q_1-q_2|=1]$ giving the terms of order $n^{-1/2}$. 
%
Thus, 
{\colorr
\beas 
E[\Psi(\sfz,\sfx)\bbI_1]
&=&
1_{\{q_1=q_2\}}q_1!n^{-1}\sum_{j}
E\big[\Psi(\sfz)a_j(q_1)^2\big]
\\&&
{\sred +n^{-0.5}1_{\{|q_1-q_2|=1\}}\Xi_n(q_1,q_2)}
+O(n^{-1})
\eeas
{\sred 
with 
\beas
\Xi_n(q_1,q_2) &=& 
\sum_jE\big[D_{1_j}\big\{\Psi(\sfz,\sfx)q_1(q_1\wedge q_2)!a_j(q_1)a_j(q_2)\big\}\big]
\eeas
}
In particular, Condition (\ref{202004111810}) simplifies the formula as 
\beas 
E[\Psi(\sfz,\sfx)\bbI_1]
&=&
1_{\{q_1=q_2\}}q_1!n^{-1}\sum_{j}
E\big[\Psi(\sfz)a_j(q_1)^2\big]
+O(n^{-1}). 
\eeas
}

\begin{en-text}
The sum $\bbI_2$ is expressed by 
\beas 
\bbI_2 
&=&
n^{2^{-1}(q_1+q_2)-2-\nu}\sum_{{\colorr j}}(nD_{1_k}a_j(q_1)) a_k(q_2) 
\sum_\nu c_\nu(q_1,q_2-1)
I_{\bar{q}(2)-1-2\nu}(1_j^{\otimes( \bar{q}(2)-1-2\nu)}).
\eeas
We have an elementary estimate
\bea\label{202004041242}
\sup_{n\in\bbN}\sup_{j,k\in\bbJ_n}\big\|nD_{1_k}a_j(q_1)\big\|_p
&=& 
\sup_{n\in\bbN}\sup_{j,k\in\bbJ_n}\bigg\|n\int_0^1D_ta_j(q_1)1_k(t)dt\bigg\|_p
\nn\\&\leq&
\sup_{n\in\bbN}\sup_{j,k\in\bbJ_n}n\int_0^1\big\|D_ta_j(q_1)\big\|_p1_k(t)dt
\nn\\&\leq&
\sup_{n\in\bbN}\sup_{j\in\bbJ_n}\sup_{t\in[0,1]}\big\|D_ta_j(q_1)\big\|_p
\><\>\infty
\eea
for every $p>1$ under (\ref{2020031201810a}) of $[A]$. 
Since Condition (\ref{202004111810}) bans $q_1=q_2-1$, which entails 
$\bar{q}(2)-1-2\nu>0$, we have 
\bea\label{202004121717}
e(\bbI_2) &=& -1
\eea
because 
\beas
\big(0.5(q_1+q_2)-2-\nu\big)-0.5(\bar{q}(2)-1-2\nu)+1-0.5
&=&
-1.
\eeas
Therefore, $\bbI_2$ is negligible. 
\end{en-text}
{\sred 
Though the simple estimate of the exponent of $\bbI_2$ becomes 
\bea\label{202012300728}
e(\bbI_2) &=& -0.5, 
\eea
this is not enough to show $\bbI_2$ has no effect  to the random symbol. 
The sum $\bbI_2$ is expressed by 
\beas 
\bbI_2 
&=&
\bbK_{2,1}+\bbI_{2,2}
\eeas
where 
\beas 
\bbK_{2,1}
&=& 
n^{2^{-1}(q_1+q_2)-2}\sum_{j,k}(nD_{1_k}a_j(q_1)) a_k(q_2) 
I_{\bar{q}(2)-1}(1_j^{\otimes q_1}\odot 1_k^{\otimes (q_2-1)})
\eeas
and 
\beas
\bbK_{2,2}
&=& 
\sum_{\nu>0} c_\nu(q_1,q_2-1)
\nn\\&&\times
n^{2^{-1}(q_1+q_2)-2}\sum_{j,k}(nD_{1_k}a_j(q_1)) a_k(q_2) 
I_{\bar{q}(2)-1-2\nu}(1_j^{\otimes q_1-\nu}\odot 1_k^{\otimes (q_2-1-\nu)})
\langle 1_j^{\otimes\nu}, 1_k^{\otimes\nu}\rangle_{\mfh^{\otimes\nu}}.
\eeas
We have an elementary estimate
\bea\label{202004041242}
\sup_{n\in\bbN}\sup_{j,k\in\bbJ_n}\big\|nD_{1_k}a_j(q_1)\big\|_p
&=& 
\sup_{n\in\bbN}\sup_{j,k\in\bbJ_n}\bigg\|n\int_0^1D_ta_j(q_1)1_k(t)dt\bigg\|_p
\nn\\&\leq&
\sup_{n\in\bbN}\sup_{j,k\in\bbJ_n}n\int_0^1\big\|D_ta_j(q_1)\big\|_p1_k(t)dt
\nn\\&\leq&
\sup_{n\in\bbN}\sup_{j\in\bbJ_n}\sup_{t\in[0,1]}\big\|D_ta_j(q_1)\big\|_p
\><\>\infty
\eea
for every $p>1$ under (\ref{2020031201810a}) of $[A]$. 
Integration-by-parts gives 
\beas 
E\big[\Psi(\sfz,\sfx)\bbK_{2,1}\big]
&=& 
O\big(n^{(0.5\bar{q}(2)-2)+2-(\bar{q}(2)-1)}\big)
\nn\\&=&
O(n^{-0.5\bar{q}(2)+1})
\yeq 
O(n^{-1}).
\eeas
For $\bbK_{2,2}$, we obtain 
\beas
e\big(\bbK_{2,2}\big)
&=& 
-1.
\eeas
Indeed,  
\beas &&
n^{2^{-1}(q_1+q_2)-2}\sum_{j,k}(nD_{1_k}a_j(q_1)) a_k(q_2) 
I_{\bar{q}(2)-1-2\nu}(1_j^{\otimes q_1-\nu}\odot 1_k^{\otimes (q_2-1-\nu)})
\langle 1_j^{\otimes\nu}, 1_k^{\otimes\nu}\rangle_{\mfh^{\otimes\nu}}
\nn\\&=&
n^{2^{-1}(q_1+q_2)-2-\nu}\sum_j(nD_{1_j}a_j(q_1)) a_j(q_2) 
I_{\bar{q}(2)-1-2\nu}(1_j^{\otimes (\bar{q}(2)-1-2\nu)})
\eeas
for $\nu>0$, and hence 
\beas 
e\big(\bbK_{2,2}\big)
&=&
0.5\bar{q}(2)-2-\nu-0.5\big(\bar{q}(2)-1-2\nu\big)+1-0.5\sfm_1(\bar{q}(2)-1-2\nu)
\nn\\&=&
-0.5-0.5\sfm_1(\bar{q}(2)-1-2\nu)
\nn\\&\leq&-1
\eeas
since Condition (\ref{202004111810}) bans $q_1=q_2-1$, which entails 
$\bar{q}(2)-1-2\nu>0$; 
if $q_1\not=q_2-1$, then $(\bar{q}(2)-1)/2>\min\{q_1,q_2-1\}\geq\nu$. 
\begin{en-text}
, we have 
\bea\label{202004121717}
e(\bbI_2) &=& -1
\eea
because 
\beas
\big(0.5(q_1+q_2)-2-\nu\big)-0.5(\bar{q}(2)-1-2\nu)+1-0.5
&=&
-1.
\eeas
Therefore, $\bbI_2$ is negligible. 
\end{en-text}
Therefore, $\bbK_{2,2}$ is negligible in $L^p$, and 
$\bbK_{2,1}$ is negligible in the random symbol. 
Therefore, $\bbI_2$ is negligible in the random symbol. 
}

For $\bbI_3$, 
\beas
\bbI_3 
&=&
-n^{2^{-1}(q_1+q_2)-2}\sum_{j}(q_1-1)(D_{1_j}a_j(q_1)) a_k(q_2)
I_{q_1-2}(1_j^{\otimes(q_1-2)})I_{q_2-1}(1_j^{\otimes(q_2-1)})
\\&=&
-n^{2^{-1}(q_1+q_2)-2-\nu}\sum_{j}(q_1-1)(D_{1_j}a_j(q_1)) a_k(q_2)
\sum_\nu c_\nu(q_1-2,q_2-1)
I_{\bar{q}(2)-3-2\nu}(1_j^{\otimes(\bar{q}(2)-3-2\nu)})
\eeas
The exponent of $\bbI_3$ is 
\bea\label{202004121721} 
e(\bbI_3) &=&  -1
\eea
{\sred 
since $\bar{q}(2)-3-2\nu\not=0$ by 
$\{(q_1-2)+(q_2-1)\}/2>\min\{q_1-2,q_2-1\}\geq\nu$ due to Condition  (\ref{202004111810}). 
}
The sum $\bbI_3$ is negligible. 

As for the sum 
\beas 
\bbI_4 
&=&
-n^{2^{-1}(q_1+q_2)-3}\sum_{j,k}(n^2D_{1_k}D_{1_j}a_j(q_1)) a_k(q_2) I_{q_1-1}(1_j^{\otimes (q_1-1)})I_{q_2-1}(1_k^{\otimes (q_2-1)}),
\eeas
the exponent is 
\bea\label{202004191359} 
e(\bbI_4) &=& \big(0.5\bar{q}(2)-3\big)-0.5(\bar{q}(2)-2)+2-0.5\times2\yeq -1.
\eea
The sum $\bbI_4$ is also negligible. 

To investigate {\sf qTan}, we need the approximation 
\beas 
\sum_{q_1,q_2\in\calq}
E\big[\Psi(\sfz,\sfx)D_{u_n}M_n\big]
&=&
\sum_{q\in\calq}q!n^{-1}\sum_{j}E\big[\Psi(\sfz,\sfx)a_j(q)^2\big]
\\&&
{\sred 
+\sum_{q_1,q_2\in\calq}n^{-0.5}1_{\{|q_1-q_2|=1\}}\Xi_n(q_1,q_2)
}
+O(n^{-1})
\nn\\&=&
\sum_{q\in\calq}q!n^{-1}\sum_{j}E\big[\Psi(\sfz,\sfx)a_j(q)^2\big]
+O(n^{-1})
\eeas
under Condition (\ref{202004111810}). 
In other words, for 
\beas 
{\mathfrak S}^{(2,0)}_{0,n}(\tti\sfz)
&=&
\half{\sf qTan}(\tti\sfz)^2
\\&=&
\half n^{1/2}\big(D_{u_n}M_n-G_\infty\big)(\tti\sfz)^2, 
\eeas
we have 
\bea\label{202003251918}
E\big[\Psi(\sfz,\sfx){\mathfrak S}^{(2,0)}_{0,n}(\tti\sfz)
\big]
&=& 
\nn\half n^{1/2}E\bigg[\Psi(\sfz,\sfx)\bigg(n^{-1}\sum_j\sum_{q\in\calq}q!a_j(q)^2-G_\infty\bigg)
\bigg](\tti\sfz)^2
\nn\\&&
+O(n^{-1/2})
\eea
under (\ref{202004111810}). 
In particular, it follows from (\ref{202003251918}) that 
\bea\label{202003251919}
E\big[\Psi(\sfz,\sfx){\mathfrak S}^{(2,0)}_{0,n}(\tti\sfz)
\big]
&=& 
o(1),
\eea
therefore, 
\bea
{\mathfrak S}^{(2,0)}_0(\tti\sfz)
&=&
0
\eea
under (\ref{202004111810}), 
if (\ref{202003251920}) is satisfied. 

\begin{remark}\rm 
If Condition (\ref{202004111810}) was not assumed, then 
\beas
E\big[\Psi(\sfz,\sfx){\mathfrak S}^{(2,0)}_{0,n}(\tti\sfz)
\big]
&=& 
\nn\half n^{1/2}E\bigg[\Psi(\sfz,\sfx)\bigg(n^{-1}\sum_j\sum_{q\in\calq}q!a_j(q)^2-G_\infty\bigg)
\bigg](\tti\sfz)^2
\nn\\&&
{\sred 
+\half n^{-0.15}\sum_{q_1,q_2\in\calq}1_{\{|q_1-q_2|=1\}}\Xi_n(q_1,q_2)
}
\nn\\&&
+O(n^{-1/2}).
\eeas
\begin{en-text}
and 
\beas
E\big[\Psi(\sfz,\sfx){\mathfrak S}^{(2,0)}_{0,n}(\tti\sfz)
\big]
&=& 
E\bigg[\Psi(\sfz,\sfx)
\bigg(\sum_{q\in\calq:q+1\in\calq}qq!
\int_0^1 a(s,q)a(s,q+1)ds\bigg)\bigg](\tti\sfz)^2
\nn\\&&
+o(1).
\eeas
\end{en-text}
This suggests the random symbol ${\mathfrak S}^{(2,0)}_{0,n}(\tti\sfz)$ does not vanish 
though it vanished in all the examples in Nualart and Yoshida \cite{nualart2019asymptotic}; 
as a matter of fact, those examples satisfied (\ref{202004111810}). 
\end{remark}
\begin{en-text}
$G_\infty$ is given by 
\beas
G_\infty &=& 
\sum_{q_1\in\calq}
q_1!\int_0^1
a_t(q_1)^2dt
\eeas
{\colorr
\beas 
E[\Psi(\sfz,\sfx)\bbI_1]
&=&
1_{\{q_1=q_2\}}q_1!n^{-1}\sum_{j}
E\big[\Psi(\sfz)a_j(q_1)^2\big]
\\&&
+1_{\{|q_1-q_2|=1\}}n^{-1.5}\sum_jE[\Psi(\sfz,\sfx)q_1(q_1\wedge q_2)!a_j(q_1)a_j(q_2)]
+O(n^{-1})
\eeas
}
\end{en-text}



\subsection{Other random symbols}\label{202005021344}
\subsubsection{${\mathfrak S}_n^{(0,1)}$ and ${\mathfrak S}_{1,n}^{(1,1)}$}
In the setting of Theorem \ref{0112121225}, 
\beas 
{\mathfrak S}^{(0,1)}\yeq0&\text{and}& {\mathfrak S}^{(1,1)}_1\yeq0
\eeas
since ${\mathfrak S}_n^{(0,1)}=0$ and ${\mathfrak S}_{1,n}^{(1,1)}=0$.

\subsubsection{${\mathfrak S}_n^{(1,1)}(\tti\sfz,\tti\sfx)$}\label{20204211726}
By definition 
\beas 
{\mathfrak S}_n^{(1,1)}(\tti\sfz,\tti\sfx)
&=& 
n^{1/2}D_{u_n}X_\infty(\tti\sfz)(\tti\sfx)
\nn\\&=&
\sum_{q\in\calq}\sum_{j=1}^n n^{0.5q-1}(nD_{1_j}X_\infty)a_j(q)I_{q-1}(1_j^{\otimes(q-1)})(\tti\sfz)(\tti\sfx)
\eeas
Then 
\beas 
E\big[\Psi(\sfz,\sfx){\mathfrak S}_n^{(1,1)}(\tti\sfz,\tti\sfx)\big]
&=&
\sum_{q\in\calq}n^{0.5q} \sum_j E\bigg[\Psi(\sfz,\sfx)
(D_{1_j}X_\infty)a_j(q)I_{q-1}(1_j^{\otimes(q-1)})
\bigg](\tti\sfz)(\tti\sfx)
\nn\\&=&
\sum_{q\in\calq}n^{0.5q} \sum_j E\bigg[D^{q-1}_{1_j^{\otimes(q-1)}}\bigg\{
\Psi(\sfz,\sfx)
(D_{1_j}X_\infty)a_j(q)\bigg\}
\bigg](\tti\sfz)(\tti\sfx).
\eeas
Then the exponent of $n$ for this expectation for $q\in\calq$ is 
\beas
0.5q+1(\sum_j)-(q-1)(D^{q-1}_{1_j^{\otimes(q-1}})-1(D_{1_j})
&=&
-0.5q+1. 
\eeas
Since $q\geq2$, the effect of this term to the asymptotic expansion exists only when $q=2$. 
Therefore, 
\bea\label{202004151552}
E\big[\Psi(\sfz,\sfx){\mathfrak S}_n^{(1,1)}(\tti\sfz,\tti\sfx)\big]
&=&
n\sum_j E\bigg[D_{1_j}\bigg\{
\Psi(\sfz,\sfx)
(D_{1_j}X_\infty)a_j(2)\bigg\}
\bigg](\tti\sfz)(\tti\sfx)+O(n^{-0.5})
\nn\\&=&
n \sum_j E\bigg[
\Psi(\sfz,\sfx)
(D_{1_j}D_{1_j}X_\infty)a_j(2)
\bigg](\tti\sfz)(\tti\sfx)+O(n^{-0.5})
\nn\\&&
+n \sum_j E\bigg[
\Psi(\sfz,\sfx)
(D_{1_j}X_\infty)D_{1_j}a_j(2)
\bigg](\tti\sfz)(\tti\sfx)+O(n^{-0.5})
\nn\\&&
+n \sum_j E\bigg[
\Psi(\sfz,\sfx)\big\{2^{-1}D_{1_j}G_\infty(\tti\sfz)^2+D_{1_j}X_\infty\tti\sfx\big\}
(D_{1_j}X_\infty)a_j(2)
\bigg](\tti\sfz)(\tti\sfx)
\nn\\&&
+O(n^{-0.5})
\eea

Like the discussion around (\ref{202003181406}), for any $r\geq1$, 
by path-wise in-$L^r([0,1],dt)$-approximation to the function $t\mapsto D_tG_\infty$ by a continuous function, 
we see 
\beas 
E\bigg[\int_0^1\sum_j1_j(s)\big|nD_{1_j}G_\infty-D_sG_\infty\big|^rds\bigg] &\to& 0.
\eeas
Similarly, 
\beas 
E\bigg[\int_0^1\sum_j1_j(s)\big|nD_{1_j}X_\infty-D_sX_\infty\big|^rds\bigg] &\to& 0.
\eeas
Moreover, 
\beas 
E\bigg[\int_{[0,1]}\sum_j 
1_{I_j}(t)\big|n^2D_{1_j}D_{1_j}X_\infty -\ddot{X}(t)\big|^rdt\bigg]
&\to&
0
\eeas
as $n\to\infty$. 
Therefore, from (\ref{202004151552}), 
\beas 
E\big[\Psi(\sfz,\sfx){\mathfrak S}_n^{(1,1)}(\tti\sfz,\tti\sfx)\big]
&\to&
E\big[\Psi(\sfz,\sfx){\mathfrak S}^{(1,1)}(\tti\sfz,\tti\sfx)\big] 
\eeas
as $n\to\infty$, where ${\mathfrak S}^{(1,1)}(\tti\sfz,\tti\sfx)$ is defined by 
(\ref{202004161636}) on p.\pageref{202004161636}.

\subsubsection{${\mathfrak S}_n^{(1,0)}(\tti\sfz,\tti\sfx)$}
In Theorem \ref{0112121225}, $N_n$ is given by (\ref{202004151750}), and 
\beas 
{\mathfrak S}^{(1,0)}_n(\tti\sfz,\tti\sfx)
&=&
{\mathfrak S}^{(1,0)}_n(\tti\sfz)
\>\equiv\> 
N_n\tti\sfz
\nn\\&=&
\sum_{q\in\calq}n^{q/2}\sum_{j=1}^n(D_{1_j}a_j(q) )I_{q-1}(1_j^{\otimes (q-1)})\tti\sfz.
\eeas
Then 
\beas 
E\big[\Psi(\sfz,\sfx){\mathfrak S}^{(1,0)}_n(\tti\sfz)\big] 
&=&
\sum_{q\in\calq}n^{q/2}\sum_{j=1}^nE\bigg[\Psi(\sfz,\sfx)
(D_{1_j}a_j(q) )I_{q-1}(1_j^{\otimes (q-1)})\bigg]\tti\sfz
\nn\\&=&
\sum_{q\in\calq}n^{q/2}\sum_{j=1}^n
E\bigg[D^{q-1}_{1_j^{\otimes(q-1)}}\big\{\Psi(\sfz,\sfx)
D_{1_j}a_j(q) \big\}\bigg]\tti\sfz
\nn\\&=&
1_{\{2\in\calq\}}\>n\sum_{j=1}^n
E\bigg[D_{1_j}\big\{\Psi(\sfz,\sfx)D_{1_j}a_j(2) \big\}\bigg]\tti\sfz+o(1)
\nn\\&=&
1_{\{2\in\calq\}}
E\bigg[\Psi(\sfz,\sfx)\>n\sum_{j=1}^n\big\{2^{-1}D_{1_j}G_\infty(\tti\sfz)^2+D_{1_j}X_\infty(\tti\sfx)\big\}
D_{1_j}a_j(2) \bigg]\tti\sfz
\nn\\&&
+1_{\{2\in\calq\}}
E\bigg[\Psi(\sfz,\sfx)\>n\sum_{j=1}^nD_{1_j}D_{1_j}a_j(2)\bigg]\tti\sfz+o(1)
\eeas
Therefore, 
\beas 
E\big[\Psi(\sfz,\sfx){\mathfrak S}^{(1,0)}_n(\tti\sfz)\big]
&\to&
E\big[\Psi(\sfz,\sfx){\mathfrak S}^{(1,0)}(\tti\sfz)\big]
\eeas
as $n\to\infty$, where ${\mathfrak S}^{(1,0)}(\tti\sfz)$ is defined by (\ref{202004161637}).

\subsubsection{${\mathfrak S}_{1,n}^{(2,0)}(\tti\sfz,\tti\sfx)$}
In Theorem \ref{0112121225}, $N_n$ is given by (\ref{202004151750}), and 
the exponent of $N_n$ is $e(N_n)=0$. 
Moreover, $q-1\not\in\calq$ under Condition (\ref{202004111810}), 
which gives 
\beas 
e(D_{u_n}N_n) &\leq& -0.5
\eeas
by Proposition \ref{202004091035}. 
In particular, 
\beas
E\big[\Psi(\sfz,\sfx){\mathfrak S}_{1,n}^{(2,0)}\big] 
&=& 
o(1)
\eeas
as $n\to\infty$. Therefore, 
\beas 
{\mathfrak S}_1^{(2,0)} &=& 0.
\eeas

\subsection{Condition $[B]$}
We will verify Condition $[B]$ for $\psi_n=1$, under Condition $[A]$. 

\subsubsection{$[B]$ (i)}\label{202004220324} 
Condition $[B]$ (i) is verified with (\ref{2020031201810a}) 
of $[A]$ for 
$\psi_n=1$ and any pair $(p,p_1)$ of indices satisfying $5p^{-1}+p_1^{-1}\leq1$. 

\subsubsection{Condition (\ref{202003301101}) of $[B]$ (ii)}\label{202004220327}
For later use, we show a stronger estimate. 
\begin{lemma}\label{202004031134} 
Suppose that (\ref{2020031201810a}) of $[A]$ is satisfied. 
Then 
\beas 
\|D^iu_n\|_{p} &=& O(1)
\eeas
as $n\to\infty$ for every $i\in\bbZ_+$ and $p>1$. 
\end{lemma}
\proof 
For $q\in\calq$ and $i\in\bbZ_+$, 
\beas 
D^iu_n(q) 
&=&
D^in^{2^{-1}(q-1)}\sum_{j=1}^na_j(q) I_{q-1}(1_j^{\otimes (q-1)})1_j
\\&=&
\sum_{\alpha=0}^iC_{i,\alpha}
\widetilde{\bbI}(i,\alpha,q)_n
\eeas
where $\widetilde{\bbI}(i,\alpha,q)_n$ is the symmetrization of 
\beas
\bbI(i,\alpha,q)_n
&=&
n^{2^{-1}(q-1)}\sum_{j=1}^n
I_{q-1-\alpha}(1_j^{\otimes (q-1-\alpha)})\big(D^{i-\alpha}a_j(q)\big)\otimes1_j^{\otimes(\alpha+1)}
\eeas
where $C_{i,\alpha}=(q-1)!/(q-1-\alpha)!\>1_{\{q-1-\alpha\geq0\}}$. 
Therefore, 
\beas 
\big|\bbI(i,\alpha,q)_n\big|_{\mH^{\otimes(\alpha+1)}}^2
&=&
n^{q-\alpha-2}\sum_{j=1}^n\big|D^{i-\alpha}a_j(q)\big|_{\mH^{\otimes(i-\alpha)}}^2
 I_{q-1-\alpha}(1_j^{\otimes (q-1-\alpha)})^2
\eeas
Since the hypercontractivity yields 
\beas 
\sup_{j,n}\big\| I_{q-1-\alpha}(1_j^{\otimes (q-1-\alpha)})^2\big\|_p &=& O(n^{-(q-1-\alpha)})
\eeas
for every $p>1$, 
we obtain 
\beas 
\sup_n\big\|\big|\bbI(i,\alpha,q)_n\big|_{\mH^{\otimes(\alpha+1)}}\big\|_p &<& \infty
\eeas
for every $p>1$, which concludes the proof. 
[This is a direct proof. On the other hand, we could switch it on the way by using 
Proposition \ref{202004121624}.] 
\qed\halflineskip

\subsubsection{Condition (\ref{202003301102}) of $[B]$ (ii)}
\begin{lemma}\label{202004041254}
Suppose that (\ref{2020031201810a}) and (\ref{aviii}) of $[A]$ 
{\sred except for (\ref{202003251920})} 
hold and that 
\bea\label{202004041258}
\bigg\|n^{-1}\sum_{q\in\calq}\sum_{j=1}^nq!a_j(q)^2 - G_\infty\bigg\|_p
&=& 
O(n^{-1/2})
\eea
as $n\to\infty$ for every $p>1$; 
{\sred in particular, (\ref{202003251920}) is sufficient for (\ref{202004041258})}. 
Then 
\beas
\big\|G_n^{(2)}\big\|_p &=& O(n^{-1/2})
\eeas
as $n\to\infty$ for every $p>1$. 
\end{lemma}
\proof 
From the decomposition (\ref{20200317316}) on p.\pageref{20200317316}, we have 
\bea\label{202004041259} 
G_n^{(2)} 
&=& 
D_{u_n}M_n-G_\infty
\yeq 
\sum_{q_1,q_2\in\calq}\big(\bbI_1+\bbI_2+\bbI_3+\bbI_4\big)-G_\infty
\eea
where 
$\bbI_i$ ($i=1,...,4$) are given by (\ref{202004041122})-(\ref{202004041125}).
%
The variable $\bbI_1$ has a decomposition 
$\bbI_1=\bbJ_{1,1}+\bbJ_{1,2}$ as in (\ref{202004121705}) 
with $\bbJ_{1,1}$ of (\ref{202004121706}) and $\bbJ_{1,2}$ of (\ref{202004121707}). 
From these representations, 
\bea\label{202004141102}
e(\bbJ_{1,1}) \yeq 0 &\text{and}& e(\bbJ_{1,2}) \yeq -0.5.
\eea
\begin{en-text}
By the product formula, \koko
\bea\label{202004041227}
\bbI_1 &=& 
n^{2^{-1}(q_1+q_2)-2}\sum_jq_1a_j(q_1)a_j(q_2) I_{q_1-1}(1_j^{\otimes(q_1-1)}) I_{q_2-1}(1_j^{\otimes(q_2-1)})
\nn\\&=&
n^{2^{-1}(q_1+q_2)-2}\sum_jq_1a_j(q_1)a_j(q_2) 
\sum_\nu c_\nu(q_1-1,q_2-1)I_{q_1+q_2-2-2\nu}I\big(1_j^{\otimes(q_1+q_2-2-2\nu)}\big)n^{-\nu}
\nn\\&&
\eea
\end{en-text}
Then 
\bea\label{202004041238}
\bbI_1
-1_{\{q_1=q_2\}}n^{-1}\sum_jq_1!a_j(q_1)^2
&=& O_{L^\inftym}(n^{-1/2}). 
\eea
\begin{en-text}
Indeed, Theorem \ref{0110221242} shows the exponent of the term 
on the right-hand side of (\ref{202004041227}) is 
\beas
0.5(q_1+q_2)-2-0.5(q_1+q_2-2-2\nu-1)-\nu
&=&
-0.5
\eeas
whenever $q_1+q_2-2-2\nu>0$, and the only exception is the case $\nu=q_1-1=q_2-1$ if possible. 
\end{en-text}
\begin{en-text}
We have an elementary estimate
\bea\label{202004041242}
\sup_{n\in\bbN}\sup_{j\in\bbJ_n}\big\|nD_{1_j}a_j\big\|_p
&=& 
\sup_{n\in\bbN}\sup_{j\in\bbJ_n}\bigg\|n\int_0^1D_ta_j1_j(t)dt\bigg\|_p
\yleq
\sup_{n\in\bbN}\sup_{j\in\bbJ_n}n\int_0^1\big\|D_ta_j\big\|_p1_j(t)dt
\nn\\&\leq&
\sup_{n\in\bbN}\sup_{j\in\bbJ_n}\sup_{t\in[0,1]}\big\|D_ta_j\big\|_p
\><\>\infty
\eea
for every $p>1$ under (\ref{2020031201810a}) of $[A]$. 
\end{en-text}
\begin{en-text}
Theorem \ref{0110221242} with the estimate (\ref{202004041242}) gives 
\bea\label{202004041239}
\big\|\bbI_2\big\|_p
&=&
\bigg\|n^{2^{-1}(q_1+q_2)-2}\sum_{j,k}\big\{n(D_{1_k}a_j(q_1)) a_k(q_2)\big\} I_{q_1}(1_j^{\otimes q_1})I_{q_2-1}(1_k^{\otimes(q_2-1)})\bigg\|_p
\nn\\&=&
O(n^{-1})
\eea
for every $p>1$. 
\end{en-text}
{\sred 
Equality (\ref{202012300728}) gives 
\bea\label{202004041239}
\bbI_2
&=&
O_{L^\inftym}(n^{-0.5}).
\eea
}
\begin{en-text}
In the same way by Theorem \ref{0110221242} with the estimate (\ref{202004041242}) gives
\bea\label{202004041246}
\big\|\bbI_3\big\|_p
&=&
O(n^{-1/2})
\eea
for every $p>1$; a simple $L^p$-estimate does not work. 
\end{en-text}
In the same way by (\ref{202004121721}), 
\bea\label{202004041246}
\bbI_3
&=&
O_{L^\inftym}(n^{-1}). 
\eea
Here we used Condition (\ref{202004111810}).

The sum $\bbI_4$ can be estimated from (\ref{202004121725}) as 
\bea\label{202004041252}
\bbI_4
&=&
O_{L^\inftym}(n^{-1}). 
\eea
Now we conclude the proof by combining 
(\ref{202004041259}), (\ref{202004041238}), (\ref{202004041239}), (\ref{202004041246}),
(\ref{202004041252}) and {\sred (\ref{202004041258})}.
\qed\halflineskip

\begin{lemma}\label{202004041305}
Suppose that (\ref{2020031201810a}) and (\ref{aviii}) of $[A]$ are satisfied 
Then 
\bea\label{202004041628}
\big\|G_n^{(3)}\big\|_p &=& O(n^{-1/2})
\eea
\end{lemma}
as $n\to\infty$ for every $p>1$. 
\proof 
Recall 
\bea\label{202004131421}
G_n^{(3)}=D_{u_n}G_\infty=\sum_{q\in\calq}D_{u_n(q)}G_\infty. 
\eea
Routinely, 
\bea\label{202004131314}
D_{u_n(q)}G_\infty
&=& 
n^{2^{-1}(q-1)-1}\sum_j(nD_{1_j}G_\infty) a_j(q)I_{q-1}(1_j^{\otimes(q-1)})
\eea
for $q\in\calq$. 
Then 
$e(D_{u_n(q)}G_\infty)=-0.5$,
therefore
\bea\label{202004131307}
e(G_n^{(3)})&=&-0.5. 
\eea
This gives the result. 
\qed\halflineskip

\subsubsection{Condition (\ref{202003301103}) of $[B]$ (ii)}
\begin{lemma}\label{202004041400}
Suppose that (\ref{2020031201810a}) and (\ref{aviii}) of $[A]$ are satisfied.  
Then 
\bea\label{202004041626}
\big\|D_{u_n}G_n^{(2)}\big\|_p &=& O(n^{-1/2})
\eea
as $n\to\infty$ for every $p>1$. 
\end{lemma}
\proof 
We have the decomposition (\ref{202004041259}) of $G_n^{(2)}$ and (\ref{202004121705}), 
and 
the exponents (\ref{202004141102}), (\ref{202004041239}), 
(\ref{202004041246}) and (\ref{202004041252}). 
Therefore, $(G^{(2)}_n)_{n\in\bbN}\in\call$ 
{\sred 
by expressing $\sum_{q_1,q_2\in\calq}\bbJ_{1,1}-G_\infty$ as 
\beas 
\sum_{q_1,q_2\in\calq}\bbJ_{1,1}-G_\infty
&=& 
n^{-1}\sum_{j=1}^n\big(q!a_j(q)^2-G_\infty\big)I_0(1_j^{\otimes0})
\eeas
having the exponent $0$,} and 
then Proposition \ref{202004091035} shows $(D_{u_n}G^{(2)}_n)_{n\in\bbN}\in\call$ 
and $e(D_{u_n}G^{(2)}_n)\leq-1/2$, which completes the proof. 
\qed\halflineskip
\begin{en-text}
We have 
\bea\label{202004041625}
D_{u_n}G_n^{(2)}
&=& 
\sum_{(q_1,q_2,q_3)\in\calq^3}D_{u_n(q_3)}D_{u_n(q_2)}M_n(q_1)
-G_n^{(3)}
\eea
To prove the lemma, we do not need to consider difference between 
$G_n^{(3)}$ and any term coming from 
$\sum_{(q_1,q_2,q_3)\in\calq^3}D_{u_n(q_3)}D_{u_n(q_2)}M_n(q_1)$. 
Fix $(q_1,q_2,q_3,q)\in\calq^4$. 
We have 
\bea\label{202004041617}
D_{u_n(q_3)}D_{u_n(q_2)}M_n(q_1)
&=&
D_{u_n(q_3)}\bbI_1+D_{u_n(q_3)}\bbI_2+D_{u_n(q_3)}\bbI_3+D_{u_n(q_3)}\bbI_4
\nn\\&=&
\bbI_{1,1}+\bbI_{1,1}'+\bbI_{1,2}
+\bbI_{2,1}+\bbI_{2,2}+\bbI_{2,3}
\nn\\&&
+D_{u_n(q_3)}\bbI_3
+D_{u_n(q_3)}\bbI_4
\eea

With the representations 
(\ref{202004041456}) , 
Theorem \ref{0110221242} gives 
\bea\label{202004041619}
\bbI_{1,1}\yeq O_{L^\inftym}(n^{-1/2})&\text{and}& \bbI_{1,1}' \yeq O_{L^\inftym}(n^{-1/2})
\eea
Applying Theorem \ref{0110221242} to $\{\cdots\}$ below, we see that 
\bea\label{202004041538}
 n^{0.5q_3-2.5}
 \sum_j\bigg\{\sum_\ell q_1(nD_{1_\ell}(a_j(q_1)a_j(q_2)))a_\ell(q_3)
I_{q_3-1}(1_\ell^{\otimes(q_3-1)})\bigg\}
&=& 
O_{L^\inftym}(n^{-1/2}), 
\eea
so that 
\bea\label{202004041624}
\bbI_{1,2} &=& O_{L^\inftym}(n^{-1/2})
\eea
from (\ref{202004041510}). 

By applying Theorem \ref{0110221242} to the first sum on the right-hand side of (\ref{202004041539}) 
like (\ref{202004041538}), and also to the second sum, we see the former 
is $O_{L^\inftym}(n^{-1/2})$ and the latter is $O_{L^\inftym}(n^{-1})$. 
Therefore, 
\bea\label{202004041542}
\bbI_{2,1} &=& O_{L^\inftym}(n^{-1/2}).
\eea
In case the principal part of (\ref{20200404552}) exists, the exponent becomes 
\beas 
0.5(\bar{q}(3)-3)-1+1(\sum_k)-1(D_{1_k})-0.5(q_1-1)-\nu
&=& 
-0.5
\eeas
since $\nu=0.5(q_2+q_3-3)$. 
In any case, 
\bea\label{202004041607}
\bbI_{2,2} &=& O_{L^\inftym}(n^{-1/2})
\eea
 because of (\ref{20200404552}).

Consequently, we obtain (\ref{202004041626}) 
by combining (\ref{202004041625}), 
(\ref{202004041617}), (\ref{202004041619}), (\ref{202004041624}), 
(\ref{202004041542}), (\ref{202004041607}), 
(\ref{202004041609}) on p.\pageref{202004041609}, 
(\ref{202003171314}) on p.\pageref{202003171314}, 
(\ref{202004041615}) on p.\pageref{202004041615} and 
(\ref{202004041628}) of Lemma \ref{202004041305}. 
\qed\halflineskip
\end{en-text}

\subsubsection{Condition (\ref{202003301104}) of $[B]$ (ii)}
We use (\ref{ai}) and (\ref{aviii}) of $[A]$. 
In view of (\ref{202004131314}) and (\ref{202004131307}), 
since $q-1\not\in\calq$ under Condition (\ref{202004111810}), 
Proposition \ref{202004091035} concludes that $D_{u_n}G_n^{(3)}\in\call$ and 
\beas
e(D_{u_n}G_n^{(3)}) &\leq& -1
\eeas
which implies (\ref{202003301104}). 

\subsubsection{Condition (\ref{202003301105}) of $[B]$ (ii)}
We use (\ref{ai}) and (\ref{aviii}) of $[A]$. 
From (\ref{202004041259}) and (\ref{202004131405}), we have 
\bea\label{202004131329} 
D_{u_n}G_n^{(2)} 
&=& 
D_{u_n}D_{u_n}M_n-D_{u_n}G_\infty
\nn\\&=&
\sum_{q_1,q_2,q_3\in\calq}\check{\bbI}_1+\tilde{\bbI}-G_n^{(3)}
\eea
where
\beas 
\tilde{\bbI}
&=&
\sum_{q_1,q_2,q_3\in\calq}\big(\hat{\bbI}_1+D_{u_n(q_3)}\bbI_2+D_{u_n(q_3)}\bbI_3+D_{u_n(q_3)}\bbI_4\big).
\eeas

\begin{en-text}
We apply Proposition \ref{202004091035} 
with the properties (\ref{202004131408}), (\ref{202004121717}), (\ref{202004121721}) and (\ref{202004131414}) to obtain 
\beas
e\big(\>\tilde{\bbI}\>\big)
&\leq&
-1
\eeas
and inductively 
\bea\label{202004131447}
e\big(D_{u_n}\tilde{\bbI}\big)
&\leq&
-1.
\eea
\end{en-text}
{\sred 
We apply Proposition \ref{202004091035} 
with the properties (\ref{202004131408}), 
(\ref{202004121721}) and (\ref{202004131414}) to obtain 
\beas
e\big(\hat{\bbI}_1+D_{u_n(q_3)}\bbI_3+D_{u_n(q_3)}\bbI_4\big)
&\leq&
-1
\eeas
and inductively 
\bea\label{202004131447}
e\big(D_{u_n}(\hat{\bbI}_1+D_{u_n(q_3)}\bbI_3+D_{u_n(q_3)}\bbI_4)\big)
&\leq&
-1.
\eea
Moreover, since $e(\bbI_2)=-0.5$ from (\ref{202012300728}), 
Corollary \ref{202012292154} applied under Condition (\ref{202004111810}) concludes 
\bea\label{202012300812}
e\big(D_{u_n}\bbI_2\big) 
&\leq& 
-1. 
\eea
Therefore, 
\bea\label{202004131447}
e\big(D_{u_n}\tilde{\bbI}\>\big)
&\leq&
-1.
\eea
}
Since $G_n^{(3)}$ is expressed by (\ref{202004131421}) and (\ref{202004131314}), 
Proposition \ref{202004091035} yields 
\bea\label{202004131448}
e(D_{u_n}G_n^{(3)}) &\leq& -1
\eea
from (\ref{202004131307}) 
under Condition (\ref{202004111810}).

The sum ${\colorr\check{\bbI}_1}$ consists of the three terms, as in (\ref{202004131406}). 
The component $\bbI_{1,2,1}$ has the expression (\ref{202004131428}) and 
$e(\bbI_{1,2,1})=-0.5$ according to (\ref{202004131429}). 
Then, under Condition (\ref{202004111810}), it follows from Proposition \ref{202004091035}
that 
\bea\label{2020014131433}
e(D_{u_n}\bbI_{1,2,1}) &\leq& -1. 
\eea
The component $\bbI_{1,1,1}$ has the representation (\ref{202004131435}) 
(multiplied by $I_0(1_j^{\otimes0})$) 
and the exponent $-0.5$ by (\ref{202004131436}). 
Since $0\not\in\calq$, Proposition \ref{202004091035} tells 
\bea\label{202004131450}
e(D_{u_n}\bbI_{1,1,1}) &\leq& -1.
\eea
In the same fashion, we can say 
\bea\label{202004131451}
e(D_{u_n}\bbI_{1,1,1}') &\leq& -1.
\eea
from (\ref{202004131443}) and (\ref{202004131444}). 
Thus, from (\ref{202004131329}), (\ref{202004131447}), (\ref{202004131448}), 
(\ref{202004131406}), (\ref{2020014131433}), (\ref{202004131450}) 
and (\ref{202004131451}), we conclude that 
\bea\label{2020014131453}
e(D_{u_n}D_{u_n}G_n^{(2)})
&=& 
-1, 
\eea
and in particular that 
\bea\label{202004131500}
D_{u_n}D_{u_n}G_n^{(2)} &=& O_{L^\inftym}(n^{-1}).
\eea

For $G_n^{(3)}$, we obtain 
\bea\label{2020014131456}
e(D_{u_n}D_{u_n}G_n^{(3)}) &\leq& -1. 
\eea
from (\ref{202004131448}), and in particular, 
\bea\label{202004131501}
D_{u_n}D_{u_n}G_n^{(3)} &=& O_{L^\inftym}(n^{-1}).
\eea
So Condition (\ref{202003301105}) has been verified 
by (\ref{202004131500}) and (\ref{202004131501}).

\subsubsection{Conditions (\ref{202003301106}), (\ref{202003301107}) and (\ref{202003301108}) 
 of $[B]$ (ii)}\label{202004220330}
We use (\ref{ai}) of $[A]$ and the first half of (\ref{aix}) of $[A]$. 
In the same way of deriving (\ref{202004131307}), 
(\ref{202004131448}) and (\ref{2020014131456}), we see 
\bea\label{202004131515}
e(D_{u_n(q)}X_\infty)&=&-0.5.
\eea
and 
\bea\label{202004131519}
e(D_{u_n(q)}^iX_\infty)&=&-1\quad (i=2,3)
\eea
Properties (\ref{202003301106}), (\ref{202003301107}) and (\ref{202003301108}) 
follow from 
(\ref{202004131515}) and (\ref{202004131519}).

\subsubsection{Condition (\ref{202003301109}), (\ref{202003301110}) and (\ref{202003301111})
 of $[B]$ (ii)}
 We use (\ref{ai}) of $[A]$. 
In Theorem \ref{0112121225}, $X_n=X_\infty$, therefore $\dotx_n=0$. 
The veriable $N_n$ is given by 
\bea\label{202004131526} 
N_n
&=& 
\sum_{q\in\calq}n^{q/2}\sum_{j=1}^n(D_{1_j}a_j(q) )I_{q-1}(1_j^{\otimes (q-1)})
\eea
Then the exponent of $N_n$ based on the expression (\ref{202004131526}) is 
\bea\label{202014131536} 
e(N_n) &=& 0
\eea
by definition of the exponent. 
Moreover, Proposition \ref{202004091035} ensures 
\bea\label{202004131548}
e(D_{u_n}N_n) &\leq& -0.5
\eea
under Condition (\ref{202004111810}). 
Inductively from (\ref{202004131548}), 
\bea\label{202014131537}
e(D_{u_n}^iN_n) &\leq& -0.5\quad(i=2,3)
\eea
Now (\ref{202014131536}) and (\ref{202014131537}) verify 
Properties (\ref{202003301109}), (\ref{202003301110}) and (\ref{202003301111}). 

\subsubsection{Condition (\ref{202003301112}) of $[B]$ (ii)}
Condition (\ref{202003301112}) holds obviously because $\psi_n=1$ now.

\subsubsection{Condition $[B]$ (iii) and the full random symbol ${\mathfrak S}(\tti\sfz,\tti\sfx)$}
We have already verified Condition $[B]$ (iii) 
in Sections \ref{202004220553}, \ref{202003171756} and \ref{202005021344}. 
\begin{en-text}
{\color{gray}
Then the full random symbol for Theorem \ref{0112121225} is given by 
\beas
\mathfrak{S}(\tti\sfz,\tti,\sfx)
&=&
\mathfrak{S}^{(3,0)}(\tti\sfz,\tti\sfx) 
+\mathfrak{S}^{(1,1)}(\tti\sfz,\tti\sfx) 
+\mathfrak{S}^{(1,0)}(\tti\sfz,\tti\sfx) 
\eeas
where the random symbols 
$\mathfrak{S}^{(3,0)}(\tti\sfz,\tti\sfx)$, 
$\mathfrak{S}^{(1,1)}(\tti\sfz,\tti\sfx)$ and 
$\mathfrak{S}^{(1,0)}(\tti\sfz,\tti\sfx)$ are given by 
(\ref{202003181606}), (\ref{202004161636}) and 
(\ref{202004161637}), respectively. 
}
\end{en-text}
Now Theorem \ref{202004011550} applied to the present case 
completes the proof of Theorem \ref{0112121225}.

\subsection{An example that satisfies (\ref{2020041446})}\label{202004201231}
Suppose that $A\in C^\infty_p(\bbR)$, the space of smooth functions on $\bbR$ with derivatives of at most polynomial growth. Fix a $q\in\calq$. 
For a Wiener process $w=(w_t)_{t\in[0,1]}$, let $a_j(q)^2=A(w_{1-\tj})$. 
Let $a(t,q)^2=A(w_{1-t})$. 
Then 
\beas 
U_n &:=& \sum_{j=1}^nn^{-1}a_j(q)^2-\int_0^1a(t,q)^2dt
\nn\\&=&
\sum_{j=1}^nn^{-1}A(w_{1-\tj})-\int_0^1A(w_{1-t})dt
\nn\\&=&
\sum_{j=1}^n\int_{I_j}\big\{A(w_{1-\tj})-A(w_{1-t})\big\}dt
\eeas
Let 
\beas 
f_{n-j+1}(s)&=&f_{n-j+1,n}(s)\yeq {\sred -}n\int_{I_j}1_{[1-\tj,1-t]}(s)dt
\eeas
for $j\in\{1,...,n\}$. 
Obviously, 
$\sup_{n\in\bbN}\sup_{j=1,...,n}\sup_{s\in[0.1]}|f_j(s)|\leq1$ and 
we have 
\beas 
\int_{I_j}(w_{1-\tj}-w_{1-t})dt
&=& 
\int_{I_j}I_1\big(1_{[0,1-\tj]}-1_{[0,1-t]}\big)dt
\yeq
{\sred -}\int_{I_j}I_1\big(1_{[1-\tj,1-t]}\big)dt
\nn\\&=&
{\sred -}I_1\bigg(\int_{I_j}1_{[1-\tj,1-t]}dt\bigg)
\yeq
n^{-1}I_1(f_{n-j+1}).
\eeas
With 
\beas 
R_n &=& 
\sum_j\int_{I_j}\int_0^1(1-s)A''\big(w_{1-t}+s(w_{1-\tj}-w_{1-t})\big)ds\>(w_{1-\tj}-w_{1-t})^2dt,
\eeas
we write 
\beas
U_n 
&=& 
\sum_jA'(w_{1-\tj})\int_{I_j}(w_{1-\tj}-w_{1-t})dt+R_n
\nn\\&=&
\bbL_n+R_n
\eeas
where 
\beas 
\bbL_n &=& 
n^{-1}\sum_{j=1}^nA'(w_{1-\tj})I_1(f_{n-j+1})
\yeq
n^{-1}\sum_{j=1}^nA'(w_{\tjm})I_1(f_j).
\eeas
Thne it is not difficult to show 
\beas 
R_n &=& O_{\bbD^\infty}(n^{-1})
\eeas
as $n\to\infty$. 
Since the exponent 
\beas
e(\bbL_n) &=& -1-0.5\times1+1-0.5\times1 \yeq -1, 
\eeas
we obtain 
\beas 
\bbL_n &=& O_{\bbD^\infty}(n^{-1}). 
\eeas
Consequently, 
$U_n=O_{\bbD^\infty}(n^{-1})$ as $n\to\infty$. 
This means (\ref{2020041446}) is satisfied in the present example.

\section{Proof of Theorem \ref{202004171845}}\label{202004171914}
We will verify Condition $[C]$ to apply Theorem \ref{202004011551}. 
Condition $[A^\sharp]$ is assumed now. 

\subsection{$[C]$ (i) and (iv) {\sred (a)}}
$u_n, N_n\in\bbD^\infty$ 
under $[A]$ (\ref{ai}),
$G_\infty\in\bbD^\infty$ under $[A]$ (\ref{aviii}), and 
$X_\infty=X_n\in\bbD^\infty$ uder $[A]$ (\ref{aix}). 
Condition $[C]$ (iv) {\sred (a) is ensured by} $[A^\sharp]$ (II). 

\subsection{$[C]$ (ii)}\label{202004220440}
Since $[A^\sharp]$ is stronger than $[A]$, we can use all results proved in Section \ref{202003141657}. 
Lemma \ref{202004031134} verifies the estimate (\ref{202004011421}). 

According to 
(\ref{202004041259}), 
(\ref{202004121705}). 
(\ref{202004121706}) and 
(\ref{202004121707}), 
$G_n^{(2)}$ admits the following decomposition: 
\beas 
G_n^{(2)} 
&=& 
\widetilde{\bbJ}_1+\widetilde{\bbJ}_2
\eeas
where 
\beas 
\widetilde{\bbJ}_1
&=& 
\sum_{q_1,q_2\in\calq}\bbJ_{1,1}-G_\infty
\yeq
\sum_{q\in\calq}\bigg\{\sum_jn^{-1}{\sred q!}a_j(q)^2-{\sred q!}\int_0^1a(t,q)^2dt\bigg\}
\eeas
and 
\beas 
\widetilde{\bbJ}_2 
&=& 
\sum_{q_1,q_2\in\calq}\big(\bbJ_{1,2}+\bbI_2+\bbI_3+\bbI_4\big).
\eeas
{\sred 
We already know 
\beas 
e(\widetilde{\bbJ}_2) &=& -0.5
\eeas
from (\ref{202004141102}) on p.\pageref{202004141102}, 
(\ref{202012300728}) on p.\pageref{202012300728}, 
(\ref{202004121721}) on p.\pageref{202004121721} and 
(\ref{202004121725}) on p.\pageref{202004121725}, 
therefore, 
\bea\label{202012301910}
\widetilde{\bbJ}_2
&=& 
O_{\bbD^\infty}(n^{-0.5})
\eea
Condition (\ref{2020041446}) of $[A^\sharp]$ (\ref{aviii}$^\sharp$) 
on p.\pageref{2020041446} assumes 
\bea\label{202012301911}
\widetilde{\bbJ}_2
&=& 
O_{\bbD^\infty}(n^{-0.5-\kappa}). 
\eea
Therefore, Proposition \ref{202004121624} implies 
\beas 
G_n^{(2)} &=& O_{\bbD^\infty}(n^{-0.5}),
\eeas
and hence 
Condition (\ref{202004011422}) is verified. 
On the other hand, Condition (\ref{202004011425}) follows from 
the estimate (\ref{2020014131453}). 
}
\begin{en-text}
We already know the exponent of $\widetilde{\bbJ}_2$ is 
\beas 
e\big(\widetilde{\bbJ}_2\big)
&=&
-1
\eeas
from 
(\ref{202004141102}), 
(\ref{202004121717}), 
(\ref{202004121721}) and 
(\ref{202004191359}). 
Therefore 
\beas 
\widetilde{\bbJ}_2
&=&
O_{\bbD^\infty}(n^{-1})
\eeas
by Proposition \ref{202004121624}. 
On the other hand, 
\beas
\widetilde{\bbJ}_1
&=&
O_{\bbD^\infty}(n^{-0.5-\kappa})
\eeas
for some positive number $\kappa$ 
by (\ref{2020041446}) of $[A]$. 
This means 
\beas
G_n^{(2)}&=& O_{\bbD^\infty}(n^{-0.5-\kappa'})
\eeas
for some positive number $\kappa'$. 
So we conclude the conditions (\ref{202004011422}) and (\ref{202004011425}) are met. 
\end{en-text}
\begin{en-text}
{\sred 
We already know the exponent of $\widetilde{\bbJ}_2-\bbI_2$ is 
\beas 
e\big(\widetilde{\bbJ}_2-\bbI_2\big)
&=&
-1
\eeas
from 
(\ref{202004141102}), 
(\ref{202004121721}) and 
(\ref{202004191359}). 
Therefore 
\beas 
\widetilde{\bbJ}_2-\bbI_2
&=&
O_{\bbD^\infty}(n^{-1})
\eeas
by Proposition \ref{202004121624}. 
As for $\bbI_2$, we know (\ref{202012300728}): 
\beas 
e(\bbI_2) &=& -0.5,
\eeas
besides (\ref{202012300812}): 
\beas 
e(D_{u_n}\bbI_2) &=& -1. 
\eeas
On the other hand, 
\beas
\widetilde{\bbJ}_1
&=&
O_{\bbD^\infty}(n^{-0.5-\kappa})
\eeas
for some positive number $\kappa$ 
by (\ref{2020041446}) of $[A]$. 
These estimates entail
\beas
G_n^{(2)}&=& O_{\bbD^\infty}(n^{-0.5})
\eeas
and 
\beas
D_{u_n}G_n^{(2)}&=& O_{\bbD^\infty}(n^{-0.5-\kappa'})
\eeas
for some positive number $\kappa'$. 
So we conclude the conditions (\ref{202004011422}) and (\ref{202004011425}) are met. 
}
\end{en-text}

It follows from (\ref{202004131307}) that 
$G_n^{(3)}=O_{\bbD^\infty}(n^{-1/2})$, therefore (\ref{202004011423}) is satisfied. 
Moreover, we see $D_{u_n}G_n^{(3)}=O_{\bbD^\infty}(n^{-1})$
from (\ref{202004131448}). Thus we verified (\ref{202004011423}) and (\ref{202004011424}). 
Properties (\ref{202004011426}) and (\ref{202004011427}) follow from 
(\ref{202004131515}) and (\ref{202004131519}). 
Properties (\ref{202004011428}) and (\ref{202004011429}) follow from 
(\ref{202014131536}) and (\ref{202014131537}), respectively, 
since $\dotx_n=0$ in the present situation. 
Condition $[C]$ (ii) has been checked. 

\subsection{$[C]$ (iii)}
Now $\beta_x=2$, and $[C]$ (iii) (a) is verified by using $[A^\sharp]$ (i$^\sharp$). 
We already observed the convergences requested by $[C]$ (iii) (b).

\subsection{$[C]$ (iv) {\sred (b)}}\label{202004261803}
For $q\in\calq$, the expression (\ref{202003181159}) gives 
\bea\label{202004261126} 
DM_n(q) &=& D\delta(u_n(q))
\nn\\&=&
\bbI^{(1)}_n(q) +\bbI^{(2)}_n(q) +\bbI^{(3)}_n(q) +\bbI^{(4)}_n(q) 
\eea
where 
\bea
\bbI^{(1)}_n(q)
&=&
n^{2^{-1}(q-1)}\sum_{j=1}^n(Da_j(q)) I_q(1_j^{\otimes q}),
\eea
\bea
\bbI^{(2)}_n(q)
&=&
n^{2^{-1}(q-1)}\sum_{j=1}^na_j(q)q I_{q-1}(1_j^{\otimes (q-1)})1_j,
\eea
\bea
\bbI^{(3)}_n(q)
&=&
-n^{2^{-1}(q-1)}\sum_{j=1}^n\big(D(D_{1_j}a_j(q))\big)I_{q-1}(1_j^{\otimes (q-1)})
\eea
and 
\bea
\bbI^{(4)}_n(q)
&=&
-n^{2^{-1}(q-1)}\sum_{j=1}^n(D_{1_j}a_j(q) )(q-1)I_{q-2}(1_j^{\otimes (q-2)})1_j. 
\eea
For non-zero $v\in\mH$, 
\beas 
\bbI^{(3)}_n(q)[v]
&=& 
-n^{0.5q-1.5}\sum_{j=1}^n\big(D_{|v|_\mH^{-1}v}(nD_{1_j}a_j(q))\big)I_{q-1}(1_j^{\otimes (q-1)})|v|_\mH
\>=:\>
\dot{\bbI}^{(3)}_n(q,v)|v|_\mH
\eeas
The exponents of $\dot{\bbI}^{(3)}_n(q,v)$ is 
\beas 
e\big(\dot{\bbI}^{(3)}_n(q,v)\big)
&=& 
(0.5q-1.5)-0.5(q-1)+1-0.5 
\yeq 
-0.5. 
\eeas
Then Proposition \ref{202004121624} ensures 
\beas 
\big|D^i\dot{\bbI}^{(3)}_n(q,v)\big|_{\mH^{\otimes i}} 
&\leq& 
O_{L^\inftym}(n^{-0.5)})
\eeas
and this estimate is uniform in $v\in\mH$ 
since Proposition \ref{202004121624} is based on Theorem \ref{0110221242}. 
Therefore 
\bea\label{202004261221}
\big|D^i\bbI^{(3)}_n(q)[v]\big|_{\mH^{\otimes i}} 
&\leq& 
O_{L^\inftym}(n^{-0.5})|v|_\mH\quad(v\in\mH)
\eea
where the we can take a factor for $O_{L^\inftym}(n^{-0.5)})$ independently from $v\in\mH$. 
The inequality (\ref{202004261221}) means 
\bea\label{202004261237}
\big|D^i\bbI^{(3)}_n(q)\big|_{\mH^{\otimes (i+1)}} 
&=& 
O_{L^\inftym}(n^{-0.5)})
\eea
for $i\in\bbZ_+$.
In this connection, for $\bbI^{(1)}_n[v]$, $\bbI^{(2)}_n[v]$ and $\bbI^{(4)}_n[v]$, 
the exponents are given by 
\bea\label{202004261338} 
e\big(\bbI^{(1)}_n(q)[v]\big)\yeq0,\quad
e\big(\bbI^{(2)}_n(q)[v]\big)\yeq0,\quad\text{and}\quad
e\big(\bbI^{(4)}_n(q)[v]\big)\yeq-\half\>1_{\{q>2\}},
\eea
respectively, and those sums are not negligible. 
The exponents (\ref{202004261338}) are accounted by the following representations: 
\beas 
\bbI^{(1)}_n(q)[v]
&=&
n^{0.5(q-1)}\sum_{j=1}^n(D_{|v|_\mH^{-1}v}a_j(q)) I_q(1_j^{\otimes q})|v|_\mH,
\eeas
\beas 
\bbI^{(2)}_n(q)[v]
&=&
n^{0.5q-1}\sum_{j=1}^na_j(q)q\big\langle n^{1/2}1_j,|v|_\mH^{-1}v\big\rangle I_{q-1}(1_j^{\otimes (q-1)})|v|_\mH,
\eeas
and 
\beas 
\bbI^{(4)}_n(q)[v]
&=&
-n^{0.5q-2}\sum_{j=1}^n(nD_{1_j}a_j(q) )(q-1)\big\langle n^{1/2}1_j,|v|_\mH^{-1}v\big\rangle
I_{q-2}(1_j^{\otimes (q-2)})|v|_\mH,
\eeas
respectively. 

Let 
\beas
\bbI^*=\sum_{q\in\calq}\bbI^{(2)}_n(q)&\text{ and }& 
\bbI_n^{**}=\sum_{q\in\calq}\big(\bbI^{(1)}_n(q) +\bbI^{(3)}_n(q) +\bbI^{(4)}_n(q)\big). 
\eeas
Then 
\bea\label{202004261518}
\Delta_{(M_n,X_\infty)}
&=&
\det\left(\begin{array}{cc}
\langle DM_n,DM_n\rangle & \langle DM_n,DX_\infty\rangle\\
\langle DM_n,DX_\infty\rangle & \langle DX_\infty,DX_\infty\rangle
\end{array}\right)
\nn\\&=&
\det\left(\begin{array}{cc}
\langle \bbI^*+\bbI^{**},\bbI^*+\bbI^{**}\rangle & \langle \bbI^*+\bbI^{**},DX_\infty\rangle\\
\langle \bbI^*+\bbI^{**},DX_\infty\rangle & \langle DX_\infty,DX_\infty\rangle
\end{array}\right)
\nn\\&=&
\det\left(\begin{array}{cc}
\langle\bbI^{**},\bbI^{**}\rangle & \langle\bbI^{**},DX_\infty\rangle\\
\langle\bbI^{**},DX_\infty\rangle & \langle DX_\infty,DX_\infty\rangle
\end{array}\right)
+
\det\left(\begin{array}{cc}
\langle\bbI^{*},\bbI^{*}\rangle & \langle\bbI^{**},DX_\infty\rangle\\
0 & \langle DX_\infty,DX_\infty\rangle
\end{array}\right)
+{\mathfrak R}_n
\nn\\&&
\eea
where 
\bea\label{202004261519}
{\mathfrak R}_n
&=& 
\det\left(\begin{array}{cc}
0 & \langle\bbI^{**},DX_\infty\rangle\\
\langle\bbI^{*},DX_\infty\rangle & \langle DX_\infty,DX_\infty\rangle
\end{array}\right)
+
\det\left(\begin{array}{cc}
\langle \bbI^*+\bbI^{**},\bbI^*+\bbI^{**}\rangle & \langle\bbI^{*},DX_\infty\rangle\\
\langle\bbI^{*},DX_\infty\rangle+\langle\bbI^{**},DX_\infty\rangle & 0
\end{array}\right)
\nn\\&&
+
\det\left(\begin{array}{cc}
2\langle \bbI^*,\bbI^{**}\rangle & \langle\bbI^{**},DX_\infty\rangle\\
0 &  \langle DX_\infty,DX_\infty\rangle
\end{array}\right). 
\eea
We dared compute the determinant in general way to demonstrate this passage was dimension free 
though the variables are one-dimensional in the present problem. 

It is easy to see the components appearing in the matrices 
on the right-hand side of (\ref{202004261518}) and (\ref{202004261519})
are of $O_{\bbD^\infty}(1)$. 
Now
\bea\label{202004261613} 
\langle \bbI^*_n,DX_\infty\rangle
&=& 
O_{\bbD^\infty}(n^{-0.5})
\eea
since the exponent 
\beas 
e\big(\langle \bbI_n^{(2)}(q),DX_\infty\rangle\big)
&=& 
e\bigg(n^{0.5q-1.5}\sum_{j=1}^na_j(q)(nD_{1_j}X_\infty) q I_{q-1}(1_j^{\otimes (q-1)})\bigg)
\nn\\&=&
(0.5q-1.5)-0.5(q-1)+1-0.5
\yeq
-0.5. 
\eeas
Let $q_1,q_2\in\calq$. 
Then
\bea\label{202004261607}
\langle \bbI_n^{(1)}(q_1), \bbI_n^{(2)}(q_2)\rangle
&=&
O_{\bbD^\infty}(n^{-0.5})
\eea
since 
\beas 
e\big(\langle \bbI_n^{(1)}(q_1), \bbI_n^{(2)}(q_2)\rangle\big)
&=&
e\bigg(
n^{0.5(q_1+q_2)-2}\sum_{j,k=1}^n(nD_{1_k}a_j(q_1))a_k(q_2)q_2I_{q_1}(1_j^{\otimes q_1})I_{q_2-1}(1_k^{\otimes(q_2-1)})
\bigg)
\nn\\&=&
(0.5(q_1+q_2)-2)-0.5(q_1+q_2-1)+2-0.5\times2
\yeq 
-0.5.
\eeas
On the other hand, 
\bea\label{202004261608}
\langle \bbI_n^{(4)}(q_1), \bbI_n^{(2)}(q_2)\rangle
&=&
O_{\bbD^\infty}(n^{-0.5})
\eea
from (\ref{202004261338}) and the estimate 
\beas 
e\big(\langle \bbI_n^{(4)}(2), \bbI_n^{(2)}(q)\rangle\big)
&=&
e\bigg(
n^{0.5q-2}\sum_{j=1}^n(nD_{1_j}a_j(2))a_j(q)q I_{q-1}(1_j^{\otimes (q-1)})\bigg)
\nn\\&=&
(0.5q-2)-0.5(q-1)+1-0.5
\yeq 
-1
\eeas
for $q\in\calq$, if $2\in\calq$. 
Thus, 
\bea\label{202004261614}
\langle \bbI^*,\bbI^{**}\rangle &=& O_{\bbD^\infty}(n^{-0.5})
\eea
from (\ref{202004261237}), (\ref{202004261338}), 
(\ref{202004261607}) and (\ref{202004261608}). 
The first term on the right-hand side of (\ref{202004261518}) is nonnegative. 
So it follows from 
(\ref{202004261518}), (\ref{202004261519}), (\ref{202004261613}) and (\ref{202004261614}) that 
\bea\label{202004261615}
\Delta_{(M_n,X_\infty)}
&\geq& 
\langle\bbI^{*},\bbI^{*}\rangle  \langle DX_\infty,DX_\infty\rangle+{\mathfrak R}_n^*
\eea
for some functional ${\mathfrak R}_n^*\in\bbD^\infty$ satisfying 
\beas 
{\mathfrak R}_n^* 
&=& O_{\bbD^\infty}(n^{-0.5}).
\eeas

For $q_1,q_2\in\calq$, 
\beas 
\langle \bbI^{(2)}_n(q_1),\bbI^{(2)}_n(q_2)\rangle
&=&
n^{0.5(q_1+q_2)-2}\sum_{j=1}^na_j(q_1)q_1a_j(q_2)q_2
I_{q_1-1}(1_j^{\otimes(q_1-1)})I_{q_2-1}(1_j^{\otimes(q_2-1)})
\nn\\&=&
\sum_\nu c_\nu(q_1-1,q_2-1)
\nn\\&&\hspace{10pt}\times
n^{0.5(q_1+q_2)-2-\nu}\sum_{j=1}^na_j(q_1)q_1a_j(q_2)q_2
I_{q_1+q_2-2-2\nu}(1_j^{\otimes(q_1+q_2-2-2\nu)})
\nn\\&=&
1_{\{q_1=q_2\}}
n^{-1}\sum_{j=1}^na_j(q_1)^2q_1q_1!+{\mathfrak R}_n(q_1,q_2)
\eeas
where 
\beas 
{\mathfrak R}_n(q_1,q_2)
&=&
\sum_{\nu:q_1+q_2-2-2\nu>0} c_\nu(q_1-1,q_2-1)
\nn\\&&\hspace{10pt}\times
n^{0.5(q_1+q_2)-2-\nu}\sum_{j=1}^na_j(q_1)q_1a_j(q_2)q_2
I_{q_1+q_2-2-2\nu}(1_j^{\otimes(q_1+q_2-2-2\nu)}).
\eeas
We see 
\bea\label{202004261637}
{\mathfrak R}_n(q_1,q_2)
&=& 
O_{\bbD^\infty}(n^{-0.5})
\eea
because $e\big({\mathfrak R}_n(q_1,q_2)\big)$ is 
\beas &&
\max\bigg\{
\big(0.5(q_1+q_2)-2-\nu\big)-0.5(q_1+q_2-2-2\nu)+1-0.5\>;\>
\nn\\&&\hspace{40pt}
\nu<0.5(q_1+q_2)-1
\bigg\}
\nn\\&=&
-0.5. 
\eeas
Therefore we obtain the representation 
\bea\label{202004261645}
\langle \bbI^*, \bbI^*\rangle 
&=& 
\sum_{q\in\calq}n^{-1}\sum_{j=1}^na_j(q)^2qq!+{\mathfrak R}_n^{**}
\eea
with some functional ${\mathfrak R}_n^{**}=O_{\bbD^\infty}(n^{-0.5})$. 

Let 
\beas 
{\mathfrak G}_n 
&=& 
n^{-1}\sum_{j=1}^n\sum_{q\in\calq}q!a_j(q)^2. 
\eeas
From (\ref{202004261615}) and (\ref{202004261645}), we obtain 
\bea\label{202004261651}
\Delta_{(M_n,X_\infty)}
&\geq& 
2\>{\mathfrak G}_n \langle DX_\infty,DX_\infty\rangle+{\mathfrak R}_n^{***}
\eea
with some functional ${\mathfrak R}_n^{***}$ satisfying 
\bea\label{202004261746}
{\mathfrak R}_n^{***}=O_{\bbD^\infty}(n^{-0.5})
\eea

Let 
\beas 
s_\infty=s_n &=& G_\infty|DX_\infty|_\mH^2.
\eeas
Then $s_\infty\in\bbD^\infty$ and 
$s_\infty^{-1}\in L^\inftym$ {\sred under $[A^\sharp]$}. 
By (\ref{202004261651}), we have 
\beas 
P\big[\Delta_{(M_n,X_\infty)}< s_\infty\big]
&\leq&
P\bigg[
2|{\mathfrak G}_n-G_\infty|\>|DX_\infty|_\mH^2+|{\mathfrak R}_n^{***}|>G_\infty|DX_\infty|_\mH^2
\bigg]
\nn\\&\leq&
E\bigg[\big(G_\infty|DX_\infty|_\mH^2\big)^{-L}
\bigg\{2|{\mathfrak G}_n-G_\infty|\>|DX_\infty|_\mH^2+|{\mathfrak R}_n^{***}|\bigg\}^L
\bigg]
\eeas
for any positive number $L$. 
Then, from (\ref{202004261746}), (\ref{2020041446}) and 
$s_\infty^{-1}\in L^\inftym$ (or 
$[A^\sharp]$ (II)), we obtain 
\beas 
P\big[\Delta_{(M_n,X_\infty)}< s_\infty\big]
&=&
O(n^{-L})
\eeas
as $n\to\infty$ for every $L>0$. 
In this way, Condition $[C]$ (iv) {\sred (b)} has been verified. 
This is the last step of the proof of Theorem \ref{202004171845}.

\section{Proof for Section \ref{202004021507}}\label{202004270654}
\subsection{Proof of Theorem \ref{202004210802}}\label{202004210757}
It suffices to check Condition $[B]$ for $\psi_n=1$ and $(Z_n,X_n)$ in Section \ref{202004021507}, 
that is, $Z_n=M_n+n^{1/2}N_n$, where $N_n$ and $X_n$ are general. 
Under the conditions in $[A]$ except for $[A]$ (\ref{avii}), 
the argument in Sections \ref{202004211724}, \ref{202004220553}, \ref{202003171756} and \ref{20204211726} 
(without consideration for the specific $N_n$) 
is valid. In particular, we have the convergence of 
${\mathfrak S}_n^{(3,0)}$ to ${\mathfrak S}^{(3,0)}$, 
${\mathfrak S}_{0,n}^{(2,0)}$ to ${\mathfrak S}_0^{(2,0)}=0$ 
and 
${\mathfrak S}_n^{(1,1)}$ to ${\mathfrak S}^{(1,1)}$ 
in the sense of (\ref{202004220244}) in $[B]$ (iii). 
The convergences 
${\mathfrak S}_n^{(1,0)}$ to ${\mathfrak S}^{(1,0)}$, 
${\mathfrak S}_n^{(0,1)}$ to ${\mathfrak S}^{(0,1)}$, 
${\mathfrak S}_{1,n}^{(2,0)}$ to ${\mathfrak S}_1^{(2,0)}$ 
and
${\mathfrak S}_{1,n}^{(1,1)}$ to ${\mathfrak S}_1^{(1,1)}$ 
in the sense of (\ref{202004220244}) in $[B]$ (iii) 
are assumed in $[D]$. 
Therefore, Condition $[B]$ (iii) is satisfied. 
Condition $[B]$ (i) is obvious if we combine the observations in Section \ref{202004220324} and 
Condition $[D]$ (II). 
Sections \ref{202004220327}-\ref{202004220330} verify 
the estimates (\ref{202003301101})-(\ref{202003301108}) of $[B]$ (ii), 
and the estimates (\ref{202003301109})-(\ref{202003301111}) are assumed in $[D]$ (II). 
Therefore $[B]$ (ii) has been checked. 
Now we can apply Theorem \ref{202004011550} to prove Theorem \ref{202004210802}.

\subsection{Proof of Theorem \ref{202004211628}}\label{202004211638}
We should check the conditions in $[C]$. 
Condition $[C]$ (i) holds obviously under $[D^\sharp]$. 
Condition $[A^\sharp]$ without $[A^\sharp]$ (I) (\ref{avii}$^\sharp$), which holds under $[D^\sharp]$, 
is stronger than $[A]$ without $[A]$ (\ref{avii}), 
and $[A]$ without $[A]$ (\ref{avii}) is sufficient to follow the argument 
in Section \ref{202004220440} except for the last part concerning 
(\ref{202004011428}) and (\ref{202004011429}). 
However, the last two properties are assumed by $[D^\sharp]$ (II). 
Therefore $[C]$ (ii) is satisfied under $[D^\sharp]$. 
Condition $[C]$ (iii) (a) is satisfied with $[D^\sharp]$ (I) {\sred (i$^\sfx$}) and $[D^\sharp]$ (III). 
We already found $\overline{{\mathfrak T}}_n$ for 
${\mathfrak S}_n^{(3,0)}$, ${\mathfrak S}_{0,n}^{(2,0)}$ and ${\mathfrak S}_n^{(1,1)}$,  
and their convergence in the sense of 
{\sred $[C]$ (iii) (b),}
from the argument in Sections \ref{202004211724}, \ref{202004220553}, \ref{202003171756} and \ref{20204211726}. 
Since Condition $[D^\sharp]$ (III) takes care of the convergence of 
${\mathfrak S}_n^{(1,0)}$, ${\mathfrak S}_n^{(0,1)}$, ${\mathfrak S}_{1,n}^{(2,0)}$ 
and ${\mathfrak S}_{1,n}^{(1,1)}$, 
Condition $[C]$ (iii) is verified. 
Condition $[C]$ (iv) has been verified with $[A^\sharp]$ (II); see Section \ref{202004261803}. 
%
This completes the proof of Theorem \ref{202004211628}.

\section{Proof of Theorem \ref{202005110453}}\label{2020051437}
\subsection{Proof of Lemma \ref{202005051550}}\label{202005141700}

\begin{en-text}
For the meantime, we will consider 
\beas
\wt{\bbV}_n 
&=&
\sum_{j=1}^n\theta_j(\Delta_jX)^2
\eeas
\end{en-text}

\subsubsection{Decomposition of $(\Delta_jX)^2$}
Define $f_j$ and $g_j$ in $\mH$ by 
\beas 
f_j(t) &=& 1_j(t)n(t-\tjm)\quad(t\in\bbR_+)
\eeas
and 
\bea\label{202005060441}
g_j(t) &=& 1_j(t)n(\tj-t)\quad(t\in\bbR_+)
\eea
respectively, for $j\in\bbJ_n$. 
\begin{en-text}
The squared increment $\Delta_jX$ has the following decomposition. 
\bea\label{202005031007}
(\Delta_jX)^2
&=&\koko
b^{[1]}_\tjm I_1(1_j)+b^{[2]}_\tjm h
+2^{-1}b^{[1,1]}_\tjm I_2(1_j^{\otimes2})+6^{-1}b^{[1,1,1]}_\tjm I_3(1_j^{\otimes3})
\nn\\&&
+hb^{[1,2]}_\tjm I_1(f_j)
+hb^{[2,1]}_\tjm I_1(g_j)+R_j^{(\ref{202005030151})}
\eea
\end{en-text}

The increment $\Delta_jX$ has a decomposition 
\bea\label{202005030224}
\Delta_jX
&=&
b^{[1]}_\tjm I_1(1_j)+b^{[2]}_\tjm h
+2^{-1}b^{[1,1]}_\tjm I_2(1_j^{\otimes2})+6^{-1}b^{[1,1,1]}_\tjm I_3(1_j^{\otimes3})
\nn\\&&
+hb^{[1,2]}_\tjm I_1(f_j)
+hb^{[2,1]}_\tjm I_1(g_j)+R_j^{(\ref{202005030151})}
\eea
where 
\bea\label{202005030151}
R_j^{(\ref{202005030151})}
&=&
\int_\tjm^\tj\int_\tjm^t\int_\tjm^s\int_\tjm^r b^{[1,1,1,1]}_pdw_pdw_rdw_sdw_t
\nn\\&&
+\int_\tjm^\tj\int_\tjm^t\int_\tjm^s\int_\tjm^r b^{[1,1,1,2]}_pdp\>dw_rdw_sdw_t
\nn\\&&
+\int_\tjm^\tj\int_\tjm^t\int_\tjm^s b^{[1,1,2]}_rdrdw_sdw_t
\nn\\&&
+\int_\tjm^\tj\int_\tjm^t\int_\tjm^s b^{[1,2,1]}_rdw_rdsdw_t
+\int_\tjm^\tj\int_\tjm^t\int_\tjm^s b^{[1,2,2]}_rdrdsdw_t
\nn\\&&
+\int_\tjm^\tj\int_\tjm^t\int_\tjm^s b^{[2,1,1]}_rdw_rdw_sdt
+\int_\tjm^\tj\int_\tjm^t\int_\tjm^s b^{[2,1,2]}_rdrdw_sdt
\nn\\&&
+\int_\tjm^\tj\int_\tjm^tb^{[2,2]}_sdsdt.
\eea
It is easy to see 
\bea\label{202005030344}
R_j^{(\ref{202005030151})}
&=&
\ol{O}_{\bbD^\infty}(n^{-2}),
\eea
that is, 
\beas 
\sup_{n\in\bbN}\sup_{j\in\bbJ_n} \big\|R_j^{(\ref{202005030151})}\big\|_{s,p}
&=& 
O(n^{-2})
\eeas
as $n\to\infty$ for any $s\in\bbR$ and $p>1$. 

From (\ref{202005030224}), 
\bea\label{202005030225}
(\Delta_jX)^2 
&=&
P_j^{(\ref{202005030231})}+R_j^{(\ref{202005030232})}
\eea
where 
\bea\label{202005030231}
P_j^{(\ref{202005030231})}
&=& 
(b^{[1]}_\tjm)^2 I_1(1_j)^2+2hb^{[1]}_\tjm b^{[2]}_\tjm I_1(1_j) 
{\colorr+b^{[1]}_\tjm b^{[1,1]}_\tjm  I_1(1_j)I_2(1_j^{\otimes2})}
\nn\\&&
+2hb^{[1]}_\tjm b^{[1,2]}_\tjm I_1(1_j)I_1(f_j)
+2hb^{[1]}_\tjm b^{[2,1]}_\tjm I_1(1_j)I_1(g_j)
\nn\\&&
+h^2(b^{[2]}_\tjm)^2 
+2^{-2}(b^{[1,1]}_\tjm)^2 I_2(1_j^{\otimes2})^2
\eea
and 
\bea\label{202005030232}
R_j^{(\ref{202005030232})}
&=& 
3^{-1}b^{[1]}_\tjm b^{[1,1,1]}_\tjm I_1(1_j)I_3(1_j^{\otimes3})
+2b^{[1]}_\tjm I_1(1_j)R_j^{(\ref{202005030151})}
\nn\\&&
+hb^{[2]}_\tjm b^{[1,1]}_\tjm I_2(1_j^{\otimes2})
+3^{-1}hb^{[2]}_\tjm b^{[1,1,1]}_\tjm I_3(1_j^{\otimes3})
\nn\\&&
+2h^2b^{[2]}_\tjm b^{[1,2]}_\tjm I_1(f_j)
+2h^2b^{[2]}_\tjm b^{[2,1]}_\tjm I_1(g_j)
+2h^2b^{[2]}_\tjm R_j^{(\ref{202005030151})}
\nn\\&&
+6^{-1}b^{[1,1]}_\tjm I_2(1_j^{\otimes2}) b^{[1,1,1]}_\tjm I_3(1_j^{\otimes3})
+b^{[1,1]}_\tjm I_2(1_j^{\otimes2}) hb^{[1,2]}_\tjm I_1(f_j)
\nn\\&&
+b^{[1,1]}_\tjm I_2(1_j^{\otimes2}) hb^{[2,1]}_\tjm I_1(g_j)
+b^{[1,1]}_\tjm I_2(1_j^{\otimes2}) R_j^{(\ref{202005030151})}
\nn\\&&
+\bigg\{6^{-1}b^{[1,1,1]}_\tjm I_3(1_j^{\otimes3})
+hb^{[1,2]}_\tjm I_1(f_j)
+hb^{[2,1]}_\tjm I_1(g_j)+R_j^{(\ref{202005030151})}
\bigg\}^2.
\eea

\subsubsection{A weighted sum of squares}
Suppose that the double sequence $\theta=(\theta_{n,j})_{j\in\bbJ_n,\>n\in\bbN}$ 
in $\bbD^\infty$ satisfies 
\bea\label{202005051142}
\sup_{n\in\bbN}\sup_{j\in\bbJ_n}\sup_{s_1,...,s_i\in[0,1+\eta]}
\|D_{s_i}\cdots D_{s_1}\theta_{n,j} \|_p&<&\infty
\eea
for every $i\in\bbZ_+$ and $p>1$. 
We will write $\theta_j=\theta_{n,j}$ simply.

We see that 
\begin{en-text}
\bea\label{202005030345} &&
e\bigg(\sum_j\theta_j b^{[1]}_\tjm b^{[1,1]}_\tjm  I_1(1_j)I_2(1_j^{\otimes2})\bigg)
\nn\\&=&
\max_{\nu=0,1}\big\{(0-\nu)-0.5(3-2\nu)+1-0.5\times1_{\{3-2\nu>0\}}\big\}
\yeq -1,
\eea
\end{en-text}
\bea\label{202005030346} &&
e\bigg(\sum_j\theta_j 3^{-1}b^{[1]}_\tjm b^{[1,1,1]}_\tjm I_1(1_j)I_3(1_j^{\otimes3})\bigg)
\nn\\&=&
\max_{\nu=0,1}\big\{(0-\nu)-0.5(4-2\nu)+1-0.5\times1_{\{4-2\nu>0\}}\big\}
\yeq -1.5
\eea
with the product formula, and also that 
\bea\label{202005030347} 
e\bigg(\sum_j\theta_j hb^{[2]}_\tjm b^{[1,1]}_\tjm I_2(1_j^{\otimes2})\bigg)
&=&
-1-0.5\times2+1-0.5 
\yeq -1.5.
\eea
From (\ref{202005030232}), (\ref{202005030344}), 
(\ref{202005030346}), 
and (\ref{202005030347}), we obtain 
\bea\label{202005030348}
\sum_j\theta_j R_j^{(\ref{202005030232})}
&=& 
O_{\bbD^\infty}(n^{-1.5}).
\eea

On the other hand, 
\bea\label{202005030356}
P_j^{(\ref{202005030231})}
&=&
P_j^{(\ref{202005030357})}
+F_j^{(\ref{202005030358})}+S_j^{(\ref{202005030359})}+R_j^{(\ref{202005030413})}
\eea
where 
\bea\label{202005030357}
P_j^{(\ref{202005030357})}
&=&
(b^{[1]}_\tjm)^2 h ,
\eea
\begin{en-text}
\bea\label{202005030358}
F_j^{(\ref{202005030358})}
&=&
(b^{[1]}_\tjm)^2 I_2(1_j^{\otimes2}),
\eea
\end{en-text}
\bea\label{202005030359}
S_j^{(\ref{202005030359})}
&=&
{\colorr b^{[1]}_\tjm b^{[1,1]}_\tjm  I_3(1_j^{\otimes3})
+2hb^{[1]}_\tjm b^{[1,1]}_\tjm  I_1(1_j)}
\nn\\&&
+2hb^{[1]}_\tjm b^{[2]}_\tjm I_1(1_j) 
+2hb^{[1]}_\tjm b^{[1,2]}_\tjm \langle1_j,f_j\rangle
+2hb^{[1]}_\tjm b^{[2,1]}_\tjm \langle1_j,g_j\rangle
\nn\\&&
+h^2(b^{[2]}_\tjm)^2 
+2^{-1}h^2(b^{[1,1]}_\tjm)^2
\eea
and 
\bea\label{202005030413}
R_j^{(\ref{202005030413})}
&=&
2hb^{[1]}_\tjm b^{[1,2]}_\tjm I_2(1_j\odot f_j)
+2hb^{[1]}_\tjm b^{[2,1]}_\tjm I_2(1_j\odot g_j)
\nn\\&&
+2^{-2}(b^{[1,1]}_\tjm)^2 I_4(1_j^{\otimes4})
+h(b^{[1,1]}_\tjm)^2 I_2(1_j^{\otimes2}). 
\eea
Simply 
\bea\label{202005040206}
F_j^{(\ref{202005030358})}
&=& 
\ol{O}_{\bbD^\infty}(n^{-1})
\eea
and 
\bea\label{2020050402067}
S_j^{(\ref{202005030359})}
&=& 
\ol{O}_{\bbD^\infty}(n^{-1.5})
\eea

Since
\beas
e\bigg(\sum_j\theta_jR_j^{(\ref{202005030413})}\bigg)
&=& 
\max\big\{-1.5,-1.5,-1.5,-1.5\big\}
\yeq -1.5, 
\eeas
we obtain 
\bea\label{202005030428}
\sum_j\theta_jR_j^{(\ref{202005030413})}
&=& 
O_{\bbD^\infty}(n^{-1.5}). 
\eea
Related with this, we also see 
\bea\label{202005031744} 
e\bigg(\sum_j\theta_jF_j^{(\ref{202005030358})}\bigg)
&=&
-0.5
\eea
and 
\bea\label{202005031745} 
e\bigg(\sum_j\theta_jS_j^{(\ref{202005030359})}\bigg)
&=&
\max\{-1,-1,-1,-1,-1,-1,-1
\} 
\yeq -1.
\eea

We observe 
\bea\label{202005041132}
\sum_j\theta_j\int_\tjm^\tj\beta_tdt-\sum_j\theta_j\beta_\tjm h
&=&
\sum_j  \theta_j\beta^{[1]}_\tjm hI_1(g_j)
+\sum_j  \theta_j 2^{-1}h^2{\colorr \beta^{[2]}_\tjm}+R_n^{(\ref{202005031358})}
\eea
where 
\bea\label{202005031358}
R_n^{(\ref{202005031358})}
&=& 
\sum_j \theta_j\int_\tjm^\tj (t_j-s)(\beta^{[1]}_s-\beta^{[1]}_\tjm)dw_s
+\sum_j \theta_j\int_\tjm^\tj\int_\tjm^t(\beta^{[2]}_s-\beta^{[2]}_\tjm)dsdt.
\nn\\&&
\eea
Indeed, 
\beas
\int_\tjm^\tj(\beta_t-\beta_\tjm)dt
&=&
\int_\tjm^\tj \bigg\{\int_\tjm^t\beta^{[1]}_sdw_s+\int_\tjm^t\beta^{[2]}_sds\bigg\}dt
\nn\\&=&
\int_\tjm^\tj (t_j-s)\beta^{[1]}_sdw_s
+\int_\tjm^\tj\int_\tjm^t\beta^{[2]}_sdsdt
\nn\\&=&
\beta^{[1]}_\tjm\int_\tjm^\tj (t_j-s)dw_s
+\int_\tjm^\tj (t_j-s)(\beta^{[1]}_s-\beta^{[1]}_\tjm)dw_s
\nn\\&&
+\beta^{[2]}_\tjm\int_\tjm^\tj\int_\tjm^tdsdt
+\int_\tjm^\tj\int_\tjm^t(\beta^{[2]}_s-\beta^{[2]}_\tjm)dsdt.
\eeas
The sequence 
\beas&&
\sum_j \theta_j\int_\tjm^\tj (t_j-s)(\beta^{[1]}_s-\beta^{[1]}_\tjm)dw_s
\nn\\&=&
\sum_j \theta_j\beta^{[1,1]}_\tjm\int_\tjm^\tj (t_j-s)\int_\tjm^sdw_rdw_s
+\sum_j \theta_j\int_\tjm^\tj (t_j-s)\int_\tjm^s(\beta^{[1,1]}_r-\beta^{[1,1]}_\tjm)dw_rdw_s
\nn\\&&
+\sum_j \theta_j\int_\tjm^\tj (t_j-s)\int_\tjm^s\beta^{[1,2]}_rdrdw_s
\nn\\&=&
O_{\bbD^\infty}(n^{-1.5})
\eeas
since the exponent of the first term on the right-hand side is $-1.5$. 
Therefore 
\bea\label{202005041233}
R_n^{(\ref{202005031358})}
&=& 
O_{\bbD^\infty}(n^{-1.5}).
\eea

Let 
\beas 
\bbV_n(\theta) &=& \sum_{j=1}^n\theta_j(\Delta_jX)^2
\eeas
for the sequence $\theta$ satisfying (\ref{202005051142}). 
For a generic random variable $\calv$, 
From (\ref{202005030225}), (\ref{202005041132}) and (\ref{202005030356}), 
the difference $\bbV_n - \calv$ admits the expansion 
\bea\label{202005041239}
\bbV_n(\theta) - \calv
&=& 
\bigg(\sum_{j=1}^n\theta_j(\Delta_jX)^2-\sum_{j=1}^n\theta_j\int_\tjm^\tj\beta_tdt\bigg)
+\bigg(\sum_{j=1}^n\theta_j\int_\tjm^\tj\beta_tdt -\calv\bigg)
\nn\\&=& 
\bigg(\sum_{j=1}^n\theta_jP_j^{(\ref{202005030231})}
-\sum_{j=1}^n\theta_j\beta_\tjm h
-\sum_{j=1}^n  h\theta_j\beta^{[1]}_\tjm I_1(g_j)
-\sum_{j=1}^n  \theta_j 2^{-1}h^2{\colorr \beta^{[2]}_\tjm}
\bigg)
\nn\\&& 
+\bigg(\sum_{j=1}^n\theta_j\int_\tjm^\tj\beta_tdt -\calv\bigg)
+\bigg(\sum_{j=1}^n\theta_j R_j^{(\ref{202005030232})}-R_n^{(\ref{202005031358})}\bigg)
\begin{en-text}
\nn\\&=& 
\sum_j\theta_jF_j^{(\ref{202005030358})}
+\bigg(\sum_j\theta_jS_j^{(\ref{202005030359})}
-\sum_j  h\theta_j\beta^{[1]}_\tjm I_1(g_j)
-\sum_j  \theta_j 2^{-1}h^2
\bigg)
\nn\\&& 
+\bigg(\sum_j\theta_j\int_\tjm^\tj\beta_tdt -\bbV_\infty\bigg)
\nn\\&& 
+\bigg(\sum_{j=1}^n\theta_jR_j^{(\ref{202005030232})}-R_n^{(\ref{202005031358})}
+\sum_j\theta_jR_j^{(\ref{202005030413})}\bigg)
\end{en-text}
\nn\\&=& 
\sum_{j=1}^n\theta_jF_j^{(\ref{202005030358})}
+\sum_{j=1}^n\theta_jS_j^{(\ref{202005041155})}
+\bigg(\sum_{j=1}^n\theta_j\int_\tjm^\tj\beta_tdt -\calv\bigg)
+R_n^{(\ref{202005041156})}
\eea
for 
\beas
S_j^{(\ref{202005041155})}
&=&
S_j^{(\ref{202005030359})}-h\beta^{[1]}_\tjm I_1(g_j)- 2^{-1}h^2{\colorr \beta^{[2]}_\tjm}
\eeas
and 
\bea\label{202005041156}
R_n^{(\ref{202005041156})}
&=&
\sum_{j=1}^n\theta_jR_j^{(\ref{202005030232})}-R_n^{(\ref{202005031358})}
+\sum_{j=1}^n\theta_jR_j^{(\ref{202005030413})}.
\eea
We can say 
\bea\label{202005041828}
\sum_{j=1}^n\theta_jF_j^{(\ref{202005030358})}
&=&
O_{\bbD^\infty}(n^{{\colorr -0.5}})
\eea
and 
\bea\label{202005041742}
\sum_{j=1}^n\theta_jS_j^{(\ref{202005041155})}
&=& 
O_{\bbD^\infty}(n^{-1})
\eea
by using (\ref{202005031744}) and (\ref{202005031745}), respectively, 
and also 
\bea\label{202005041236}
R_n^{(\ref{202005041156})}
&=& 
O_{\bbD^\infty}(n^{-1.5})
\eea
from (\ref{202005030348}), (\ref{202005041233}) and (\ref{202005030428}).

\subsubsection{Expansion of $L_{n,j}$ and $U_n$}
For each $j\in\bbJ_n$, applying the formula (\ref{202005041239}) to 
\beas 
\theta_k &=& 
{\colorb
\eta_{n,j}^{-1}1_{\{k\in K_j\}}
}
\eeas
and $\calv=L_{\infty,\tjm}$, we obtain 
\bea\label{202005041312}
L_{n,j} - L_{\infty,\tjm}
&=& 
\sum_{{\colorb k\in K_j}}{\colorb\eta_{n,j}}^{-1}F_k^{(\ref{202005030358})}
+\sum_{{\colorb k\in K_j}}{\colorb\eta_{n,j}}^{-1}S_k^{(\ref{202005041155})}
{\sred +\>\cale_j}
+R_{n,j}^{(\ref{202005041326})}
\eea
where 
{\sred 
\beas
\cale_j
&=& 
\sum_{k\in K_j}\eta_{n,j}^{-1}\int_\tkm^\tk\beta_tdt -L_{\infty,\tjm}
\eeas
and
}
$R_{n,j}^{(\ref{202005041326})}\in\bbD^\infty$ such that 
uniformly in $j\in\bbJ_n$, 
\bea\label{202005041326}
R_{n,j}^{(\ref{202005041326})}
&=& 
\ol{O}_{\bbD^\infty}(n^{-1.5}).
\eea 
{\sred 
[We remark that  
\beas 
\cale_j
&=&
\sum_{k\in K_j}\eta_{n,j}^{-1}\int_\tkm^\tk\beta_tdt -L_{\infty,\tjm}
\nn\\&=&
\sum_{k\in K_j}\eta_{n,j}^{-1}\int_\tkm^\tk(\beta_t-\beta_\tjm)dt 
-\eta_{\infty,\tjm}^{-1}\int_{(\tjm-\lambda)\vee0}^{(\tjm+\lambda)\wedge1}(\beta_t-\beta_\tjm)dt 
\nn\\&=&
\sum_{k\in K_j}\big(\eta_{n,j}^{-1}-\eta_{\infty,\tjm}^{-1}\big)\int_\tkm^\tk(\beta_t-\beta_\tjm)dt 
\nn\\&&
+\eta_{\infty,\tjm}^{-1}\bigg\{
-\int_{(\tjm+\lfloor n\lambda\rfloor/n)\wedge1}^{(\tjm+\lambda)\wedge1}(\beta_t-\beta_\tjm)dt
+\int_{(\tjm-(\lfloor n\lambda\rfloor-1)/n)\vee0}^{(\tjm-\lambda)\vee0}(\beta_t-\beta_\tjm)dt
\bigg\}
\nn\\&=&
\ol{O}_{\bbD^\infty}(n^{-1})
\eeas
since $\liminf_{n\to\infty}\inf_{j}\eta_{n,j}>0$ and $\inf_{t}\eta_{\infty,t}>0$. 
Thus, $\cale_j$ is not negligible in $L^p$ and in $\bbD^\infty$.] 
}

\begin{en-text}
$\theta_j=\eta^{-1}$ and $\theta_j=1$, respectively, 
we see that each of $L_{n,j}-L_{\infty,\tjm}$ and $U_n-U_\infty$ 
admits a stochastic expansion similar to (\ref{202005041239}) as follows: 
\bea\label{202005041312}
L_{n,j}-L_{\infty,\tjm}
&=& 
\sum_{k=j}^{j+\lfloor n\eta\rfloor}\eta^{-1}F_k^{(\ref{202005030358})}
+R_j^{(\ref{202005041326})}
\eea
with 
\bea\label{202005041326}
R_j^{(\ref{202005041326})}
&=&
\ol{O}_{\bbD^\infty}(n^{-1})
\eea
since 
$R_n^{(\ref{202005041326})}$ consists of 
$\sum_{k=1}^{j+\lfloor n\eta\rfloor}\eta^{-1}S_j^{(\ref{202005041155})}$ 
admitting an estimate similar to (\ref{202005041742}), 
\beas 
\sum_{k=j}^{j+\lfloor n\eta\rfloor}\eta^{-1}\int_{t_{k-1}}^{t_k}\beta_tdt -L_{\infty,\tjm}
&=&
O_{\bbD^\infty}(n^{-1})
\eeas
and a smaller term corresponding to $R_n^{(\ref{202005041156})}$. 
\end{en-text}
Similarly, (\ref{202005041239}) applied to $\theta_j=1$ and $U_\infty$ gives 
\bea\label{202005041757}
U_n-U_\infty
&=&
\sum_{k=1}^nF_k^{(\ref{202005030358})}
+\sum_{k=1}^nS_k^{(\ref{202005041155})}
+R_n^{(\ref{202005041758})}
\eea
with 
\bea\label{202005041758}
R_n^{(\ref{202005041758})}
&=&
O_{\bbD^\infty}(n^{-1.5})
\eea
thanks to the estimate (\ref{202005041236}).

\begin{en-text}
\beas 
\ol{\theta}_t &=& 
\psi\left(U_\infty^{-1}\right)\varphi\bigg(U_\infty^{-1}\eta^{-1}\int_t^{t+\eta}\sigma_s^2ds\bigg)
\yeq 
\varphi\bigg(U_\infty^{-1}\eta^{-1}\int_t^{t+\eta}\sigma_s^2ds\bigg). 
\eeas
\end{en-text}

\subsubsection{Expansion of $\Phi(U_n,L_{n,j})$}
Let 
\bea\label{202005051431}
\theta_j &=& \theta_{n,j} \yeq 
\Phi\left(U_n,L_{n,j}\right). 
\eea
Then Condition (\ref{202005051142}) is satisfied. 
\begin{en-text}
By definition, 
$\theta_j=\varphi(U_n^{-1}L_j)$ if $U_n\geq c_0$, and 
$\theta_j=0$ if $U_n\leq c_0/2$. 
\end{en-text}
Recall $\Psi_{j,k}(x,y)$ defined by (\ref{202005120042}) and 
$\Xi_{j,k,\ell}(x,y)$ defined by (\ref{202005131723}). 
\begin{en-text}
Let 
\bea\label{202005120042}
\Psi_{j,k}(x,y) 
&=& 
\Psi_{n,j,k}(x,y) 
\yeq
{\colorb1_{\{k\in\bbJ_n\}}}\partial_1\Phi(x,y)
+{\colorb1_{\{k\in K_j\}}}\partial_2\Phi(x,y){\colorb\eta_{n,j}}^{-1}.
\eea
Let 
\bea\label{202005131723}
\Xi_{j,k,\ell}(x,y) 
&=& 
\Xi_{n,j,k,\ell}(x,y) 
\nn\\&=&
\half{\colorb1_{\{k,\ell\in\bbJ_n\}}}\partial_1^2\Phi(x,y)
+\partial_1\partial_2\Phi(x,y){\colorb1_{\{k\in\bbJ_n\}}1_{\{\ell\in K_j\}}\eta_{n,j}^{-1}}
\nn\\&&
+\half\partial_2^2\Phi(x,y)
{\colorb1_{\{k,\ell\in K_j\}}\eta_{n,j}^{-2}}
\eea
\end{en-text}
We have 
\bea
\theta_j-\ol{\theta}_\tjm
&=&
\partial\Phi(U_\infty,L_{\infty,\tjm})\big[(U_n-U_\infty,L_{n,j}-L_{\infty,\tjm})\big]
\nn\\&&
+\half\partial^2\Phi(U_\infty,L_{\infty,\tjm})\big[(U_n-U_\infty,L_{n,j}-L_{\infty,\tjm})^{\otimes2}\big]
+R^{(\ref{202005041847})}_{j}
\nn\\&=&
\sum_{k=1}^n\Psi_{j,k}(U_\infty,L_{\infty,\tjm})F_k^{(\ref{202005030358})}
+
\sum_{k=1}^n\Psi_{j,k}(U_\infty,L_{\infty,\tjm})S_k^{(\ref{202005041155})}
\nn\\&&
\begin{en-text}
+\half\partial_1^2\Phi(U_\infty,L_{\infty,\tjm})\sum_{k,\ell=1}^n
F_k^{(\ref{202005030358})}F_\ell^{(\ref{202005030358})}
\nn\\&&
+\partial_1\partial_2\Phi(U_\infty,L_{\infty,\tjm})\sum_{k,\ell=1}^n1_{\{j\leq\ell\leq\nu_j\}}\eta_n^{-1}
F_k^{(\ref{202005030358})}F_\ell^{(\ref{202005030358})}
\nn\\&&
+\half\partial_2^2\Phi(U_\infty,L_{\infty,\tjm})\sum_{k,\ell=1}^n
1_{\{j\leq k\leq\nu_j\}}1_{\{j\leq\ell\leq\nu_j\}}\eta_n^{-2}
F_k^{(\ref{202005030358})}F_\ell^{(\ref{202005030358})}
\end{en-text}
\nn\\&&
{\sred 
+\partial_2\Phi(U_\infty,L_{\infty,\tjm})\cale_j
}
\nn\\&&
+\sum_{k,\ell=1}^n\Xi_{j,k,\ell}(U_\infty,L_{\infty,\tjm})F_k^{(\ref{202005030358})}F_\ell^{(\ref{202005030358})}
+R^{(\ref{202005041911})}_{j},
\eea
where 
\bea\label{202005041847}
R^{(\ref{202005041847})}_{j}
&=&
\ol{O}_{\bbD^\infty}(n^{-1.5})
\eea
and 
\bea\label{202005041911}
R^{(\ref{202005041911})}_{j}
&=&
\ol{O}_{\bbD^\infty}(n^{-1.5})
\eea
by (\ref{202005041758}), (\ref{202005041326}) and (\ref{202005041847}). 
In particular, 
\beas 
\theta_j-\ol{\theta}_\tjm
&=&
\ol{O}_{\bbD^\infty}(n^{-0.5})
\eeas
as $n\to\infty$. 
\begin{en-text}
Here we used (\ref{202005041828}) and a similar estimate, 
as well as 
(\ref{202005041312}), 
(\ref{202005041326}), 
(\ref{202005041757}) 
and (\ref{202005041758}). 
\end{en-text}
\begin{en-text}
We have
\bea\label{202005041947}
R_n^{(\ref{202005041947})} 
&:=& 
\sum_{j=1}^n\ol{\theta}_\tjm\int_\tjm^\tj\beta_tdt-\bbV_\infty
\nn\\&=&
\sum_{j=1}^n\int_\tjm^\tj\big\{\Phi(U_\infty,L_{\infty,\tjm})-\Phi(U_\infty,L_{\infty,t})\big\}\beta_tdt
\nn\\&=&
O_{\bbD^\infty}(n^{-1.5})
\eea
since 
{\colorb 
\beas 
L_{\infty,\tjm}-L_{\infty,t}
&=&
\eta_{\infty,\tjm}^{-1}\int_{(\tjm-\lambda)\vee0}^{(\tjm+\lambda)\wedge1}\beta_sds
-\eta_{\infty,t}^{-1}\int_{(t-\lambda)\vee0}^{(t+\lambda)\wedge1}\beta_sds
\nn\\&=&
\eta_{\infty,\tjm}^{-1}\int_{(\tjm-\lambda)\vee0}^{(\tjm+\lambda)\wedge1}(\beta_s-\beta_\tjm)ds
-\eta_{\infty,t}^{-1}\int_{(t-\lambda)\vee0}^{(t+\lambda)\wedge1}(\beta_s-\beta_\tjm)ds
\nn\\&=&
\big(\eta_{\infty,\tjm}^{-1}-\eta_{\infty,t}^{-1}\big)\int_{(\tjm-\lambda)\vee0}^{(\tjm+\lambda)\wedge1}(\beta_s-\beta_\tjm)ds
\nn\\&&
-\eta_{\infty,t}^{-1}\bigg\{\int_{(\tjm+\lambda)\wedge1}^{(t+\lambda)\wedge1}(\beta_s-\beta_\tjm)ds
-\int_{(\tjm-\lambda)\vee0}^{(t-\lambda)\vee0}(\beta_s-\beta_\tjm)ds\bigg\}
\nn\\&=&
O_{\bbD^\infty}(n^{-1.5})
\eeas
}
for $t\in I_j$ 
by $\beta_s-\beta_\tjm=O_{\bbD^\infty}(n^{-0.5})$ uniformly in $s\in I_j$ and $j\in\bbJ_n$. 
\end{en-text}

\subsubsection{Expansion of $\bbV_n$}
{\sred 
Recall that the target of the estimator is 
\beas 
\ol{\bbV}_n
&=& 
\int_0^1\Phi(U_n,L_{n,j})\sigma_t^2dt
\yeq
\sum_{j\in\bbJ_n}\int_\tjm^\tj\theta_j\beta_tdt
\eeas
for $\theta=(\theta_j)$ of (\ref{202005051431}).
}
We apply (\ref{202005041239}) to $\theta=(\theta_j)$ 
and {\sred$\calv=\ol{\bbV}_n$} 
to obtain 
{\sred 
\bea\label{202005041852}
\bbV_n - \ol{\bbV}_n
&=& 
\sum_{j=1}^n\theta_jF_j^{(\ref{202005030358})}
+\sum_{j=1}^n\theta_jS_j^{(\ref{202005041155})}
+R_n^{(\ref{202005041156})}
\nn\\&=&
{\colorr
\sum_{j=1}^n\ol{\theta}_\tjm F_j^{(\ref{202005030358})}
}
+\sum_{j=1}^n\ol{\theta}_\tjm S_j^{(\ref{202005041155})}
\nn\\&& 
+{\colorr\sum_{j=1}^n\sum_{k=1}^{{\colorb n}}}
\Psi_{j,k}(U_\infty,L_{\infty,\tjm})F_j^{(\ref{202005030358})}F_k^{(\ref{202005030358})}
+{\colorr\sum_{j=1}^n\sum_{k=1}^{{\colorb n}}}
\Psi_{j,k}(U_\infty,L_{\infty,\tjm})F_j^{(\ref{202005030358})}S_k^{(\ref{202005041155})}
\nn\\&&
+\sum_{j=1}^nF_j^{(\ref{202005030358})}\partial_2\Phi(U_\infty,L_{\infty,\tjm})\cale_j
\nn\\&&
+{\colorr\sum_{j=1}^n\sum_{k,\ell=1}^{{\colorb n}}}
\Xi_{j,k,\ell}(U_\infty,L_{\infty,\tjm})F_j^{(\ref{202005030358})}F_k^{(\ref{202005030358})}F_\ell^{(\ref{202005030358})}
+R_n^{(\ref{202005050041})}
\eea
\begin{en-text}
\bea\label{202005041852}
\bbV_n - \bbV_\infty
&=& 
\sum_{j=1}^n\theta_jF_j^{(\ref{202005030358})}
+\sum_{j=1}^n\theta_jS_j^{(\ref{202005041155})}
+\sum_{j=1}^n(\theta_j-\ol{\theta}_\tjm)\int_\tjm^\tj\beta_tdt 
+R_n^{(\ref{202005041947})} 
+R_n^{(\ref{202005041156})}
\nn\\&=&
{\colorr
\sum_{j=1}^n\ol{\theta}_\tjm F_j^{(\ref{202005030358})}
}
+\sum_{j=1}^n\ol{\theta}_\tjm S_j^{(\ref{202005041155})}
\nn\\&& 
+{\colorr\sum_{j=1}^n\sum_{k=1}^{{\colorb n}}}
\Psi_{j,k}(U_\infty,L_{\infty,\tjm})F_j^{(\ref{202005030358})}F_k^{(\ref{202005030358})}
+{\colorr\sum_{j=1}^n\sum_{k=1}^{{\colorb n}}}
\Psi_{j,k}(U_\infty,L_{\infty,\tjm})F_j^{(\ref{202005030358})}S_k^{(\ref{202005041155})}
\nn\\&&
+{\colorr\sum_{j=1}^n\sum_{k,\ell=1}^{{\colorb n}}}
\Xi_{j,k,\ell}(U_\infty,L_{\infty,\tjm})F_j^{(\ref{202005030358})}F_k^{(\ref{202005030358})}F_\ell^{(\ref{202005030358})}
\nn\\&& 
+
{\colorr
{\colorr\sum_{j=1}^n\sum_{k=1}^{{\colorb n}}}\Psi_{j,k}(U_\infty,L_{\infty,\tjm})\int_\tjm^\tj\beta_tdt F_k^{(\ref{202005030358})}
}
+
{\colorr\sum_{j=1}^n\sum_{k=1}^{{\colorb n}}}\Psi_{j,k}(U_\infty,L_{\infty,\tjm})\int_\tjm^\tj\beta_tdt S_k^{(\ref{202005041155})}
\nn\\&&
+{\colorr\sum_{j=1}^n\sum_{k,\ell=1}^{{\colorb n}}}\Xi_{j,k,\ell}(U_\infty,L_{\infty,\tjm})\int_\tjm^\tj\beta_tdt F_k^{(\ref{202005030358})}F_\ell^{(\ref{202005030358})}
+R_n^{(\ref{202005050041})}
\eea
\end{en-text}
where $R_n^{(\ref{202005050041})}\in\bbD^\infty$ satisfying 
\bea\label{202005050041}
R_n^{(\ref{202005050041})}
&=&
O_{\bbD^\infty}(n^{-1.5}). 
\eea
}
The estimate (\ref{202005050041}) used 
(\ref{202005041236}). 
\begin{en-text}
Here we used 
(\ref{202005041742}) with $n^{0.5}(\theta_j-\ol{\theta}_\tjm)$ for $\theta_j$ 
to make the second term on the right-hand side of (\ref{202005041852}). 
\end{en-text}

{\colorr 
We have 
\bea\label{202005110356}
e\bigg(\sum_{j=1}^n\sum_{k=1}^{{\colorb n}}\Psi_{j,k}(U_\infty,L_{\infty,\tjm})F_j^{(\ref{202005030358})}S_k^{(\ref{202005041155})}\bigg)
&=&
\max\big\{0-0.5\times5+2-0.5\times2,
\nn\\&&\hspace{30pt}
-1-0.5\times3+2-0.5\times2,
\nn\\&&\hspace{30pt}
(1-2)-0.5\times2+1-0.5
\big\}
\nn\\&=&
-1.5
\eea
since this sum consists of the three types of sums involving $I_3$, $I_1$ and $h^2$ respectively. 
Therefore, this sum is $O_{\bbD^\infty}(n^{-1.5})$ and negligible in the asymptotic expansion. 
Similarly, 
\bea\label{202005110407}
e\bigg(\sum_{j=1}^n\sum_{k,\ell=1}^{{\colorb n}}\Xi_{j,k,\ell}(U_\infty,L_{\infty,\tjm})F_j^{(\ref{202005030358})}F_k^{(\ref{202005030358})}F_\ell^{(\ref{202005030358})}\bigg)
&=&
0-0.5\times6+3-0.5\times3
\nn\\&=&
-1.5
\eea
Therefore, this sum is also $O_{\bbD^\infty}(n^{-1.5})$ and negligible in the asymptotic expansion. 
{\sred 
Moreover, 
\bea\label{202012312245}
e\bigg(\sum_{j=1}^nF_j^{(\ref{202005030358})}\partial_2\Phi(U_\infty,L_{\infty,\tjm})\cale_j\bigg)
&=&
-1.5
\eea
}

{\sred
Let 
\beas
\Theta_j
&=&
\ol{\theta}_\tjm 
\eeas
equivalently (\ref{202005110334}). 
\begin{en-text}
\bea\label{202005110334}
\Theta_j
&=&
{\colorb 
\Phi(U_\infty,L_{\infty,\tjm })
+\sum_{k=1}^n\partial_1\Phi(U_\infty,L_{\infty,\tkm})\int_\tkm^\tk\beta_rdr 
}
\nn\\&&\hspace{10pt}
{\colorb
+\sum_{k=1}^n1_{\{ j\in K_k\}}
\partial_2\Phi(U_\infty,L_{\infty,\tkm})\eta_{n,k}^{-1}\int_\tkm^\tk\beta_rdr
}
\eea
}
\begin{en-text}
{\colorr
Let 
\beas
\Theta_j
&=&
\ol{\theta}_\tjm 
+\sum_{k=1}^n\Psi_{k,j}(U_\infty,L_{\infty,\tkm})\int_\tkm^\tk\beta_tdt, 
\eeas
equivalently (\ref{202005110334}). 
\begin{en-text}
\bea\label{202005110334}
\Theta_j
&=&
{\colorb 
\Phi(U_\infty,L_{\infty,\tjm })
+\sum_{k=1}^n\partial_1\Phi(U_\infty,L_{\infty,\tkm})\int_\tkm^\tk\beta_rdr 
}
\nn\\&&\hspace{10pt}
{\colorb
+\sum_{k=1}^n1_{\{ j\in K_k\}}
\partial_2\Phi(U_\infty,L_{\infty,\tkm})\eta_{n,k}^{-1}\int_\tkm^\tk\beta_rdr
}
\eea
}
\end{en-text}
%
\begin{en-text}
In this way, (\ref{202005041852}) can be simplified as 
\bea\label{202005110412}
\bbV_n - \bbV_\infty
&=& 
{\colorr
\sum_{j=1}^{{\colorb n}}\Theta_j F_j^{(\ref{202005030358})}
}
+\sum_{j=1}^{{\colorb n}}\Theta_j S_j^{(\ref{202005041155})}
\nn\\&& 
+\sum_{j=1}^n\sum_{k=1}^{{\colorb n}}\Psi_{j,k}(U_\infty,L_{\infty,\tjm})F_j^{(\ref{202005030358})}F_k^{(\ref{202005030358})}
\nn\\&& 
+\sum_{j=1}^n\sum_{k,\ell=1}^{{\colorb n}}\Xi_{j,k,\ell}(U_\infty,L_{\infty,\tjm})\int_\tjm^\tj\beta_tdt F_k^{(\ref{202005030358})}F_\ell^{(\ref{202005030358})}
+R_n^{(\ref{202005110413})}
\eea
with the residual term 
\bea\label{202005110413}
R_n^{(\ref{202005110413})}
&=&
O_{\bbD^\infty}(n^{-1.5}).
\eea
\end{en-text}
In this way, (\ref{202005041852}) can be simplified as 
\bea\label{202005110412}
\bbV_n - \ol{\bbV}_n
&=& 
{\colorr
\sum_{j=1}^{{\colorb n}}\Theta_j F_j^{(\ref{202005030358})}
}
+\sum_{j=1}^{{\colorb n}}\Theta_j S_j^{(\ref{202005041155})}
\nn\\&& 
+\sum_{j=1}^n\sum_{k=1}^{{\colorb n}}\Psi_{j,k}(U_\infty,L_{\infty,\tjm})F_j^{(\ref{202005030358})}F_k^{(\ref{202005030358})}
+R_n^{(\ref{202005110413})}
\eea
with the residual term 
\bea\label{202005110413}
R_n^{(\ref{202005110413})}
&=&
O_{\bbD^\infty}(n^{-1.5}).
\eea
}

%
%
%
{\colorr 
We consider $M_n=\delta(u_n)$ for $u_n$ defined by 
\bea\label{202005042346}
u_n 
&=& 
n^{1/2}\sum_{j=1}^{{\colorb n}}\Theta_j(b^{[1]}_\tjm)^2I_1(1_j)1_j.
\eea
Then 
\bea\label{202005042341}
M_n
&=&
\delta(u_n) 
\nn\\&=& 
n^{1/2}\sum_{j=1}^{{\colorb n}}\Theta_j (b^{[1]}_\tjm)^2I_1(1_j)^2
-n^{1/2}\sum_{j=1}^{{\colorb n}}\big\langle D\big(\Theta_j(b^{[1]}_\tjm)^2I_1(1_j)\big),1_j\big\rangle
\nn\\&=&
n^{1/2}\sum_{j=1}^{{\colorb n}}\Theta_j F_j^{(\ref{202005030358})}
-n^{1/2}\sum_{j=1}^{{\colorb n}} \big(D_{1_j}(\Theta_j(b^{[1]}_\tjm)^2)\big) I_1(1_j)
\eea
from (\ref{202005030358}).
Therefore we obtained (\ref{202005050001}). 
The estimate (\ref{202005110436}) follows from (\ref{202005110413}). 

}

\subsection{Proof of Theorem \ref{202005110453}
}\label{202005141659}
\subsubsection{Decomposition of $Z_n$}
For the reference variable $X_\infty$ in the general theory, 
let us consider $X_1$, the value of the process $X=(X_t)_{t\in\bbR_+}$ evaluated at $t=1$. 
It is known that $X_1\in\bbD^\infty$ and $\Delta_{X_1}$ is non-degenerate if 
the stochastic differential equation is uniformly elliptic. 
In this situation, $\calq=\{2\}$ and 
\bea\label{202005112246}
a_j(2) 
&=& 
{\colorr\Theta_j}(b^{[1]}_\tjm)^2
\yeq
{\colorr\Theta_j}\beta_\tjm.
\eea

The density of the Malliavin derivative $DX_\tau$ is expressed explicitly by 
\beas 
D_sX_\tau 
&=& 
Y_\tau Y_s^{-1}\sigma_s1_{\{s\leq\tau\}}
\eeas
for $\tau\in\bbR_+$, 
where $Y=(Y_t)_{t\in\bbR_+}$ is the unique solution to the variational equation
\beas 
\left\{\begin{array}{ccl}
dY_t &=& \sigma'(X_t)Y_tdw_t+b'(X_t)Y_tdt \\
Y_0 &=& 1.
\end{array}\right.
\eeas
Similarly, 
\beas 
D_tD_sX_\tau
&=& 
Y_\tau Y_s^{-1}\sigma_s\sigma'_t1_{\{s<t\leq\tau\}}
+
Y_\tau Y_t^{-1}\sigma_t\sigma'_s1_{\{t<s\leq\tau\}}
\nn\\&=&
Y_\tau Y_{s\wedge t}^{-1}\sigma_{s\wedge t}\sigma'_{s\vee t}1_{\{s,t\leq\tau\}}.
\eeas
So it is natural to write 
\beas 
D_tD_tX_\tau
&=& 
Y_\tau Y_t^{-1}\sigma_t\sigma'_t1_{\{t\leq\tau\}}.
\eeas
The derivatives of $U_\infty$ and $L_{\infty,t}$ up to the second order 
are expressed with those of $X_\tau$ for $\tau\in\bbR_+$. 
For example, $D_tD_tU_\infty$ is well defined through the expression of $D_tD_tX_\tau$. 
In particular, these derivatives are stable in Lebesgue measure in limiting operations.

{\sred 
From (\ref{202005112246}) and (\ref{202005110334}), 
\bea\label{202005112315}
a_j(2) 
&=& 
\ol{\theta}_\tjm \beta_\tjm
\eea
}
{\colorr 
\begin{en-text}
From (\ref{202005112246}) and (\ref{202005110334}), 
\bea\label{202005112315}
a_j(2) 
&=& 
\ol{\theta}_\tjm \beta_\tjm
+\sum_{k=1}^n\Psi_{k,j}(U_\infty,L_{\infty,\tkm})\int_\tkm^\tk\beta_rdr\beta_\tjm
\eea
\end{en-text}
To specify the terms in the asymptotic expansion formula, it is necessary to investigate 
certain projections of the Malliavin derivatives of $a_j(2)$. 
%
The first two projections concerning {\sred $a_j(2)$} 
are 
\beas 
D_{1_k}{\sred a_j(2)}
&=&
\partial_1\Phi(U_\infty,L_{\infty,\tjm})\beta_\tjm D_{1_k}U_\infty
+\partial_2\Phi(U_\infty,L_{\infty,\tjm})\beta_\tjm D_{1_k}L_{\infty,\tjm}
\nn\\&&
+\Phi(U_\infty,L_{\infty,\tjm})\beta'_\tjm D_{1_k}X_\tjm
\eeas
and 
\beas 
D_{1_k}D_{1_k}{\sred a_j(2)}
&=&
\partial_1^2\Phi(U_\infty,L_{\infty,\tjm})\beta_\tjm (D_{1_k}U_\infty)^2
\nn\\&&
+
2\partial_1\partial_2\Phi(U_\infty,L_{\infty,\tjm})\beta_\tjm (D_{1_k}U_\infty)(D_{1_k}L_{\infty,\tjm})
\nn\\&&
+
2\partial_1\Phi(U_\infty,L_{\infty,\tjm})\beta'_\tjm(D_{1_k}X_\tjm)(D_{1_k}U_\infty)
\nn\\&&
+
\partial_1\Phi(U_\infty,L_{\infty,\tjm})\beta_\tjm D_{1_k}D_{1_k}U_\infty
\nn\\&&
+\partial_2^2\Phi(U_\infty,L_{\infty,\tjm})\beta_\tjm (D_{1_k}L_{\infty,\tjm})^2
\nn\\&&
+2\partial_2\Phi(U_\infty,L_{\infty,\tjm})\beta'_\tjm(D_{1_k}X_\tjm)( D_{1_k}L_{\infty,\tjm})
\nn\\&&
+\partial_2\Phi(U_\infty,L_{\infty,\tjm})\beta_\tjm D_{1_k}D_{1_k}L_{\infty,\tjm}
\nn\\&&
+\Phi(U_\infty,L_{\infty,\tjm})\beta''_\tjm (D_{1_k}X_\tjm)^2
\nn\\&&
+\Phi(U_\infty,L_{\infty,\tjm})\beta'_\tjm D_{1_k}D_{1_k}X_\tjm. 
\eeas

For Condition $[A]$ (\ref{aiii}), the process $a(s,2)$ is 
corresponding to ${\sred a_j(2)}$. 
Similarly, for Condition $[A]$ (\ref{aiv}), 
the process 
\beas 
{\sred \dota(t,s,2) }
&=& 
\partial_1\Phi(U_\infty,L_{\infty,s})\beta_sD_tU_\infty
+\partial_2\Phi(U_\infty,L_{\infty,s})\beta_sD_tL_{\infty,s}
\nn\\&&
+\Phi(U_\infty,L_{\infty,s})\beta'_sD_tX_s
\eeas
appears. 
In Condition $[A]$ (\ref{av}), 
\beas 
{\sred \dot{a}(t,2)}
&=& 
\partial_1\Phi(U_\infty,L_{\infty,t})\beta_tD_tU_\infty
+\partial_2\Phi(U_\infty,L_{\infty,t})\beta_tD_tL_{\infty,t}
\nn\\&&
+\Phi(U_\infty,L_{\infty,t})\beta'_t D_tX_t
\eeas
for {\sred$a_j(2)$}.  
For Condition $[A]$ (\ref{avi}), 
\beas 
{\sred\ddota(t,s,2)} 
&=& 
\partial_1^2\Phi(U_\infty,L_{\infty,s})\beta_s (D_tU_\infty)^2
\nn\\&&
+
2\partial_1\partial_2\Phi(U_\infty,L_{\infty,s})\beta_s (D_tU_\infty)(D_tL_{\infty,s})
\nn\\&&
+
2\partial_1\Phi(U_\infty,L_{\infty,s})\beta'_s(D_tX_s)(D_tU_\infty)
\nn\\&&
+
\partial_1\Phi(U_\infty,L_{\infty,s})\beta_s D_tD_tU_\infty
\nn\\&&
+\partial_2^2\Phi(U_\infty,L_{\infty,s})\beta_s (D_tL_{\infty,s})^2
\nn\\&&
+2\partial_2\Phi(U_\infty,L_{\infty,s})\beta'_s(D_tX_\tjm)( D_tL_{\infty,s})
\nn\\&&
+\partial_2\Phi(U_\infty,L_{\infty,s})\beta_s D_tD_tL_{\infty,s}
\nn\\&&
+\Phi(U_\infty,L_{\infty,s})\beta''_s (D_tX_s)^2
\nn\\&&
+\Phi(U_\infty,L_{\infty,s})\beta'_s D_tD_tX_s. 
\eeas
As already noted, these functionals appearing as the limits are well defined and have an explicit expression, 
respectively. 
We could go further with involved expressions though we change the way 
of description below. 
}

{\colorr
To simplify the notation, we will use the process $(\mba_s)_{s\in\bbR_+}$ defined 
by (\ref{202005120000}) on p.\pageref{202005120000}. 
\begin{en-text}
we introduce the process $(\mba_s)_{s\in\bbR_+}$ defined by 
{\colorb 
\bea\label{202005120000} 
\mba_s 
&=&
\bigg\{\Phi(U_\infty,L_s)
+\int_0^1\big(\partial_1\Phi(U_\infty,L_{\infty,r})\big)\beta_rdr
\nn\\&&\hspace{10pt}
+\int_{(s-\lambda)\vee0}^{(s+\lambda)\wedge1}\big(\partial_2\Phi(U_\infty,L_{\infty,r})\big)\eta_{\infty,r}^{-1}\beta_rdr\bigg\}\beta_s.
\eea
}
\end{en-text}
From (\ref{202005112246}) and 
(\ref{202005110334}), 
we have 
{\colorb
\beas
a_j(2) 
&=&
\Theta_j\beta_\tjm
{\sred \yeq\>}
\begin{en-text}
\nn\\&=&
\bigg\{
\ol{\theta}_\tjm 
+\sum_{k=1}^n\Psi_{k,j}(U_\infty,L_{\infty,\tkm})\int_\tkm^\tk\beta_rdr 
\bigg\}\beta_\tjm
\nn\\&=&
\end{en-text}
\begin{en-text}
\bigg\{
\Phi(U_\infty,L_{\infty,\tjm })
+\sum_{k=1}^n\partial_1\Phi(U_\infty,L_{\infty,\tkm})\int_\tkm^\tk\beta_rdr 
\nn\\&&\hspace{10pt}
+\sum_{k=1}^n1_{\{ j\in K_k\}}
\partial_2\Phi(U_\infty,L_{\infty,\tkm})\eta_{n,k}^{-1}\int_\tkm^\tk\beta_rdr
\bigg\}\beta_\tjm.
\end{en-text}
{\sred 
\Phi(U_\infty,L_{\infty,\tjm })
\beta_\tjm.
}
\eeas
}
\begin{en-text}
\beas
\Psi_{k,j}(x,y) 
&=& 
{\colorr1_{\{j\leq n\}}}\partial_1\Phi(x,y)+1_{\{k\leq j\leq\nu_k\rfloor\}}\partial_2\Phi(x,y)\eta_n^{-1}.
\eeas
\end{en-text}
Then 
\bea\label{202005111302}
\sup_{j\in\bbJ_n}\sup_{s\in I_j}
\sup_{t_1,...,t_i\in[0,1]}
\big\|D_{t_i}\cdots D_{t_1}a_j(2)-D_{t_i}\cdots D_{t_1}\mba_s\big\|_p &=& O(n^{-0.5})
\eea
for every $i\in\bbZ_+$ and $p>1$. 
Thanks to (\ref{202005111302}), we can investigate  
functionals related with the Malliavin derivatives of $a_j(2)$.

\subsubsection{Asymptotic expansion of $Z^{\sf o}_n$}
Now we apply Theorem \ref{202004211628} on p.\pageref{202004211628} 
to $Z_n^{\sf o}$ of (\ref{202005120238}) in place of $Z_n$. 
Condition $[A]$ (\ref{ai}) is routine. 
Condition $[A]$ (\ref{aiii}) is satisfied for 
\bea\label{202005120002}
a(s,2)&=&\mba_s.
\eea
Condition $[A]$ (\ref{aiv}) holds for 
\bea\label{202005120003}
\dota(t,s,2) &=&D_t\mba_s.
\eea
Condition $[A]$ (\ref{av}) holds for 
\bea\label{202005120004}
\dot{a}(t,2) 
&=& 
D_t\mba_t
\eea
Condition $[A]$ (\ref{avi}) holds for 
\bea\label{202005120005}
\ddota(t,s,2) 
&=& 
D_tD_t\mba_s
\eea
Condition $[A]$ (\ref{aix}) holds for 
\bea\label{202005120008}
\ddot{X}_\infty(t) 
&=& 
D_tD_tX_1.
\eea
In the present problem, we have set $G_\infty$ in (\ref{202005141652}): 
\beas
G_\infty
=
\int_0^12\>\mba_t^2\>dt. 
\eeas
Let 
{\sred 
\bea\label{202005140526} 
\ol{\Theta}_s 
&=& 
\Phi(U_\infty,L_s).
\eea
}
\begin{en-text}
\bea\label{202005140526} 
\ol{\Theta}_s 
&=& 
\Phi(U_\infty,L_s)
+\int_0^1\big(\partial_1\Phi(U_\infty,L_{\infty,r})\big)\beta_rdr
\nn\\&&\hspace{10pt}
+{\colorb \int_{(s-\lambda)\vee0}^{(s+\lambda)\wedge1}}\big(\partial_2\Phi(U_\infty,L_{\infty,r})\big)
{\colorb\eta_{\infty,r}^{-1}}\beta_rdr.
\eea
\end{en-text}
Then 
\bea\label{202005140527} 
\Theta_j-\ol{\Theta}_s 
&=& 
O_{\bbD^\infty}(n^{-1})
\eea
uniformly in $s\in I_j$ and $j\leq n$
since $|\eta_n-\eta|\leq n^{-1}$ and 
\beas 
L_{\infty,t_{k-1}}-L_{\infty,r} &=& O_{\bbD^\infty}(n^{-1})
\eeas
uniformly in $r\in I_k$ and $k\leq n$.

We have 
\beas&&
n^{-1}\sum_{j=1}^n2a_j(2)^2-G_\infty 
\nn\\&=& 
n^{-1}\sum_{j=1}^n2\Theta_j^2\beta_\tjm^2
-\int_0^1 2\>\mba_t^2\>dt
\nn\\&=&
\sum_{j=1}^n\int_\tjm^\tj 
2\big\{\Theta_j^2-(\ol{\Theta}_t)^2\big\}
\beta_t^2dt
+
\sum_{j=1}^n\int_\tjm^\tj 
2\Theta_j^2\big(\beta_\tjm^2-\beta_t^2\big)dt
\nn\\&=&
O_{\bbD^\infty}(n^{-1})
-\sum_{j=1}^n n^{-1}
2\Theta_j^2
(\beta^2)^{[1]}_\tjm I_1(g_j)
\nn\\&&
-\sum_{j=1}^n\int_\tjm^\tj 2\Theta_j^2\bigg(
\int_\tjm^t\big\{(\beta^2)^{[1]}_s-(\beta^2)^{[1]}_\tjm\big\}dw_s+\int_\tjm^t(\beta^2)^{[2]}_sds
\bigg)dt
\eeas
for $g_j$ defined in (\ref{202005060441}) on \pageref{202005060441}. 
We know the exponent of the second term on the right-hand side of the above equation is $-1$. 
Therefore 
\beas 
n^{-1}\sum_{j=1}^n2a_j(2)^2-G_\infty
&=& 
O_{\bbD^\infty}(n^{-1}), 
\eeas
which shows $[A^\sharp]$ (\ref{aviii}$^\sharp$). 
\underline{$[A^\sharp]$ (II) is now assumed. }
[This condition is not always satisfied. Our setting admits 
$\Phi(x,y)=x^{-1}$ (e.g. it is possible to consider when $\sigma$ is uniformly non-degenerate). ]
$[D^\sharp]$ (\ref{ai}$^{\sf x}$) is obvious. 

The functional $N_n^{\sf o}$ of (\ref{202005120741}) satisfies 
\bea\label{202005120749}
e\big(D_{u_n}^i N_n^{\sf o}\big) &=& -0.5
\eea
for $i=1$
due to Propositioin \ref{202004091035} applied to the expressions 
(\ref{202005120741}) and (\ref{202005041155}), and it entails 
the estimate (\ref{202005120749}) for $i=2$. 
Thus, Conditions (\ref{202004011428}) and (\ref{202004011429}) are fulfilled 
for $N_n^o$ 
since $\dotx_n=0$. Condition (II) of $[D^\sharp]$ has been verified. 
Related with these estimates, we conclude 
\beas 
{\mathfrak S}^{(0,1)}(\tti\sfx) &=& 0,\\
{\mathfrak S}^{(1,1)}_1(\tti\sfz,\tti\sfx) &=& 0
\eeas
for $\dotx_n=0$, and 
\beas 
{\mathfrak S}^{(2,0)}_1(\tti\sfz) &=& 0
\eeas
for $N_n^{\sf o}$. 

For ${\mathfrak S}^{(1,0)}_n(\tti\sfz)=N_n^{\sf o}\tti\sfz$, 
the integration-by-parts is used in $E[\Psi(\sfz,\sfx)N_n^{\sf o}]$ written by 
\bea\label{202005120909}
E\big[\Psi(\sfz,\sfx)N_n^{\sf o}\tti\sfz\big]
&=&
E\bigg[\Psi(\sfz,\sfx)n\sum_{j=1}^n \big(D_{1_j}({\colorr\Theta_j}(b^{[1]}_\tjm)^2)\big) I_1(1_j)
+\Psi(\sfz,\sfx)n\sum_{j=1}^n\Theta_j S_j^{(\ref{202005041155})}\bigg]\tti\sfz
\nn\\&=&
E\bigg[\Psi(\sfz,\sfx)n\sum_{j=1}^n \big(D_{1_j}({\colorr\Theta_j}(b^{[1]}_\tjm)^2)\big) I_1(1_j)
\nn\\&&\hspace{20pt}
+\Psi(\sfz,\sfx)n\sum_{j=1}^n\Theta_j \bigg\{
{\colorr b^{[1]}_\tjm b^{[1,1]}_\tjm  I_3(1_j^{\otimes3})
+2hb^{[1]}_\tjm b^{[1,1]}_\tjm  I_1(1_j)}
\nn\\&&\hspace{20pt}
+2hb^{[1]}_\tjm b^{[2]}_\tjm I_1(1_j) 
+h^2b^{[1]}_\tjm b^{[1,2]}_\tjm 
+h^2b^{[1]}_\tjm b^{[2,1]}_\tjm 
\nn\\&&\hspace{20pt}
+h^2(b^{[2]}_\tjm)^2 
+2^{-1}h^2(b^{[1,1]}_\tjm)^2
{\coloro -h\beta^{[1]}_\tjm I_1(g_j)}- 2^{-1}h^2{\colorr \beta^{[2]}_\tjm}
\bigg\}\bigg]\tti\sfz.
\nn\\&&
\eea
Obviously, the term involving the multiple integrals of order three vanishes in the limit by the IBP. 
Define the following three random symbols: 
\beas 
\mbbb^{(1)}_t(\tti\sfz,\tti\sfx)
&=& 
D_tD_t\mba_t(\tti\sfz)
+2^{-1}(D_tG_\infty)D_t\mba_t(\tti\sfz)^3+(D_tX_\infty)D_t\mba_t(\tti\sfz)(\tti\sfx),
\eeas
\beas 
\mbbb^{(2)}_t(\tti\sfz,\tti\sfx)
&=& 
D_t\hat{\mbbb}^{(2)}_t(\tti\sfz)
+2^{-1}(D_tG_\infty)\hat{\mbbb}^{(2)}_t(\tti\sfz)^3+(D_tX_\infty)\hat{\mbbb}^{(2)}_t(\tti\sfz)(\tti\sfx),
\eeas
where 
\beas 
\hat{\mbbb}^{(2)}_t
&=&
\ol{\Theta}_t\big\{
2hb^{[1]}_t b^{[1,1]}_t+2hb^{[1]}_t b^{[2]}_t
 -2^{-1}h\beta^{[1]}_t
 \big\}
\eeas
and 
\beas 
\mbbb^{(3)}_t(\tti\sfz,\tti\sfx)
&=& 
\big\{b^{[1]}_t b^{[1,2]}_t
+b^{[1]}_t b^{[2,1]}_t
+(b^{[2]}_t)^2 
+2^{-1}(b^{[1,1]}_t)^2
- 2^{-1}{\colorr \beta^{[2]}_t}\big\}(\tti\sfz).
\eeas
Let 
\bea\label{202005120907}
{\mathfrak S}^{(1,0)}(\tti\sfz,\tti\sfx)
&=&
\int_0^1\sum_{i=1}^3\mbbb^{(i)}_t(\tti\sfz,\tti\sfx)dt. 
\eea
Then 
\bea\label{202005120910}
E\big[\Psi(\sfz,\sfx)N_n^{\sf o}\tti\sfz\big]
&\to& 
E\big[\Psi(\sfz,\sfx){\mathfrak S}^{(1,0)}(\tti\sfz,\tti\sfx)\big]
\eea
as $n\to\infty$ 
from (\ref{202005120909}). Finally this proved $[D^\sharp]$ (III). 

Applying Theorem \ref{202004211628} to $Z_n^{\sf o}$, we obtain 
the asymptotic expansion for $(Z_n^{\sf o},X_1)$ as 
\bea\label{202005120916}
\sup_{f\in\cale(M,\gamma)}
\bigg| E\big[f(Z_n^{\sf o},X_1)\big] 
-\int_{\bbR^2}f(z,x)p_n(z,x)dzdx\bigg|
&=& 
o(n^{-1/2})
\eea
as $n\to\infty$ for every $(M,\gamma)\in(0,\infty)^2$. 
Here the density $p_n(z,x)$ is defined by 
\bea\label{202005120933} 
p_n(z,x) &=& 
E\bigg[\mathfrak{S}_n(\partial_z,\partial_x)^*
\bigg\{\phi(z;0,G_\infty)\delta_x(X_1)\bigg\}\bigg]
\eea
for 
\beas 
\mathfrak{S}_n
&=&
1+n^{-1/2}\mathfrak{S}
\eeas
with the random symbol $\mathfrak{S}$ given by 
\bea\label{202005120920}
\mathfrak{S}(\tti\sfz,\tti\sfx)
&=&
\mathfrak{S}^{(3,0)}(\tti\sfz,\tti\sfx) 
+\mathfrak{S}^{(1,1)}(\tti\sfz,\tti\sfx) 
+\mathfrak{S}^{(1,0)}(\tti\sfz,\tti\sfx) 
\eea
with the components 
$\mathfrak{S}^{(3,0)}(\tti\sfz,\tti\sfx)$  
and $\mathfrak{S}^{(1,1)}(\tti\sfz,\tti\sfx)$  
are generally given by (\ref{202003181606}) 
and 
(\ref{202004161636}), respectively. 
In the present situation they are described by 
\bea\label{202005120926}
{\mathfrak S}^{(3,0)}(\tti\sfz,\tti\sfx)
&=& 
\frac{4}{3}
\int_0^1 \mba_t^3 
dt\>(\tti\sfz)^3
+2
\int_{[0,1]^2}\mba^{(3,0)}(t,s)dsdt\>(\tti\sfz)^3
\nn\\&&
+2
\int_{t=0}^1
\bigg[\bigg(
{\sred 2^{-1}}
D_tG_\infty
(\tti \sfz)^5+D_tX_\infty(\tti\sfz)^3(\tti\sfx)\bigg)
\int_0^1(D_t\mba_s)ds \bigg]\mba_tdt, 
\eea
where the random field $\mba^{(3,0)}(t,s)$ comes from 
$a^{(3,0)}(t,s,q_1,q_2)$ defined by (\ref{202003181607}) is at present given by 
\bea\label{202005121036} 
\mba^{(3,0)}(t,s)
&=&
(D_tD_t\mba_s){\sred \mba_s}\mba_t
+(D_t\mba_s)^2\mba_t
+(D_t\mba_s)\mba_sD_t\mba_t,
\eea
and by 
\bea\label{202005120928}
{\mathfrak S}^{(1,1)}(\tti\sfz,\tti\sfx)
&=& 
\int_0^1D_tD_tX_1
\>\mba_t\>dt(\tti\sfz)(\tti\sfx)
\nn\\&&\hspace{10pt}
+\int_0^1(D_tX_1)D_t\mba_t\>dt(\tti\sfz)(\tti\sfx)
\nn\\&&\hspace{10pt}
+\half\int_0^1(D_tG_\infty)(D_tX_1)\mba_t\>dt(\tti\sfz)^3(\tti\sfx)
\nn\\&&\hspace{10pt}
+\int_0^1(D_tX_1)^2\mba_t\>dt(\tti\sfz)(\tti\sfx)^2. 
\eea
}

\begin{en-text}
Condition $[A]$ (\ref{aiv}) holds for 
\beas 
\dota(t,s,2) 
&=& 
\partial_1\Phi(U_\infty,L_{\infty,s})(b^{[1]}_s)^2D_tU_\infty
+\partial_2\Phi(U_\infty,L_{\infty,s})(b^{[1]}_s)^2D_tL_{\infty,s}
\nn\\&&
+\Phi(U_\infty,L_{\infty,s})\beta'_sD_tX_s.
\eeas

Condition $[A]$ (\ref{av}) holds for 
\beas 
\dot{a}(t,2) 
&=& 
\partial_1\Phi(U_\infty,L_{\infty,t})(b^{[1]}_t)^2D_tU_\infty
+\partial_2\Phi(U_\infty,L_{\infty,t})(b^{[1]}_t)^2D_tL_{\infty,t}
\eeas

Condition $[A]$ (\ref{avi}) holds for 
\beas 
\ddota(t,s,2) 
&=& 
\partial_1^2\Phi(U_\infty,L_{\infty,s})\beta_s (D_tU_\infty)^2
\nn\\&&
+
2\partial_1\partial_2\Phi(U_\infty,L_{\infty,s})\beta_s (D_tU_\infty)(D_tL_{\infty,s})
\nn\\&&
+
2\partial_1\Phi(U_\infty,L_{\infty,s})\beta'_s(D_tX_s)(D_tU_\infty)
\nn\\&&
+
\partial_1\Phi(U_\infty,L_{\infty,s})\beta_s D_tD_tU_\infty
\nn\\&&
+\partial_2^2\Phi(U_\infty,L_{\infty,s})\beta_s (D_tL_{\infty,s})^2
\nn\\&&
+2\partial_2\Phi(U_\infty,L_{\infty,s})\beta'_s(D_tX_\tjm)( D_tL_{\infty,s})
\nn\\&&
+\partial_2\Phi(U_\infty,L_{\infty,s})\beta_s D_tD_tL_{\infty,s}
\nn\\&&
+\Phi(U_\infty,L_{\infty,s})\beta''_s (D_tX_s)^2
\nn\\&&
+\Phi(U_\infty,L_{\infty,s})\beta'_s D_tD_tX_s. 
\eeas

Condition $[A]$ (\ref{aix}) holds for 
\beas 
\ddot{X}_\infty(t) 
&=& 
D_tD_tX_1.
\eeas

In the present problem, we set 
\beas 
G_\infty
=
\int_0^12(\ol{\theta}_t)^2\sigma(X_t)^4dt. 
\eeas
\beas&&
n^{-1}\sum_{j=1}^n2a_j(2)^2-G_\infty 
\nn\\&=& 
n^{-1}\sum_{j=1}^n2\big(\ol{\theta}_\tjm\big)^2\beta_\tjm^2
-\int_0^1 2(\ol{\theta}_t)^2\beta_t^2dt
\nn\\&=&
\sum_{j=1}^n\int_\tjm^\tj 
2\big\{\big(\ol{\theta}_\tjm\big)^2-(\ol{\theta}_t)^2\big\}
\beta_t^2dt
+
\sum_{j=1}^n\int_\tjm^\tj 
2\ol{\theta}_\tjm^2\big(\beta_\tjm^2-\beta_t^2\big)dt
\nn\\&=&
O_{\bbD^\infty}(n^{-1})
-\sum_{j=1}^n n^{-1}
2\ol{\theta}_\tjm^2
(\beta^2)^{[1]}_\tjm I_1(g_j)
\nn\\&&
-\sum_{j=1}^n\int_\tjm^\tj 2\ol{\theta}_\tjm^2\bigg(
\int_\tjm^t\big\{(\beta^2)^{[1]}_s-(\beta^2)^{[1]}_\tjm\big\}dw_s+\int_\tjm^t(\beta^2)^{[2]}_sds
\bigg)dt
\eeas
for $g_j$ defined in (\ref{202005060441}) on \pageref{202005060441}. 
We know the exponent of the second term on the right-hand side of the above equation is $-1$. 
Therefore 
\beas 
n^{-1}\sum_{j=1}^n2a_j(2)^2-G_\infty
&=& 
O_{\bbD^\infty}(n^{-1}), 
\eeas
which shows $[A^\sharp]$ (\ref{aviii}$^\sharp$). 

$[A^\sharp]$ (II) is now obvious. 

$[D^\sharp]$ (\ref{ai}$^{\sf x}$) is obvious.

We will apply the perturbation method because Condition (\ref{202004011429}) can 
not be verified easily. 
Remark that the first two terms in the expression (\ref{202005061538}) of $N_n$ 
have the property that the exponent of 
their projection by $D_{u_n}$ decreases by $0.5$.

${\mathbb a}$ kokoko
${\mathbb \Theta}$

\end{en-text}

\subsubsection{Perturbation method and asymptotic expansion of $Z_n$}\label{202005150029}
{\colorr
The variable aimed at is $Z_n$ having the decomposition (\ref{202005120237}): 
\beas
Z_n &=& Z_n^{\sf o} +n^{-1/2}N_n^{\sf x}, 
\eeas
We need convergence of the joint law $(M_n,N_n^{\sf x})$ 
to apply the perturbation method. 

Though the following limit theorem is easily generalized to more sophisticated cases, 
we only consider functionals 
\beas 
\calx_n&=&n^{0.5(q_0-1)}\sum_{j=1}^n A_n(j) I_{q_0}(1_j^{\otimes q_0})
\eeas
and 
\beas 
\caly_n 
&=& 
n^{0.5(q\cdot m-|m|)}
\sum_{j_{1,1},...j_{1,m_1}=1}^n\cdots \sum_{j_{\nu,1},...j_{\nu,m_\nu}=1}^n 
B_n(j_{1,1},...,j_{\nu,m_\nu})
\nn\\&&\hspace{90pt}\times
I_{q_1}(1_{j_{1,1}}^{\otimes q_1})\cdots I_{q_1}(1_{j_{1,m_1}}^{\otimes q_1})
\cdots
I_{q_\nu}(1_{j_{\nu,1}}^{\otimes q_\nu})\cdots I_{q_\nu}(1_{j_{\nu,m_\nu}}^{\otimes q_\nu})
\eeas
where $\nu\in\bbN$, $m_1,...,m_\nu\geq1$, 
$q_0\geq2$, $2\leq q_1<q_3<\cdots<q_m$, 
$q\cdot m=m_1q_1+\cdots+m_\nu q_\nu$ and $|m|=m_1+\cdots+m_\nu$.

\bd\im[[E\!\!]] 
{\bf (i)} The coefficients $A_n(j)$ and $B_n(j_{1,1},...,j_{1,m_1},...,j_{\nu,1},...,j_{\nu,m_\nu})$ satisfy 
$[S]$ on p.\pageref{202005131639}, respectively. 
\bd\im[(ii)]
There exist some continuous 
random fields $[0,1]\ni s\mapsto A(s)\in\bbD^\infty$ and 
$[0,1]^{q\cdot m}\ni (t_{1,1},...,t_{1,m_1},...,t_{\nu,1},...,t_{\nu,m_\nu})\mapsto
B(t_{1,1},...,t_{1,m_1},...,t_{\nu,1},...,t_{\nu,m_\nu})\in\bbD^\infty$ 
such that 
\beas 
\lim_{n\to\infty}\sup_{s\in I_j\atop j\in\bbJ_n}\in\|A_n(j)-A(s)\|_{k,p}&=&0
\eeas
and
\beas 
\lim_{n\to\infty}\sup_{}
\big\|B_n(j_{1,m_1},...,j_{\nu,m_\nu})-B(t_{1,1},...,t_{1,m_1},...,t_{\nu,1},...,t_{\nu,m_\nu})\big\|_{k,p}
&=&
0, 
\eeas
for every $k\in\bbR$ and $p>1$, 
where the supremum is taken over 
\beas &&
(t_{1,1},...,t_{1,m_1},...,t_{\nu,1},...,t_{\nu,m_\nu})
\in I_{j_{1,1}}\times\cdots\times I_{j_{1,m_1}}\times\cdots\times I_{j_{\nu,1}}\times\cdots\times I_{j_{\nu,m_\nu}},
\\&&
j_{1,1},...,j_{1,m_1},...,j_{\nu,1},...,j_{\nu,m_\nu}\in\bbJ_n. 
\eeas
\ed
\ed

In the following lemma, $d_s$ stands for the $\calf$-stable convergence. 
\begin{lemma}\label{202005132036}
 Suppose that Condition $[E]$ is satisfied. Then 
there exist Wiener processes $W^{(0)}$, $W^{(1)}$, ..., $W^{(\nu)}$ independent of $\calf$ such that 
\beas 
\big(\calx_n,\caly_n\big) &\to^{d_s}& \big(\calx_\infty,\caly_\infty\big)
\eeas
as $n\to\infty$, where 
\beas
\calx_\infty &=& \int_0^1 (q_0!)^{1/2}A(s)dW^{(0)}_s
\eeas
and 
\bea\label{2020051319122} 
\caly_\infty &=& \int_{[0,1]^{q\cdot m}} 
(q_1!)^{m_1/2}\cdots(q_\nu!)^{m_\nu/2}
B\big(t_{1,1},...,t_{1,m_1},...,t_{\nu,1},...,t_{\nu,m_\nu}\big)
\nn\\&&\hspace{50pt}
\cdot
dW^{(1)}_{t_{1,1}}\cdots dW^{(1)}_{t_{1,m_1}}\cdots dW^{(\nu)}_{t_{\nu,1}}\cdots dW^{(\nu)}_{t_{\nu,m_\nu}}
\eea
Moreover, 
$W^{(1)}$, ..., $W^{(\nu)}$ are independent. 
$W^{(0)}=W^{(q_i)}$ if 
$q_i=q_0$ for some $i\in\{1,..,\nu\}$, 
otherwise, $W^{(0)}$, $W^{(1)}$, ..., $W^{(\nu)}$ are independent. 
The integral $\caly_\infty$ includes the diagonal elements. 
\end{lemma}
\proof 
$L^p$-estimate of the functional of type $\caly_n$ is discussed in Section \ref{202003141640}, 
and it is observed that $\caly_n$ is near to $0$ in $\bbD^\infty$ if 
the coefficients are uniformly close to $0$ in $\bbD^\infty$. 
This implies that for any $\ep>0$, there exists $K\in\bbN$ such that 
\beas 
\limsup_{n\to\infty}\big\|\caly_n-\wt{\caly}_n\big\|_2 &<& \ep
\eeas
for
\beas 
\wt{\caly}_n 
&=& 
n^{0.5(q\cdot m-|m|)}
\sum_{j_{1,1},...j_{1,m_1}=1}^n\cdots \sum_{j_{\nu,1},...j_{\nu,m_\nu}=1}^n 
B(S(j_{1,1},...,j_{\nu,m_\nu}))
\nn\\&&\hspace{90pt}\times
I_{q_1}(1_{j_{1,1}}^{\otimes q_1})\cdots I_{q_1}(1_{j_{1,m_1}}^{\otimes q_m})
\cdots
I_{q_\nu}(1_{j_{\nu,1}}^{\otimes q_\nu})\cdots I_{q_\nu}(1_{j_{\nu,m_\nu}}^{\otimes q_\nu})
\eeas
where $S(j_{1,m_1},...,j_{\nu,m_\nu})$ is the nearest lower grid point in $\{0/K,...,(K-1)/K\}^{q\cdot m}$ 
from $(j_{1,1},...,j_{\nu,m_\nu})/n$. 
Therefore, by the Bernstein blocks, 
the problem is reduced to the stable central limit theorem 
in sub-cubes made with vertices in  $\{0/K,...,(K-1)/K\}^{q\cdot m}$. 
Then the limit of $\wt{\caly}_n$ becomes an elementary integral of $B$ with respect to 
the suitably scaled Wiener processes appearing as the stable limit. 
And it is near to $\caly_\infty$. 
\qed\halflineskip

%
In view of (\ref{202005042341}), (\ref{202005140526}), (\ref{202005140527}) and (\ref{202005030358}), 
let $\calx_n$ be the principal part of $M_n$, that is, 
\beas
\calx_n
&=&
n^{1/2}\sum_{j=1}^{{\colorb n}}\ol{\Theta}_\tjm \beta_\tjm I_2(1_j^{\otimes2})
\eeas
satisfies 
\beas 
M_n-\calx_n &\to^p& 0. 
\eeas

Recall (\ref{202005120742}): 
{\sred 
\beas
N_n^{\sf x}
&=&
n\sum_{j,k=1}^n\Psi_{j,k}(U_\infty,L_{\infty,\tjm})F_j^{(\ref{202005030358})}F_k^{(\ref{202005030358})}
.
\eeas
}
\begin{en-text}
\beas
N_n^{\sf x}
&=&
n\sum_{j,k=1}^n\Psi_{j,k}(U_\infty,L_{\infty,\tjm})F_j^{(\ref{202005030358})}F_k^{(\ref{202005030358})}
\nn\\&& 
+n\sum_{j,k,\ell=1}^n\Xi_{j,k,\ell}(U_\infty,L_{\infty,\tjm})\int_\tjm^\tj\beta_tdt F_k^{(\ref{202005030358})}F_\ell^{(\ref{202005030358})}.
\eeas
\end{en-text}
Define $\Lambda_{s,t}$ by 
{\sred\beas 
\Lambda_{s,t}
&=&
\partial_1\Phi(U_\infty,L_{\infty,s})\beta_s\beta_t
+1_{\{(s-\lambda)\vee0\leq t\leq (s+\lambda)\wedge1\}}\eta_{\infty,s}^{-1}\partial_2\Phi(U_\infty,L_{\infty,s})
\beta_s\beta_t
.
\eeas
}
\begin{en-text}
\beas 
\Lambda_{s,t}
&=&
\partial_1\Phi(U_\infty,L_{\infty,s})\beta_s\beta_t
+1_{\{(s-\lambda)\vee0\leq t\leq (s+\lambda)\wedge1\}}\eta_{\infty,s}^{-1}\partial_2\Phi(U_\infty,L_{\infty,s})
\beta_s\beta_t
\nn\\&&
+\int_0^1\bigg\{
\half\partial_1^2\Phi(U_\infty,L_{\infty,r})
+\partial_1\partial_2\Phi(U_\infty,L_{\infty,r})1_{\{(r-\lambda)\vee0\leq t\leq(r+\lambda)\wedge1\}}
\eta_{\infty,r}^{-1}
\nn\\&&\hspace{40pt}
+\half\partial_2^2\Phi(U_\infty,L_{\infty,r})1_{\{(r-\lambda)\vee0\leq s,t\leq(r+\lambda)\wedge1\}}
\eta_{\infty,r}^{-2}
\bigg\}\beta_rdr\beta_s\beta_t.
\eeas
\end{en-text}
Then 
\beas 
\caly_n 
&:=& 
n\sum_{j,k=1}^n\Lambda_{\tjm,\tkm}I_2(1_j^{\otimes2})I_2(1_k^{\otimes2})
\eeas
satisfies 
\beas 
N_n^{\sf x}-\caly_n &\to^p& 0
\eeas
By Lemma \ref{202005132036}, we obtain  
\beas
(Z_n^o,\caly_n)
&=&
(\calx_n,\caly_n) + o_p(1)
\nn\\&\to^{d_s}&
(\calx_\infty,\caly_\infty) 
\nn\\&=&
\bigg(\int_{[0,1]}\sqrt{2}\>\ol{\Theta}_s\beta_s\>dW_s,\>
\int_{[0.1]^2}2\Lambda_{s,t}\cdot dW_sdW_t\bigg)
\eeas
for some standard Wiener process $W=(W_s)_{s\in[0,1]}$ independent of $\calf$. 
Remark that 
\bea\label{202005140626} 
\int_{[0.1]^2}2\Lambda_{s,t}\cdot dW_sdW_t
&=& 
\int_0^1\int_0^t{\sred4}\wt{\Lambda}_{s,t}\> dW_sdW_t
+\int_0^12\Lambda_{s,s}ds
\eea
where the first term on the right-hand side of (\ref{202005140626}) is 
{\sred equal to} 
a double Wiener integral 
and $\wt{\Lambda}_{s,t}$ is the symmetrization of $\Lambda_{s,t}$. 

We have 
\beas 
E\big[\caly_\infty|(\calx_\infty,X_1)=(z,x)\big]
&=&
E\bigg[E\big[\caly_\infty|(\calx_\infty,X_1)=(z,x),\calf\big]\bigg|(\calx_\infty,X_1)=(z,x)\bigg]
\nn\\&=&
E\big[{\mathfrak C}(\omega,z)\big| ({\sred \calx_\infty},X_1)=(z,x)\big]
\eeas
where
\beas 
{\mathfrak C}(z)
&=&
{\mathfrak C}(\omega,z)
\nn\\&=&
\int_0^1\int_0^t{\sred8\>}G_\infty^{-2}\ol{\Theta}_t\beta_t\wt{\Lambda}_{s,t}
\ol{\Theta}_s\beta_sdsdt\>(z^2-G_\infty)
+\int_0^12\Lambda_{s,s}ds
\nn\\&=:&
{\mathfrak C}_3(z^2-G_\infty)+{\mathfrak C}_1.
\eeas
{\sred 
See e.g. Yoshida \cite{Yoshida1992} for a formula of the above conditional expectation.}
\begin{en-text}
\beas &=&
E\bigg[\int_0^1\int_0^t2G_\infty^{-2}\ol{\Theta}_t\beta_t\wt{\Lambda}_{s,t}
\ol{\Theta}_s\beta_sdsdt\>(z^2-G_\infty)
\nn\\&&\hspace{40pt}
+\int_0^12\Lambda_{s,s}ds
\bigg| ({\sred \calx_\infty},X_1)=(z,x)\bigg]
\eeas
\end{en-text}
In 
Theorem \ref{202005171835}. 
we take $(Z_n,X_1)$, $(Z_n^{\sf o},X_1)$, $(N_n^{\sf x},0)$ and $n^{-1/2}$ 
for $\bbS_T$, $\bbX_T$, $\bbY_T$ and $r_T$, respectively. 
Then the asymptotic expansion of $Z_n$ is adjusted by the function 
\beas 
{\mathfrak g}(z,x)
&=& 
-\partial_{(z,x)}\cdot\left\{E\left[
\left(\begin{array}{c} \caly_\infty \\ 0 \end{array}    \right)
 \bigg|  \begin{pmatrix} {\sred \calx_\infty} \\ X_1 \end{pmatrix} 
 =\begin{pmatrix} z \\ x \end{pmatrix}
 \right]p^{({\sred \calx_\infty},X_1)}(z,x)\right\}
 \nn\y&=&
 -\partial_z\big\{E\big[{\mathfrak C}(\omega,z)\big|({\sred \calx_\infty},X_1)=(z,x)\big]p^{({\sred \calx_\infty},X_1)}(z,x)\big\}
  \nn\y&=&
 -\partial_z\big\{E\big[{\mathfrak C}(\omega,z)\delta_z({\sred \calx_\infty})\delta_x(X_1)\big]\big\}
  \nn\y&=&
 -\partial_z\bigg\{E\bigg[{\mathfrak C}(\omega,z)
 E\big[\delta_z({\sred \calx_\infty})|\calf\big]\delta_x(X_1)\bigg]\bigg\}
   \nn\y&=&
 -\partial_z\bigg\{E\bigg[{\mathfrak C}(\omega,z)
 \phi(z;0,G_\infty)\delta_x(X_1)\bigg]\bigg\}
 \nn\y&=&
 -\partial_z\bigg\{E\bigg[\big({\mathfrak C}_3(z^2-G_\infty)+{\mathfrak C}_1\big)
 \phi(z;0,G_\infty)\delta_x(X_1)\bigg]\bigg\}
 \nn\y&=&
( -\partial_z)^3\bigg\{E\bigg[{\mathfrak C}_3{\sred G_\infty^2}
 \phi(z;0,G_\infty)\delta_x(X_1)\bigg]\bigg\}
 \nn\y&&
 -\partial_z\bigg\{E\bigg[{\mathfrak C}_1
 \phi(z;0,G_\infty)\delta_x(X_1)\bigg]\bigg\}. 
\eeas
{\sred Here we used a suitable representation of random variables on a probability space.} 
Define the random symbol ${\mathfrak G}$ by 
\beas 
{\mathfrak G}(\tti\sfz) &=& {\mathfrak G}(\tti\sfz,\tti\sfx) 
\yeq
{\mathfrak C}_3{\sred G_\infty^2}(\tti\sfz)^3+{\mathfrak C}_1(\tti\sfz). 
\eeas
Consequently, 
\beas 
{\mathfrak g}(z,x)
&=& 
E\big[{\mathfrak G}(\partial_z,\partial_x)^*\big\{ \phi(z;0,G_\infty)\delta_x(X_1)\big\}\big].
\eeas
and 
the asymptotic expansion of $\call\{(Z_n,X_1)\}$ is given by the formula 
\beas 
{\mathfrak p}_n(z,x) 
&=& 
p_n(z,x)+n^{-1/2}{\mathfrak g}(z,x). 
\eeas
Equivalently, 
the asymptotic expansion of $\call\{(Z_n,X_1)\}$ is given by the random symbol 
\beas 
{\mathfrak A}_n
&=& 
1+n^{-1/2}{\mathfrak A}
\eeas
where ${\mathfrak A}={\mathfrak S}+{\mathfrak G}$. 
This completes the proof of Theorem \ref{202005110453}. 
}

\section{Estimate of a multilinear form of multiple Wiener integrals}\label{202003141640}
\subsection{Basic estimates}
Let $\bbW=\{\bbW(h)\}_{h\in\mH}$ be a Gaussian process on a real separable Hilbert space $\mH$. 
Let ${\bf q}=(q_1,...,q_m)\in\bbN^m$, for $m\in\bbN=\{1,2,...\}$.  
We will consider a sum of multilinear forms of multiple Wiener integrals
\bea\label{202005150105} 
\bbI_n &=& \sum_{j_1,...,j_m=1}^n
A_n(j_1,...,j_m)I_{q_1}(f_{j_1}^{(1)})\cdots I_{q_m}(f_{j_m}^{(m)})
\eea
for $f^{(i)}_{j_i}\in\mH^{\tilde{\otimes} q_i}$, the $q_i$-th symmetric tensors, $i=1,...,m$. 
The random coefficients $A_n(j_1,...,j_m)$ can be anticipative. 
We already met many such functionals in the previous sections. 
An example of $\bbI_n$ is 
\begin{en-text}
\beas 
\bbI_n 
&=&
\sum_{j=1}^n\sum_{k=1}^n\sum_{\ell=1}^n
\bigg((D_ta_\tjm)1_j(s)a_{t_{k-1}}1_k(t)a_{t_{\ell-1}}1_\ell(s)\bigg)
\big(B_{t_j}-B_{t_{j-1}}\big)\big(B_t-B_{t_{k-1}}\big)\big(B_s-B_{t_{\ell-1}}\big)
\eeas
for $s,t\in[0,T]$. 
\end{en-text}
\beas 
\bbI_n 
&=&
\sum_{j=1}^n\sum_{k=1}^n\sum_{\ell=1}^n
\bigg((D_sD_ta_{t_{j-1}}a_{t_{k-1}}1_k(t)a_{t_{\ell-1}}1_\ell(s)\bigg)
\\&&
\hspace{50pt}\times
\big((B_\tj-B_\tjm)^2-(\tj-\tjm)\big)
\big(B_t-B_{t_{k-1}}\big)\big(B_s-B_{t_{\ell-1}}\big)
\eeas
for $s,t\in[0,T]$, a fractional Brownian motion $B=(B_t)$ dwelling in $\bbW$ and 
$(a_r)\subset\bbD^\infty$.  
The monomial $\bbI_n^{2k}$ is also an example of $\bbI_n$, $k\in\bbN$. 
The coefficient 
$A_n(j_1,...,j_m)$ can have the factor $\Psi_\infty$ appearing in the context of the theory of asymptotic expansion. 
To introduce the notion of exponent, an $L^p$-estimate for $\bbI_n$ and the generated algebra 
played an essential role. 
This section will give an $L^p$-estimate of $\bbI_n$.

 Suppose that 
Denote $\bar{q}=\bar{q}_1=\sum_{i=1}^mq_i$ and $\bar{r}=\bar{r}_{1} = r_1+\cdots+r_{m-2}+r_{m-1}$.  
Suppose that 
\bea\label{202005150138}
A_n(j_1,...,j_m)&\in&\bbD^{\bar{q},2}. 
\eea
The symmetrized $r$-contraction of a tensor product is denoted by $\tilde{\otimes}_r$. 
for $(r_1,...,r_{m-1})\in\bbZ_+^{m-1}$. 
We have the following estimate for $\bbI_n$ defined by (\ref{202005150105}). 
\begin{proposition}\label{0110210701}
Suppose that Condition (\ref{202005150138}) is satisfied. Then 
\beas 
\big|E[\bbI_n]\big|
&\leq&
\sum_{j_1,...,j_m=1}^n\big|E[\bbI_n]\big|
\\&\leq&
C(q_1,...,q_m)\sum_{r_1,...,r_{m-1}}
\sum_{j_1,...,j_m=1}^n
\bigg|E\big[\big\langle D^{\bar{q}_{1}-2\bar{r}_{1}}A_n(j_1,...,j_m),
f_{j_{1}}^{(1)}\tilde{\otimes}_{r_{1}}\cdots\tilde{\otimes}_{r_{m-1}}f_{j_m}^{(m)}
\big\rangle_{\mH^{\otimes(\bar{q}_{1}-2\bar{r}_{1})}}\big]\bigg|
\eeas
for some constant $C(q_1,...,q_m)$ depending only on $q_1,...,q_m$. 
\end{proposition}
\proof 
Repeatedly apply the product formula to 
$I_{q_1}(f_{j_1}^{(1)})\cdots I_{q_m}(f_{j_m}^{(m)})$ 
to reach a linear combination of Skorohod integrals of order $\bar{q}_1-2\bar{r}_1$. 
The sum $\sum_{r_1,...,r_{m-1}}$ is finite. 
Next use duality for $\delta^{\bar{q}_1-2\bar{r}_1}$. 
\qed\halflineskip
\begin{en-text}
ex. Case $m=2$. 
\beas 
E[\bbI_n]
&=&
\sum_{r=0}^{q_{1}\wedge q_2}r!
\left(\begin{array}{c}q_{1}\\r\end{array}\right)
\left(\begin{array}{c}q_2\\r\end{array}\right)
\sum_{(j_1,j_2)\in\{1,...,n\}^2}
E\big[A_n(j_1,j_2)I_{q_{1}+q_2-2r}(f_{j_{1}}^{(1)}\tilde{\otimes}_rf_{j_2}^{(2)})\big]
\\&=&
\sum_{r=0}^{q_{1}\wedge q_2}r!
\left(\begin{array}{c}q_{1}\\r\end{array}\right)
\left(\begin{array}{c}q_2\\r\end{array}\right)
\sum_{(j_1,j_2)\in\{1,...,n\}^2}
E\big[\big\langle D^{q_1+q_2-2r}A_n(j_1,j_2), f_{j_{1}}^{(1)}\tilde{\otimes}_rf_{j_2}^{(2)}\big\rangle_{\mH^{\otimes(q_1+q_2-2r)}}\big]
\eeas
Therefore, 
\beas 
\big|E[\bbI_n]\big|
&\leq&
C(q_1,q_2)\sum_r
\sum_{(j_1,j_2)\in\{1,...,n\}^2}
\bigg|E\big[\big\langle D^{q_1+q_2-2r}A_n(j_1,j_2), f_{j_{1}}^{(1)}\tilde{\otimes}_rf_{j_2}^{(2)}\big\rangle_{\mH^{\otimes(q_1+q_2-2r)}}\big]\bigg|
\eeas
\end{en-text}

Let 
\beas 
\bbI_n(j_1,...,j_m) 
&=& 
A_n(j_1,...,j_m)I_{q_1}(f_{j_1}^{(1)})\cdots I_{q_m}(f_{j_m}^{(m)})
\eeas
\begin{corollary}\label{0110221510}
Suppose that Condition (\ref{202005150138}) is satisfied. Then 
\beas 
\big|E[\bbI_n(j_1,...,j_m) ]\big|
&\leq&
C(q_1,...,q_m)
\\&&\times
\sum_{r_1,...,r_{m-1}}\bigg|E\big[\big\langle D^{\bar{q}_{1}-2\bar{r}_{1}}A_n(j_1,...,j_m),
f_{j_{1}}^{(1)}\tilde{\otimes}_{r_{1}}\cdots\tilde{\otimes}_{r_{m-1}}f_{j_m}^{(m)}
\big\rangle_{\mH^{\otimes(\bar{q}_{1}-2\bar{r}_{1})}}\big]\bigg|
\eeas
for some constant $C(q_1,...,q_m)$ depending only on $q_1,...,q_m$. 
\end{corollary}

\begin{en-text}
product formula
\beas &&
I_{q_1}(f_{j_1}^{(1)})\cdots I_{q_m}(f_{j_m}^{(m)})
\\&=&
\sum_{r_{m-1}=0}^{\wedge q_{m-1}}c(q_{m-1},q_m,r_{m-1})
I_{q_1}(f_{j_1}^{(1)})\cdots I_{q_{m-2}}(f_{j_{m-2}}^{(m-2)})
I_{\bar{q}_{m-1}-2r_{m-1}}(f_{j_{m-1}}^{(m-1)}\tilde{\otimes}_{r_{m-1}}f_{j_m}^{(m)})
\\&=&
\sum_{r_{m-1}=0}^{\wedge q_{m-1}}
\sum_{r_{m-2}=0}^{\wedge q_{m-2}(\bar{r}_{m-1})}
c(q_{m-1},q_m,r_{m-1})c(q_{m-2},\bar{q}_{m-1}-2\bar{r}_{m-1},r_{m-2})
\\&&\times
I_{q_1}(f_{j_1}^{(1)})\cdots I_{q_{m-3}}(f_{j_{m-3}}^{(m-3)})
I_{\bar{q}_{m-2}-2\bar{r}_{m-2}}(f_{j_{m-2}}^{(m-2)}\tilde{\otimes}_{r_{m-2}}f_{j_{m-1}}^{(m-1)}\tilde{\otimes}_{r_{m-1}}f_{j_m}^{(m)})
\\&=&
\sum_{r_{m-1}=0}^{\wedge q_{m-1}}
\sum_{r_{m-2}=0}^{\wedge q_{m-2}(\bar{r}_{m-1})}
\sum_{r_{m-3}=0}^{\wedge q_{m-3}(\bar{r}_{m-2})}
c(q_{m-1},q_m,r_{m-1})
c(q_{m-2},\bar{q}_{m-1}-2\bar{r}_{m-1},r_{m-2})
\\&&\times
c(q_{m-3},\bar{q}_{m-2}-2\bar{r}_{m-2},r_{m-3})
\\&&\times
I_{q_1}(f_{j_1}^{(1)})\cdots I_{q_{m-4}}(f_{j_{m-4}}^{(m-4)})
I_{\bar{q}_{m-3}-2\bar{r}_{m-3}}(f_{j_{m-3}}^{(m-3)}\tilde{\otimes}_{r_{m-3}}\cdots\tilde{\otimes}_{r_{m-1}}f_{j_m}^{(m)})
\\&=&
\sum_{r_{m-1}=0}^{\wedge q_{m-1}}
\sum_{r_{m-2}=0}^{\wedge q_{m-2}(\bar{r}_{m-1})}
\sum_{r_{m-3}=0}^{\wedge q_{m-3}(\bar{r}_{m-2})}
\cdots
\sum_{r_{1}=0}^{\wedge q_{1}(\bar{r}_{2})}
c(q_{m-1},q_m,r_{m-1})
c(q_{m-2},\bar{q}_{m-1}-2\bar{r}_{m-1},r_{m-2})
\\&&\times
c(q_{m-3},\bar{q}_{m-2}-2\bar{r}_{m-2},r_{m-3})
\cdots
c(q_{1},\bar{q}_{2}-2\bar{r}_{2},r_{1})
I_{\bar{q}_{1}-2\bar{r}_{1}}(f_{j_{1}}^{(1)}\tilde{\otimes}_{r_{1}}\cdots\tilde{\otimes}_{r_{m-1}}f_{j_m}^{(m)})
\eeas

Therefore
\beas &&
E\big[A_n(j_1,...,j_m)I_{q_1}(f_{j_1}^{(1)})\cdots I_{q_m}(f_{j_m}^{(m)})\big]
\\&=&
\sum_{r_{m-1}=0}^{\wedge q_{m-1}}
\sum_{r_{m-2}=0}^{\wedge q_{m-2}(\bar{r}_{m-1})}
\sum_{r_{m-3}=0}^{\wedge q_{m-3}(\bar{r}_{m-2})}
\cdots
\sum_{r_{1}=0}^{\wedge q_{1}(\bar{r}_{2})}
c(q_{m-1},q_m,r_{m-1})
c(q_{m-2},\bar{q}_{m-1}-2\bar{r}_{m-1},r_{m-2})
\\&&\times
c(q_{m-3},\bar{q}_{m-2}-2\bar{r}_{m-2},r_{m-3})
\cdots
c(q_{1},\bar{q}_{2}-2\bar{r}_{2},r_{1})
\\&&\times
E\big[A_n(j_1,...,j_m)
I_{\bar{q}_{1}-2\bar{r}_{1}}(f_{j_{1}}^{(1)}\tilde{\otimes}_{r_{1}}\cdots\tilde{\otimes}_{r_{m-1}}f_{j_m}^{(m)})\big]
\\&=&
\sum_{r_{m-1}=0}^{\wedge q_{m-1}}
\sum_{r_{m-2}=0}^{\wedge q_{m-2}(\bar{r}_{m-1})}
\sum_{r_{m-3}=0}^{\wedge q_{m-3}(\bar{r}_{m-2})}
\cdots
\sum_{r_{1}=0}^{\wedge q_{1}(\bar{r}_{2})}
c(q_{m-1},q_m,r_{m-1})
c(q_{m-2},\bar{q}_{m-1}-2\bar{r}_{m-1},r_{m-2})
\\&&\times
c(q_{m-3},\bar{q}_{m-2}-2\bar{r}_{m-2},r_{m-3})
\cdots
c(q_{1},\bar{q}_{2}-2\bar{r}_{2},r_{1})
\\&&\times
E\big[\big\langle D^{\bar{q}_{1}-2\bar{r}_{1}}A_n(j_1,...,j_m),
f_{j_{1}}^{(1)}\tilde{\otimes}_{r_{1}}\cdots\tilde{\otimes}_{r_{m-1}}f_{j_m}^{(m)}
\big\rangle_{\mH^{\otimes(\bar{q}_{1}-2\bar{r}_{1})}}\big]
\eeas
\end{en-text}
\begin{en-text}
\bi
\im notation. \koko
\beas 
c(a,b,r) &=& 
r!
\left(\begin{array}{c}a\\r\end{array}\right)
\left(\begin{array}{c}b\\r\end{array}\right)
\eeas

\beas 
\bar{q}_i &=& q_i+q_{i+1}+\cdots+q_m
\eeas

\beas 
\bar{r}_{m-1} &=& r_{m-1},
\\
\bar{r}_{m-2} &=& r_{m-2}+r_{m-1},
\\&\cdots&
\\
\bar{r}_{1} &=& r_1+\cdots+r_{m-2}+r_{m-1}. 
\eeas

\beas
\wedge q_{m-1} &=& q_{m-1}\wedge q_m,
\\
\wedge q_{m-2}(\bar{r}_{m-1}) &=& q_{m-2}\wedge (\bar{q}_{m-1}-2\bar{r}_{m-1}),
\\
\wedge q_{m-3}(\bar{r}_{m-2}) &=& q_{m-3}\wedge (\bar{q}_{m-2}-2\bar{r}_{m-2}),
\\&\cdots&
\\
\wedge q_{1}(\bar{r}_{2}) &=& q_{1}\wedge (\bar{q}_{2}-2\bar{r}_{2})
\eeas

\im 

\im 
\ei
\end{en-text}

\subsection{$L^p$-estimate for a multilinear form of multiple Wiener integrals for discretization}
Consider a double sequence of numbers $t_j=t^n_j$ ($i=1,...n$, $n\in\bbN$) such that 
\bea\label{202005170230} 
0= t_0<t_1<\cdots<t_n, &&\max_{j=0,1,...,n}(\tj-\tjm)\leq T/n
\eea
for a constant $T>0$. 
For example, $(t_j)$ is a double sequence such that 
$0= t_0<t_1<\cdots<t_n\leq T_0$ and that $\max_{j=1,...,n}(\tj-\tjm)\leq2T_0/n$ for $T=2T_0$. 

Let $I_j=I^n_j=1_{(\tjm,\tj]}$.\footnote{$I_j$ denotes the $j$-th interval. 
The mapping of the $q$-th multiple integral is denoted by $I_q(\cdot)$ but 
there is no fear of confusion thanks to the argument of $I_q$. } 
From now on we let $\mH=L^2([0,T],dt)$. 
We will consider a double sequence $f_{j_i}^{(i)}=f_{j_i}^{(i),n}\in\mH^{\tilde{\otimes}q_j}$ such that 
$\text{supp}(f_{j_i}^{(i)})\subset I_j^{\>q_i}$ and that 
\bea\label{0110231323}
\max_{i=1,...,m}\>\sup_{n\in\bbN}\sup_{j_i=1,...,n}
\sup_{t_1,...,t_{q_i}\in[0,T]}|f_{j_i}^{(i)}(t_1,...,t_{q_i})|&\leq&C_1
\eea
for a constant $C_1$. 
Let $k\in\bbN$. 
Suppose that $A_n(j_1,...,j_m)\in\bbD^{2k\bar{q},4k}
$ (the derivatives up to $2k$-th order are in $L^{4k}$) and 
\bea\label{0110231324}
\max_{i=0,1,...,2k\bar{q}}\>
\sup_{n\in\bbN}\sup_{(j_1,...,j_m)\in\{1,...,n\}^m}
\sup_{t_1,...,t_i\in[0,T]}
\big\|D^i_{t_1,...,t_i}A_n(j_1,...,j_m)\big\|_{4k}&\leq&C_0
\eea
for some constant $C_0$. 

In the following theorem, we apply the duality of the Skorohod integrals, as it was used in 
Proposition \ref{0110210701} and Corollary \ref{0110221510}, to
\beas 
\check{\bbI}_n({\bf j})
&=&
A_n(j_{1,1},...,j_{1,m})A_n(j_{2,1},...,j_{2,m})\cdots A_n(j_{2k,1},...,j_{2k,m})
\\&&\qquad\qquad\times
I_{q_1}(f_{j_{1,1}}^{(1)})\cdots I_{q_m}(f_{j_{1,m}}^{(m)})
I_{q_1}(f_{j_{2,1}}^{(1)})\cdots I_{q_m}(f_{j_{2,m}}^{(m)})
\cdots
I_{q_1}(f_{j_{2k,1}}^{(1)})\cdots I_{q_m}(f_{j_{2k,m}}^{(m)})
\eeas
for ${\bf j}=\big( j_{1,1},...,j_{1,m},\cdots,j_{2k,1},...,j_{2k,m}\big)\in\{1,...,n\}^{2km}$. 
The exponent introduced in Section \ref{202004270705} is based on the following estimate. 
\halflineskip
\begin{theorem}\label{0110221242}
There exists a constant $C_2=C\big(q_1,...,q_m,k,T,C_0,C_1\big)$ such that 
\beas 
\sum_{{\bf j}\in\{1,...,n\}^{2km}}
\sup\!\!\>^*\>
\big|E\big[\check{\bbI}_n({\bf j})\big]\big|
&\leq&
C_2^{2k}\>n^{-k(\bar{q}-m)}
\eeas
for all $n\in\bbN$, 
where $\sup^*$ is the supremum taken over the all elements 
$f^*_*$ and $A_n(*)$ satisfying (\ref{0110231323}) and (\ref{0110231323}), 
respectively. 
In particular, 
\beas 
\big\|\bbI_n\big\|_{2k}
&\leq&
C_2n^{-2^{-1}(\bar{q}-m)}
\eeas
for all $n\in\bbN$. 
{\colorr More precisely, one can take 
\beas 
C_2&=& C(q_1,...,q_m,k,T)C_0C_1^m
\eeas
for some constant $C(q_1,...,q_m,k,T)$. }
\end{theorem}
\proof 
Let $\bar{Q}_1 = 2k\bar{q}$. 
We will estimate the sum 
\beas 
S({\bf R})&:=&\sum_{{\bf j}}\sup\!\!\>^*\big|\bbE_n({\bf j}, {\bf R})\big|
\eeas
for
\beas 
\bbE_n({\bf j}, {\bf R})
&=& 
E\bigg[\bigg\langle D^{\bar{Q}_1-2\bar{R}_1}
B_n({\bf j}),\>
\Phi({\bf j},\R)
\bigg\rangle_{\mH^{\otimes (\bar{Q}_1-2\bar{R}_1)}}\bigg],
\eeas
where 
\beas 
B_n({\bf j})
&=&
A_n(j_{1,1},...,j_{1,m})A_n(j_{2,1},...,j_{2,m})\cdots A_n(j_{2k,1},...,j_{2k,m}), 
\eeas
\beas 
\R &=& 
\big(R_{1,1},...,R_{1,m},R_{2,1},...,R_{2,m},...,R_{2k,1},...,R_{2k,m-1}\big)
\eeas
\beas 
\Phi({\bf j},\R)
&=&
f_{j_{1,1}}^{(1)}\otimes_{R_{1,1}}\cdots\otimes_{R_{1,m-1}}f_{j_{1,m}}^{(m)}
\otimes_{R_{1,m}}
f_{j_{2,1}}^{(1)}\otimes_{R_{2,1}}\cdots\otimes_{R_{2,m-1}}f_{j_{2,m}}^{(m)}
\otimes_{R_{2,m}}\cdots
\\&&\qquad
\otimes_{R_{2k-1,m}}
f_{j_{2k,1}}^{(1)}
\otimes_{R_{2k,1}}\cdots\otimes_{R_{2k,m-1}}f_{j_{2k,m}}^{(m)}
\eeas
and 
\beas 
\bar{R}_1 &=& 
R_{1,1}+\cdots+R_{1,m}+R_{2,1}+\cdots+R_{2,m}+\cdots+
R_{2k,1}+\cdots+R_{2k,m-1}. 
\eeas
\begin{en-text}
We release $\tilde{\Phi}(j_{a,b},R_{a,b};a=1,...,2k;b=1,...,m)$ 
from the symmetrization and consider 
\beas &&
\tilde{\Phi}(j_{a,b},R_{a,b};a=1,...,2k;b=1,...,m)
\\&=&
f_{j_{1,1}}^{(1)}\tilde{\otimes}_{R_{1,1}}\cdots\tilde{\otimes}_{R_{1,m-1}}f_{j_{1,m}}^{(m)}
\tilde{\otimes}_{R_{1,m}}
f_{j_{2,1}}^{(1)}\tilde{\otimes}_{R_{2,1}}\cdots\tilde{\otimes}_{R_{2,m-1}}f_{j_{2,m}}^{(m)}
\tilde{\otimes}\cdots\tilde{\otimes}
f_{j_{2k,1}}^{(1)}\tilde{\otimes}_{R_{2k,1}}\cdots\tilde{\otimes}_{R_{2k,m-1}}f_{j_{2k,m}}^{(m)}.
\eeas
\end{en-text}
Remark that each nonnegative integer $R_{a,b}$ is bounded by 
a suitable number so that each contraction makes sense. 
Each contraction has implicitly chosen a set of arguments 
in the factor in front of it and the all factors after it. 
Later this set will vary between all possible cases, but 
we fix an ${\bf R}$ for the meantime. 

We call the factor $f^{(i)}_{j_i}$ is {\it single} if 
all its arguments are free from any contraction specifically allocated according to ${\bf R}$. 
A single factor $f^{(i)}_{j_i}$ has $q_i$ arguments $(t_1,...,t_{q_i})$, and 
this factor contributions to $S$ by 
a factor of order $O(n^{-(q_i-1)})$ after summation 
$\sum_{j_i=1}^n\big|\bbE_n({\bf j}, {\bf R})\big|$, due to 
the specific support of $f^{(i)}_{j_i}$. 

We consider connection between two $f^{(*)}_*$s by a contraction. 
In a connected component $\calc$ (having at least two $f^{(*)}_*$s), suppose that ${\sf r}$ $(\geq1)$ contractions exist. 
Let ${\sf q}$ be the sum of the numbers of arguments of ${\sf m}$ functions, 
$f^{(i_1)}_{j_1},...,f^{(i_{\sf m})}_{ 
j_{\sf m}} $ say, consisting of $\calc$. 
We have ${\sf q}=q_{i_1}+\cdots+q_{i_{{\sf m}}}$ 
since each $f^{(i)}_{j_i}$ has $q_i$ arguments. 
Connected two arguments of (different) $f^{(*)}_*$s in $\calc$ give 
a factor of order $O(n^{-1})$ (before summing up with respect to involving $j$'s). 
So, $2{\sf r}$ arguments give totally a factor of order $O(n^{-{\sf r}})$. 
On the other hand, the rest of the arguments appearing in $\calc$ gives 
a factor of order $O(n^{-({\sf q}-2{\sf r})})$. 
The multiple sum with respect to $j$'s in $\calc$ becomes a single sum 
by the orthogonality,
because $\calc$ is connected. 
Thus, after summation with respect $j_1=\cdots=j_{{\sf m}}$, 
the component $\calc$ contributes to $S({\bf R})$ with a factor of the order 
\bea\label{0110231459} 
n\times n^{-({\sf q}-2{\sf r})}\times n^{-{\sf r}}
\yeq 
n^{1+{\sf r}-{\sf m}-\sum_{i=i_1,...,i_{\sf m}}(q_i-1)}
\ygeq
n^{-\sum_{i=i_1,...,i_{\sf m}}(q_i-1)}. 
\eea
Here we used the inequality $1+{\sf r}\geq {\sf m}$. 
Inequality (\ref{0110231459}) suggests we should consider 
the allocations of contractions that do not involve any single $f^{(*)}_{*}$, 
in order to find the dominating term. 

We consider an allocation of contractions according to ${\sf R}$ that  
has 
$c$ connected components with $c\geq1$ and has no single $f^{(*)}_*$. 
Then 
\beas 
S({\bf R})
&=&
O(n^{c+\bar{R}_1-2km-(\bar{Q}_1-2km)})
\yeq
O(n^{c+\bar{R}_1-\bar{Q}_1}). 
\eeas
Since the contraction takes place between different two functions $f^{(*)}_{*}$, 
each connected component has at least two functions $f^{(*)}_{*}$, therefore, 
the maximum $c$ is $km$. 
The maximum $\bar{R}_1$ is $k\bar{q}$. 
(Indeed, there is a pattern of contractions for which $c=km$ and $\bar{R}_1=k\bar{q}$.) 
Consequently, 
\beas 
S({\bf R})
&=&
O(n^{km+k\bar{q}-\bar{Q}_1})
\yeq
O(n^{-k(\bar{q}-m)})
\eeas
for any pattern of allocation of contractions. 
The existence of a constant $C\big(q_1,...,q_m,k,T\big)
$ is 
obvious if one observes the above argument precisely. 
This completes the proof. 
\qed\halflineskip

\subsection{$L^p$-estimate of an integrated functional on a compact group}
Let $m\in\bbN$ and $q_1,...,q_m\in\bbN$. 
For each $i\in\{1,...,m\}$, let $G^{(i)}$ be a compact group with a Haar measure $\mu^{(i)}$. 
Let $G=G_1\times \cdots G_m$ and $\mu=\mu_1\times\cdots\times\mu_m$. 
%
Let $e^{(i)}_{n,j}:G^{(i)}\to\mH^{\tilde{\otimes}q_i}$ for $j\in\bbJ_n=\{1,...,n\}$, $n\in\bbN$ and $i=1,...,m$. 
%
Suppose that $A_n(j_1,...,j_m;\mbx)$ is a random variable 
for each $(j_1,...,j_m)\in\bbJ_n^m$ and $\mbx=(x^{(1)},...,x^{(m)})\in G$. 
We assume that 
the mappings $A_n(j_1,...,j_m;\cdot)$ and $I_{q_i}\big(e^{(i)}_{j_i}(\cdot)\big)$ are 
measurable random fields on $G$ and $G^{(i)}$, respectively, 
and 
\beas 
\int_G 
\bigg|A_n(j_1,...,j_m;\mbx)\prod_{i=1}^m
I_{q_i}\big(e_{j_i}^{(i)}(x^{(i)})\big)\bigg|\mu(dx)
&<&
\infty\quad a.s.
\eeas
\begin{en-text}
$e^{(i)}_{n,j_1}$ is Bochner integrable with respect to $\mu^{(i)}$, 
the mapping $G^{(i)}\ni x^{(i)}\to I_{q_i}\big(e^{(i)}_{n,j}(x^{(i)})\big)$ is integrable with respect to $\mu^{(i)}$, 
and that 
\beas 
I_{q_i}\bigg(\int_{G^{(i)}}e^{(i)}_{n,j}(x^{(i)})\mu^{(i)}(dx^{(i)})\bigg) 
&=& 
\int_{G^{(i)}}I_{q_1}\big(e^{(i)}_{n,j}(x^{(i)})\big)\mu^{(i)}(dx^{(i)})\koko
\eeas
\end{en-text}
\begin{en-text}
\beas 
\bbH_n 
&=& 
\int_G
\sum_{j_1,...,j_m=1}^n 
A_n(j_1,...,j_m;\mbt)I_{q_1}\big(e_{j_1}^{(1)}(x^{(1)})\big)\cdots I_{q_m}\big(e_{j_m}^{(m)}(x^{(m)})\big)
\mu(dx)
\eeas
\end{en-text}
Let 
\bea\label{202005160359} 
\dot{\bbH}_n(\mbx) 
&=& 
\sum_{j_1,...,j_m=1}^n 
A_n(j_1,...,j_m;\mbx)I_{q_1}\big(e_{j_1}^{(1)}(x^{(1)})\big)\cdots I_{q_m}\big(e_{j_m}^{(m)}(x^{(m)})\big).
\eea
We are interested in the $L^p$-estimate of the integral 
\beas
\bbH_n 
&=& 
\int_G\dot{\bbH}_n(\mbx)\mu(d\mbx). 
\eeas

Assume the localization condition that for every $j_i\in\bbJ_n$, 
\bea\label{202005171254} 
e_{n,j_i}^{(i)}(x^{(i)}) &=& e_{n,j_i}^{(i)}(x^{(i)})1_{\chi_n^{(i)}}(L_{g^{(i)}_{n,j_i}}x^{(i)}) 
\qquad (x^{(i)}\in G^{(i)})
\eea
for some $g^{(i)}_{n,j_i}\in G^{(i)}$ and some measurable set $\chi^{(i)}_n$ in $G^{(i)}$, 
where $L_{g^{(i)}}$ is the left transformation. 
\begin{en-text}
 such that 
\beas 
\mu^{(i)}(\chi^{(i)}_n)  &\leq & Cn^{-1}
\eeas
for each $i\in\{1,...,m\}$. 
\end{en-text}
[For example, when $G^{(i)}$ is $\bbR/\bbZ$, $\chi^{(i)}_n=[0,c/n]$, a fixed interval of bigger than the maximum length of $(\tjm,\tj]$. 
$e^{(i)}_{n,j}(x^{(i)})$ can still include an indicator function. ] 

Then
\beas 
\bbH_n
&=& 
\begin{en-text}
\int_G 
\sum_{j_1,...,j_m=1}^nA_n(j_1,...,j_m;\mbx)\prod_{i=1}^m
I_{q_i}\big(e_{j_i}^{(i)}(x^{(i)})\big)\mu(dx)
\nn\\&=&
\end{en-text}
\int_{G_m}\cdots\int_{G_1} 
\sum_{j_1,...,j_m=1}^nA_n(j_1,...,j_m;\mbx)
\prod_{i=1}^mI_{q_i}\big(e_{n,j_i}^{(i)}(x^{(i)})\big)
\>\mu^{(1)}(dx^{(1)})\cdots\mu^{(m)}(dx^{(m)})
\nn\\&=&
\int_{G_m}\cdots\int_{G_1} 
\sum_{j_1,...,j_m=1}^nA_n(j_1,...,j_m;\mbx)
\nn\\&&\hspace{30pt}\times
\prod_{i=1}^m\big\{I_{q_i}\big(e_{j_i}^{(i)}(x^{(i)})\big)1_{\chi_n^{(i)}}(L_{g^{(i)}_{n,j_i}}x^{(i)}) \big\}
\>\mu^{(1)}(dx^{(1)})\cdots\mu^{(m)}(dx^{(m)})
\nn\\&=&
\begin{en-text}
\int_{G_m}\cdots\int_{G_1} 
\sum_{j_1,...,j_m=1}^nA_n\big(j_1,...,j_m;(L_{{(g^{(i)}_{n,j_i})^{-1}}}x^{(i)})_{i=1}^m\big)
\nn\\&&\hspace{30pt}\times
\prod_{i=1}^mI_{q_i}\big(e_{j_i}^{(i)}(L_{{(g^{(i)}_{n,j_i})^{-1}}}x^{(i)})\big)
\prod_{i=1}^m1_{\chi_n^{(i)}}(x^{(i)}) 
\>\mu^{(1)}(dx^{(1)})\cdots\mu^{(m)}(dx^{(m)})
\nn\\&=&
\end{en-text}
\int_G{\stackrel{\circ}{\bbH}}_n(\mbx)
\prod_{i=1}^m1_{\chi_n^{(i)}}(x^{(i)}) \>\mu(d\mbx)
\eeas
where 
\bea\label{202005190511} 
{\stackrel{\circ}{\bbH}}_n(\mbx)
&=& 
\sum_{j_1,...,j_m=1}^nA_n\big(j_1,...,j_m;(L_{{(g^{(i)}_{n,j_i})^{-1}}}x^{(i)})_{i=1}^m\big)
\prod_{i=1}^mI_{q_i}\big(e_{j_i}^{(i)}(L_{{(g^{(i)}_{n,j_i})^{-1}}}x^{(i)})\big).
\eea
Therefore, 
\begin{theorem}\label{202005190527}
Suppose that Condition (\ref{202005171254}) is satisfied. Then 
\bea\label{202005190606}
\bbH_n
&=& 
\int_G{\stackrel{\circ}{\bbH}}_n(\mbx)
\prod_{i=1}^m1_{\chi_n^{(i)}}(x^{(i)}) \>\mu(d\mbx)
\eea
and 
\bea\label{202005171343} 
\big\|\bbH_n\big\|_p 
&\leq&
\sup_{\mbx\in G}\big\|{\stackrel{\circ}{\bbH}}_n(\mbx)\big\|_p\>
\mu\bigg(\prod_{i=1}^m\chi_n^{(i)}\bigg)
\eea
for $p\geq1$. 
\end{theorem}
\halflineskip

\begin{en-text}
Extend the function $[0,T]^M\ni\mbt\to A_n(j_1,...,j_m;\mbt)$ to $\bbR^m$ 
(e.g. give $0$ to $A_n(j_1,...,j_m;\cdot)$ outside of the rectangle). 
Let $\varphi_\ep^{(i)}$ be a mollifier on the space of $t^{(i)}$, 
and let $\varphi_\ep=\otimes_{i=1}^m\varphi_\ep^{(i)}$, i.e., we assume that 
\beas 
A_n(j_1,...,j_m;\cdot)*\varphi_\ep &\to& A_n(j_1,...,j_m;\cdot) 
\eeas
in $L^1(\bbR^M,d\mbt)$ as $\ep\down0$. 

Define $\bbH_{n,\ep}(\mbt)$ by 
\beas 
\bbH_{n,\ep}(\mbt)
&=& 
\int \sum_{j_1,...,j_m=1}^n  
\big(A_{n,\ep}(j_1,...,j_m;\cdot)*\varphi_\ep\big)(\mbt) 
I_{q_1}\big(e_{j_1}^{(1)}(t^{(1)})\big)\cdots I_{q_m}\big(e_{j_m}^{(m)}(t^{(m)})\big)\>d\mbt
\nn\\&=&
\int \sum_{j_1,...,j_m=1}^n 
\int A_n(j_1,...,j_m;\mbs)\varphi_\ep(\mbt-\mbs)d\mbs\>
I_{q_1}\big(e_{j_1}^{(1)}(t^{(1)})\big)\cdots I_{q_m}\big(e_{j_m}^{(m)}(t^{(m)})\big)\>d\mbt
\nn\\&=&
\int \sum_{j_1,...,j_m=1}^n 
\int A_n(j_1,...,j_m;\mbs)\varphi_\ep(\mbt-\mbs)d\mbs\>
I_{q_1}\big(e_{j_1}^{(1)}(t^{(1)})\big)\cdots I_{q_m}\big(e_{j_m}^{(m)}(t^{(m)})\big)\>d\mbt
\eeas
\end{en-text}

Here is an example of applications of Theorem \ref{0110221242}. 
Let $\mH=L^2([0,1],dt)$. 
A functional 
that appeared in Nualart and Yoshida \cite{nualart2019asymptotic} on p.32 is 
\beas 
\cali_1 &=& 
\int_0^1\int_0^1	 \sqrt{n}\sum_{j=1}^n {\colorr(}D_sD_ta_{t_{j-1}}{\colorr)}
{\>\colorr {\sf q}}_j
\\&&\hspace{40pt}\times 
\sqrt{n}\sum_{k=1}^na_{t_{k-1}}(B_t-B_{t_{k-1}})1_k(t)dt
\\&&\hspace{40pt}\times
\sqrt{n}\sum_{\ell=1}^na_{t_{\ell-1}}(B_s-B_{t_{\ell-1}})1_{\ell}(s)ds,
\eeas
where $\tj=j/n$, 
$1_j=1_{(\tjm,\tj]}$, 
${\colorr{\sf q}}_j=(B_\tj-B_\tjm)^2-n^{-1}$, and 
$a_t=a(B_t)$ with $a\in C^\infty(\bbR)$ is of moderate growth and $B=(B_t)$ is a Brownian motion. 

Let $G=\{0\}\times\bbT\times\bbT$ for $\bbT=\bbR/\bbZ$. The element $t\in\bbT$ is naturally embedded into $[0,1]$, 
so $B_t$ makes sense for $t\in\bbT$. 
\begin{en-text}
Let 
\beas 
\dot{\bbH}_n(x^{(1)},x^{(2)},x^{(3)})
&=&
n^{3/2} \sum_{j,k,\ell=1}^n
\big(D_sD_ta_{t_{j-1}}\big)a_{t_{k-1}}a_{t_{\ell-1}}\>{\colorr{\sf q}}_j(B_t-B_{t_{k-1}})(B_s-B_{t_{\ell-1}})
1_k(t)1_{\ell}(s).
\eeas
\end{en-text}
For (\ref{202005160359}), let 
\beas
\dot{\bbH}_n(x^{(1)},x^{(2)},x^{(3)})
&=&
n^{1.5}\sum_{j,k,\ell=1}^n\big(D_{x^{(3)}}D_{x^{(2)}}a_{t_{j-1}}\big)a_{t_{k-1}}a_{t_{\ell-1}}
e_{n,j}^{(1)}(x^{(1)})e_{n,k}^{(2)}(x^{(2)})e_{n,\ell}^{(3)}(x^{(3)})
\eeas
with
\beas
e_{n,j}^{(1)}(x^{(1)}) &=& I_2(1_j^{\otimes2})\qquad (x^{(1)}\equiv0),\\
e_{n,k}^{(2)}(x^{(2)}) &=& I_1\big(1_k(x^{(2)})1_{(\tkm,x^{(2)}]}(\cdot)\big), \\
e_{n,\ell}^{(3)}(x^{(3)}) &=& I_1\big(1_\ell(x^{(3)})1_{(t_{\ell-1},x^{(3)}]}(\cdot)\big). 
\eeas
Then 
\beas 
\cali_1 &=& \bbH_n
\yeq 
\int_G \dot{\bbH}_n(x^{(1)},x^{(2)},x^{(3)})\mu(dx^{(1)}dx^{(2)}dx^{(3)}). 
\eeas

For $(\chi^{(1)}_n,g^{(1)}_{n,j})=(\{0\},0)$, 
$(\chi^{(2)}_n,g^{(2)}_{n,k})=([0,1/n],-\tkm)$ and 
$(\chi^{(3)}_n,g^{(3)}_{n,k})=([0,1/n],-t_{\ell-1})$, 
Condition (\ref{202005171254}) is satisfied because 
\beas
e_{n,j}^{(1)}(t) &=& I_2(1_j^{\otimes2})1_{\{0\}}(L_0x^{(1)})\quad (x^{(1)}\equiv0),
\eeas
\beas 
e_{n,k}^{(2)}(x^{(2)}) &=& I_1\big(1_k(x^{(2)})1_{(\tkm,x^{(2)}]}(\cdot)\big)1_k(x^{(2)})
\nn\\&=&
 I_1\big(1_k(x^{(2)})1_{(\tkm,x^{(2)}]}(\cdot)\big)1_{[0,1/n]}(x^{(2)}-\tkm)
\eeas
and 
\beas
e_{n,\ell}^{(3)}(x^{(3)}) &=& 
 I_1\big(1_\ell(x^{(3)})1_{(t_{\ell-1},x^{(3)}]}(\cdot)\big)1_{[0,1/n]}(x^{(3)}-t_{\ell-1}). 
\eeas

Now 
\beas
{\stackrel{\circ}{\bbH}}_n(x^{(1)},x^{(2)},x^{(3)})
&=& 
\begin{en-text}
n^{3/2} \sum_{j,k,\ell=1}^n
\big(D_{x^{(3)}-t_{\ell-1}}D_{x^{(1)}-\tkm}a_{t_{j-1}}\big)a_{t_{k-1}}a_{t_{\ell-1}}
1_1(x^{(2)})1_1(x^{(3)})
\nn\\&&\hspace{30pt}\times
{\colorr{\sf q}}_j(B_{x^{(1)}+\tkm}-B_{t_{k-1}})(B_{x^{(3)}+t_{\ell-1}}-B_{t_{\ell-1}})
\nn\\&=&
\end{en-text}
n^{3/2} \sum_{j,k,\ell=1}^n
\big(D_{t_{\ell-1}+x^{(3)}}D_{\tkm+x^{(2)}}a_{t_{j-1}}\big)a_{t_{k-1}}a_{t_{\ell-1}}
\nn\\&&\hspace{30pt}\times
I_2(1_j^{\otimes2})\>
I_1\big(1_{[\tkm,\tkm+x^{(2)}]}1_1(x^{(2)})\big)\>
I_1\big(1_{[\tkm,\tkm+x^{(3)}]}1_1(x^{(3)})\big).
\eeas
Then we obtain 
$\cali_1=O_{L^\inftym}(n^{-1})$ as $n\to\infty$ from (\ref{202005171343}) 
since $\mu\big(\prod_{i=1}^3\chi_n^{(i)}\big)=O(n^{-2})$ and 
${\stackrel{\circ}{\bbH}}_n(\mbx)=O_{L^\inftym}(n)$ uniformly in $\mbx\in G$ 
due to 
the exponent $e\big({\stackrel{\circ}{\bbH}}_n(x^{(1)},x^{(2)},x^{(3)})\big)=1$ 
by Theorem \ref{0110221242}. 
\begin{en-text}
Let $p\in\bbN$. 
Then \koko
\beas 
E\big[\cali_1^{2p}\big]
&=&
n^{3p}
\sum_{j_1,k_1,\ell_1=1}^n\cdots \sum_{j_{2p},k_{2p},\ell_{2p}=1}^n
\int_{[0,1]^{4p}}E\big[\bbI_n(s_1,t_1,j_1,k_1,\ell_1)\cdots\bbI_n(s_{2p},t_{2p},j_{2p},k_{2p},\ell_{2p})\big]\>
\\&&\hspace{150pt}\times
1_{k_1}(t_1)1_{\ell_1}(s_1)\cdots1_{k_{2p}}(t_{2p})1_{\ell_{2p}}(s_{2p})
ds_1dt_1\cdots ds_{2p}dt_{2p}
\\&\leq&
n^{3p}
\sum_{j_1,k_1,\ell_1=1}^n\cdots \sum_{j_{2p},k_{2p},\ell_{2p}=1}^n
\bigg\{
\sup_{s_1,t_1,...,s_{2p},t_{2p}}E\big[\bbI_n(s_1,t_1,j_1,k_1,\ell_1)\cdots\bbI_n(s_{2p},t_{2p},j_{2p},k_{2p},\ell_{2p})\big]
\\&&\hspace{100pt}\times
\int_{[0,1]^{4p}}1_{k_1}(t_1)1_{\ell_1}(s_1)\cdots1_{k_{2p}}(t_{2p})1_{\ell_{2p}}(s_{2p})
ds_1dt_1\cdots ds_{2p}dt_{2p}\bigg\}
\\&=&
n^{-p}
\sum_{j_1,k_1,\ell_1=1}^n\cdots \sum_{j_{2p},k_{2p},\ell_{2p}=1}^n
\sup_{s_1,t_1,...,s_{2p},t_{2p}}E\big[\bbI_n(s_1,t_1,j_1,k_1,\ell_1)\cdots\bbI_n(s_{2p},t_{2p},j_{2p},k_{2p},\ell_{2p})\big]
\\&=&
n^{-p}\times O(n^{-p})\yeq O(n^{-2p}). 
\eeas
Therefore, 
\beas 
\|\cali_1\|_{2p} &=& O(n^{-1}).
\eeas
\end{en-text}
\begin{en-text}
Then 
\beas 
\|\cali_1\|_p
&\leq&
n^{1.5}\sum_{j,k,\ell=1}^n\int_0^1\int_0^1
\bigg\|\big(D_sD_ta_{t_{j-1}}\big)a_{t_{k-1}}a_{t_{\ell-1}}\>q_j(B_t-B_{t_{k-1}})(B_s-B_{t_{\ell-1}})
1_k(t)1_{\ell}(s)
\bigg\|_p
1_k(t)1_{\ell}(s)dtds
\\&\leq&
n^{1.5}\sum_{j,k,\ell=1}^n\bigg\{
\sup_{s,t\in[0,1]}\bigg\|\big(D_sD_ta_{t_{j-1}}\big)a_{t_{k-1}}a_{t_{\ell-1}}
\>q_j\>(B_t-B_{t_{k-1}})1_k(t)\>(B_s-B_{t_{\ell-1}})1_{\ell}(s)
\bigg\|_p
\\&&\hspace{60pt}\times
\int_0^1\int_0^1
1_k(t)1_{\ell}(s)dtds\bigg\}
\\&=&
O(n^{1.5-0.5-2})\yeq O(n^{-1}). 
\eeas
\end{en-text}
The order $O(n^{-1})$ is different from the order  
$E[\Psi({\sf z}, {\sf x})\cali_1]=O(n^{-1.5})$, however 
these bounds do not conflict with each other 
because the bound obtained here is an $L^p$-bound. 
As for $L^p$-estimate of the derivative of the quasi-torsion, 
$\|D_{u_n}D_{u_n}D_{u_n}M_n\|_p=o(r_n)$ follows from this estimate 
in the martingale expansion approach. 
\begin{en-text}
In the present situation, by the orthogonality due to 
the specific support of $f^{(i)}_{j_i}$, 

the summations in $E[\check{\bbI}_n]$ (or in $\check{\bbI}_n$) 
become $c$ summations corresponding to $c$ connected components. 
Since each contraction gives a factor of order $O(n^{-1})$, 
we get a factor of order $O(n^{-\bar{R}_1})$ from the whole contractions. 
Moreover, the inner product in the expression of $\bbE_n$ gives 
a factor of order $O(n^{-(\bar{Q}_1-2\bar{R}_1)})$. 
Therefore, 
\beas
|\bbE_n| 
&=&
O(n^{c+\bar{R}_1-\bar{Q}_1}).
\eeas
Since the contraction takes place, if exists, between different factors $f^{(*)}_{*}$, 
the maximum $c$ is $km$. 
The maximum $\bar{R}_1$ is $k\bar{q}$. 
(Indeed, there is a pattern of contractions for which $c=km$ and $\bar{R}_1=k\bar{q}$.) 
Therefore
\beas 
c+\bar{R}_1-\bar{Q}_1
&=& 
c+\bar{R}_1-2k\bar{q}
\yleq 
km-k\bar{q},
\eeas
and hence we generally conclude 
\beas
|\bbE_n| 
&=&
O(n^{-k(\bar{q}-m)}).
\eeas
This completes the proof. 
\qed\halflineskip
\end{en-text}
\begin{en-text}
In fact, one contraction of an argument of $f^{(i)}_{j_i}$ and 
and argument of $f^{(i')}_{j_{i'}}$ 
gives order $n^{-1}$ 
\xout{and erases $\sum_{j_{i'}}$.} 
The number of the remaining sums is $2km-\bar{r}_1$. 
The number of the remaining arguments is totally $2k\bar{q}-2\bar{r}_1$, 
that gives $n^{-(2k\bar{q}-2\bar{r}_1)}$ by integration-by-parts. 
Therefore the order of a term of $E[\check{\bbI}_n]$ involving 
$\bar{r}$-contractions is \koko
\beas 
n^{2km\xout{-\bar{r}_1}}\times n^{-\bar{r}_1}\times n^{-(2k\bar{q}-2\bar{r}_1)}
&=&
n^{-2k(\bar{q}-m-2^{-1}\bar{r}_1)}\leq n^{-2k(2^{-1}\bar{q}-m)}.
\eeas
Thus, 
\beas 
E\big[\bbI_n^{2k}\big] 
&\leq&
C\big(q_1,...,q_m,k,T,C_0,C_1\big)^{2k}n^{-2k(\bar{q}-m)}.
\eeas
\end{en-text}

\begin{en-text}
\begin{theorem}\label{0110210316}
\beas 
\sup_{n\in\bbN}\big\|\bbI_n\big\|_{2k} &\leq& C(T,k,q,C_0,C_1)\yleq\infty.
\eeas
The bound on the right-hand side depends on $A_n(j_1,...,j_m)$ and $(f_j)$ 
only through the constants $C_0$ and $C_1$. 
\end{theorem}
\proof
We will consider the case $k=1$ since it is sufficient to show the general case. 
\beas 
\bbI_n^2 
&=&
\sum_{j_1,...,j_m,j_{m+1},...,j_{2m}=1}^n
F_n(j_1,...,j_m,j_{m+1},...,j_{2m})B(f_{j_1})\cdots B(f_{j_m})
B(f_{j_{m+1}})\cdots B(f_{j_{2m}})
\eeas
where 
\beas 
F_n(j_1,...,j_m,j_{m+1},...,j_{2m})
&=& 
A_n(j_1,...,j_m)A(j_{m+1},...,j_{2m}).
\eeas
\end{en-text}

\subsection{A bound for a polynomial of multiple Wiener integrals
}
Let ${\bf p}=(p_1,...,p_m)\in\bbN^m$ and ${\bf q}=(q_1,...,q_m)\in\bbN^m$. 
We will assume (\ref{0110231323}) and Condition that 
\bea\label{0205171525}
\max_{i=0,1,...,s}\>
\sup_{n\in\bbN}\sup_{(j_1,...,j_m)\in\{1,...,n\}^m}
\sup_{t_1,...,t_i\in[0,T]}
\big\|D^i_{t_1,...,t_i}A_n(j_1,...,j_m)\big\|_{4k}&\leq&C_2
\eea
for $s=2k\sum_{i=1}p_iq_i$ and some constant $C_2$. 
In this section, we will consider the functional 
\beas 
\bbK_n &=& \sum_{j_1,...,j_m=1}^n
A_n(j_1,...,j_m)I_{q_1}(f_{j_1}^{(1)})^{p_1}\cdots I_{q_m}(f_{j_m}^{(m)})^{p_m}. 
\eeas

\begin{theorem}\label{0205171708} 
Suppose that Conditions (\ref{0110231323}) and (\ref{0205171525}) are fulfilled. 
Then there exists a constant $C$ depending on $(C_0,C_2,{\bf p}, {\bf q})$ such that 
\beas 
\|\bbK_n\|_{2k} &\leq& C\>n^{-\half\xi}
\eeas
for 
\beas 
\xi
&=&
{\bf p}\cdot{\bf q}-m-\#\big\{i\in\{1,...,m\};\>p_i\geq2\text{ and }\>p_iq_i\text{ is even}\big\}.
\eeas
\end{theorem}
\proof 
Without loss of generality, we may assume 
with $p_1,...,p_{m-\ell}\geq2$ and $p_{m-\ell+1}=\cdots=p_m=1$; 
in particular, all $P_1,...,p_m\geq2$ when $\ell=0$. 
By the product formula applied to each factor $I_{q_i}(f^{(i)}_{j_i})^{p_i}$, we obtain 
\beas 
\bbK_n 
&=& 
\sum_{\nu_1,...,\nu_\ell\in\bbZ_+}
\bigg\{\prod_{i'=1}^\ell c_{\nu_{i'}}[q_{i'},p_{i'}]\bigg\}
\\&&\times 
\sum_{j_1,...,j_m=1}^n n^{-\sum_{i=1}^\ell\nu_i}
\bigg\{A_n(j_1,...,j_m)
\bigg(\prod_{i=1}^\ell I_{p_iq_i-2\nu_i}\big(n^{\nu_i}\tilde{\otimes}^{p_iq_i}_{\nu_i}f^{(i)}_{j_i}\big)\bigg)
\bigg(\prod_{i=\ell+1}^m I_{q_i}(f_{j_i}^{(i)})\bigg)
\bigg\}
\eeas
The desired estimate follows from the product formula and the estimate (\ref{202004051415}) in Section \ref{202004270705} based on Theorem \ref{0110221242}. 
Indeed, 
the exponent 
\beas 
e(\bbK_n) 
&=&
\max_{(\nu_1,,.,\nu_m)\in\bbZ_+^\ell}
\bigg\{-\sum_{i=1}^\ell\nu_i-0.5\bigg(\sum_{i=1}^\ell(p_iq_i-2\nu_i)+\sum_{i=\ell+1}^mp_iq_i\bigg)
+m
\nn\\&&\hspace{80pt}
-0.5\bigg(\#\big\{i\in\{1,...,\ell\};\>p_iq_i-2\nu_i>0\big\}+(m-\ell)\bigg)
\bigg\}
\nn\\&=&
-0.5{\bf p}\cdot{\bf q}+0.5m+0.5\#\big\{i\in\{1,...,\ell\};\>p_iq_i\text{ is  even}\big\}. 
\eeas
\qed\halflineskip
\begin{en-text}
For the meantime, we assume $p_i\geq2$ for all $i\in\{1,...,m\}$.

\koko

We denote ${\bf p}\cdot{\bf q}=\sum_{i=1}^mp_iq_i$. 
We use the product formula. For $f\in\calh^{\tilde{\otimes}q}$, 
\beas 
I_q(f)^p
&=& 
\sum_{\nu\in\bbZ_+}c_\nu[q,p]I_{pq-2\nu}(\tilde{\otimes}^{pq}_{2\nu}f),
\eeas
where $c_\nu[q,p]$ are universal constants; in particular, $c_\nu[q,p]=0$ unless $0\leq pq-2\nu\leq pq$.

\bi

\im 
\im 
\beas 
\bbK_n &=& \sum_{j_1,...,j_m=1}^n\sum_{\nu_1,...,\nu_m\in\bbZ_+}
A_n(j_1,...,j_m)\prod_{i=1}^m\bigg\{c_{\nu_i}(q_i,p_i)
I_{p_iq_i-2\nu_i}(\tilde{\otimes}^{p_iq_i}_{2\nu_i}f^{(i)}_{j_i})\bigg\}
\eeas
\beas 
\bbK_n &=& \sum_{\nu_1,...,\nu_m\in\bbZ_+}\bigg\{\prod_{i'=1}^mc_{\nu_{i'}}(q_{i'},p_{i'})\bigg\}\sum_{j_1,...,j_m=1}^n
A_n(j_1,...,j_m)
\prod_{i=1}^m\bigg\{I_{p_iq_i-2\nu_i}(\tilde{\otimes}^{p_iq_i}_{2\nu_i}f^{(i)}_{j_i})\bigg\}
\eeas
\im $\pi(0)=\{0\}$, $\pi(1)=\bbN$
\beas 
\bbK_n 
&=& 
\sum_{\nu_1,...,\nu_m\in\bbZ_+}\bigg\{\prod_{i'=1}^mc_{\nu_{i'}}(q_{i'},p_{i'})\bigg\}\sum_{j_1,...,j_m=1}^n
A_n(j_1,...,j_m)
\\&&\times 
\prod_{i=1}^m\bigg\{\sum_{\ep_i\in\{0,1\}}
1_{\{p_iq_i-2\nu_i\in \pi(\ep_i)\}}I_{p_iq_i-2\nu_i}(\tilde{\otimes}^{p_iq_i}_{2\nu_i}f^{(i)}_{j_i})\bigg\}
\\&=& 
\sum_{\nu_1,...,\nu_m\in\bbZ_+}\bigg[\bigg\{\prod_{i'=1}^mc_{\nu_{i'}}(q_{i'},p_{i'})\bigg\}
\\&&\times 
\sum_{\ep_1,...,\ep_m\in\{0,1\}}\sum_{j_1,...,j_m=1}^nA_n(j_1,...,j_m)
\prod_{i=1}^m\bigg\{
1_{\{p_iq_i-2\nu_i\in \pi(\ep_i)\}}I_{p_iq_i-2\nu_i}(\tilde{\otimes}^{p_iq_i}_{2\nu_i}f^{(i)}_{j_i})\bigg\}\bigg]
\eeas

\im 
\beas 
\|\bbK_n\|_{2k} 
&\leq&
C'\sum_{\nu_1,...,\nu_m\in\bbZ_+}\bigg[\bigg\{\prod_{i=1}^mc_{\nu_{i}}(q_{i},p_{i})\bigg\}
\\&&\times 
n^{-\sum_{i=1}^m\nu_i1_{\{p_iq_i-2\nu_i=0\}}}
\bigg\|\sum_{j_1,...,j_m=1}^nA_n(j_1,...,j_m)
\prod_{i\in\{1,...,m\}:p_iq_i-2\nu_i\geq1}^m
I_{p_iq_i-2\nu_i}(\tilde{\otimes}^{p_iq_i}_{2\nu_i}f^{(i)}_{j_i})\bigg\|_{2k}\bigg]
\\&\leq&
C''\sum_{\nu_1:0\leq\nu_1\leq p_1q_1/2}\cdots\sum_{\nu_m:0\leq\nu_m\leq p_mq_m/2}
\bigg[\bigg\{\prod_{i=1}^mc_{\nu_{i}}(q_{i},p_{i})\bigg\}
\\&&\times 
n^{-\sum_{i=1}^m\nu_i1_{\{p_iq_i-2\nu_i=0\}}}n^{\sum_{i=1}^m1_{\{p_iq_i-2\nu_i=0\}}}
n^{\sum_{i=1}^m[-2^{-1}(p_iq_i-2\nu_i-1)-\nu_i]1_{\{p_iq_i-2\nu_i\geq1\}}}
\eeas
by Theorem \ref{0110221242}, 
where 
the factor $n^{\sum_{i=1}^m1_{\{p_iq_i-2\nu_i=0\}}}$ comes from the sum in $\sum_{j_1,...,j_m=1}^n$ 
that goes out of the norme $\|\cdot\|_{2k}$, and 
the last $-\nu_i$ comes from $\nu_i$ contractions in $\tilde{\otimes}^{p_iq_i}_{2\nu_i}f^{(i)}_{j_i}$. 
The exponent of $n$ is totally 
\beas &&
-\sum_{i=1}^m(\nu_i-1)1_{\{p_iq_i-2\nu_i=0\}}+\sum_{i=1}^m[-2^{-1}(p_iq_i-2\nu_i-1)-\nu_i]1_{\{p_iq_i-2\nu_i\geq1\}}
\\&=&
-\sum_{i=1}^m(\nu_i-1)1_{\{p_iq_i-2\nu_i=0\}}-2^{-1}\sum_{i=1}^m(p_iq_i-1)1_{\{p_iq_i-2\nu_i\geq1\}}
\\&=&
-\sum_{i=1}^m(\nu_i-1)1_{\{p_iq_i-2\nu_i=0\}}-2^{-1}\sum_{i=1}^m(p_iq_i-1)\big(1-1_{\{p_iq_i-2\nu_i=0\}}\big)
\\&=&
\sum_{i=1}^m\big\{2^{-1}(p_iq_i-1)-\nu_i+1\big\}1_{\{p_iq_i-2\nu_i=0\}}-2^{-1}\sum_{i=1}^m(p_iq_i-1)
\\&=&
2^{-1}\sum_{i=1}^m1_{\{p_iq_i-2\nu_i=0\}}-2^{-1}\sum_{i=1}^m(p_iq_i-1)
\\&\leq&
-2^{-1}\sum_{i=1}^m(p_iq_i-1)+2^{-1}\#\big\{i\in\{1,...,m\};\>p_iq_i\text{ is even}\big\}.
\eeas

\im Conclusively, 
\begin{theorem} Let $q_i\geq1$ and $p_i\geq2$ for all $i\in\{1,...,m\}$. Then 
\beas 
\|\bbK_n\|_{2k} &\leq& C\>n^{-\half\xi}
\eeas
for 
\beas 
\xi
&=&
{\bf p}\cdot{\bf q}-m-\#\big\{i\in\{1,...,m\};\>p_iq_i\text{ is even}\big\}.
\eeas
\end{theorem}

\im 
In the above proof, if $p_i=1$, then we did not need to apply the product formula to $I_{q_i}$. 
When $p_i=1$, the parity of $p_iq_i$ would unnecessarily increase the exponent. 
We should distinguish the factors depending on $p_i=1$ or $p_i\geq2$.

\im Try once again. 
\im Let $m\geq2$, $1\leq\ell< m$. 
\im A functional
\beas 
\bbK_n &=& \sum_{j_1,...,j_m=1}^n
A_n(j_1,...,j_m)I_{q_1}(f_{j_1}^{(1)})^{p_1}\cdots I_{q_m}(f_{j_m}^{(m)})^{p_m}
\eeas
with $p_1,...,p_{m-\ell}\geq2$ and 
$p_{m-\ell+1},...,p_m=1$. (Only the last $\ell$ exponents are $1$.)

\im product formula. For $f\in\calh^{\tilde{\otimes}q}$, 
\beas 
I_q(f)^p
&=& 
\sum_{\nu\in\bbZ_+}c_\nu(q,p)I_{pq-2\nu}(\tilde{\otimes}^{pq}_{2\nu}f)
\eeas
Remark that $c_\nu(q,p)=0$ unless $0\leq pq-2\nu\leq pq$. 

\im 
\beas 
\bbK_n &=& \sum_{j_1,...,j_m=1}^n\sum_{\nu_1,...,\nu_{m-\ell}\in\bbZ_+}
A_n(j_1,...,j_m)
\bigg[\prod_{i=1}^{m-\ell}\bigg\{c_{\nu_i}(q_i,p_i)
I_{p_iq_i-2\nu_i}(\tilde{\otimes}^{p_iq_i}_{2\nu_i}f^{(i)}_{j_i})\bigg\}\bigg]
\bigg[\prod_{i=m-\ell+1}^mI_{q_i}(f^{(i)}_{j_i})\bigg]
\eeas

\im $\pi(0)=\{0\}$, $\pi(1)=\bbN$
\beas 
\bbK_n 
&=& 
\sum_{\nu_1,...,\nu_{m-\ell}\in\bbZ_+}\bigg\{\prod_{i'=1}^{m-\ell}c_{\nu_{i'}}(q_{i'},p_{i'})\bigg\}\sum_{j_1,...,j_m=1}^n
A_n(j_1,...,j_m)
\\&&\times 
\prod_{i=1}^{m-\ell}\bigg\{\sum_{\ep_i\in\{0,1\}}
1_{\{p_iq_i-2\nu_i\in \pi(\ep_i)\}}I_{p_iq_i-2\nu_i}(\tilde{\otimes}^{p_iq_i}_{2\nu_i}f^{(i)}_{j_i})\bigg\}
\times\prod_{i=m-\ell+1}^mI_{q_i}(f^{(i)}_{j_i})
\\&=& 
\sum_{\nu_1,...,\nu_{m-\ell}\in\bbZ_+}\bigg[\bigg\{\prod_{i'=1}^{m-\ell}c_{\nu_{i'}}(q_{i'},p_{i'})\bigg\}
\\&&\times 
\sum_{\ep_1,...,\ep_{m-\ell}\in\{0,1\}}\sum_{j_1,...,j_m=1}^nA_n(j_1,...,j_m)
\prod_{i=1}^{m-\ell}\bigg\{
1_{\{p_iq_i-2\nu_i\in \pi(\ep_i)\}}I_{p_iq_i-2\nu_i}(\tilde{\otimes}^{p_iq_i}_{2\nu_i}f^{(i)}_{j_i})\bigg\}\bigg]
\\&&\times
\prod_{i=m-\ell+1}^mI_{q_i}(f^{(i)}_{j_i})
\eeas

\im 
\beas 
\|\bbK_n\|_{2k} 
&\leq&
C'\sum_{\nu_1,...,\nu_m\in\bbZ_+}\bigg[\bigg\{\prod_{i=1}^{m-\ell}c_{\nu_{i}}(q_{i},p_{i})\bigg\}
\\&&\times 
n^{-\sum_{i=1}^{m-\ell}\nu_i1_{\{p_iq_i-2\nu_i=0\}}}
\bigg\|\sum_{j_1,...,j_m=1}^nA_n(j_1,...,j_m)
\prod_{i\in\{1,...,m-\ell\}:p_iq_i-2\nu_i\geq1}
I_{p_iq_i-2\nu_i}(\tilde{\otimes}^{p_iq_i}_{2\nu_i}f^{(i)}_{j_i})
\\&&\hspace{150pt}\times
\prod_{i=m-\ell+1}^mI_{q_i}(f^{(i)}_{j_i})\bigg\|_{2k}\bigg]
\\&\leq&
C''\sum_{\nu_1:0\leq\nu_1\leq p_1q_1/2}\cdots\sum_{\nu_m:0\leq\nu_m\leq p_mq_m/2}
\bigg\{\prod_{i=1}^{m-\ell}c_{\nu_{i}}(q_{i},p_{i})\bigg\}
\\&&\times 
n^{-\sum_{i=1}^{m-\ell}\nu_i1_{\{p_iq_i-2\nu_i=0\}}}n^{\sum_{i=1}^{m-\ell}1_{\{p_iq_i-2\nu_i=0\}}}
n^{-\sum_{i=1}^{m-\ell}\nu_i1_{\{p_iq_i-2\nu_i\geq1\}}}
\\&&
\times
n^{\sum_{i=1}^{m-\ell}[-2^{-1}(p_iq_i-2\nu_i-1)]1_{\{p_iq_i-2\nu_i\geq1\}}-2^{-1}\sum_{i=m-\ell+1}^m(p_iq_i-1)}
\eeas
by Theorem \ref{0110221242}, 
where 
the factor $n^{\sum_{i=1}^{m-\ell}1_{\{p_iq_i-2\nu_i=0\}}}$ comes from the sum in $\sum_{j_1,...,j_m=1}^n$ 
that goes out of the norme $\|\cdot\|_{2k}$, and 
the factor $n^{-\sum_{i=1}^{m-\ell}\nu_i}$ comes from $\nu_i$ contractions in $\tilde{\otimes}^{p_iq_i}_{2\nu_i}f^{(i)}_{j_i}$. 

The exponent of $n$ is totally 
\beas &&
-\sum_{i=1}^{m-\ell}(\nu_i-1)1_{\{p_iq_i-2\nu_i=0\}}+\sum_{i=1}^{m-\ell}[-2^{-1}(p_iq_i-2\nu_i-1)-\nu_i]1_{\{p_iq_i-2\nu_i\geq1\}}-2^{-1}\sum_{i=m-\ell+1}^m(p_iq_i-1)
\\&=&
-\sum_{i=1}^{m-\ell}(\nu_i-1)1_{\{p_iq_i-2\nu_i=0\}}-2^{-1}\sum_{i=1}^{m-\ell}(p_iq_i-1)1_{\{p_iq_i-2\nu_i\geq1\}}
-2^{-1}\sum_{i=m-\ell+1}^m(p_iq_i-1)
\\&=&
-\sum_{i=1}^{m-\ell}(\nu_i-1)1_{\{p_iq_i-2\nu_i=0\}}-2^{-1}\sum_{i=1}^{m-\ell}(p_iq_i-1)\big(1-1_{\{p_iq_i-2\nu_i=0\}}\big)
-2^{-1}\sum_{i=m-\ell+1}^m(p_iq_i-1)
\\&=&
\sum_{i=1}^{m-\ell}\big\{2^{-1}(p_iq_i-1)-\nu_i+1\big\}1_{\{p_iq_i-2\nu_i=0\}}-2^{-1}\sum_{i=1}^{m-\ell}(p_iq_i-1)
-2^{-1}\sum_{i=m-\ell+1}^m(p_iq_i-1)
\\&=&
2^{-1}\sum_{i=1}^{m-\ell}1_{\{p_iq_i-2\nu_i=0\}}-2^{-1}\sum_{i=1}^{m-\ell}(p_iq_i-1)
-2^{-1}\sum_{i=m-\ell+1}^m(p_iq_i-1)
\\&\leq&
-2^{-1}\sum_{i=1}^{m}(p_iq_i-1)+2^{-1}\#\big\{i\in\{1,...,m-\ell\};\>p_iq_i\text{ is even}\big\}
-2^{-1}\sum_{i=m-\ell+1}^m(p_iq_i-1).
\eeas

\im In this way, we obtain 
\begin{theorem}\label{0111301152} Let $1\leq\ell<m$. 
Let $q_i\geq1$ for $i\in\{1,...,m\}$, and let $p_i\geq2$ for $i\in\{1,...,m-\ell\}$ and 
$p_i=1$ for $i=\{m-\ell+1,...,m\}$. Then 
\beas 
\|\bbK_n\|_{2k} &\leq& C\>n^{-\half\xi}
\eeas
for 
\beas 
\xi
&=&
{\bf p}\cdot{\bf q}-m-\#\big\{i\in\{1,...,m-\ell\};\>p_iq_i\text{ is even}\big\}.
\eeas
\end{theorem}


Combining the above two results, we obtain 
\begin{theorem}\label{0112061610} 
Let $q_i\geq1$ for $i\in\{1,...,m\}$. 
Then 
\beas 
\|\bbK_n\|_{2k} &\leq& C\>n^{-\half\xi}
\eeas
for 
\beas 
\xi
&=&
{\bf p}\cdot{\bf q}-m-\#\big\{i\in\{1,...,m\};\>p_i\geq2\text{ and }\>p_iq_i\text{ is even}\big\}.
\eeas
\end{theorem}
{\colorr \proof This follows from the product formula and the estimate (\ref{202004051415}) in Section \ref{202004270705} based on Theorem \ref{0110221242}. }
\ei
\end{en-text}

\bibliographystyle{spmpsci}      
\bibliography{bibtex-20180615-20191212-20200312-20200418-20201101}   

\end{document}